%% file: whittle_likelihood_revision.tex
\documentclass[12pt]{article}
\usepackage[utf8]{inputenc}
\usepackage{graphicx}
\usepackage{amsmath, amssymb,theorem,verbatim}
\usepackage{hyperref}
\usepackage{apacite}
\usepackage{textcomp} 
\usepackage{setspace}
\usepackage{mathtools}
\usepackage[round]{natbib}
\usepackage{multirow} 
\usepackage{pdflscape}
\usepackage[usenames,dvipsnames]{xcolor}
\usepackage{makecell}

\numberwithin{equation}{section}

\newtheorem{theorem}{Theorem}[section]
\newtheorem{lemma}{Lemma}[section] 
\newtheorem{assumption}{Assumption}[section] 
\newtheorem{defin}{Definition}[section] 
\newtheorem{corollary}{Corollary}[section] 
\newtheorem{example}{Example}[section]
\newtheorem{remark}{Remark}[section] 

\newcommand{\Pcon}{\stackrel{\mathcal{P}}{\rightarrow}}

\newcommand{\R}{\operatorname{Re}}
\newcommand{\I}{\operatorname{Im}}

\newcommand{\Xunder}{\underline{X}}

\newcommand{\cov}{\mathrm{cov}}
\newcommand{\spa}{\mathrm{sp}}
\newcommand{\var}{\mathrm{var}}
\newcommand{\cum}{\mathrm{cum}}
\newcommand{\Ex}{\mathbb{E}}
\newcommand{\diag}{\mathrm{diag}}

\newcommand{\svdots}{\raisebox{3pt}{$\scalebox{.75}{\vdots}$}}
\newcommand{\sddots}{\raisebox{3pt}{$\scalebox{.75}{$\ddots$}$}}

\title{Reconciling the Gaussian and Whittle Likelihood with an
application to estimation in the frequency domain}
\author{Suhasini Subba Rao\footnote{Texas A\&M University, College
    Station, Texas, TX 77845, U.S.A.} \hspace{1mm} and Junho
  Yang\footnote{Authors ordered alphabetically}} 
\date{\today}

\begin{document}
\maketitle

\begin{abstract}

In time series analysis there is an apparent dichotomy between time
and frequency domain methods. The aim of this paper is to draw connections 
between frequency and time domain methods. Our focus will be on reconciling
the Gaussian likelihood and the Whittle likelihood.
We derive an exact, interpretable, bound between the
Gaussian and Whittle likelihood of a second order stationary time
series. The
derivation is based on obtaining the transformation which is
biorthogonal to the discrete Fourier transform of the time
series. Such a transformation yields a new decomposition for the 
inverse of a Toeplitz matrix and enables the representation of the Gaussian
likelihood within the frequency domain. 
We show that the
difference between the Gaussian and Whittle likelihood is due
to the omission of the best linear predictions outside the domain of
observation in the periodogram associated with the Whittle
likelihood. Based on this result, we obtain an approximation for the difference between the
Gaussian and Whittle likelihoods in terms of the best
fitting, finite order autoregressive parameters. These approximations are
used to define two new frequency domain
quasi-likelihoods criteria. We show that
these new criteria can yield a better approximation of the spectral
divergence criterion,  as compared to both the Gaussian and Whittle
likelihoods.
In simulations, we show that the proposed estimators have satisfactory finite sample properties. 

\vspace{2mm}

\end{abstract}

\input{main}
\bibliography{bib_pred}
\bibliographystyle{plainnat}

\newpage

\appendix

\input{appendix_main_proof}

\input{appendix_baxter}

\input{appendix_consistency}

\input{appendix_bias}

\input{appendix_simulations}

\end{document}

%% file: main.tex
\section{Introduction}

In his seminal work, Whittle (1951, 1953) \nocite{p:whi-51,p:whi-53} 
introduced the Whittle likelihood as an approximation of the Gaussian likelihood.
A decade later, the asymptotic sampling properties of moving average
models fitted using the Whittle likelihood were derived in \cite{p:wal-64}. Subsequently, the Whittle
likelihood has become a popular method for parameter estimation of
various stationary
time series (both long and short memory) and spatial models. The
Whittle likelihood is computationally a very attractive method for
estimation. Despite the considerable improvements in technology, interest in the Whittle likelihood has not
abated. The Whittle likelihood has gained further traction
as a quasi-likelihood (or as an information criterion, see \cite{p:par-83}) between the periodogram and the
spectral density. Several diverse applications of the Whittle likelihood can be
found in \cite{p:fox-taq-86}, \cite{p:dah-87} (for spatial processes), \cite{p:rob-95},
\cite{p:dah-00}, \cite{p:che-00}, \cite{p:gir-01}, \cite{p:cho-04}, \cite{p:abd-07}, \cite{p:sha-wu-07},
\cite{b:gir-12} (long memory time series and local Whittle methods),
  \cite{p:pan-13}, \cite{p:kir-19} (Bayesian spectral methods),
 and \cite{p:del-eic-19} (functional time series), to name but
a few. 

Despite its advantages, it is well known that 
for small samples the Whittle likelihood 
can give rise to estimators with a substantial bias (see
\cite{b:pri-81} and \cite{p:dah-88}). \cite{p:dah-88} shows that the
finite sample bias in the periodogram impacts
the performance of the Whittle likelihood.  Motivated by this discrepancy,
\cite{p:olh-19} proposes the debiased Whittle likelihood, which fits
directly to the expectation of the periodogram rather than the
limiting spectral density. Alternatively,
\cite{p:dah-88} shows that the tapered periodogram is
better at capturing the features in the spectral density, such as
peaks, than the regular periodogram. He uses this as the basis of 
the tapered Whittle likelihood. Empirical studies show that the tapered Whittle
likelihood yields a smaller bias than the regular Whittle
likelihood. As a theoretical justification, Dahlhaus (1988,
1990)\nocite{p:dah-88,p:dah-90} uses an alternative asymptotic framework to show 
that tapering yields a good approximation to the inverse of the Toeplitz matrix.
It is worth
mentioning that within the time domain, several authors, including
Shaman (1975, 1976), \nocite{p:sha-75,p:sha-76} \cite{p:bha-82}
and \cite{p:cou-82}, have studied approximations to the
inverse of the Toeplitz matrix. These results can be used to
approximate the Gaussian likelihood. 

However, as far as we are aware, there are no results which explain
what is lost when using the Whittle likelihood rather than the Gaussian likelihood.  
The objective of this paper is to address some of these issues. The benefits of
such insights are not only of theoretical interest but also lead to 
the development of computationally simple frequency domain
methods which are comparable with the Gaussian
likelihood. 

We first recall the definition of the Gaussian and Whittle
likelihood. Our aim is to fit a parametric second order stationary
model with spectral density $f_{\theta}(\omega)$ and corresponding
autocovariance function $\{c_{f_\theta}(r)\}_{r\in \mathbb{Z}}$
to the observed time series $\{X_{t}\}_{t=1}^{n}$.  The
(quasi) log-Gaussian likelihood is proportional to
\begin{equation}
\label{eq:likeG}
\mathcal{L}_{n}(\theta;\Xunder_{n}) =
n^{-1}\left(\Xunder_{n}^{\prime}\Gamma_{n}(f_{\theta})^{-1}\Xunder_{n} + \log |\Gamma_{n}(f_{\theta})|\right)
\end{equation}
where $\Gamma_{n}(f_{\theta})_{s,t} = c_{f_\theta}(s-t)$, $|A|$ denotes the
determinant of the matrix $A$ and $\Xunder_{n}^{\prime} = (X_1,\ldots,X_{n})$. 
In contrast, the Whittle likelihood is a ``spectral divergence'' between the periodogram and the candidate
spectral density. There are two subtly different methods for
defining this contrast, one is with an integral the other is to use the
Riemann sum. In this paper, we focus on the Whittle likelihood defined
in terms of the Riemann sum over the fundamental frequencies
\begin{equation}
K_{n}(\theta;\Xunder_{n}) = n^{-1}\sum_{k=1}^{n}\left(\frac{|J_{n}(\omega_{k,n})|^{2}}{f_{\theta}(\omega_{k,n})}
+ \log f_{\theta}(\omega_{k,n})\right) \quad \omega_{k,n} = \frac{2\pi k}{n},
\end{equation}
where $J_{n}(\omega_{k,n}) =
n^{-1/2}\sum_{t=1}^{n}X_{t}e^{it\omega_{k,n}}$ is 
the discrete Fourier transform (DFT) of the observed time series. To compare
the Gaussian and Whittle likelihood, we rewrite the Whittle likelihood in 
matrix form. We define the $n\times n$ circulant matrix 
$C_{n}(f_\theta)$ with entries $(C_{n}(f_\theta))_{s,t}=
n^{-1}\sum_{k=1}^{n}f_{\theta}(\omega_{k,n}) e^{-i(s-t)\omega_{k,n}}$. 
The Whittle likelihood $K_{n}(\theta;\Xunder_{n})$ can be written as 
\begin{equation}
\label{eq:like1}
K_{n}(\theta;\Xunder_{n}) = n^{-1}\left( \underline{X}_{n}^{\prime}C_{n}(f_{\theta}^{-1})\underline{X}_{n}+
\sum_{k=1}^{n}\log f_{\theta}(\omega_{k,n}) \right).
\end{equation}
To obtain an exact expression for
$\Gamma_{n}(f_{\theta})^{-1} -C_{n}(f_{\theta}^{-1})$ and 
$\underline{X}^{\prime}_{n}[\Gamma_{n}(f_{\theta})^{-1}-C_{n}(f_{\theta}^{-1})]\Xunder_n$,
we focus
on the DFT of the time series. The idea is to obtain the  linear transformation of the
observed time series $\{X_{t}\}_{t=1}^{n}$ which is biorthogonal to the regular DFT,
$\{J_{n}(\omega_{k,n})\}_{k=1}^{n}$.
The biorthogonal transform, when coupled with the regular
DFT, exactly decorrelates the time series. In Section \ref{sec:DFT}, we show that the
biorthogonal transform corresponding to the regular DFT contains the regular
DFT \emph{plus} the Fourier transform of the best linear predictors of the
time series outside the domain of observation. Since this
transformation completes the information not found in the regular DFT, 
 we call it the complete DFT. It is common to use the Cholesky
 decomposition to decompose the inverse of a Toeplitz matrix.
An interesting aspect of the biorthogonal transformation is that it
provides an alternative decomposition of the inverse of a Toeplitz matrix.

In Section \ref{sec:gauss}, we show that 
the complete DFT, together with the regular DFT, allows us to rewrite the Gaussian
likelihood within the frequency domain (which, as far as we are aware,
is new). Further, it is well known that the
Whittle likelihood has a bias due to the boundary effect. By rewriting
the Gaussian likelihood within the frequency domain we show that the
Gaussian likelihood avoids the boundary effect problem by predicting
the time series outside the domain of observation. Precisely,
the approximation
error between the Gaussian and Whittle likelihood is due to the
omission of these linear predictors in the regular DFT.  From this
result, we observe that the greater
the persistence in the time series model (which corresponds to a more
peaked spectral density) the larger the loss in
approximating the complete DFT with the regular DFT. 
In order to obtain a better approximation of the Gaussian likelihood
in the frequency domain, it is of interest to approximate the
difference of the two likelihoods $\mathcal{L}_{n}(\theta;\Xunder_{n})-K_{n}(\theta;\Xunder_{n})$.
For autoregressive processes of finite order, we obtain an analytic
expression for the difference in the two likelihoods in terms of the
AR parameters (see equation (\ref{eq:LKdiff})). For general second
order stationary models, the expression is more complex. In Section
\ref{sec:approx}, we obtain an approximation for 
$\mathcal{L}_{n}(\theta;\Xunder_{n})-K_{n}(\theta;\Xunder_{n})$
in terms of the infinite order (causal/minimum phase) autoregressive factorisation of
$f_{\theta}(\omega)=\sigma^2|1-\sum_{j=1}^{\infty}\phi_{j}e^{-ij\omega}|^{-2}$. 
We show that this approximation is the first order term in a series
expansion of the inverse of the Toeplitz matrix, $\Gamma_{n}(f)^{-1}$. More precisely, in Section \ref{sec:higherexpan},
 we show that $\Gamma_{n}(f)^{-1}$ can be expressed in terms of
$C_{n}(f_{\theta})^{-1}$ plus a polynomial-type series expansion of the
AR$(\infty)$ coefficients. 

In Section \ref{sec:newlikelihood}, we obtain an approximation for the
difference
$\mathcal{L}_{n}(\theta;\Xunder_{n})-K_{n}(\theta;\Xunder_{n})$ in
terms of a finite order autoregressive process. We use this to
define two spectral divergence criteria which are ``almost''
unbiased estimators of the spectral divergence between the true
(underlying spectral) density and the parametric spectral density. We
use these criteria to define two new frequency domain estimators. In
Section \ref{sec:sample}, we obtain the asymptotic sampling properties
of the new likelihood estimators including the asymptotic bias and variance.
Finally, in Section \ref{sec:sim}, we illustrate and compare the
proposed frequency domain estimators through some simulations.
We study the performance of the estimation scheme when the
parametric model is both correctly specified and 
misspecified. 

The proofs can be found in the Supplementary material. The main proofs can be found
in Appendix \ref{sec:proofbio}, \ref{sec:proofs2},
\ref{sec:consistent} and 
\ref{sec:bias}. Baxter type inequalities for derivatives of
finite predictors can be found in Appendix \ref{sec:baxter}. These
results are used to obtain an approximation for the difference between the derivatives
of the Gaussian and Whittle likelihood. 
In Appendix \ref{sec:bias} we derive an expression
for the asymptotic bias of the Gaussian, Whittle likelihoods, and the new
frequency domain likelihoods, described above. In Appendix
\ref{sec:sim_appendix}, \ref{sec:sim-long} and \ref{sec:alternative} we present
additional simulations.

\section{The Gaussian likelihood in the frequency domain}

\subsection{Preliminaries}

In this section, we introduce most of the notation used in the paper, it can be skipped on first
reading. 
To reduce notation, 
we omit the symbol $\Xunder_{n}$ in the Gaussian and Whittle
likelihood.
Moreover, since the focus in this paper will be on the first terms in the
Gaussian and Whittle likelihoods we use 
$\mathcal{L}_{n}(\theta) $ and $K_{n}(\theta)
$ to denote only these terms:
\begin{equation}
\label{eq:L1L2}
\mathcal{L}_{n}(\theta) =
n^{-1}\Xunder_{n}^{\prime}\Gamma_{n}(f_{\theta})^{-1}\Xunder_{n}
 \quad \textrm{and} \quad
K_{n}(\theta) =
  n^{-1}\underline{X}_{n}^{\prime}C_{n}(f_{\theta}^{-1})\underline{X}_{n}.
\end{equation}
Let $A^{*}$ denote the conjugate transpose of the matrix $A$.  
We recall that the circulant matrix $C_{n}(g)$ can be written as
$C_{n}(g)= F_{n}^{*}\Delta_{n}(g)F_{n}$, 
where  $\Delta_{n}(g) =
\diag(g(\omega_{1,n}),\ldots,g(\omega_{n,n}))$ (diagonal matrix) and
$F_{n}$ is the $n \times n$ DFT matrix with entries $(F_{n})_{k,t}
=n^{-1/2} e^{it\omega_{k,n}}$.  
We recall that the eigenvalues and  the corresponding eigenvectors of any
circulant matrix $C_{n}(g)$ are $\{g(\omega_{k,n})\}_{k=1}^{n}$ and
$\{\underline{e}_{k,n}^{\prime} = (e^{ik\omega_{1,n}},\ldots,e^{ik\omega_{n,n}})\}_{k=1}^{n}$ respectively.

In general, we assume that $\Ex[X_{t}]=0$ (as it makes the derivations
cleaner).
We use $\{c_{f}(r)\}_{r \in \mathbb{Z}}$ to denote an autocovariance function and 
$f(\omega) =
\sum_{r\in\mathbb{Z}}c_{f}(r)e^{ir\omega}$ its corresponding
spectral density. Sometimes, it will be necessary to make explicit the
true underlying
covariance (equivalently the spectral density) of the process. In this
case, we use the notation $\cov_{f}(X_{t},X_{t+r}) = \Ex_{f}[X_{t}X_{t+r}] =c_{f}(r)$.
Next we define the norms we will use. 
Suppose $A$ is a $n \times n$ square matrix, let $\|A\|_{p} =
(\sum_{i,j=1}^{n}|a_{i,j}|^{p})^{1/p}$ be an entrywise $p$-norm for $p\geq 1$, and $\|A\|_{spec}$ denote the
spectral norm. Let $\|X\|_{\Ex,p} = \left(\Ex |X|^{p}\right)^{1/p}$, where
$X$ is a random variable.  
For the $2\pi$-periodic square integrable function $g$ 
with $g(\omega) = \sum_{r\in \mathbb{Z}}g_{r}e^{ir\omega}$, we use the 
sub-multiplicative norm
$\|g\|_{K} = \sum_{r\in \mathbb{Z}}(2^{K}+|r|^{K})|g_{r}|$. Note that
if $\sum_{j=0}^{K+2}
\sup_{\omega}|g^{(j)}(\omega)|<\infty$ then
$\|g\|_{K}<\infty$, where $g^{(j)}(\cdot)$ denotes
the $j$th derivative of  $g$. 

Suppose $f, g:[0,2\pi]\rightarrow \mathbb{R}$ are
bounded functions, that are strictly larger than zero and 
are symmetric about $\pi$. By using the classical factorisation results
in \cite{p:sze-21} and \cite{p:bax-62} we can write $f(\cdot) =
\sigma_{f}^{2}|\psi_{f}(\cdot)|^{2}= \sigma^{2}_{f}|\phi_{f}(\cdot)|^{-2}$,
where $\phi_{f}(\omega) = 1 -\sum_{j=1}^{\infty}\phi_{j}(f)e^{-ij\omega}$
and $\psi_{f}(\omega) =
1+\sum_{j=1}^{\infty}\psi_{j}(f)e^{-ij\omega}$, the terms $\sigma_g$, $\phi_g(\cdot)$, and $\psi_g(\cdot)$ are defined similarly.
We use these expansions in Sections \ref{sec:approx} and
\ref{sec:newlikelihood}, where we require
the following notation
\begin{eqnarray*}
 \rho_{n,K}(f) &=& \sum_{r=n+1}^{\infty}|r^{K}\phi_{r}(f)|, \\
A_{K}(f,g) &=& 2\sigma_{g}^{-2}\|\psi_{f}\|_{0}\|\phi_{g}\|_0^{2}\|\phi_{f}\|_{K}, \\
\textrm{and} \qquad  C_{f,K} &=& \frac{3-\varepsilon}{1-\varepsilon}\left\|\phi_f\right\|_{K}^{2}
\left\|\psi_f\right\|_{K}^{2}
\end{eqnarray*} 
for some $0<\varepsilon<1$. 

For postive sequences $\{a_n\}$ and $\{b_n\}$, we denote $a_n \sim b_n$ if there exist
$0<C_1 \leq C_2 <\infty$ such that $C_1 \leq a_n/b_n \leq C_2$ for all $n$.
Lastly, we denote $\R$ and $\I$ as 
the real and imaginary part of a complex variable respectively.


\subsection{Motivation}

\begin{figure}[h!]
\begin{center}
\includegraphics[scale = 0.4]{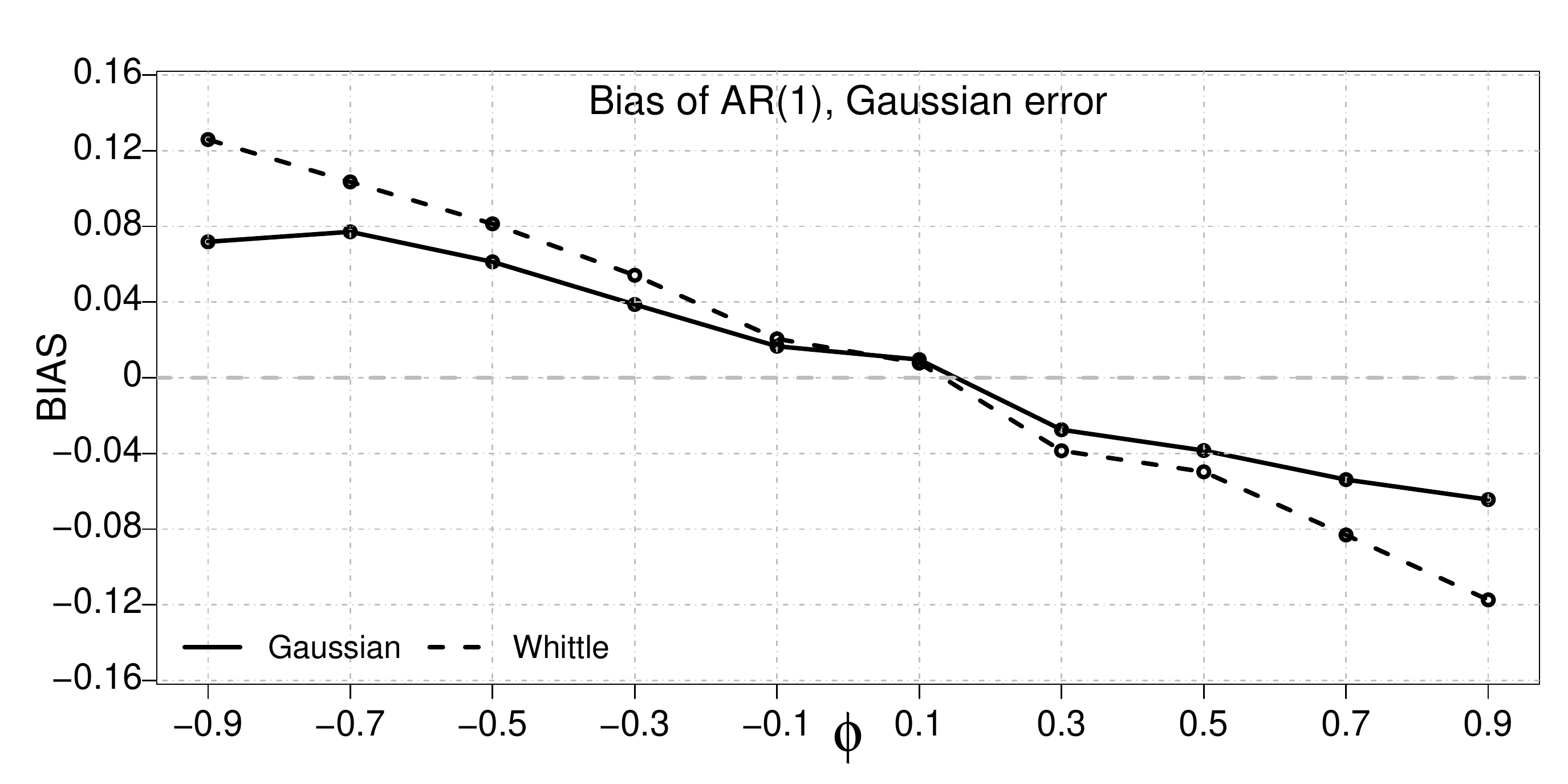}
\caption{{\it The model $X_{t} = \phi X_{t-1} + \varepsilon_{t}$ with independent
  standard normal errors is
  simulated. The bias of the estimator of $\phi$ based on sample size
  $n=20$ over   1000 replications.} \label{fig:AR} }
\end{center}
\end{figure}

In order to motivate our approach, we first study the difference in
the bias of the AR$(1)$ parameter estimator using both the Gaussian and
Whittle likelihood. In Figure \ref{fig:AR}, we plot the bias in the
estimator of $\phi$ in the AR$(1)$ model  $X_{t} = \phi
X_{t-1}+\varepsilon_{t}$ for different values of $\phi$ (based on
sample size $n=20$). We observe that the difference between the bias
of the two estimators increases as $|\phi|$
approaches one. Further, the Gaussian likelihood clearly has a
smaller bias than the Whittle likelihood (which is more pronounced
when $|\phi|$ is close to one). 
Let $\{X_{t}\}_{t=1}^{n}$ denote the observed time series. 
Straightforward calculations (based on expressions for
$\Gamma_{n}(f_\phi)^{-1}$ and $C_{n}(f_{\phi}^{-1})$) 
show that the difference between the
Gaussian and Whittle likelihoods for an AR$(1)$ model is
\begin{equation}\label{eq:diffAR1}
\mathcal{L}_{n}(\phi) - K_{n}(\phi) = 
n^{-1}\left[2\phi X_{1}X_{n} - \phi^2 (X_{1}^2+X_{n}^2)\right]
\end{equation}
Thus we observe that the closer $|\phi|$ is to one, the larger the expected
difference between the likelihoods. Using (\ref{eq:diffAR1}) and the
Bartlett correction (see \cite{p:bar-53} and \cite{p:cox-68}, it is
possible to obtain an asymptotic expression for the difference in the
biases (see also Appendix \ref{sec:AR1bias})
Generalisations of this result to higher order AR$(p)$
models may also be possible using the analytic expression for the inverse of the Toeplitz
matrix corresponding to an AR$(p)$ model derived in \cite{p:sid-58}
and \cite{p:gal-74}.

However, for more general models, such as the MA$(q)$ or ARMA$(p,q)$
models, using brute force calculations
for deriving the difference $\mathcal{L}_{n}(\theta) - K_{n}(\theta)$
and its derivatives is extremely difficult. Furthermore,  
such results do not 
offer any insight on how the Gaussian and Whittle likelihood are
related, nor what is ``lost'' when going from the Gaussian
likelihood to the Whittle likelihood. In the remainder of this section,
we derive an exact expression for the Gaussian likelihood in the frequency
domain. Using these derivations, we obtain a simple expression for the
difference between the Whittle and Gaussian likelihood for AR$(p)$
models. In subsequent sections, we obtain approximations for this
difference for general time series models.

\subsection{The biorthogonal transform to the discrete Fourier transform}\label{sec:DFT}


In order to obtain an exact bound, we start with the Whittle
likelihood and recall that 
the DFT of the time
series plays a fundamental role in its formulation. With this in mind, our approach is based on
deriving the transformation $\{Z_{k,n}\}_{k=1}^{n}\subset
\spa(\underline{X}_{n})$ (where  $\spa(\Xunder_{n})$ 
denotes the linear space over a complex field spanned by
$\Xunder_{n}=\{X_{t}\}_{t=1}^{n}$), which is biorthogonal to
$\{J_{n}(\omega_{k,n})\}_{k=1}^{n}$. That is, we derive a transformation
$\{Z_{k,n}\}_{k=1}^{n}$ which when coupled with $\{J_{n}(\omega_{k,n})\}_{k=1}^{n}$
satisfies the following condition
\begin{equation*}
\cov_{f}\left(Z_{k_{1},n},J_{n}(\omega_{k_2,n}) \right) = f(\omega_{k_1}) \delta_{k_{1},k_{2}}
\end{equation*}
where $\delta_{k_1,k_2}=1$ if $k_1=k_2$ (and zero otherwise).
Since $\underline{Z}_{n}^{\prime} = (Z_{1,n},\ldots,Z_{n,n})\in
\spa(\Xunder_{n})^{n}$, there exists an $n\times n$ complex matrix
$U_{n}$, such that $\underline{Z}_{n}=U_{n}\underline{X}_{n}$. Since 
$(J_{n}(\omega_{k,1}),\ldots,J_{n}(\omega_{n,n}))^{\prime} =
F_{n}\underline{X}_{n}$,
the biorthogonality of
$U_{n}\underline{X}_n$ and $F_{n}\underline{X}_{n}$ gives 
$\cov_{f}\left( U_{n}\underline{X}_{n},F_{n}\underline{X}_{n}\right)
= \Delta_{n}(f)$.
The benefit of biorthogonality is that it leads to the following
simple identity on the inverse of the variance matrix.

\begin{lemma}
Suppose that $U_{n}$ and $V_{n}$ are invertible
matrices which are biorthogonal 
with respect to the variance matrix $\var(\underline{X}_n)$. That is
$\cov(U_{n}\underline{X}_{n}, V_{n}\underline{X}_{n}) =
\Delta_{n}$, where $\Delta_{n}$ is a diagonal matrix. Then 
 \begin{equation}
\label{eq:UXUX}
\var(\underline{X}_{n})^{-1} = V_{n}^{*}\Delta_{n}^{-1}U_{n}.
 \end{equation}
\end{lemma}
\noindent PROOF. It follows immediately from $\cov(U_{n}\underline{X}_{n}, V_{n}\underline{X}_{n}) 
= U_{n}\var(\underline{X}_{n})V_{n}^{*} = \Delta_{n} $ and
$\var(\underline{X}_{n}) =
U_{n}^{-1}\Delta_{n}(V_{n}^{*})^{-1}$. \hfill $\Box$

\vspace{1em}

\noindent To understand how $U_{n}\underline{X}_{n}$ is related to
$F_{n}\underline{X}_{n}$ we rewrite $U_{n} = F_{n}+ D_{n}(f_{})$.
We show in the following theorem that $D_{n}(f)$ has a specific form with an
intuitive interpretation. In
order to develop these ideas, we use methods from linear
prediction. In particular, we define the best linear predictor of 
$X_{\tau}$ for $\tau \leq 0$ and $\tau  >n$ given
$\{X_{t}\}_{t=1}^{n}$ as 
\begin{equation}
\label{eq:Xtaun}
\widehat{X}_{\tau,n} = \sum_{t=1}^{n}\phi_{t,n}(\tau;f)X_{t},
\end{equation}
where $\{\phi_{t,n}(\tau;f)\}_{t=1}^{n}$ are the coefficients which minimize the $L_{2}$-distance
$\Ex_{f}[X_{\tau} - \sum_{t=1}^{n}\phi_{t,n}(\tau;f)X_{t}]^{2}$.
Using this notation we obtain the following theorem.

\begin{theorem}[The biorthogonal transform]\label{theorem:bio}
Let $\{X_{t}\}$ be a second order stationary, zero mean time series with spectral density
$f_{}$ which is bounded
away from zero and whose autocovariance satisfies 
$\sum_{r\in \mathbb{Z}}|rc_{f}(r)|<\infty$. Let
$\widehat{X}_{\tau,n}$ denote the best linear predictors of $X_{\tau}$
as defined in (\ref{eq:Xtaun}) and $\{\phi_{t,n}(\tau;f)\}_{t=1}^{n}$ the
corresponding coefficients. Then 
\begin{equation}
\label{eq:UX}
\cov_{f}\big( (F_{n} +
  D_{n}(f_{}))\underline{X}_{n},F_{n}\underline{X}_{n}\big) = \Delta_{n}(f_{}),
\end{equation} 
where $D_{n}(f_{})$ has entries
\begin{equation}
\label{eq:Dftheta}
D_{n}(f_{})_{k,t}=
     n^{-1/2}\sum_{\tau\leq 0}\left(\phi_{t,n}(\tau;f)e^{i\tau\omega_{k,n}}
+\phi_{n+1-t,n}(\tau;f)e^{-i(\tau-1)\omega_{k,n}}\right),
\end{equation} for $1\leq k,t\leq n$.
And, entrywise $1\leq k_1,k_2\leq n$, we have  
\begin{equation}
\label{eq:DFTPred}
\cov_{f}\left(\widetilde{J}_{n}(\omega_{k_1,n};f_{}), J_{n}(\omega_{k_2,n})\right)
  = f_{}(\omega_{k_1,n})\delta_{k_1,k_2}
\end{equation}
where 
$\widetilde{J}_{n} (\omega;f_{}) = J_{n} (\omega)+
\widehat{J}_{n} (\omega;f_{})$ and
\begin{equation} 
\label{eq:28}
\widehat{J}_{n} (\omega;f_{}) = n^{-1/2}\sum_{\tau \leq 0}
  \widehat{X}_{\tau, n} e^{i \tau \omega} +
n^{-1/2}\sum_{\tau > n} \widehat{X}_{\tau, n} e^{i \tau \omega}.
\end{equation}
\end{theorem}
\noindent PROOF. See Appendix \ref{sec:proofbio} (note that identity
(\ref{eq:DFTPred}) can be directly verified using results on best
linear predictors). \hfill $\Box$

\vspace{2mm}

\begin{corollary}[Inverse Toeplitz identity]\label{corollary:inverse}
Let $\Gamma_{n}(f)$ denote an $n\times n$ Toeplitz matrix generated by the
spectral density $f$. Then equations (\ref{eq:UXUX}) and (\ref{eq:UX})  yield the following identity
\begin{eqnarray}
\label{eq:UX1}
\Gamma_{n}(f)^{-1} = F_{n}^{*}\Delta_{n}(f^{-1})(F_{n}+D_{n}(f)),
\end{eqnarray}
where $D_{n}(f)$ is defined in (\ref{eq:Dftheta}). Observe that two spectral density functions $f_{1}(\omega)$ and
$f_{2}(\omega)$ with the same autocovariance up to lag $(n-1)$, $\{c(r)\}_{r=0}^{n-1}$,
 can give rise to two different representations 
\begin{equation*}
\Gamma_{n}(f_{1})^{-1}
= F_{n}^{*}\Delta_{n}(f_1^{-1})(F_{n}+D_{n}(f_1)) = F_{n}^{*}\Delta_{n}(f_2^{-1})(F_{n}+D_{n}(f_2)) = \Gamma_{n}(f_{2})^{-1}.
\end{equation*}
\end{corollary}

\begin{figure}[h!]
\begin{center}
\includegraphics[scale = 0.35]{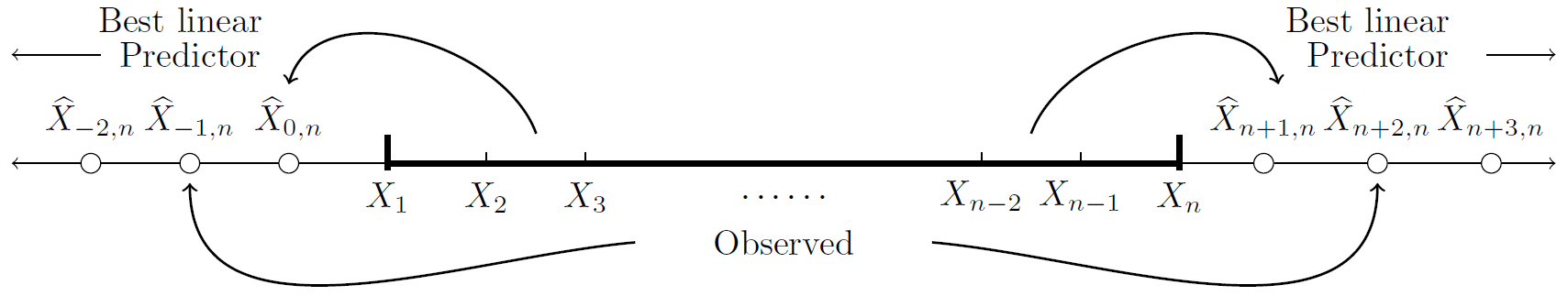}
\caption{ \textit{$\widetilde{J}_{n}(\omega;f)$ is the Fourier transform
  over both the observed time series and its predictors outside this domain. \label{fig:1}} }
\end{center}
\end{figure}

\noindent  What we observe is that the biorthogonal transformation
$(F_{n}+D_{n}(f))\underline{X}_{n}$
extends the domain of observation by predicting
outside the boundary. A visualisation of the observations and the
predictors that are involved in the construction of
$\widetilde{J}_{n}(\omega;f)$ is given in Figure
\ref{fig:1}.

It is quite surprising that only a small modification of the
regular DFT leads to its biorthogonal transformation. Furthermore, the
contribution of the additional DFT term is
$\widehat{J}_{n}(\omega_{k,n};f) = O_{p}(n^{-1/2})$. This is
why the regular DFT satisfies the well known ``near'' orthogonal property
\begin{equation*}
\cov_{f}(J_{n}(\omega_{k_1,n}),J_{n}(\omega_{k_2,n})) =
  f_{}(\omega_{k_1})\delta_{k_1,k_2} + O(n^{-1}),
\end{equation*}
see \cite{p:lah-03} and \cite{b:bri-01}.
For future reference, we will use the following definitions. 

\begin{defin} We refer to $\widehat{J}_{n}
(\omega;f)$ as the predictive DFT (as it is the Fourier
transform of all the linear predictors), noting that basic algebra
yields the expression
 \begin{equation}
\label{eq:DFTgen}
\widehat{J}_{n}(\omega_{};f_{}) =
  n^{-1/2}\sum_{t=1}^{n}X_{t}\sum_{\tau\leq 0}(\phi_{t,n}(\tau;f_{})e^{i\tau\omega_{}}
+e^{in\omega}\phi_{n+1-t,n}(\tau;f_{})e^{-i(\tau-1)\omega_{}}). 
\end{equation}
Note that when $\omega = \omega_{k,n}$, 
the term $e^{in\omega}$ in (\ref{eq:DFTgen}) vanishes.
Further, we refer to $\widetilde{J}_{n}
(\omega;f)$ as the complete DFT (as it contains the classical
DFT of the time series together with the predictive DFT). Note that
both $\widetilde{J}_{n}(\omega;f)$ and $\widehat{J}_{n}
(\omega;f)$ are functions of $f$ since they involve the
spectral density $f(\cdot)$, unlike the regular DFT which is model-free. 
\end{defin}

\begin{example}[The AR$(1)$ process] \label{example:AR1}
Suppose that $X_{t}$ has an AR$(1)$ representation $X_{t} = \phi
X_{t-1}+\varepsilon_{t}$ ($|\phi|<1$). Then the best linear
predictors are simply a function of the observations at the two endpoints. 
That is for $\tau\leq 0$, $\widehat{X}_{\tau,n} =
\phi^{|\tau|+1}X_{1}$ and for $\tau>n$ $\widehat{X}_{\tau,n} =
\phi^{\tau-n}X_{n}$. 
 An illustration is given in Figure \ref{fig:AR1example}.

\begin{figure}[!ht] 
\begin{center}
\includegraphics[scale = 0.4]{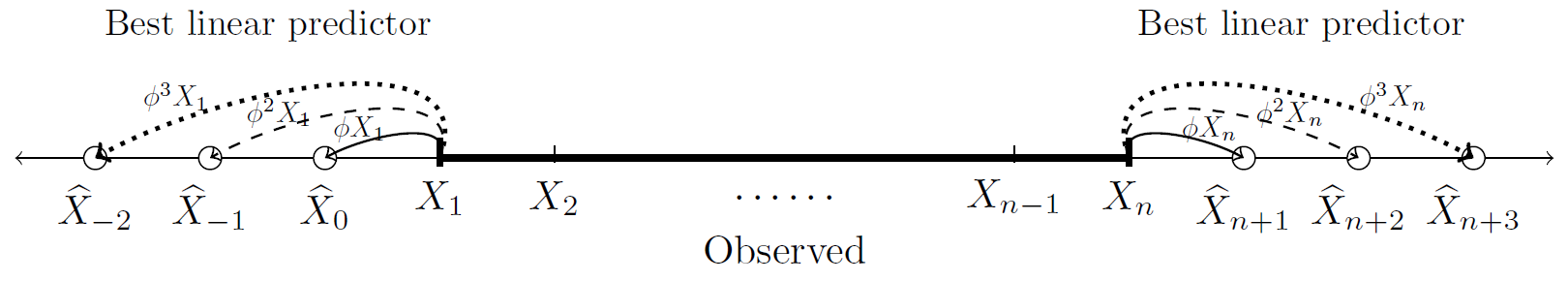}
\end{center}
\caption{{\it The past and future best linear predictors based on a AR(1) model.}\label{fig:AR1example}}
\end{figure}
\noindent Then the predictive DFT for the AR$(1)$ model is
\begin{equation*}
\widehat{J}_{n}(\omega;f_{\phi}) = \frac{\phi}{\sqrt{n}}\left(
\frac{1}{\phi(\omega)} X_{1}+
\frac{e^{i(n+1)\omega}}{ \overline{\phi(\omega)}}X_{n}\right) \quad \textrm{where}\quad\phi(\omega) = 1-\phi e^{-i\omega}.
\end{equation*} 
In other words, a small adjustment of the boundary leads to
$\widetilde{J}_{n}(\omega;f_{\phi})\overline{J_{n}(\omega)}$ being an
unbiased estimator of $f(\omega) = \sigma^{2}|\phi(\omega)|^{-2}$.
\end{example}

\begin{remark}
Biorthogonality of random variables is rarely used in
  statistics. An interesting exception is
\cite{p:kas-09}. They apply the notion of biorthogonality to problems
in prediction. In particular they consider 
the biorthogonal transform of 
$\underline{X}_{n}$, which is the random vector $\underline{\widetilde{X}}_{n} =
\Gamma_{n}(f)^{-1}\underline{X}_{n}$ (since
$\cov_{f}(\underline{\widetilde{X}}_{n},\underline{X}_{n})=I_{n}$). They obtain an
expression for the entries of $\underline{\widetilde{X}}_{n}$ in terms
of the Cholesky
decomposition of $\Gamma_{n}(f)^{-1}$. However, there is an 
interesting duality between $\underline{\widetilde{X}}_{n}$ and $\underline{\widetilde{J}}_{n} =
(\widetilde{J}_{n}(\omega_{1,n};f),\ldots,\widetilde{J}_{n}(\omega_{n,n};f))^{\prime}$.
In particular, applying identity (\ref{eq:UX1}) to 
the DFT of $\underline{\widetilde{X}}_{n}$ gives
\begin{equation*}
F_{n}\underline{\widetilde{X}}_{n}= F_{n}\Gamma_{n}(f)^{-1}\underline{X}_{n}=
 \Delta_{n}(f^{-1})\underline{\widetilde{J}}_{n}.
\end{equation*}
This shows that the DFT of the biorthogonal transform of
$\underline{X}_{n}$ is the \emph{standardized} complete DFT. Conversely, the
inverse DFT of the standardized complete DFT gives the biorthogonal
transform to the original time series, where the entries of
$\underline{\widetilde{X}}_{n}$ are 
\begin{equation*}
\widetilde{X}_{j,n} =
  \frac{1}{\sqrt{n}}\sum_{k=1}^{n}\frac{\widetilde{J}_{n}(\omega_{k,n};f)}{f_{}(\omega_{k,n})}
e^{-ij\omega_{k,n}}.
\end{equation*}
\end{remark}

\begin{remark}[Connection to the orthogonal increment process]
Suppose that $Z(\omega)$ is the orthogonal increment process
associated with the stationary time series $\{X_{t}\}$ and $f$
the corresponding spectral density. If $\{X_{t}\}$ is a Gaussian
time series, then we have 
\begin{equation*}
\widehat{X}_{\tau,n}= \Ex\left[X_{\tau}|\underline{X}_{n}\right] =
\frac{1}{2\pi}\int_{0}^{2\pi} e^{-i\omega\tau}\Ex [Z(d\omega)|\underline{X}_{n}]
  =   \frac{\sqrt{n}}{2\pi}\int_{0}^{2\pi}e^{-i\omega\tau}\widetilde{J}_{n}(\omega;f) d\omega.
\end{equation*} 
\end{remark}

\subsection{The Gaussian likelihood in the frequency domain}\label{sec:gauss}
In the following theorem, we exploit the biorthogonality between the
regular DFT and the complete DFT to yield an \emph{exact} ``frequency domain'' representation for
the Gaussian likelihood.  We use the notation
defined in Theorem \ref{theorem:bio}.

\begin{theorem}[A frequency domain representation of the Gaussian likelihood]\label{theorem:gauss}
Suppose the spectral density $f_{\theta}$ is bounded away from zero,
and the corresponding autocovariance is such that
$\sum_{r}|rc_{f_\theta}(r)|<\infty$. Let 
$\mathcal{L}_{n}(\theta)$ 
and $K_{n}(\theta)$
be defined as in
(\ref{eq:L1L2}). Then we have 
\begin{equation}
\label{eq:GL}
\mathcal{L}_{n}(\theta) =  \frac{1}{n}\Xunder_{n}^{\prime}\Gamma_{n}(f_\theta)^{-1}\Xunder_{n}
= \frac{1}{n}\sum_{k=1}^{n}\frac{\widetilde{J}_{n}(\omega_{k,n};f_{\theta})\overline{J_{n}(\omega_{k,n})}}{f_{\theta}(\omega_{k,n})}.
\end{equation}
Further
\begin{equation}
\label{eq:GammaCn}
\Gamma_{n}(f_\theta)^{-1} - C_{n}(f_\theta^{-1}) = F_{n}^{*}\Delta_{n}(f_\theta^{-1})D_{n}(f_\theta).
\end{equation}
This yields the difference between the Gaussian and Whittle likelihood
\begin{eqnarray}
\mathcal{L}_{n}(\theta) - K_{n}(\theta) &=& 
n^{-1}\Xunder_{n}^{\prime}\left[\Gamma_{n}(f_\theta)^{-1}-C_{n}(f_{\theta}^{-1})\right]\Xunder_{n} \nonumber \\
 &=&
\frac{1}{n}\sum_{k=1}^{n}\frac{\widehat{J}_{n}(\omega_{k,n};f_{\theta})
\overline{J_{n}(\omega_{k,n})}}{f_{\theta}(\omega_{k,n})}.
\label{eq:GLWL}
\end{eqnarray}
\end{theorem}
\noindent PROOF. (\ref{eq:GammaCn}) follows immediately from Corollary \ref{corollary:inverse}.
Next, we note that  $F_{n}\underline{X}_{n}=
\underline{J}_{n}$ and $(F_{n}+D_{n}(f_{\theta}))\underline{X}_{n}=
\underline{\widetilde{J}}_{n}$,
thus we immediately obtain equation (\ref{eq:GL}), and since
$\widetilde{J}_{n}(\omega_{k,n};f_{\theta}) = J_{n}(\omega_{k,n})+\widehat{J}_{n}(\omega_{k,n};f_{\theta})$,
it proves (\ref{eq:GLWL}). 
\hfill $\Box$

\vspace{2mm}

\noindent From the above theorem, we observe that the Gaussian
likelihood is the Whittle likelihood plus an additional ``correction''
\begin{equation*}
\mathcal{L}_{n}(\theta)  =
 \underbrace{\frac{1}{n}\sum_{k=1}^{n}
\frac{|J_{n}(\omega_{k,n})|^{2}}{f_{\theta}(\omega_{k,n})}}_{=K_{n}(\theta)} +
  \frac{1}{n}\sum_{k=1}^{n}
\frac{\widehat{J}_{n}(\omega_{k,n};f_{\theta})
\overline{J_{n}(\omega_{k,n})}}{f_{\theta}(\omega_{k,n})}.
\end{equation*}
To summarize, 
the Gaussian likelihood compensates for the well known boundary effect
in the Whittle likelihood, by predicting outside the domain of observation. 
The Whittle likelihood estimator selects
the spectral density $f_\theta$ which best fits the periodogram. On
the other hand, since $\Ex_{f_\theta}[\widetilde{J}_{n}(\omega_{k,n};f_\theta) \overline{J_{n}(\omega_{k,n})}]=f_{\theta}(\omega_{k,n})$,
the Gaussian likelihood estimator selects the spectral density which
best fits $\widetilde{J}_{n}(\omega_{k,n};f_\theta) \overline{J_{n}(\omega_{k,n})}$
by simultaneously predicting and fitting. 
Therefore, the ``larger'' the level of ``persistence'' in the time series, the greater the predictive DFT 
$\widehat{J}_{n}(\omega_{k,n};f_\theta)$, and subsequently the larger the
approximation error between the two likelihoods. This fits with
the insights of \cite{p:dah-88}, who shows that the more peaked the spectral
density the greater the leakage effect in the Whittle likelihood,
leading to a large finite sample bias.

In the remainder of this section and the subsequent section, we study 
the difference between the two likelihoods and corresponding
matrices. This will allow us to develop methods that better capture the Gaussian likelihood within the
frequency domain. By using Theorem \ref{theorem:gauss}, we have 
\begin{equation*}
\mathcal{L}_{n}(\theta) - K_{n}(\theta) = \frac{1}{n}\sum_{k=1}^{n}\frac{\widehat{J}_{n}(\omega_{k,n};f_{\theta})
\overline{J_{n}(\omega_{k,n})}}{f_{\theta}(\omega_{k,n})} = 
n^{-1}\underline{X}_{n}^{\prime}F_{n}^{*}\Delta_{n}(f_{\theta}^{-1})D_{n}(f_\theta)\underline{X}_{n},
\end{equation*}
where the entries of $F_{n}^{*}\Delta_{n}(f_{\theta}^{-1})D_{n}(f_\theta)$ are
\begin{eqnarray}
 &&(F_{n}^{*}\Delta_{n}(f_{\theta}^{-1})D_{n}(f_\theta) )_{s,t}  \nonumber \\
&& \qquad = \sum_{\tau \leq 0}\left[\phi_{t,n}(\tau;f_\theta)
    G_{1,n}(s,\tau;f_\theta)+\phi_{n+1-t,n}(\tau;f_\theta) G_{2,n}(s,\tau;f_\theta)\right]
\label{eq:DDeltaD}
\end{eqnarray}
with
\begin{eqnarray*}
 G_{1,n}(s,\tau;f_\theta) &=& \frac{1}{n}\sum_{k=1}^{n}\frac{1}{f_{\theta}(\omega_{k,n})}
e^{i(\tau-s)\omega_{k,n}}  = \sum_{a\in \mathbb{Z}}K_{f_\theta^{-1}}(\tau - s +an) \\
G_{2,n}(s,\tau;f_\theta) &=&  \frac{1}{n}\sum_{k=1}^{n}\frac{1}{f_{\theta}(\omega_{k,n})}
e^{-i(\tau+s-1)\omega_{k,n}} = \sum_{a\in \mathbb{Z}}K_{f_\theta^{-1}}(\tau + s -1 +an)
\end{eqnarray*}
and $K_{f_\theta^{-1}}(r) = \int_{0}^{2\pi}
f_{\theta}(\omega)^{-1}e^{ir\omega}d\omega$. 
We observe that for $1<<t<<n$, $\phi_{t,n}(\tau;f_\theta)$ and
$\phi_{n+1-t,n}(\tau;f_\theta)$ will be ``small'' as compared with $t$
close to one or $n$. The same is true for $G_{1,n}(s,\tau;f_\theta)$
and $G_{2,n}(s,\tau;f_\theta)$ when $1<<s<<n$. 
Thus the entries of  $F_{n}^{*}\Delta_{n}(f_{\theta}^{-1})D_{n}(f_\theta)$
will be ``small'' far from the four corners of the matrix. In contrast,
the entries of $F_{n}^{*}\Delta_{n}(f_{\theta}^{-1})D_{n}(f_\theta)$
will be largest at the four corners at the matrix.
This can be clearly seen in the following theorem, where we consider
the special case of AR$(p)$ models. 
We showed in Example \ref{example:AR1} that for AR$(1)$ processes, the
predictive DFT has a simple form. In the following theorem, we obtain
an analogous result for AR$(p)$ models (where $p \leq n$). 

\begin{theorem}[Finite order autoregressive models]\label{lemma:ARp1}
Suppose that $f_{\theta}(\omega) =\sigma^{2}|\phi_{p}(\omega)|^{-2}$
where $\phi_{p}(\omega) = 
1-\sum_{u=1}^{p}\phi_{u}e^{-iu\omega}$ (the roots of the corresponding
characteristic polynomial lie outside the unit circle) and $p\leq n$. 
The predictive DFT has
the analytic form 
\begin{eqnarray}
&& \widehat{J}_{n}(\omega;f_\theta) = \nonumber \\
&& \qquad 
\frac{n^{-1/2}}{\phi_{p}(\omega)} \sum_{\ell=1}^{p}X_{\ell}\sum_{s=0}^{p-\ell}\phi_{\ell+s}e^{-is\omega}+
e^{in\omega} \frac{n^{-1/2}}{ \overline{\phi_{p}(\omega)}} \sum_{\ell=1}^{p}X_{n+1-\ell}\sum_{s=0}^{p-\ell}
\phi_{\ell+s}e^{i(s+1)\omega}.
\label{eq:JAR}
\end{eqnarray}
If $p\leq n/2$, then  
$D_{n}(f_{\theta})$ is a rank $2p$ matrix where 
\begin{eqnarray}
&& D_{n}(f_{\theta}) \nonumber \\
&& \qquad  = n^{-1/2} \left( 
\begin{smallmatrix}
\phi_{1,p}(\omega_{1,n}) & \ldots & \phi_{p,p}(\omega_{1,n})
& 0 & \ldots & 0 & e^{i\omega_{1,n}} \overline{\phi_{p,p}(\omega_{1,n})} & \ldots &
e^{i\omega_{1,n}} \overline{\phi_{1,p}(\omega_{1,n})} \\
\phi_{1,p}(\omega_{2,n}) & \ldots &
\phi_{p,p}(\omega_{2,n})
& 0 & \ldots & 0 &  e^{i\omega_{2,n}} \overline{\phi_{p,p}(\omega_{2,n})} & \ldots &
 e^{i\omega_{2,n}} \overline{\phi_{1,p}(\omega_{2,n})} \\
\svdots & \sddots & \svdots & \svdots & \sddots & \svdots &  \svdots & \sddots & \svdots \\
\phi_{1,p}(\omega_{n,n}) & \ldots &
\phi_{p,p}(\omega_{n,n})
& 0 & \ldots & 0 &  e^{i\omega_{n,n}} \overline{\phi_{p,p}(\omega_{n,n})} & \ldots &
 e^{i\omega_{n,n}} \overline{\phi_{1,p}(\omega_{n,n})}\\
\end{smallmatrix} \right)
\label{eq:DDD}
\end{eqnarray}
and $\phi_{j,p}(\omega) = \phi_{p}(\omega)^{-1}
\sum_{s=0}^{p-j}\phi_{j+s}e^{-is\omega}$.  
Note, if $n/2< p\leq n$, then the entries of $D_{n}(f_\theta)$ will overlap.
Let $\widetilde{\phi}_{0}
= 1$ and for $1\leq s \leq p$,  
$\widetilde{\phi}_{s} = -\phi_{s}$ (zero otherwise), then if $1\leq p\leq
n/2$ we have 
\begin{eqnarray}
&& \left(\Gamma_{n}(f_\theta)^{-1} - C_{n}(f_\theta^{-1})\right)_{s,t} =
 (F_{n}^{*}\Delta_{n}(f_{\theta}^{-1})D_{n}(f_\theta) )_{s,t} \nonumber \\
&& \qquad = \left\{ \begin{array}{cl}
\sigma^{-2} \sum_{\ell=0}^{p-t}\phi_{\ell+t}\widetilde{\phi}_{(\ell+s) \bmod n}
  & 1\leq t \leq p \\
\sigma^{-2} \sum_{\ell=1}^{p-(n-t)}
  \phi_{\ell + (n-t)}\widetilde{\phi}_{(\ell -s) \bmod n} &   n-p+1 \leq t \leq n\\
0 &  otherwise
\end{array}.
\right.
\label{eq:DDDmatrix}
\end{eqnarray}
\end{theorem}
\noindent PROOF. In Appendix \ref{sec:proofbio}. \hfill $\Box$

\vspace{2mm}

Theorem \ref{lemma:ARp1} shows that for AR$(p)$ models, the
predictive DFT only involves the $p$ observations on each side of the
observational boundary $X_{1},\ldots,X_{p}$ and
$X_{n-p+1},\ldots,X_{n}$, where the coefficients in the prediction are a
linear combination of the AR parameters (excluding the denominator $\phi_{p}(\omega)$).
The well known result (see \cite{p:sid-58} and
  \cite{p:sha-75}, equation (10)) that 
$F_{n}^{*}\Delta_{n}(f_{\theta}^{-1})D_{n}(f_\theta)$ is non-zero only
at the $(p\times p)$ submatrices located in the four corners of
$F_{n}^{*}\Delta_{n}(f_{\theta}^{-1})D_{n}(f_\theta)$ follows from
equation (\ref{eq:DDDmatrix}).

By using (\ref{eq:JAR}) we obtain an analytic expression 
for the Gaussian likelihood of the AR$(p)$ model in terms of the 
autoregressive coefficients. In particular, the 
Gaussian likelihood (written in the frequency domain) corresponding to
the AR$(p)$ model $X_{t} = \sum_{j=1}^{p}\phi_{j}X_{t-j}+\varepsilon_{t}$  is 
\begin{eqnarray}
&& \mathcal{L}_{n}(\phi) =\frac{\sigma^{-2}}{n}\sum_{k=1}^{n}|J_{n}(\omega_{k,n})|^{2}|\phi_p(\omega_{k,n})|^{2} \nonumber \\
&& \quad  +  \frac{\sigma^{-2}}{n}\sum_{\ell=1}^{p}X_{\ell}\sum_{s=0}^{p-\ell}\phi_{\ell+s}\left(X_{(-s)\bmod
    n}-\sum_{j=1}^{p}\phi_{j}X_{(j-s)\bmod n}\right)\nonumber\\
&& \quad +
  \frac{\sigma^{-2}}{n}\sum_{\ell=1}^{p}X_{n+1-\ell}\sum_{s=0}^{p-\ell}\phi_{\ell+s}\left(X_{(s+1)\bmod
  n}-\sum_{j=1}^{p}\phi_{j}X_{(s+1-j)\bmod n}\right),
\label{eq:LKdiff}
\end{eqnarray} 
where $\phi = (\phi_1, ..., \phi_{p})^{\prime}$ and $\phi_p(\omega) = 1-\sum_{j=1}^{p}\phi_{j}e^{-ij\omega}$. A
proof of the above identity can be found in Appendix \ref{sec:proofbio}.
Equation (\ref{eq:LKdiff}) offers a simple
representation of the Gaussian likelihood in terms of a
Whittle likelihood plus an additional term in terms of the AR$(p)$ coefficients.



\section{Frequency domain approximations of the Gaussian likelihood}\label{sec:approx}

In Theorem \ref{theorem:gauss} we rewrote the Gaussian likelihood 
within the frequency domain. This allowed us to obtain an
expression for the difference between the Gaussian
and Whittle likelihoods for AR$(p)$ models (see (\ref{eq:LKdiff})). 
This is possible  because the predictive DFT $\widehat{J}_{n}(\cdot;f_\theta)$
has a simple analytic form.

It would be of interest to generalize this result to general time series models.
However, for infinite order autoregressive models, the predictions
across the boundary and the predictive DFT given in
(\ref{eq:DFTgen}) do not have a simple, analytic form. 
In Section \ref{sec:firstorder} we show that we can obtain an 
approximation of the predictive DFT in terms of the 
AR$(\infty)$ coefficients corresponding to $f_\theta$. In turn, this allows us
to obtain an approximation for $\Gamma_{n}(f_\theta)^{-1} -
C_{n}(f_\theta^{-1})$, which is 
analogous to equation (\ref{eq:DDDmatrix}) for AR$(p)$ models. 
Such a result proves to be very useful from both a theoretical and 
practical perspective. Theoretically, we use this result to show that the difference between the Whittle and Gaussian
likelihood is of order $O(n^{-1})$.  Furthermore, in Section \ref{sec:higherexpan} we show that the
approximation described in Section \ref{sec:firstorder} is the
first order term of a polynomial-type series expansion of
$\Gamma_{n}(f_\theta)^{-1}$ in terms of the AR$(\infty)$ parameters.
From a practical perspective, the approximations are
used in Section \ref{sec:newlikelihood} to
motivate alternative quasi-likelihoods defined within the frequency
domain.

First, we require the following set of
assumptions on the spectral density $f_{\theta}$.
\begin{assumption}\label{assum:A}
\begin{itemize}
\item[(i)] The spectral density $f$ is bounded away from zero. 
\item[(ii)] For some $K> 1$, the autocovariance function is such that
$\sum_{r\in \mathbb{Z}}|r^{K}c_{f}(r)|<\infty$.
\end{itemize}
\end{assumption}
Under the above assumptions, we can write $f(\omega) = \sigma^{2}|\psi(\omega;f)|^{2}  =
\sigma^{2}|\phi(\omega;f)|^{-2}$ where 
\begin{equation}
\label{eq:thetaphi}
\psi(\omega;f)=
1+\sum_{j=1}^{\infty}\psi_{j}(f)e^{-ij\omega} \textrm{ and } \phi(\omega;f) = 1 -
\sum_{j=1}^{\infty}\phi_{j}(f)e^{-ij\omega}. 
\end{equation}
Further, under Assumption \ref{assum:A} we have
$\sum_{r=1}^{\infty}|r^{K}\psi_{r}(f)|$ and
$\sum_{r=1}^{\infty}|r^{K}\phi_{r}(f)|$ are both finite (see
\cite{p:kre-11}). Thus if $f$ satisfies Assumption \ref{assum:A} with some $K> 1$,
 then  $\|\psi_f\|_K<\infty$ and
 $\|\phi_f\|_K<\infty$. 


\subsection{The first order approximation} \label{sec:firstorder}

In order to obtain a result analogous  to 
Theorem \ref{lemma:ARp1}, we 
replace $\phi_{s,n}(\tau;f_\theta)$ in $D_{n}(f_{\theta})$ with
$\phi_{s}(\tau;f_\theta)$ which are the coefficients of the
best linear predictor of $X_{\tau}$ (for $\tau \leq 0$) given
$\{X_{t}\}_{t=1}^{\infty}$ i.e.
$\widehat{X}_{\tau} = \sum_{t=1}^{\infty}\phi_{t}(\tau;f_\theta)X_{t}$.
This gives the matrix $D_{\infty,n}(f_\theta)$, where 
\begin{equation*}
(D_{\infty,n}(f_{\theta}))_{k,t}  =
     n^{-1/2}\sum_{\tau\leq0}\left(\phi_{t}(\tau;f_\theta)e^{i\tau\omega_{k,n}}
+\phi_{n+1-t}(\tau;f_\theta)e^{-i(\tau-1)\omega_{k,n}}\right).
\end{equation*}
It can be shown that for $1\leq k,t \leq n$,
\begin{equation}
\label{eq:widetildeDalt}
(D_{\infty,n}(f_{\theta}))_{k,t} = 
n^{-1/2} \frac{\phi_{t}^{\infty}(\omega_{k,n};f_\theta)}{\phi(\omega_{k,n};f_\theta)} +
     n^{-1/2}e^{i\omega_{k,n}}\frac{\overline{\phi_{n+1-t}^{\infty}(\omega_{k,n};f_\theta)}}{\overline{\phi(\omega_{k,n};f_\theta)}},
\end{equation}
where  $\phi_{t}^{\infty}(\omega_{};f_\theta) =
\sum_{s=0}^{\infty}\phi_{t+s}(f_\theta)e^{-is\omega_{}}$. 
The proof of the above identity can be found in 
Appendix \ref{sec:approxproofs}. 
Using the above we can show that
$(D_{\infty,n}(f_\theta)\underline{X}_{n})_{k}=
\widehat{J}_{\infty,n}(\omega_{k,n};f_\theta)$ where 
\begin{eqnarray}
&& \widehat{J}_{\infty,n}(\omega_{};f_\theta) \nonumber \\
&& \quad = 
\frac{n^{-1/2}}{\phi(\omega;f_\theta)} \sum_{t=1}^{n}X_{t}\phi_{t}^{\infty}(\omega;f_\theta)
+e^{i(n+1)\omega}
\frac{n^{-1/2}}{ \overline{\phi(\omega;f_\theta)}} \sum_{t=1}^{n}X_{n+1-t}\overline{\phi_{t}^{\infty}(\omega;f_\theta)}.
\label{eq:JARIn}
\end{eqnarray}
We show below that 
$\widehat{J}_{\infty,n}(\omega_{k,n};f_\theta)$ is an approximation of
$\widehat{J}_{n}(\omega_{k,n};f_\theta)$.


\begin{theorem}[An AR$(\infty)$ approximation for general processes]\label{theorem:approx}
Suppose $f$ satisfies Assumption \ref{assum:A}, $f_\theta$ is bounded
away from zero and 
$\|f_\theta\|_{0}<\infty$ (with $f_\theta(\omega) = \sigma^{2}_\theta |\phi_{\theta}(\omega)|^{-2}$).
Let $D_{n}(f)$, $D_{\infty,n}(f)$ and  $\widehat{J}_{\infty,n}(\omega_{k,n};f)$ be defined as
in (\ref{eq:Dftheta}) and (\ref{eq:widetildeDalt}) and
(\ref{eq:JARIn}) respectively. Then we have 
\begin{eqnarray}
&& \underline{X}_{n}^{\prime}F_{n}^{*}\Delta_{n}(f_{\theta}^{-1})\left(D_{n}(f)-
 D_{\infty,n}(f)\right)\underline{X}_{n} \nonumber \\
&& \quad = 
\sum_{k=1}^{n} 
\frac{\overline{J_{n}(\omega_{k,n})}}{f_{\theta}(\omega_{k,n})}
\big(
\widehat{J}_{n}(\omega_{k,n};f) -
\widehat{J}_{\infty,n}(\omega_{k,n};f) \big)
\label{eq:approx1AA}
\end{eqnarray}
and 
\begin{equation}
\label{eq:approx1A}
\left\|F_{n}^{*}\Delta_{n}(f_{\theta}^{-1})\left(D_{n}(f)-
 D_{\infty,n}(f)\right)
  \right\|_{1} \leq
  \frac{C_{f,0}\rho_{n,K}(f)}{n^{K-1}}A_{K}(f,f_\theta).
\end{equation} 
Further, if $\{X_{t}\}$ is a time
series where $\sup_{t} \|X_{t}\|_{\Ex,2q} = \|X\|_{\Ex,2q}<\infty$ (for some
$q>1$), then 
\begin{eqnarray}
&& n^{-1}\left\|\underline{X}_{n}^{\prime}F_{n}^{*}\Delta_{n}(f_{\theta}^{-1})\left(D_{n}(f)-
 D_{\infty,n}(f)\right)\underline{X}_{n}
  \right\|_{\Ex,q} \nonumber \\
&& \quad \leq
  \frac{C_{f,0}\rho_{n,K}(f)}{n^{K}}A_{K}(f,f_\theta)\|X\|_{\Ex,2q}^{2}.
\label{eq:approx1B}
\end{eqnarray}
\end{theorem}
\noindent PROOF. See Appendix \ref{sec:approxproofs}. \hfill $\Box$ 

\vspace{1em}
\noindent We mention that we state the above theorem in the general case that the spectral density $f$ is
used to construct the predictors $D_{n}(f)$. It does not necessarily have to be the same as
$f_\theta$. This is to allow generalisations of the Whittle and
Gaussian likelihoods, which we discuss in Section \ref{sec:newlikelihood}. 

Applying the above theorem to the Gaussian likelihood gives an
approximation which is analogous to (\ref{eq:LKdiff})
\begin{eqnarray}
\mathcal{L}_{n}(\theta) &=&  
K_{n}(\theta)
+\frac{1}{n}\sum_{k=1}^{n}\frac{\widehat{J}_{\infty,n}(\omega_{k,n};f_{\theta})\overline{J_{n}(\omega_{k,n})} }
{f_\theta(\omega_{k,n})} +O_{p}(n^{-K}) \nonumber\\
&=&  K_{n}(\theta) + \frac{1}{n}\sum_{s,t=1}^{n}X_{s}
X_{t}\frac{1}{n}\sum_{k=1}^{n} e^{-is\omega_{k,n}}\varphi_{t,n}(\omega_{k,n};f_\theta) + O_{p}(n^{-K}),
\label{eq:GaussApprox}
\end{eqnarray}
where 
$\varphi_{t,n}(\omega;f_\theta) = \sigma^{-2}
\left[ \overline{\phi(\omega;f_\theta)}\phi_{t}^{\infty}(\omega;f_\theta)    +
  e^{i\omega}\phi(\omega;f_\theta)\overline{\phi_{n+1-t}^{\infty}(\omega;f_\theta)}\right]$.
The above approximation shows that if the 
autocovariance function, corresponding to $f_\theta$ decays sufficiently fast (in the sense that $\sum_{r\in
  \mathbb{Z}}|r^{K}c_{f_\theta}(r)|<\infty$ for some $K> 1$). Then replacing the
finite predictions with the predictors using the infinite past (or future)
gives a  close approximation of the Gaussian likelihood. 

\begin{remark} Following from the above, the 
entrywise difference between the two matrices is approximately
\begin{equation*}
(\Gamma_{n}(f_\theta)^{-1}-C_{n}(f_{\theta}^{-1}))_{s,t} \approx
(F_{n}^{*}\Delta_{n}(f_{\theta}^{-1})D_{\infty, n}(f_\theta))_{s,t} =
\frac{1}{n}\sum_{k=1}^{n} e^{-is\omega_{k,n}}\varphi_{t,n}(\omega_{k,n};f_\theta),
\end{equation*} 
thus giving an analytic approximation to (\ref{eq:DDeltaD}).
\end{remark}

In the following theorem, we obtain a bound between the Gaussian and
Whittle likelihood. 
\begin{theorem}[The difference in the likelihoods]\label{theorem:Bound}
Suppose $f_{\theta}$ satisfies  Assumption \ref{assum:A}.
Let $D_{n}(f_{\theta})$ and $D_{\infty,n}(f_{\theta})$ be defined as
in (\ref{eq:Dftheta}) and (\ref{eq:widetildeDalt}) respectively.
Then we have
\begin{equation}
\label{eq:Bound1}
\left\|  F_{n}^{*}\Delta_{n}(f_{\theta}^{-1})D_{\infty,n}(f_{\theta})
  \right\|_{1} \leq A_{1}(f_\theta,f_\theta)
\end{equation}
and 
\begin{equation}
\label{eq:Bound2}
\left\|\Gamma_{n}(f_{\theta})^{-1} - C_{n}(f_{\theta}^{-1})
  \right\|_{1} 
\leq  \left( A_{1}(f_\theta,f_\theta) + \frac{C_{f_\theta,0}\rho_{n,K}(f_\theta)}{n^{K-1}}A_{K}(f_\theta,f_\theta)\right).
\end{equation}
Further, if $\{X_{t}\}$ is a time
series where $\sup_{t} \|X_{t}\|_{\Ex,2q} = \|X\|_{\Ex,2q}<\infty$ (for some
$q>1$), then  
\begin{eqnarray}
\|\mathcal{L}_{n}(\theta)-K_{n}(\theta)\|_{\Ex,q} 
 \leq n^{-1}\left( A_{1}(f_\theta,f_\theta) + \frac{C_{f_\theta,0}\rho_{n,K}(f_\theta)}{n^{K-1}}A_{K}(f_\theta,f_\theta)\right)\|X\|_{\Ex,2q}^{2}.
\label{eq:BD1}
\end{eqnarray}
\end{theorem}
\noindent PROOF. See Appendix \ref{sec:approxproofs}. \hfill $\Box$ 

\vspace{2mm}
\noindent

The above result shows that under the stated conditions
\begin{equation*}
n^{-1}\left\|\Gamma_{n}(f_{\theta})^{-1} - C_{n}(f_{\theta}^{-1})
  \right\|_{1} = O(n^{-1}),
\end{equation*}
and the difference between 
 the Whittle and Gaussian likelihoods is of
 order $O(n^{-1})$. 
We conclude this section by obtaining a higher order expansion of
$\Gamma_{n}(f_{\theta})^{-1}$.

\subsection{A series expansion}\label{sec:higherexpan}

Theorem \ref{theorem:approx} gives an approximation of the predictive
DFT $\widehat{J}_{n}(\omega;f)$ in terms of
$\widehat{J}_{\infty,n}(\omega;f)$, which is comprised of the
AR$(\infty)$ coefficients corresponding to $f$. In the following lemma
we show that it is possible to obtain a series expansion of
$\widehat{J}_{n}(\omega;f)$ and $\Gamma_{n}(f)^{-1}-C_{n}(f)^{-1}$ in
terms of the products of AR$(\infty)$ coefficients. The proof of the
results in this section hinge on  applying von Neumann's alternative
projection theorem to stationary time series. This technique was first
developed for time series in \cite{p:ino-06}. We make use of Theorem 2.5,
\cite{p:ino-06}, where an expression for the coefficients of the
finite predictors  $\phi_{t,n}(\tau)$ is given.

We define the function
$\zeta_{t,n}^{(1)}(\omega;f)=\phi_{t}^{\infty}(\omega)$ and for $s\geq
2$ 
\begin{eqnarray}
\label{eq:Dinftys}
\zeta_{t,n}^{(s)}(\omega;f) &=&
  \frac{1}{(2\pi)^{s-1}}\int_{[0,2\pi]^{s}}
\left(\prod_{a=1}^{s-1}\phi(\lambda_{a+1};f)^{-1}\Phi_{n}(\lambda_{a},\lambda_{a+1})\right)\times \nonumber\\
&& \bigg(\phi_{t}^{\infty}(\lambda_s;f)\delta_{ s\equiv 1 (\bmod 2)} 
+\phi_{n+1-t}^{\infty}(\lambda_s;f)\delta_{ s\equiv 0(\bmod 2)}\bigg)\delta_{\lambda_1=\omega}d\underline{\lambda}_s,
\end{eqnarray}
where $d\underline{\lambda}_s = d\lambda_1 \cdots d\lambda_s$ denotes the $s$-dimensional Lebesgue measure, 
\begin{eqnarray*}
\Phi_{n}(\lambda_{1},\lambda_{2}) 
=  \sum_{u=0}^{\infty} \phi_{n+1+u}^{\infty}(\lambda_1;f) e^{iu\lambda_2}
= \sum_{u=0}^{\infty} \sum_{s=0}^{\infty} \phi_{n+1+u+s}(f) e^{-is\lambda_1} e^{iu\lambda_2}
\end{eqnarray*}
and  $\delta$ denotes the indicator variable. 
In the following lemma, we show that $\zeta_{t,n}^{(s)}(\omega;f)$
plays the same role as $\phi_{t}^{\infty}(\omega;f)$ in the predictive
DFT approximation given in equation (\ref{eq:JARIn}). It will be used to
approximate $\widehat{J}_{n}(\omega;f)$ to a greater degree of accuracy.

\begin{theorem}\label{thm:higherorder}
Suppose $f$ satisfies Assumption \ref{assum:A}, where
$f(\omega) = \sigma^{2}|\phi(\omega;f)|^{-2}$.  
Let $D_{n}(f)$ and $\zeta_{j,n}^{(s)}(\omega;f)$ be defined as in 
(\ref{eq:Dftheta}) and 
(\ref{eq:Dinftys}) respectively. Define the $s$-order predictive DFT
\begin{eqnarray*}
\widehat{J}_{n}^{(s)}(\omega;f) =
  \frac{n^{-1/2}}{\phi(\omega;f)}\sum_{t=1}^{n}X_{t}\zeta_{t,n}^{(s)}(\omega;f)
  + e^{i(n+1)\omega}\frac{n^{-1/2}}{\overline{\phi(\omega;f)}}\sum_{t=1}^{n}X_{n+1-t}\overline{\zeta_{t,n}^{(s)}(\omega;f)}.
\end{eqnarray*}
Then 
\begin{eqnarray}
\label{eq:Japprox1}
\widehat{J}_{n}(\omega,f) = \sum_{s=1}^{\infty}\widehat{J}_{n}^{(s)}(\omega;f) 
\end{eqnarray}
and $D_{n}(f) = \sum_{s=1}^{\infty}D_{n}^{(s)}(f)$, where 
\begin{eqnarray}
\label{eq:Dapprox1}
(D_{n}^{(s)}(f))_{k,t} =  n^{-1/2} \frac{\zeta_{t,n}^{(s)}(\omega_{k,n};f)}{\phi(\omega_{k,n};f)} +
     n^{-1/2}e^{i\omega_{k,n}}\frac{\overline{\zeta_{n+1-t,n}^{(s)}(\omega_{k,n};f)}}{\overline{\phi(\omega_{k,n};f)}}.
\end{eqnarray}
Further, for a sufficiently large $n$ we have 
\begin{eqnarray}
\label{eq:Jhatapprox2}
\widehat{J}_{n}(\omega;f) =
                              \sum_{s=1}^{m}\widehat{J}_{n}^{(s)}(\omega;f) + O_{p}\left(\frac{1}{n^{m(K-1)+1/2}}\right).
\end{eqnarray}
\end{theorem}
PROOF. See Appendix \ref{sec:higherproof}. \hfill $\Box$ 

\vspace{1em}

\noindent In the case $s=1$,  it is straightforward to show that
\begin{eqnarray*}
\widehat{J}_{n}^{(1)}(\omega;f)=\widehat{J}_{\infty,n}(\omega;f)
\quad \text{and} \quad
D_{n}^{(1)}(f)=D_{\infty,n}(f).
\end{eqnarray*} 
Therefore, the first term in the
expansion of $\widehat{J}_{n}(\omega,f)$ and $D_{n}(f)$
is the AR$(\infty)$ approximation
$\widehat{J}_{\infty,n}(\omega;f)$ and $D_{\infty,n}(f)$ respectively. 
We mention, that it is simple to check that if $f$ corresponds to an AR$(p)$ spectral density for some $p \leq n$, then 
$\widehat{J}_{n}^{(s)}(\omega;f) = 0$ for all $s\geq 2$.
For general spectral densities, the higher order expansion gives a
higher order approximation of $\widehat{J}_{n}(\omega;f)$ and
$\Gamma_{n}(f)^{-1}$ in terms of products of the AR$(\infty)$ coefficients. Using the above result we have the expansions
\begin{eqnarray*}
\Gamma_{n}(f)^{-1} = C_{n}(f^{-1}) +\sum_{s=1}^{\infty} F_{n}^{*}\Delta_{n}(f^{-1})D_{n}^{(s)}(f)
\end{eqnarray*}
and 
\begin{eqnarray*}
\mathcal{L}_{n}(\theta) = K_{n}(\theta) +
\sum_{s=1}^{\infty}\frac{1}{n}\sum_{k=1}^{n}\frac{\widehat{J}^{(s)}_{n}(\omega_{k,n};f_{\theta})
\overline{J_{n}(\omega_{k,n})}}{f_{\theta}(\omega_{k,n})}.
\end{eqnarray*}
It is interesting to note that $\zeta_{t,n}^{(s)}(\omega;f)$ can be
evaluated recursively using
\begin{eqnarray}
&&\zeta_{t,n}^{(s+2)}(\omega;f) \nonumber \\
&&\quad =
  \frac{1}{(2\pi)^{2}}\int_{[0,2\pi]^{2}}\phi(y_{1};f)^{-1}\phi(y_{2};f)^{-1}\Phi_{n}(\omega,y_1)
\Phi_{n}(y_1,y_2) \zeta_{t,n}^{(s)}(y_2;f)dy_{1}dy_{2}.
\label{eq:recursionZeta}
\end{eqnarray}
In a similar vein, both the $s$-order predictive DFT
$\widehat{J}_{n}^{(s)}(\omega;f)$  and $D_{n}^{(s)}(f)$ can be
evaluated recursively using a recursion similar to the above (see
Appendix \ref{sec:higherproof} for the details). 

The above results show that it is possible to  obtain an analytic expression for
$\widehat{J}_{n}(\omega;f)$ and $\Gamma_{n}(f)^{-1}$
in terms of the products of the AR$(\infty)$ coefficients. This expression for the inverse of
a Toeplitz matrix may have
applications outside time series. However, from the perspective of
estimation,  the first order approximation $\widehat{J}_{n}^{(1)}(\omega;f)=\widehat{J}_{\infty,n}(\omega;f)$ is
sufficient. We discuss some applications in the
next section. 

\section{New frequency domain  quasi-likelihoods}\label{sec:newlikelihood}

In this section, we apply the approximations from the previous
section to define two new spectral divergence criteria.

To motivate the criteria, we recall from Theorem \ref{theorem:gauss} that the Gaussian likelihood can be
written as a contrast between
$\widetilde{J}_{n}(\omega;f_\theta)\overline{J_{n}(\omega)}$
and $f_{\theta}(\omega)$. The resulting estimator is based on
simultaneously predicting and fitting the spectral density. In the
case that the model is correctly specified, in the sense there exists
a $\theta\in \Theta$ where $f=f_{\theta}$ (and $f$ is the true
spectral density). Then
\begin{equation*}
\Ex_{f_\theta}[\widetilde{J}_{n}(\omega;f_\theta)\overline{J_{n}(\omega)}]
  = f_{\theta}(\omega)
\end{equation*}
and the Gaussian criterion has a clear interpretation. However, if the model is
misspecified (which for real data is likely), 
$\Ex_{f}[\widetilde{J}_{n}(\omega;f_\theta)\overline{J_{n}(\omega)}]$
has no clear interpretation. Instead, to understand what the Gaussian
likelihood is estimating, we use that
$\Ex_{f}[\widehat{J}_{n}(\omega;f_\theta)\overline{J_{n}(\omega)}] =
O(n^{-1})$, which leads to the approximation
$\Ex_{f}[\widetilde{J}_{n}(\omega;f_\theta)\overline{J_{n}(\omega)}] =  f(\omega) + O(n^{-1})$. 
From this, we observe that the expected negative log Gaussian likelihood is
\begin{equation*}
n^{-1}\Ex_f[\Xunder_{n}^{\prime}\Gamma_{n}(f_{\theta})^{-1}\Xunder_{n}] +
  n^{-1}\log |\Gamma_{n}(f_\theta)|=
 I(f,f_\theta)+ O(n^{-1}),
\end{equation*}
where 
\begin{equation}
\label{eq:Ifthetaf}
I_{n}(f;f_\theta) = \frac{1}{n}\sum_{k=1}^{n}\left(\frac{f(\omega_{k,n})}{f_{\theta}(\omega_{k,n})}
  + \log f_{\theta}(\omega_{k,n})\right).
\end{equation}
Since $I_{n}(f;f_\theta)$ is the spectral divergence between the true spectral $f$
density and parametric spectral density $f_\theta$, asymptotically the
misspecified Gaussian likelihood estimator has a meaningful
interpretation. However, there is still a finite sample bias in 
the Gaussian likelihood of order $O(n^{-1})$. This can have a knock-on effect,
 by increasing the finite sample bias in the resulting Gaussian likelihood estimator. 
To remedy this, in the following section, we obtain a frequency
domain criterion which approximates the spectral divergence
$I_{n}(f;f_\theta)$ to a greater degree of accuracy. This may 
lead to estimators which may give a more accurate fit of the underlying
spectral density. We should emphasis at this point, that reducing the
bias in the likelihood, does not necessarily translate to a provable
reduction in the bias of the resulting estimators (this is discussed
further in Section \ref{sec:asymprop}).

It is worth noting that, strictly, the spectral divergence is defined
as
\\* $n^{-1} \sum_{k=1}^{n}\left(\frac{f(\omega_{k,n})}{f_{\theta}(\omega_{k,n})}
  - \log\frac{f(\omega_{k,n})}{f_{\theta}(\omega_{k,n})}-1\right)$.
It is zero when $f_\theta =f$ and positive for other values of
$f_\theta$. But since $-\log f-1$ does not depend on $\theta$ we
ignore this term.

\subsection{The boundary corrected Whittle likelihood}\label{sec:boundary}

In order to address some of
the issues raised above, we recall from Theorem \ref{theorem:bio} that 
\\
$\Ex_{f}[\widetilde{J}_{n}(\omega;f)\overline{J_{n}(\omega)} ]=
f(\omega)$.  In other words, by predicting 
over the boundary using the (unobserved) spectral density which \emph{generates} the data,
the ``complete periodogram'' 
$\widetilde{J}_{n}(\omega;f)\overline{J_{n}(\omega)} $
is an inconsistent but \emph{unbiased} of the true spectral density $f$. 
This motivates the (infeasible) boundary corrected Whittle likelihood
\begin{equation}
\label{eq:Winf}
W_{n}(\theta) = \frac{1}{n}\sum_{k=1}^{n}
\frac{\widetilde{J}_{n}(\omega_{k,n};f)\overline{J_{n}(\omega_{k,n})}}{f_{\theta}(\omega_{k,n})}
+\frac{1}{n}\sum_{k=1}^{n}
\log f_{\theta}(\omega_{k,n}).
\end{equation}
Thus, if $\{X_{t}\}$ is a second order stationary time series with spectral
density $f$, then we have $\Ex_{f}[W_{n}(\theta) ] = I_{n}(f;f_\theta)$. 

Of course $f$ and thus 
$\widetilde{J}_{n}(\omega_{k,n};f)$ are unknown. 
However, we recall that $\widetilde{J}_{n}(\omega_{k,n};f)$
is comprised of the best linear predictors based on the unobserved
time series. The coefficients of the best linear predictors can be
replaced with the $h$-step ahead predictors evaluated with the best fitting
autoregressive parameters of order $p$ (the so called plug-in
estimators; see \cite{p:bha-96} and \cite{p:kle-19}). This is equivalent to
replacing $f$ in $\widetilde{J}_{n}(\omega_{k,n};f)$ with the spectral
density function corresponding to the best
fitting AR$(p)$ process $\widetilde{J}_{n}(\omega_{k,n};f_p)$, where
an analytic form is given in (\ref{eq:JAR}). Since we have replaced
$f$ with $f_p$, the ``periodogram''
$\widetilde{J}_{n}(\omega_{k,n};f_p)\overline{J_{n}(\omega_{k,n})}$
does have a bias, but it is considerably smaller than the bias of the
usual periodogram. In particular, it follows from the proof of Lemma
\ref{lemma:12}, below, that 
\begin{equation*}
\Ex_{f}[ \widetilde{J}_{n}(\omega_{k,n};f_p) \overline{J_{n}(\omega_{k,n})}]
  = f(\omega_{k,n}) + O\left(\frac{1}{np^{K-1}}\right).
\end{equation*}
The above result leads to an approximation of the boundary corrected Whittle likelihood
\begin{equation}
W_{p,n}(\theta) = 
\frac{1}{n}\sum_{k=1}^{n}\frac{\widetilde{J}_{n}(\omega_{k,n};f_p)\overline{J_{n}(\omega_{k,n})}}{f_{\theta}(\omega_{k,n})}+
\frac{1}{n}\sum_{k=1}^{n}
\log f_{\theta}(\omega_{k,n}).
\end{equation}
In the following lemma, we obtain a bound between the ``ideal''
boundary corrected Whittle likelihood $W_{n}(\theta)$ and 
$W_{p,n}(\theta)$. 

\begin{lemma}\label{lemma:12}
Suppose $f$ satisfies Assumption \ref{assum:A}, $f_\theta$ is
bounded away from zero and $\|f_\theta\|_{0}<\infty$. 
Let $\{a_{j}(p)\}$ denote the coefficients of the best fitting
AR$(p)$ model corresponding to the spectral density $f$ and define
$f_{p}(\omega) =|1-\sum_{j=1}^{p}a_{j}(p)e^{-ij\omega}|^{-2}$. 
Suppose $1\leq p<n$, then we have 
\begin{eqnarray}
&& \left\|F_{n}^{*}\Delta_{n}(f_{\theta}^{-1})\left(D_{n}(f_{})
  - D_{n}(f_p)\right)\right\|_{1} \nonumber \\
&& \quad \leq \rho_{p,K}(f) A_{K}(f,f_\theta)  \left( \frac{(C_{f,1}+1)}{p^{K-1}} +
\frac{2(C_{f,1}+1)^2 }{p^{K}}\|\psi_f\|_{0} \|\phi_f\|_{1} +
   \frac{C_{f,0}}{n^{K-1}} \right).
\label{eq:BD12}
\end{eqnarray}
Further, if $\{X_{t}\}$ is a time
series where $\sup_{t} \|X_{t}\|_{\Ex,2q} = \|X\|_{\Ex,2q}<\infty$ (for some
$q>1$), then
\begin{eqnarray}
 \| W_{n}(\theta) - W_{p,n}(\theta) \|_{\Ex,q} &\leq& 
\rho_{p,K}(f) A_{K}(f,f_\theta) \times  \nonumber \\
&& \left( \frac{(C_{f,1}+1)}{np^{K-1}} +
\frac{2(C_{f,1}+1)^2 }{np^{K}}\|\psi_f\|_{0} \|\phi_f\|_{1} + \frac{C_{f,0}}{n^{K}} \right)
\|X\|_{\Ex,2q}^{2}.
\label{eq:BD2}
\end{eqnarray}
\end{lemma}
 \noindent PROOF. See Appendix \ref{sec:approxproofs}. \hfill $\Box$ 

\begin{remark}
We briefly discuss what the above bounds mean for different types of
spectral densities $f$.
\begin{itemize}
\item[(i)] Suppose $f$ is the spectral density of 
a finite order AR$(p_0)$. If $p\geq p_0$, then 
\\*
$\left\|F_{n}^{*}\Delta_{n}(f_{\theta}^{-1})\left(D_{n}(f_{})
  - D_{n}(f_p)\right)\right\|_{1}=0$ and $ \| W_{n}(\theta) -
W_{p,n}(\theta) \|_{\Ex,q} =0$. On the other hand, if $p<p_0$ we replace the 
$p^{K}$ and $p^{K-1}$ terms in Lemma \ref{lemma:12} with 
$\sum_{j=p+1}^{p_0}|\phi_{j}|$ and $\sum_{j=p+1}^{p_0}|j\phi_{j}|$
respectively, where $\{\phi_{j}\}_{j=1}^{p}$ are the AR$(p)$
coefficients corresponding to $f$.
\item[(ii)] If the autocovariances corresponding to $f$ decay 
  geometrically fast to zero (for example an ARMA processes), then
  for some $0\leq \rho<1$ we have
\begin{eqnarray}
 \| W_{n}(\theta) - W_{p,n}(\theta) \|_{\Ex,q} = 
O\left(\frac{\rho^{p}}{n}+\rho^{n}\right).
\end{eqnarray}
\item[(iii)] If the autocovariances corresponding to $f$ decay to zero
  at a polynomial rate
  with $\sum_{r}|r^{K}c(r)|<\infty$, then
\begin{eqnarray}
 \| W_{n}(\theta) - W_{p,n}(\theta) \|_{\Ex,q} = 
O\left(\frac{1}{np^{K-1}}\right).
\end{eqnarray} 
Roughly speaking, the faster the rate of decay of the autocovariance
function, the ``closer'' $W_{p,n}(\theta)$ will be to $W_{n}(\theta)$
for a given $p$. 
\end{itemize}
\end{remark}

\vspace{2mm}
\noindent It follows from the lemma above that if $1\leq p<n$, $\Ex_{f}[W_{p,n}(\theta)]
=I_{n}(f;f_\theta) + O((np^{K-1})^{-1})$  and 
\begin{equation*}
W_{p,n}(\theta) = W_{n}(\theta) + O_{p}\left(\frac{1}{np^{K-1}}\right).
\end{equation*}
Thus if $p\rightarrow \infty$ as $n\rightarrow
\infty$, then $W_{p,n}(\theta)$ 
yields a better approximation to the ``ideal'' $W_{n}(\theta)$
than both the Whittle and the Gaussian likelihood. 

Since $f$ is unknown, $f_{p}$ is also unknown. But $f_{p}$ is easily
estimated from the data. We use the Yule-Walker estimator to fit an AR$(p)$ process to the observed time series,
where we select the order $p$ using the AIC. 
We denote this estimator as
$\widehat{\underline{\phi}}_{p}$ and the corresponding spectral density as 
$\widehat{f}_{p}$. Using this we define 
$\widehat{J}_{n}(\omega_{k,n};\widehat{f}_{p})$ where 
\begin{eqnarray*}
&& \widehat{J}_{n}(\omega;\widehat{f}_{p})= \\
&& \quad  
\frac{n^{-1/2}}{\widehat{\phi}_{p}(\omega)} \sum_{\ell=1}^{p}X_{\ell}\sum_{s=0}^{p-\ell}\widehat{\phi}_{\ell+s,p}e^{-is\omega}+
e^{in\omega}\frac{n^{-1/2}}{\overline{\widehat{\phi}_{p}(\omega)}} \sum_{\ell=1}^{p}X_{n+1-\ell}\sum_{s=0}^{p-\ell}
\widehat{\phi}_{\ell+s,p}e^{i(s+1)\omega},
\end{eqnarray*}
and $\widehat{\phi}_{p}(\omega) = 
1-\sum_{u=1}^{p}\widehat{\phi}_{u,p}e^{-iu\omega}$. This estimator
allows us to replace $\widehat{J}_{n}(\omega_{k,n};f_{p})$ in $W_{p,n}(\theta)$
with $\widehat{J}_{n}(\omega_{k,n};\widehat{f}_{p})$ to give the
``observed'' boundary corrected Whittle likelihood
\begin{equation} 
\label{eq:GF}
\widehat{W}_{p,n}(\theta)
=  \frac{1}{n}\sum_{k=1}^{n}
\frac{\widetilde{J}_{n}(\omega_{k,n};\widehat{f}_p)\overline{J_{n}(\omega_{k,n})}}{f_{\theta}(\omega_{k,n})}+
\frac{1}{n}\sum_{k=1}^{n}
\log f_{\theta}(\omega_{k,n}).
\end{equation}
We use as an
estimator of $\theta$, $\widehat{\theta}_{n} = \arg\min
\widehat{W}_{p,n}(\theta)$. 
It is worth bearing in mind that
\begin{equation*}
\I \frac{\widetilde{J}_{n}(\omega_{k,n};\widehat{f}_p)\overline{J_{n}(\omega_{k,n})}}{f_{\theta}(\omega_{k,n})}
  =-\I \frac{\widetilde{J}_{n}(\omega_{n-k,n};\widehat{f}_p)\overline{J_{n}(\omega_{n-k,n})}}{f_{\theta}(\omega_{n-k,n})}
\end{equation*}
thus $\widehat{W}_{p,n}(\theta)$ is real for all $\theta$. However,
due to rounding errors it is prudent to use
$\R\widehat{W}_{p,n}(\theta)$ in the minimisation algorithm. Sometimes 
$\R\widetilde{J}_{n}(\omega_{k,n};\widehat{f}_p)\overline{J_{n}(\omega_{k,n})}$
can be negative, when this arises we threshold it to be positive
(the method we use is given in Section \ref{sec:sim}).
 
In this paper, we focus on estimating
$\widehat{J}_{n}(\omega_{k,n};f_p)$ using the Yule-Walker
estimator. However, as pointed out by two referees, other estimators
could be used.  These may, in certain
situations, give better results. For example, in the case that $f$ has
a more peaked spectral density (corresponding to AR parameters close
to the unit circle) it may be better to replace the Yule-Walker
estimator with the tapered Yule-Walker estimator (as described in
\cite{p:dah-88} and \cite{p:zha-92}) or the Burg estimator. We show in Appendix \ref{sec:alternative}, that
using the tapered Yule-Walker estimator tends to give better results
for peaked spectral density functions. Alternatively one could
directly estimate $\widehat{J}_{\infty,n}(\omega_{k,n};f)$, where we
use a non-parametric spectral density estimator of $f$. This is
described in greater detail in Appendix \ref{sec:alternative} together with the results
of some simulations.

\subsection{The hybrid Whittle likelihood}

The simulations in Section \ref{sec:sim} suggest that the
boundary corrected Whittle likelihood estimator (defined in (\ref{eq:GF}))
yields an estimator with a smaller bias than the regular Whittle
likelihood. However, the bias of the tapered Whittle likelihood (and
often the Gaussian likelihood) is in some cases lower. The tapered
Whittle likelihood (first proposed in \cite{p:dah-88}) 
gives a better resolution at the peaks in the
spectral density. It also ``softens'' the observed domain of
observation.  With this in mind, we propose the hybrid 
Whittle likelihood which incorporates the notion of tapering.
 
Suppose $\underline{h}_{n} = \{h_{t,n}\}_{t=1}^{n}$ is a data taper, where
the weights $\{h_{t,n}\}$ are non-negative and 
$\sum_{t=1}^{n}h_{t,n} = n$. We
define the tapered DFT as 
\begin{equation*}
J_{n,\underline{h}_{n}}(\omega_{k,n})= n^{-1/2} \sum_{t=1}^{n} h_{t,n} X_{t} e^{it\omega_{k,n}}.
\end{equation*}
Suppose $f$ is the best fitting spectral density function. Using
that $\sum_{t=1}^{n}h_{t,n} = n$ and 
\\ $\cov_f(X_{t},\widehat{X}_{\tau,n})
= c_f(t-\tau)$ we have
\begin{equation}
\label{eq:taper}
\Ex_f [\widetilde{J}_{n}(\omega;f)\overline{J_{n,\underline{h}_{n}}(\omega)}]
  = f(\omega),
\end{equation}
which is analogous to the non-tapered result
$\Ex_{f}[\widetilde{J}_{n}(\omega;f)\overline{J_{n}(\omega)}] =
f(\omega)$.  Based on the above result we define the
infeasible hybrid Whittle likelihood which combines the regular DFT of the tapered time series and  the
complete DFT (which is not tapered)
\begin{equation}
\label{eq:GFinf}
H_{n}(\theta) = \frac{1}{n}\sum_{k=1}^{n}
\frac{\widetilde{J}_{n}(\omega_{k,n};f)\overline{J_{n,\underline{h}_{n}}(\omega_{k,n})}}{f_{\theta}(\omega_{k,n})}
+\frac{1}{n}\sum_{k=1}^{n}
\log f_{\theta}(\omega_{k,n}).
\end{equation}

\begin{figure}[h!]
\begin{center}
\includegraphics[scale = 0.5]{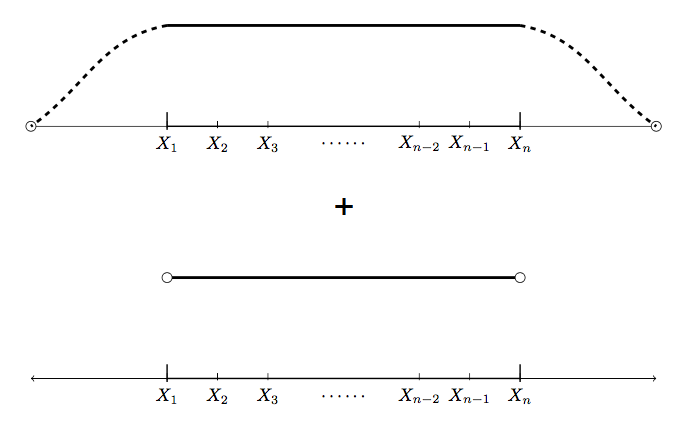}~~~~
 \includegraphics[scale = 0.5]{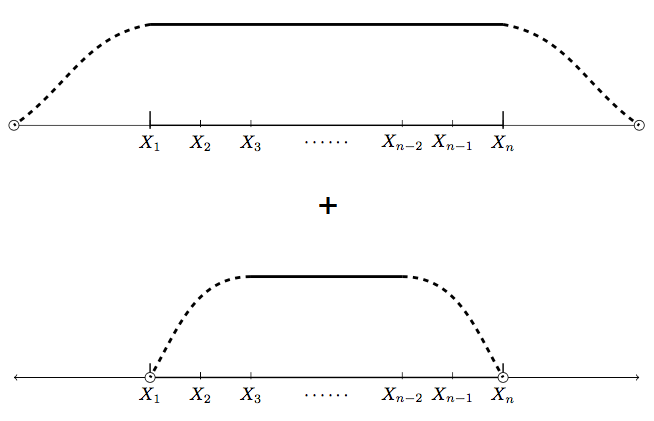}
\end{center}
\caption{{\it Left: The estimated complete DFT and the regular DFT which
  yields the boundary corrected Whittle likelihood. Right: The estimated complete DFT and the tapered DFT which
  forms the hybrid Whittle likelihood.} \label{fig:diff}}
\end{figure}

Using (\ref{eq:taper}), it can be shown that $\Ex_{f}[H_{n}(\theta) ] =
I_{n}(f;f_\theta)$. Thus $H_{n}(\theta)$ is an unbiased estimator of
$I_{n}(f;f_\theta)$.
Clearly, it is not possible to estimate
$\theta$ using the (unobserved) criterion $H_{n}(\theta)$. Instead we
replace $\widetilde{J}_{n}(\omega_{k,n};f)$ with its estimator
$\widetilde{J}_{n}(\omega_{k,n};\widehat{f}_{p})$ and define 
\begin{equation}
\label{eq:Hyb}
\widehat{H}_{p,n}(\theta) =
  \frac{1}{n}\sum_{k=1}^{n}\frac{\widetilde{J}_{n}(\omega_{k,n};
  \widehat{f}_{p})\overline{J_{n,\underline{h}_{n}}(\omega_{k,n})} }{f_{\theta}(\omega_{k,n})}+
\frac{1}{n}\sum_{k=1}^{n} \log f_{\theta}(\omega_{k,n}).
\end{equation} 
We then use as an estimator of $\theta$, $\widehat{\theta}_n =
\arg\min \widehat{H}_{p,n}(\theta)$. An illustration which visualises and compares the 
boundary corrected Whittle likelihood and hybrid Whittle likelihood is given
in Figure \ref{fig:diff}.
 
\section{The sampling properties of the hybrid Whittle likelihood}\label{sec:sample}

In this section, we study the sampling properties of the 
boundary corrected and hybrid Whittle likelihood. Our focus will be on 
the hybrid Whittle likelihood as it includes the boundary corrected likelihood as a special case,
  when  $h_{t,n}=1$ for $1\leq t \leq n$. 
In \cite{p:dsy-20} we study the sampling properties of 
the estimated complete periodogram
$\widetilde{J}_{n}(\omega;\widehat{f}_p)\overline{J_{n,\underline{h}_{n}}(\omega)}$.
Using these results and the results in Appendix \ref{sec:consistent}
and \ref{sec:bias}, we obtain the
bias and variance of the
boundary corrected and hybrid Whittle likelihood. 

Suppose we fit the spectral density $f_{\theta}(\omega)$ 
(where $\theta$ is an unknown $d$-dimension parameter vector) to the stationary
 time series $\{X_{t}\}_{t=1}^{n}$ whose true spectral density is
 $f$. The best fitting spectral density is $f_{\theta_n}$, where $\theta_{n} = \arg\min
 I_{n}(f;f_\theta)$. Let 
$\widehat{\theta}_{n}=(\widehat{\theta}_{1,n},\ldots,\widehat{\theta}_{d,n})$ be its
 estimator, where  $\widehat{\theta}_{n} =
\arg\min \widehat{H}_{p,n}(\theta)$. 

\subsection{Assumptions}

To derive the sampling properties
of $\widehat{\theta}_{n}$ we assume the data taper has the following form
\begin{eqnarray}
\label{eq:Htn}
h_{t,n} =  c_{n} h_{n}(t/n), 
\end{eqnarray} 
where $h_{n}:[0,1]\rightarrow \mathbb{R}$ is a sequence of positive functions
that satisfy the taper assumptions in Section 5, \cite{p:dah-88}
and $c_n = n/H_{1,n}$ with $H_{q,n} = \sum_{t=1}^{n}
h_{n}(t/n)^{q}$. We will assume $\sup_{t,n}h_{t,n} <\infty$, using
this it is straightforward to show 
$H_{2,n}/H_{1,n}^2 \sim
n^{-1}$. Under this condition, the hybrid 
Whittle is $n^{1/2}$--consistency and the equivalence result in
Theorem \ref{theorem:equivalence} holds. 
This assumption is used in
\cite{p:dah-83} and in practice one often
assumes that a fixed percentage of the data is tapered. A relaxation
of the condition $H_{2,n}/H_{1,n}^2 \sim n^{-1}$ will lead to a change
of rate in Theorem \ref{theorem:equivalence}.


\begin{assumption}[Assumptions on the parameter space]\label{assum:B}
\begin{itemize}
\item[(i)] The parameter space $\Theta \subset \mathbb{R}^{d}$ is compact,
$0<\inf_{\theta\in \Theta}\inf_{\omega} f_{\theta}(\omega)
\leq \sup_{\theta\in \Theta}\sup_{\omega} f_{\theta}(\omega)<\infty$
and $\theta_{n}$ lies in the interior of $\Theta$.
\item[(ii)] The one-step
  ahead prediction error $\sigma^{2}=\exp((2\pi)^{-1}\int_{0}^{2\pi} \log
f_\theta(\omega)d\omega)$ is not a function of the parameter $\theta$. 

\item[(iii)] Let $\{\phi_{j}(f_\theta)\}$ and $\{\psi_{j}(f_\theta)\}$
  denote the AR$(\infty)$  and  MA$(\infty)$ coefficients
  corresponding to the spectral density
  $f_\theta$ respectively. Then for all $\theta \in \Theta$ and $0\leq s \leq \kappa$ (for some $\kappa \geq 4$), we have 
\begin{equation*}
(a)\sup_{\theta\in
  \Theta}\sum_{j=1}^{\infty}\|j^{K}\nabla_{\theta}^{s}\phi_{j}(f_\theta)\|_1<\infty
\qquad 
(b) \sup_{\theta\in \Theta}\sum_{j=1}^{\infty}\|j^{K}\nabla_{\theta}^{s}\psi_{j}(f_\theta)\|_1<\infty,\quad 
\end{equation*} 
where $K>3/2$, 
$\nabla_{\theta}^{a}g(f_{\theta})$ is the $a$th order partial derivative of
$g$ with respect to $\theta$, and
$\|\nabla_{\theta}^{a}g(f_{\theta})\|_{1}$ denotes the absolute sum of all the
partial derivatives in $\nabla_{\theta}^{a}g(f_{\theta})$.
\end{itemize}
\end{assumption}
We use Assumption \ref{assum:B}(ii, iii)  to show 
that the $n^{-1}\sum_{k=1}^{n}\log f_{\theta}(\omega_{k,n})$ term in 
boundary corrected and hybrid Whittle likelihoods are
negligible with respect the other bias terms.
This allows us to simplify some of the bias expansions. 
Without Assumption \ref{assum:B}(ii, iii-a) the asymptotic bias of the new-frequency domain likelihood estimators 
would contain some additional terms. 
Assumption \ref{assum:B}(iii-b) is used to bound the $s$th derivative of the
spectral density.

\begin{assumption}[Assumptions on the time series]\label{assum:TS}
\begin{itemize}
\item[(i)] $\{X_{t}\}$ is a stationary time series. 
Let $\kappa_{\ell}(t_{1},\ldots,t_{\ell-1})$ denote the
joint cumulant  
$\cum(X_{0},X_{t_{1}},\ldots,X_{t_{\ell-1}})$.

Then for all $1\leq j\leq \ell\leq 12$, 
\begin{eqnarray*}
\sum_{t_{1},\ldots,t_{\ell-1}}|(1+t_{j})\kappa_{\ell}(t_{1},\ldots,t_{\ell-1})|<\infty.
\end{eqnarray*}
\item[(ii)] The spectral density of $\{X_t\}$ is such that the spectral density $f$ is bounded away from zero
and for some $K > 1$, the autocovariance function satisfies
$\sum_{r\in \mathbb{Z}}|r^{K}c_{f}(r)|<\infty$.
\item[(iii)] $I(\theta_{n})$ is invertible where 
\begin{eqnarray}
\label{eq:Ithetan}
I(\theta) = -\frac{1}{2\pi}\int_{0}^{2\pi}[\nabla_{\theta}^{2}f_{\theta}(\omega)^{-1}]f(\omega)d\omega.
\end{eqnarray}
\end{itemize} 
\end{assumption}
We require Assumption \ref{assum:TS}(i), when $\ell=4$ and $6$ to obtain a bound for the
expectation of the terms in the bias expansions and $\ell=12$ to show
equivalence between the feasible estimator based on
$\widehat{H}_{p,n}(\theta)$ and its
infeasible counterparts $H_{n}(\theta)$. Under Assumption \ref{assum:TS}(i,ii), Theorem 3.1
in \cite{p:dsy-20}, we can show that
\begin{eqnarray*}
\widehat{H}_{p,n}(\theta) = H_{n}(\theta) + O_{p}\left(\frac{p^{3}}{n^{3/2}}+\frac{1}{np^{K-1}} \right).
\end{eqnarray*}
Under Assumption \ref{assum:B}(i,iii) the above error 
is uniform over the parameter space. If the model is an AR$(p_0)$
and $p_0 \leq p$, then the term $O((np^{K-1})^{-1})$ in the above
disappears. 

To obtain a bound for the mean and variance of 
$\widehat{\theta}_{n}=(\widehat{\theta}_{1,n},\ldots,\widehat{\theta}_{d,n})$
we require the following quantities. Let 
\begin{eqnarray}
V(g,h) &=& \frac{2}{2\pi}\int_{0}^{2\pi}g(\omega)h(\omega)f(\omega)^{2}d\omega \nonumber \\
&&  +
  \frac{1}{(2\pi)^2  }\int_{0}^{2\pi}\int_{0}^{2\pi}g(\omega_1)h(\omega_2)f_{4}(\omega_1,-\omega_1,\omega_{2}) d\omega_{1}d\omega_{2} \nonumber\\
\text{and} \qquad \qquad J(g) &=&
                                  \frac{1}{2\pi}\int_{0}^{2\pi}g(\omega)
                                  f(\omega)   d\omega, \label{eq:VVV}
\end{eqnarray}
where $f_{4}$ denotes the fourth order cumulant density of the time
series $\{X_{t}\}$. We denote the $(s,r)$th element of
  $I(\theta_{n})^{-1}$ (where $I(\theta_{n})$ is defined in
  (\ref{eq:Ithetan})) as $I^{(s,r)}$, and define 
\begin{eqnarray}
G_{r}(\theta) &=& \sum_{s_1,s_2=1}^{d}I^{(s_1,s_2)}V\left(\frac{\partial f_\theta^{-1}}{
   \partial \theta_{s_2}}, \frac{\partial^{2}f_\theta^{-1}}{\partial
   \theta_{s_1}\partial \theta_{r}}\right)\nonumber\\
&& + \frac{1}{2}\sum_{s_1,s_2,s_3,s_4=1}^{d}I^{(s_1,s_3)}I^{(s_2,s_4)}
V\left(\frac{\partial f_{\theta}^{-1}}{\partial \theta_{s_3}}, \frac{\partial f_{\theta}^{-1}}{\partial \theta_{s_4}}\right)
J\left(\frac{\partial^{3}f_{\theta}^{-1}}{\partial
   \theta_{s_1}\partial\theta_{s_2}\partial \theta_r} \right). \label{eq:GGG}
\end{eqnarray}

\subsection{The asymptotic sampling properties}\label{sec:asymprop}

Using the assumptions above we
obtain a bound between the feasible and infeasible estimators. 
\begin{theorem}[Equivalence of feasible and infeasible estimators]\label{theorem:equivalence}
Suppose Assumptions \ref{assum:B} and \ref{assum:TS} hold. Define the
feasible and infeasible estimators as
$\widetilde{\theta}_{n}^{}=\arg\min H_{n}(\theta)$ and 
$\widehat{\theta}_{n}^{}=\arg\min \widehat{H}_{p,n}(\theta)$
 respectively. Then for $p\geq 1$ we have 
\begin{eqnarray*}
|\widehat{\theta}_{n}^{}-\widetilde{\theta}_{n}^{}|_{1} = 
O_{p}\left(\frac{p^{3}}{n^{3/2}}+\frac{1}{np^{K-1}}\right),
\end{eqnarray*}
where $|a|_{1}=\sum_{j=1}^{d}|a_{j}|$. For the case
  $p=0$, $\widehat{\theta}_n$ is the parameter estimator based on
  the Whittle likelihood using the one-sided tapered periodogram
$J_{n}(\omega_{k,n}) \overline{J_{n,\underline{h}_n}(\omega_{k,n})}$
rather than the regular tapered periodogram. In this case,
$|\widehat{\theta}_{n}^{}-\widetilde{\theta}_{n}^{}|_{1} = 
O_{p}\left(n^{-1}\right)$.

Note if the true spectral density of the time series is that of an
AR$(p_0)$ where $p_0\leq p$, then the $O((np^{K-1})^{-1})$ term is zero. 
\end{theorem}
PROOF. In Appendix \ref{sec:consistent}. \hfill $\Box$

\vspace{2mm}
\noindent The implication of the equivalence result is if $p^{3}/n^{1/2}\rightarrow 0$ as
$p\rightarrow \infty$ and
$n\rightarrow \infty$, then
$n|\widehat{\theta}_{n}^{}-\widetilde{\theta}_{n}^{}|_{1}\rightarrow
0$ and asymptotically
the properties of the infeasible estimator (such as bias and variance)
transfer to the feasible estimator.

\subsubsection{The bias and variance of the hybrid Whittle likelihood} \label{sec:bias_and_variance}

The expressions in this section are derived under 
Assumptions \ref{assum:B} and \ref{assum:TS}. 

\vspace{2mm}

\noindent \underline{The bias} We show in Appendix \ref{sec:multi},
that the asymptotic bias (in the sense of Bartlett) for
$\widehat{\theta}_{n}=(\widehat{\theta}_{1,n},\ldots,\widehat{\theta}_{d,n})$
is  
\begin{eqnarray}
\label{eq:biasexpan}
\Ex[\widehat{\theta}_{j,n} - \theta_{j,n}] =
  \frac{H_{2,n}}{H_{1,n}^{2}}\sum_{r=1}^{d}I^{(j,r)}G_{r}(\theta_n)
+ O\left(\frac{p^{3}}{n^{3/2}} + \frac{1}{np^{K-1}}\right)
\quad 1\leq j \leq d,
\end{eqnarray}
where $I^{(j,r)}$ and $G_{r}(\theta_n)$ is defined in (\ref{eq:GGG}).
We note that if no tapering were used then $H_{2,n}/H_{1,n}^{2}= n^{-1}$.
The Gaussian and Whittle likelihood have a bias which
includes the above term (where $H_{2,n}/H_{1,n}^{2}
= n^{-1}$) plus an additional term of the form
$\sum_{r=1}^{d}I^{(j,r)}\Ex[\nabla_{\theta}L_{n}(\theta_{n})]$, where $L_{n}(\cdot)$ is the
Gaussian or Whittle likelihood (see Appendix
\ref{sec:multi} for the details). 

Theoretically,
it is unclear which criteria has the smallest bias (since the
inclusion of additional terms does not necessarily increase the
bias). However, for the hybrid Whittle likelihood estimator, a 
straightfoward ``Bartlett correction'' can be made to estimate the
 bias in (\ref{eq:biasexpan}). We briefly outline how this can be
 done. We observe that the bias is built of
$I(\cdot)$, $J(\cdot)$ and $V(\cdot,\cdot)$. Both $I(\cdot)$ and $J(\cdot)$ can easily
be estimated with their sample means. The term $V(\cdot,\cdot)$
can also be estimated  by using an adaption of 
orthogonal samples (see \cite{p:sub-18}), which we now describe. Define the random variable
\begin{eqnarray*}
h_{r}(g;f) = \frac{1}{n}\sum_{k=1}^{n}g(\omega_{k,n})
\widetilde{J}_{n}(\omega_{k+r,n};f) \overline{J_{n}(\omega_{k,n})}
\quad \textrm{for} \quad r\geq 1,
\end{eqnarray*}
where $g$ is a continuous and bounded function. Suppose $g_{1}$ and
$g_{2}$ are continuous and bounded functions. 
If $r\neq n \mathbb{Z}$, then $\Ex_{f}[h_{r}(g_{j};f)]=0$ (for $j=1$ and $2$). But interestingly, if
$r<<n$, then
$n\cov_{f}[h_{r}(g_{1};f),h_{r}(g_{2};f)]=n\Ex_{f}[h_{r}(g_{1};f)\overline{h_{r}(g_{2};f)}]
=  V(g_{1},g_{2})+O(r/n)$.
Using these results, we estimate $V(g_{1},g_{2})$ by replacing
$h_{r}(g_{j};f)$ with $h_{r}(g_{j};\widehat{f}_{p})$ and defining the
``sample covariance'' 
\begin{eqnarray*}
\widehat{V}_{M}(g_{1},g_{2}) =
  \frac{n}{M}\sum_{r=1}^{M}h_{r}(g_{1};\widehat{f}_{p})
\overline{h_{r}(g_{2};\widehat{f}_{p})}
\end{eqnarray*}
where $M<<n$. Thus, $\widehat{V}_{M}(g_{1},g_{2})$ is an estimator of
$V(g_1,g_2)$. Based on this construction,
\begin{eqnarray*}
\widehat{V}_{M}\left(\frac{\partial}{
   \partial \theta_{s_2}} f_{\widehat{\theta}_n}^{-1}, \frac{\partial^{2}}{
   \partial \theta_{s_1}\partial\theta_{r}} f_{\widehat{\theta}_n}^{-1} \right) 
\quad \textrm{and} \quad
\widehat{V}_{M}\left(\frac{\partial}{
   \partial \theta_{s_3}} f_{\widehat{\theta}_n}^{-1}, \frac{\partial}{
   \partial \theta_{s_4}} f_{\widehat{\theta}_n}^{-1} \right)
\end{eqnarray*}
are estimators of 
$V\left(\frac{\partial f_\theta^{-1}}{
   \partial \theta_{s_2}}, \frac{\partial^{2}f_\theta^{-1}}{\partial
   \theta_{s_1}\partial \theta_{r}}\right)$ and 
$V\left(\frac{\partial f_{\theta}^{-1}}{\partial \theta_{s_3}},
  \frac{\partial f_{\theta}^{-1}}{\partial \theta_{s_4}}\right)$
respectively. 
This estimation scheme yields a consistent estimate of the
bias  even when the model is misspecified. In contrast, it is unclear how a bias
correction would work for the Gaussian and Whittle likelihood under
misspecification, as they also involve the term
$\Ex_{f}[\nabla_{\theta}L_{n}(\theta_{n})]$. In the case of misspecification,
$\Ex_{f}[\nabla_{\theta}L_{n}(\theta_{n})]\neq 0$ and is of order $O(n^{-1})$.

It is worth mentioning that the asymptotic expansion in (\ref{eq:biasexpan}) does
not fully depict what we observe in the simulations in Section \ref{sec:sim}. A
theoretical comparison of the biases of both new likelihoods show
that for the boundary corrected Whittle likelihood, the bias is asymptotically
$n^{-1}\sum_{r=1}^{d}I^{(j,r)}G_{r}(\theta_n)$, whereas when tapering is
used the bias is $(H_{2,n}/H_{1,n}^{2})\sum_{r=1}^{d}I^{(j,r)}G_{r}(\theta_n)\geq n^{-1}\sum_{r=1}^{d}I^{(j,r)}G_{r}(\theta_n)$. This
would suggest that the hybrid Whittle likelihood should have a larger
bias than the boundary corrected Whittle likelihood. 
But the simulations (see Section \ref{sec:sim}) suggest this is not
necessarily true and the hybrid likelihood tends to have a smaller bias. 

\vspace{2mm}

\noindent \underline{The variance} We show in Corollary 3.1, \cite{p:dsy-20} that the inclusion of the 
prediction DFT
in the hybrid Whittle likelihood has a variance which asymptotically
is small as compared with the main Whittle term if $p^{3}/n\rightarrow 0$ as
$p,n\rightarrow \infty$ (under the condition $H_{2,n}/H_{1,n}^{2}\sim n^{-1}$)
Using this observation, standard Taylor expansion methods and 
Corollary 3.1 in \cite{p:dsy-20},
the asymptotic variance of $\widehat{\theta}_{n}$ is 
\begin{eqnarray*}
\frac{H_{1,n}^2}{H_{2,n}} \var (\widehat{\theta}_{n})
  = I(\theta_n)^{-1}V\left(\nabla_{\theta}f_{\theta}^{-1},
  \nabla_{\theta}f_{\theta}^{-1}\right)\rfloor_{\theta=\theta_n}I(\theta_{n})^{-1}
  + o(1),
\end{eqnarray*}
where $V(\cdot)$ is defined in (\ref{eq:VVV}).


\subsubsection{The role of order estimation on the rates}

The order in the AR$(p)$ approximation is selected using the AIC, where
$\widehat{p}=\arg\min \text{AIC}(p)$ with
\begin{eqnarray*}
\text{AIC}(p) = \log \widehat{\sigma}_{p,n}^{2} + \frac{2p}{n},
\end{eqnarray*}
$\widehat{\sigma}_{p,n}^{2} =
\frac{1}{n-K_{n}}\sum_{t=K_{n}}^{n}(X_{t} -
\sum_{j=1}^{p}\widehat{\phi}_{j,p}X_{t-j})^{2}$, $K_{n}$ is such that
$K_{n}^{2+\delta} \sim n$ for some $\delta>0$.
\cite{p:ing-wei-05} assume that the underlying time series is a
linear, stationary time series with an AR$(\infty)$ that satisfies
Assumption K.1$-$K.4 in \cite{p:ing-wei-05}. They show that under the
condition that the AR$(\infty)$ coefficients satisfy $(\sum_{j=p+1}^{\infty}|\phi_{j}|)^{2} = O(p^{-2K})$,
then $\widehat{p}  = O_{p}(n^{1/(1+2K)})$ (see Example 2 in \cite{p:ing-wei-05}).
Thus, if $K>5/2$, then $\widehat{p}^3/n^{1/2} \Pcon 0$ (where
$\widehat{p}  = O_{p}(n^{1/(1+2K)})$) and $\widehat{p} \Pcon \infty$ as
$n\rightarrow \infty$. These rates ensure that the difference between
the feasible and infeasible estimator is
$|\widehat{\theta}_{n}^{}-\widetilde{\theta}_{n}^{}|_{1}=o_{p}(n^{-1})$. Thus
the feasible estimator, constructed using the AIC, and the
infeasible estimator are equivalent and the bias and variance derived
above are valid for this infeasible estimator.

\subsubsection{The computational cost of the estimators}

We now discuss some of the implementation issues of the new estimators. 

The Durbin-Levinson algorithm is often used to maximize the Gaussian
likelihood. If this is employed, then the computational cost of the
algorithm is $O(n^{2})$. On the other hand, by using the FFT, the
computational cost of the Whittle likelihood is $O(n\log n)$. 

For the boundary corrected Whittle and hybrid Whittle likelihood algorithm,
there is an additional cost over the Whittle likelihood due to the
estimation of $\{\widehat{J}_{n}(\omega_{k,n};
\widehat{f}_p)\}_{k=1}^{n}$. We recall that $\widehat{f}_p$ is
constructed using the Yule-Walker estimator
$\widehat{\underline{\phi}}_p =(\widehat{\phi}_{1,p}, ...,\widehat{\phi}_{p,p})^{\prime}$
where  $p$ is selected with the AIC.
We now calculate the complexity of calculating $\{\widehat{J}_{n}(\omega_{k,n};
\widehat{f}_p)\}_{k=1}^{n}$. 

The sample autocovariances, $\{\widehat{c}_{n}(r)\}_{r=0}^{n-1}$ (which are required
in the Yule-Walker estimator) can be calculated in $O(n\log n)$ operations. 
Let $K_n$ denote the maximum order used for the evaluation of the AIC. If
we implement the Durbin-Levinson algorithm, 
then evaluating $\widehat{\underline{\phi}}_p$ for $1\leq p \leq K_{n}$
requires in total
$O(K_n^2)$ arithmetic operations. 
Given the estimated 
AR coefficients $\widehat{\underline{\phi}}_{\hat{p}}$, the predictive DFT $\{\widehat{J}_{n}(\omega_{k,n};
\widehat{f}_p)\}_{k=1}^{n}$ can be calculated in $O(\min(n \log n, n\hat{p}))$ arithmetic operations
(the details of the algorithm for optimal calculation can be found in Appendix \ref{sec:sup}).
Therefore, the overall computational cost of implementing both the
 boundary corrected Whittle and hybrid Whittle likelihood algorithms
 is  $O(n\log n + K_n^2)$.

Using \cite{p:ing-wei-05} Example 2, for consistent order selection
$K_{n}$ should be such that $K_n \sim
n^{1/(2K+1)+\varepsilon}$ for some $\varepsilon>0$ (where $K$ is
defined in Assumption \ref{assum:A}). 
Therefore, we conclude that the computational cost of the new
likelihoods is of the same order as the Whittle likelihood.

\input{simul2}

\section{Concluding remarks and discussion}\label{sec:conc}

In this paper we have derived an exact expression for the differences  
$\Gamma_{n}(f_{\theta})^{-1} -C_{n}(f_{\theta}^{-1})$ and 
$\underline{X}^{\prime}_{n}[\Gamma_{n}(f_{\theta})^{-1}-C_{n}(f_{\theta}^{-1})]\Xunder_n$.
These expressions are simple, with an intuitive interpretation, in terms
of predicting outside the boundary of observation. 
They also provide a new perspective to the Whittle likelihood as an
approximation based on a biorthogonal transformation.
We have used these expansions and approximations to define two new 
spectral divergence criteria (in the frequency domain).
Our simulations show that both new estimators (termed the boundary
corrected Whittle and hybrid Whittle) tend to outperform the Whittle
likelihood. Intriguingly, the hybrid Whittle likelihood 
tends to outperform the boundary corrected Whittle likelihood. Currently, we have
no theoretical justification for this and one future aim is to investigate these differences. 

We believe that it is possible to use a similar construction to obtain
an expression for the difference between the Gaussian likelihood of a multivariate time
series and the corresponding multivariate Whittle likelihood. 
The construction we use in this paper hinges on past and future predictions. In the
univariate set-up there is an elegant symmetry for the predictors in
the past and future. In the multivariate set-up there are some
important differences. This leads to interesting, but different
expressions for the predictive DFT. To prove analogous
results to those in this paper, we will require Baxter-type
inequalities for the multivariate framework. The
bounds derived in \cite{p:che-pou-93} and \cite{p:ino-18} may be useful
in this context. 


The emphasis of this paper is on short memory time
  series. But we conclude by briefly discussing
  extensions to long memory time series. The fundamental feature (in the frequency domain) that
distinguishes a short memory time series from a long memory time
series is that the spectral density of a long memory time series is
not bounded at the origin. However, we conjecture that the complete DFT
described in Theorem \ref{theorem:bio} can have applications within this
setting too. Suppose $f$ is the spectral density corresponding to a long
memory time series. And let $\widehat{J}_{n}(\omega_{k,n};f)$ be
defined as in (\ref{eq:28}). We assume that for $1\leq k \leq n-1$
that $\widehat{J}_{n}(\omega_{k,n};f)$ is a well defined random
variable. Under these conditions, simple calculations show that equation (\ref{eq:DFTPred}) in
Theorem \ref{theorem:bio} applies for $1\leq k_1,k_2\leq n-1$, but not when
$k=n$ (the frequency at the origin). Using this,  
a result analogous to Corollary
\ref{corollary:inverse} can be obtained for long memory time
series. 


\begin{theorem}[Inverse Toeplitz identity for long memory time series]\label{theorem:longmemory}
Suppose that the spectral density of a time series, $f$, is bounded
away from zero and bounded on
$[\varepsilon,\pi]$ for any $\varepsilon>0$ and satisfies
$f(\omega)\sim c_{f}|\omega|^{-2d}$ as $\omega \rightarrow 0$ for some
$d \in (0,1/2)$ and $c_{f}>0$. Define the $(n-1)\times n$ submatrix
$(\widetilde{D}_{n}(f))_{k,t}=(D_{n}(f))_{k,t}$ for $1\leq k \leq
(n-1)$, where $D_{n}(f)$ is defined in Theorem
\ref{thm:higherorder}, equation (\ref{eq:Dapprox1}). We assume that the entries of
$\widetilde{D}_{n}(f)$ are finite. 
Let $\widetilde{F}_{n}$
denotes $(n-1)\times n$ submatrix of $F_n$,
where 
$(\widetilde{F}_{n})_{k,t} = n^{-1/2}e^{ik\omega_{t,n}}$. Let $ \widetilde{\Delta}_{n}(f) =
\diag(f(\omega_{1,n}),\ldots,f(\omega_{n-1,n}))$. Then  
\begin{eqnarray}
\label{eq:longinverse}
(I_n - n^{-1} \textbf{1}_n \textbf{1}_n^{\prime}) \Gamma_n(f)^{-1} =
\widetilde{F}_n^{*} \widetilde{\Delta}_{n}(f)^{-1}(\widetilde{F}_{n}+\widetilde{D}_{n}(f))
\end{eqnarray}
where $\textbf{1}_n = (1, ..., 1)^{\prime}$ is an $n$-dimension vector
\end{theorem}
PROOF. See Appendix \ref{sec:proofbio}. \hfill $\Box$

\vspace{2mm}
We note that an expansion of $D_{n}(f)_{k,t}$ and
$\widehat{J}_{n}(\omega_{k,n};f)$ is given in Theorem
\ref{thm:higherorder} in terms of the AR$(\infty)$ coefficients
corresponding to the spectral density $f$.  The theorem above is contingent
on these quantities being finite. We
conjecture that this holds for invertible ARFIMA$(p,d,q)$
time series where $d\in (0,1/2)$. 
Unfortunately, a precise proof of this condition uses a different set of tools to those
developed in this paper. Thus we leave it for future research. 

The plug-in Gaussian likelihood replaces the population mean in the
Gaussian likelihood with the sample mean. 
We now apply the above result to representing the long memory plug-in Gaussian
likelihood in the frequency domain. 
We define the demeaned time series $\underline{X}_{n}^{(c)} =
\underline{X}_{n} - \overline{X} \textbf{1}_n$, where $\bar{X} =
n^{-1}\sum_{t=1}^{n}X_{t}$. Then by using Theorem \ref{theorem:longmemory} we have
\begin{eqnarray*}
\frac{1}{n}\underline{X}_{n}^{(c) \prime}\Gamma_{n}(f_\theta)^{-1}\underline{X}_{n}^{(c)}
&=& \frac{1}{n}\sum_{k=1}^{n-1}
    \frac{|J_{n}^{}(\omega_{k,n})|^{2}}{f_\theta(\omega_{k,n})}
    + \frac{1}{n}\sum_{k=1}^{n-1} \frac{\widehat{J}_{n}^{(c)}(\omega_{k,n};f_\theta)
    \overline{J_{n}^{}(\omega_{k,n})}}{f_\theta(\omega_{k,n})},
\end{eqnarray*} 
where $\widehat{J}_{n}^{(c)}(\cdot;f_\theta)$ denotes the predictive
DFT of the demeaned time series $\underline{X}_n^{(c)}$.  It would be
of interest to show that the new likelihoods defined in this paper (and a local
 frequency version of them) could be used in the analysis of
 long memory time series. In Appendix \ref{sec:sim-long} we have presented some preliminary
 simulations for long memory time series. The results are not
 conclusive, but they do suggest that the new likelihoods can, in some
 settings, reduce the bias for long memory parameter estimators. 

In summary, the notion of
biorthogonality and its application to the inversion of certain
variance matrices may be of value in future research.


\section*{Acknowledgements}
The research was partially conducted while SSR was visiting the
Universit\"at Heidelberg, SSR is extremely gratefully to the
hospitality of everyone at Mathematikon. SSR is grateful to Thomas
Hotz for several fruitful discussions. 
Finally, this paper is dedicated to SSR's Father, Tata Subba Rao, who
introduced the Whittle likelihood to (the young and rather confused)
SSR many years ago. SSR and JY gratefully acknowledge the  support of the National
Science Foundation (grant DMS-1812054). 

The authors wish to thank three anonymous referees and editors for their
insightful observations and suggestions, which substantially improved all aspects of the paper.

%% file: simul2.tex
\section{Empirical results}\label{sec:sim}

To substantiate our theoretical results, we conduct some
simulations (further simulations can be found in Appendix \ref{sec:sim_appendix},
\ref{sec:sim-long} and \ref{sec:alternative}). To compare different methods, 
we evaluate six different quasi-likelihoods: the Gaussian likelihood
(equation (\ref{eq:likeG})), 
the Whittle likelihood (equation (\ref{eq:like1})),
the boundary corrected Whittle likelihood (equation (\ref{eq:GF})),
the hybrid Whittle likelihood (equation (\ref{eq:Hyb})),
the tapered Whittle likelihood (p.810 of \cite{p:dah-88}) and the
debiased Whittle likelihood 
(equation (7) in \cite{p:olh-19}).

The tapered and hybrid Whittle likelihoods require the use of data tapers. We use a Tukey taper (also known as the cosine-bell taper) where
\begin{eqnarray*} 
h_{n}\left(t/n\right) = \left\{ \begin{array}{ll}
\frac{1}{2}[1-\cos(\pi(t-\frac{1}{2})/d)] & 1\leq t \leq d \\
1 & d+1\leq t \leq n-d \\
\frac{1}{2}[1-\cos(\pi(n-t+\frac{1}{2})/d)] & n-d+1\leq t \leq n
\end{array}.
\right. 
\end{eqnarray*}
We set the proportion of tapering at each end of the time series is
$0.1$, i.e. $d=n/10$ (the default in R). 

When evaluating the boundary corrected Whittle likelihood and
hybrid Whittle likelihood, the order $p$ is selected with the AIC and
$\widehat{f}_{p}$ is estimated using the Yule-Walker estimator.  

Unlike the Whittle, the tapered Whittle and debiased
Whittle likelihood,
$\R\widetilde{J}_{n}(\omega_{k,n};\widehat{f}_p)\overline{J_{n}(\omega_{k,n})}$
and $\R \widetilde{J}_{n}(\omega_{k,n};\widehat{f}_p)\overline{J_{n,\underline{h}_n}(\omega_{k,n})}$
can be negative. 
To avoid negative values, we
apply the thresholding function 
$f(t) = \max(t, 10^{-3})$ to $\R \widetilde{J}_{n}(\omega_{k,n};\widehat{f}_p)\overline{J_{n}(\omega_{k,n})}$
and $\R\widetilde{J}_{n}(\omega_{k,n};\widehat{f}_p)\overline{J_{n,\underline{h}_n}(\omega_{k,n})}$
over all the frequencies. Thresholding induces an additional (small)
bias to the new criteria. The proportion of times that
$\R \widetilde{J}_{n}(\omega_{k,n};\widehat{f}_p)\overline{J_{n}(\omega_{k,n})}$ drops
 below the threshold increases for spectral
 density functions with large peaks and when the spectral density is
 close to zero.
However, at least for the models that we studied in the simulations, the bias due to the thresholding 
is negligible.

All simulations are conducted over 1000 replications
with sample sizes $n=20$, $50$, and $300$.  In all the tables below
and Appendix, the bias of the estimates are in parenthesis and the
standard deviation are in
parenthesis. The ordering of the performance of the
estimators is colour coded and is based on their squared root of the mean squared error (RMSE).

\subsection{Estimation with correctly specified models} \label{sec:specified}

We first study the AR(1) and MA(1) parameter estimates when the models
are correctly specified. 
We generate two types of time series models $\Xunder_{n}$ and
$\underline{Y}_{n}$, which satisfy the following recursions 
\begin{eqnarray*}
{\bf AR(1)}:&& X_{t} = \theta X_{t-1}+e_{t}; \quad  \phi_{X}(\omega) = 1-\theta e^{-i\omega}\\
{\bf MA(1)}:&& Y_{t} = e_{t}+\theta e_{t-1}; \quad \phi_{Y}(\omega) = (1+\theta e^{-i\omega})^{-1},
\end{eqnarray*}
where
$|\theta|<1$, $\{e_{t}\}$ are independent, identically distributed
Gaussian random variables with mean 0 and variance 1.
Note that the Gaussianity of the innovations is not required to obtain the theoretical properties of the 
estimations. In Appendix \ref{sec:AR-MA-chi}, we include simulations 
when the innovations follow a standardized chi-squared distribution with two degrees of freedom. 
The results are similar to those with Gaussian innovations.
We generate the AR$(1)$ and
MA$(1)$ models with parameters $\theta = 0.1,0.3,0.5,0.7$ and
$0.9$. For the time series generated by an AR$(1)$ process, we fit an
AR$(1)$ model, similarly, for the time series generated by a MA$(1)$ process we fit a
MA$(1)$ model. 

\begin{figure}[]
\begin{center}

\textbf{AR$(1)$ model}

\vspace{1em}

\includegraphics[scale=0.35,page=1]{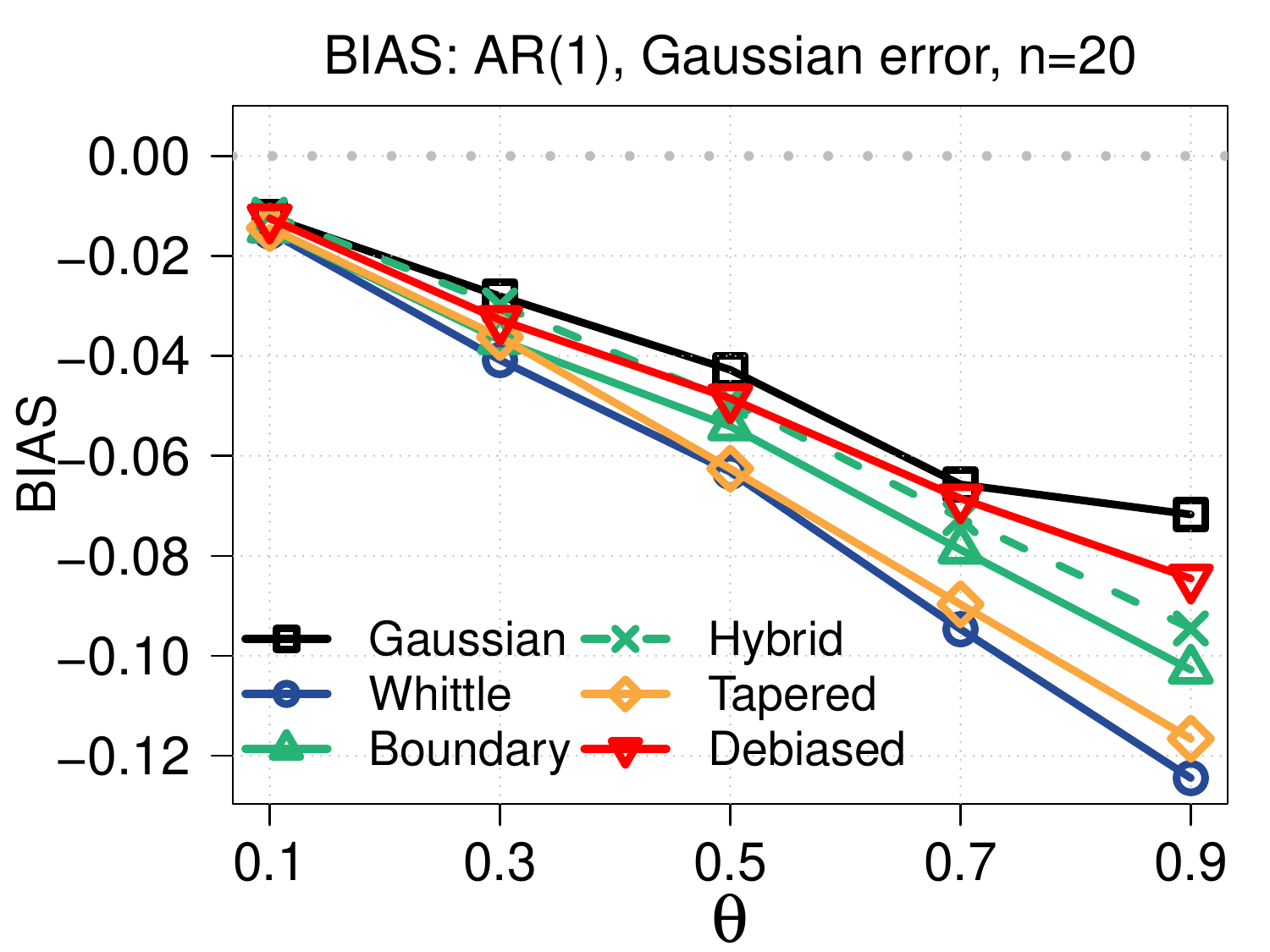}
\includegraphics[scale=0.35,page=3]{plot/AR_gaus2.pdf}
\includegraphics[scale=0.35,page=5]{plot/AR_gaus2.pdf}

\includegraphics[scale=0.35,page=2]{plot/AR_gaus2.pdf}
\includegraphics[scale=0.35,page=4]{plot/AR_gaus2.pdf}
\includegraphics[scale=0.35,page=6]{plot/AR_gaus2.pdf}

\vspace{1em}

\textbf{MA$(1)$ model}

\vspace{1em}

\includegraphics[scale=0.35,page=1]{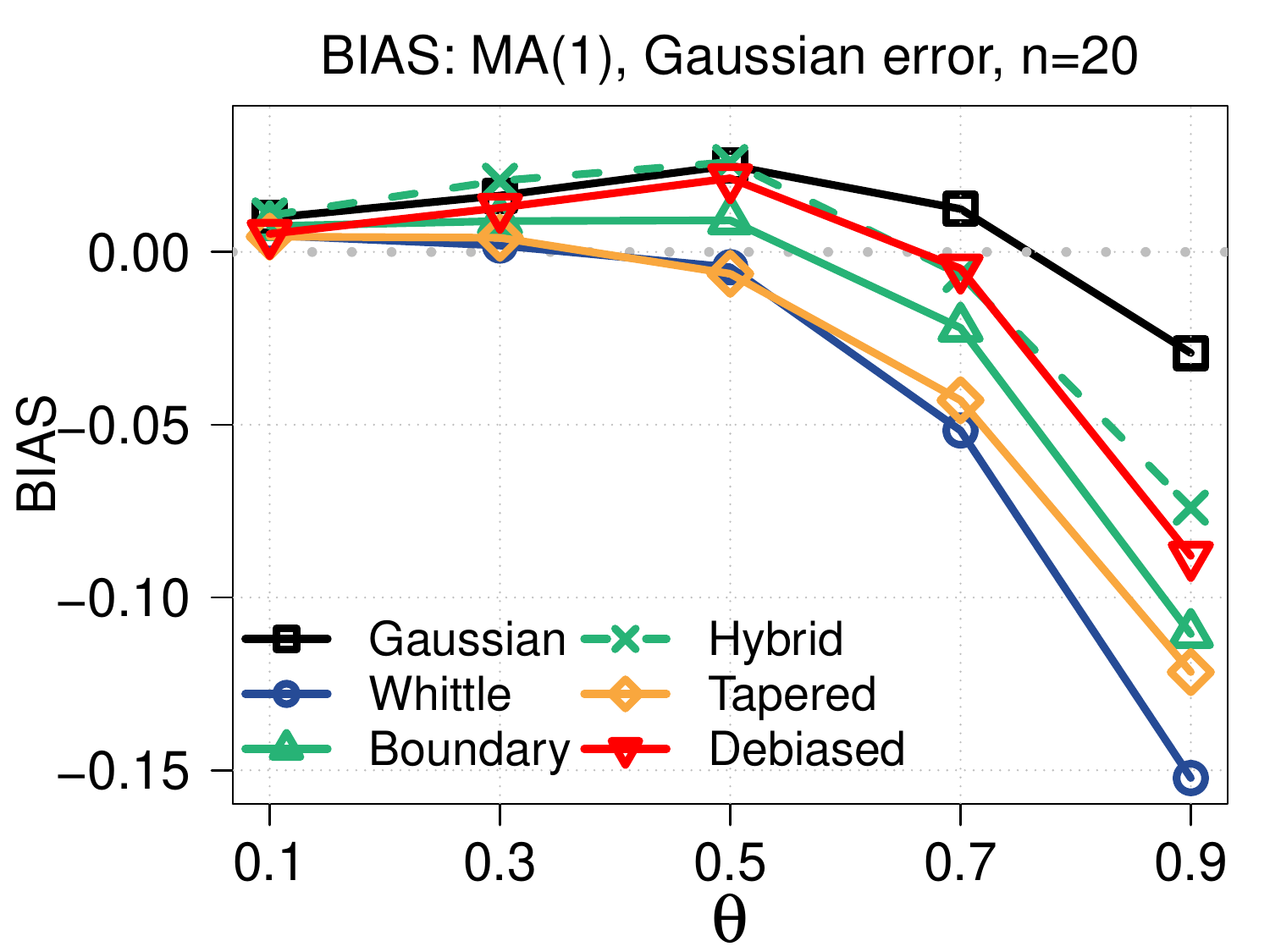}
\includegraphics[scale=0.35,page=3]{plot/MA_gaus2.pdf}
\includegraphics[scale=0.35,page=5]{plot/MA_gaus2.pdf}

\includegraphics[scale=0.35,page=2]{plot/MA_gaus2.pdf}
\includegraphics[scale=0.35,page=4]{plot/MA_gaus2.pdf}
\includegraphics[scale=0.35,page=6]{plot/MA_gaus2.pdf}

\caption{\textit{Bias (first row) and the RMSE (second row) of the parameter estimates for the Gaussian AR(1) models and Gaussian MA(1) models.  Length of the time series $n=20$(left), $50$(middle), and $300$(right). } }
\label{fig:ar1ma1} 
\end{center}
\end{figure}

For each simulation, we evaluate the six different parameter
estimators. The empirical bias and standard deviation are calculated. 
Figures \ref{fig:ar1ma1} gives the bias (first row) and the RMSE
(second row) of 
each estimated parameter $\theta$ for both AR$(1)$ and MA$(1)$ models. We focus on positive $\theta$,
similar results are obtained for negative $\theta$. The
results are also summarized in Table \ref{tab:AR} in Appendix \ref{sec:AR-MA-Gaussian}. 

For both AR$(1)$ and MA$(1)$ models, we observe a stark difference
between the bias of the Whittle likelihood estimator (blue line) and the other
five other methods, which in most cases have a lower bias.
The Gaussian likelihood performs uniformly well for both models and all sample sizes. 
Whereas, the tapered Whittle estimator performs very well for
the MA$(1)$ model but not quite as well for the AR$(1)$ model.
The debiased Whittle likelihood performs quite well for both models, especially  when the parameter values 
are small (e.g. $\theta = 0.1, 0.3$, and $0.5$).

The simulations suggest that the boundary corrected and hybrid Whittle
likelihoods (referred from now on as the new likelihoods) are 
competitive with the benchmark Gaussian likelihood for both AR$(1)$ and MA$(1)$ models.
For the AR$(1)$ model the new likelihoods tend to have the smallest or second smallest
RMSE (over all sample sizes and more so when $\phi$ is large). 
A caveat is that for the AR$(1)$ model the bias of the boundary corrected
and hybrid Whittle tends to be a little larger than the bias of the Gaussian
likelihood (especially for the smaller sample sizes). This is
interesting, because in Appendix \ref{sec:AR1bias} we show that if the 
AR$(1)$ model is correctly specified, the first order bias of the boundary corrected Whittle likelihood and
the Gaussian likelihood are the same (both are $-2\theta/n$).
The bias of the hybrid Whittle likelihood is slightly large, due
to the data taper. However, there are differences in the second order
expansions. Specifically, for the Gaussian likelihood,
it is $O(n^{-3/2})$, whereas, for the boundary corrected and
hybrid Whittle it is $O(p^3n^{-3/2})$. Indeed, the
$O(p^{3}n^{-3/2})$ term arises because of the parameter estimation in
the predictive DFT. This term is likely to dominate the $O(n^{-3/2})$ in the
Gaussian likelihood. Therefore, for small sample sizes, the second order terms
can impact the bias. It is this second order term that may be causing the
larger bias seen in the boundary corrected Whittle likelihood as
compared with the Gaussian likelihood.  

On the other hand, the bias for the MA$(1)$ model tends to be smaller
for the new likelihoods, including the benchmark
Gaussian likelihood. Surprisingly, there appears to be examples
where the new likelihood does better (in terms of RMSE) than the
Gaussian likelihood. This happens when $n \in \{50, 300\}$ for $\theta =
0.9$. This observation is noteworthy, as the computational cost of
the Gaussian likelihood is greater than the computational cost of the
new likelihoods (see Section \ref{sec:asymprop}). Thus the simulations suggest that in certain
situations the new estimator may outperform the Gaussian likelihood
at a lower computational cost. 

In summary, the new likelihoods perform well compared with the
standard methods, including the benchmark Gaussian likelihood. As
expected, for large sample sizes the performance of all the estimators
improves considerably. And for some models, the new likelihood is able
to outperform the Gaussian likelihood estimator. Though there is no
clear rule when this will happen.

\subsection{Estimation under misspecification}\label{sec:misspecifiedmodel}

\noindent Next, we turn into our attention to the case that the model
is misspecified (which is more realistic for real data). 
As we mentioned above, the estimation of the AR parameters in the
predictive DFT of the new likelihoods leads to an additional error of
order $O(p^{3}n^{-3/2})$. The more complex the model, the larger $p$
will be, leading to a larger $O(p^{3}n^{-3/2})$. To understand the
effect this may have for small sample sizes,  
in this section we fit a simple model to a relatively complex process.

\begin{figure}[ht]
\begin{center}
\includegraphics[scale=0.65]{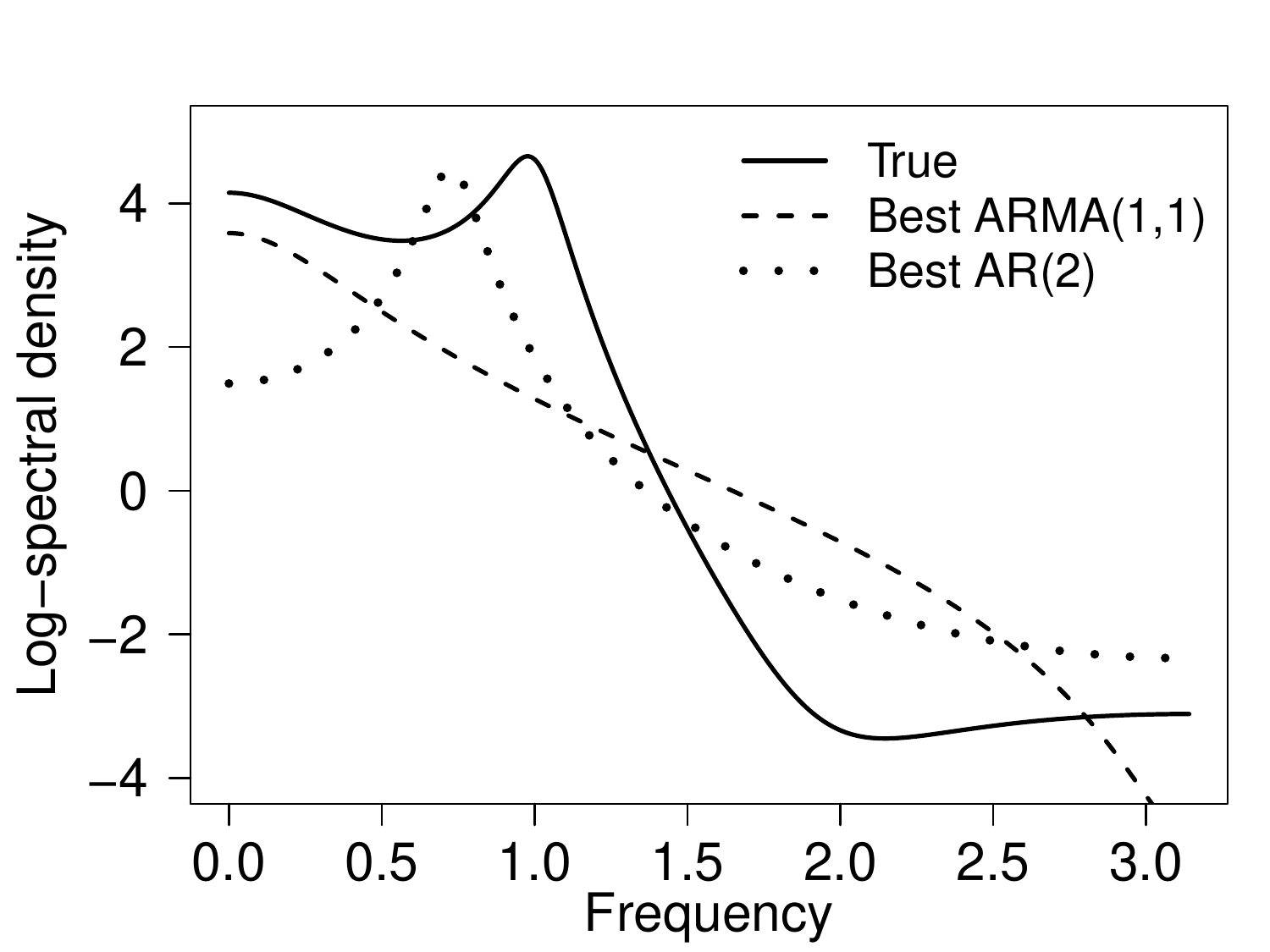}
\caption{\textit{Plot of $\log f_{Z}(\omega)$ and $\log f_{\theta^{Best}}(\omega)$; 
Theoretical ARMA(3,2) spectral density (solid), best fitting ARMA(1,1) spectral density (dashed),
and best fitting AR(2) spectral density (dotted) for $n=20$.
}}
 \label{fig:armaspec} 
\end{center}
\end{figure}

For the ``true'' data generating process we use an ARMA$(3,2)$  Gaussian time series 
with spectral density 
$f_{Z}(\omega) =
|\psi_{Z}(e^{-i\omega})|^2/|\phi_{Z}(e^{-i\omega})|^{2}$, where 
AR and MA characteristic polynomials are
\begin{equation*}
\phi_{Z}(z) = (1-0.7z)(1-0.9e^{i}z)(1-0.9e^{-i}z) \quad {\textrm{and}}\quad 
\psi_{Z}(z) = (1+0.5z+0.5z^2).
\end{equation*}
This spectral density has some
interesting characteristics: a pronounced peak, a large amount of power at the low frequencies,
and a sudden drop in power at the higher frequencies. 
We consider sample sizes $n=20,50$ and $300$, and fit a model with
fewer parameters. Specifically, we
fit two different ARMA models with the same number of unknown parameters.
The first is the ARMA(1,1) model with spectral density
\begin{equation*}
f_{\theta}(\omega) = |1+\psi_{}e^{-i\omega}|^2 |1-\phi_{}e^{-i\omega}|^{-2}
\qquad \theta=(\phi_{}, \psi_{}).
\end{equation*} 
The second is the AR(2) model with spectral density
\begin{equation*}
f_{\theta}(\omega) = |1-\phi_{1}e^{-i\omega} - \phi_{2}e^{-2i \omega}|^{-2}
\qquad \theta=(\phi_{1}, \phi_{2}).
\end{equation*} 
Figure \ref{fig:armaspec} shows the logarithm of the theoretical ARMA(3,2)
spectral density (solid line, $f_{Z}$) and the corresponding log spectral densities
 of the best fitting ARMA(1,1) (dashed line) and AR(2) (dotted line) processes
for $n=20$. The best fitting models are obtained by minimizing the spectral divergence 
$\theta^{Best} = \arg\min_{\theta \in \Theta} I_{n}(f;f_{\theta})$, 
where $I_{n}(f,f_\theta)$ is defined in  (\ref{eq:Ifthetaf}) and
$\Theta$ is the parameter space. 
The best fitting models  for $n=50$ and $300$ are similar.
We observe that neither of the misspecified models capture all of the features of the true spectral density.
The best fitting ARMA(1,1) model has a large amount of power at the
low frequencies and the power declines for the higher frequencies. 
The best fitting AR(2) model peaks around frequency 0.8, but
the power at the low frequencies is small.
Overall, the spectral divergence between the true and the best fitting
AR(2) model is smaller than the
spectral divergence between the true and the best ARMA(1,1) model.

For each simulation, we calculate the six different parameter estimators and the spectral
divergence. The result of the estimators using the six different
quasi-likelihoods is given in Table \ref{tab:arma11} (for ARMA(1,1))
and Table \ref{tab:ar2} 
(for AR(2)). 

\begin{table}[ht]
\centering
\small
\begin{tabular}{cc|ccccccc}

$n$ & Parameter & Gaussian & Whittle & {\color{blue}Boundary} & {\color{blue}Hybrid} & Tapered & Debiased \\ \hline \hline
	
\multirow{3}{*}{20} & $\phi$ & {\color{blue} $0.031(0.1)$} & -$0.095(0.16)$ & -$0.023(0.12)$
&{\color{red} -$0.006(0.1)$} & -$0.080(0.13)$ & $0.187(0.11)$ \\ 

& $\psi$ & {\color{red} $0.069(0.08)$} & -$0.172(0.18)$ & -$0.026(0.14)$ 
& {\color{blue}$0.028(0.1)$} & -$0.068(0.12)$ & $0.093(0.06)$ \\ 

& $I_{n}(f;f_\theta)$ & $1.653(0.81)$ & $1.199(1.57)$ & {\color{blue} $0.945(0.84)$} 
& $1.024(0.89)$ & {\color{red}$0.644(0.61)$} & $2.727(0.73)$\\ \hline

\multirow{3}{*}{50} & $\phi$ & {\color{red}$0.012(0.07)$} & -$0.054(0.09)$ & -$0.006(0.07)$ 
& {\color{blue}$0.004(0.07)$} & -$0.005(0.07)$ & $0.154(0.11)$ \\

& $\psi$ & $0.029(0.06)$ & -$0.116(0.12)$ & -$0.008(0.08)$ 
& {\color{blue} $0.009(0.07)$} & {\color{red} $0.011(0.06)$} & $0.093(0)$ \\ 

& $I_{n}(f;f_\theta)$ & $0.354(0.34)$ & $0.457(0.46)$ & $0.292(0.3)$ 
& {\color{blue}$0.235(0.28)$} & {\color{red} $0.225(0.26)$} & $1.202(0.34)$ \\ \hline

\multirow{3}{*}{300} & $\phi$ & {\color{red} $0.002(0.03)$} & -$0.014(0.03)$ & {\color{red} $0(0.03)$} 
& $0.001(0.03)$ & $0(0.03)$ & $0.093(0.08)$ \\ 

& $\psi$ & {\color{red}$0.005(0.03)$} & -$0.033(0.05)$ & $0.001(0.03)$ 
& {\color{red}$0.003(0.03)$} & $0.003(0.03)$ & $0.092(0.01)$ \\ 

& $I_{n}(f;f_\theta)$ & $0.027(0.05)$ & $0.064(0.09)$ & $0.029(0.05)$
& {\color{red}$0.026(0.04)$} & {\color{blue}$0.027(0.05)$} & $0.752(0.22)$ \\ \hline

\multicolumn{8}{l}{Best fitting ARMA$(1,1)$ coefficients $\theta = (\phi, \psi)$ and spectral divergence:} \\
\multicolumn{8}{l}{~~$-$ $\theta_{20}=(0.693, 0.845)$, $\theta_{50}=(0.694,0.857)$, $\theta_{300}=(0.696,0.857)$. } \\
\multicolumn{8}{l}{~~$-$ $I_{20}(f; f_{\theta}) = 3.773$, $I_{50}(f; f_{\theta}) = 3.415$, $I_{300}(f; f_{\theta}) = 3.388$.} \\
\end{tabular} 
\caption{\textit{The bias of estimated
coefficients for six different estimation methods for the Gaussian ARMA$(3,2)$
misspecified case fitting ARMA$(1,1)$ model. Standard deviations are in the parentheses. We
use {\color{red}red} to denote the smallest RMSE and 
{\color{blue}blue} to denote the second smallest RMSE.}}
\label{tab:arma11}
\end{table}

\begin{table}[ht]
\centering
\footnotesize
\begin{tabular}{cc|ccccccc}

$n$ & Parameter & Gaussian & Whittle & {\color{blue}Boundary} & {\color{blue}Hybrid} & Tapered & Debiased \\ \hline \hline
	
\multirow{3}{*}{20} & $\phi_1$ & {\color{blue}$0.028(0.14)$} & -$0.162(0.22)$ & -$0.032(0.16)$
&{\color{red} $0.003(0.14)$} & -$0.123(0.16)$ & $0.069(0.15)$ \\ 

& $\phi_2$ & {\color{red} -$0.004(0.09)$} & $0.169(0.18)$ & $0.052(0.14)$ 
& {\color{blue}$0.025(0.12)$} & $0.132(0.12)$ & -$0.034(0.11)$ \\ 

& $I_{n}(f;f_\theta)$ & {\color{red}$0.679(0.72)$} & $1.203(1.46)$ & $0.751(0.85)$ 
& {\color{blue}$0.684(0.8)$} & $0.862(0.97)$ & $0.686(0.81)$\\ \hline

\multirow{3}{*}{50} & $\phi_1$ & {\color{red}$0.019(0.09)$} & -$0.077(0.12)$ & -$0.009(0.09)$ 
& {\color{red}$0.003(0.09)$} & {\color{red}-$0.017(0.09)$} & $0.156(0.15)$ \\

& $\phi_2$ & -$0.024(0.06)$ & $0.066(0.1)$ & $0.006(0.07)$ 
& {\color{red} -$0.003(0.06)$} & {\color{blue}$0.013(0.06)$} & -$0.121(0.06)$ \\ 

& $I_{n}(f;f_\theta)$ & {\color{red}$0.275(0.33)$} & $0.382(0.45)$ & $0.283(0.37)$ 
& $0.283(0.37)$ & {\color{blue}$0.283(0.36)$} & $0.65(0.7)$ \\ \hline

\multirow{3}{*}{300} & $\phi_1$ & {\color{red} $0.004(0.04)$} & -$0.013(0.04)$ & {\color{red}$0(0.04)$}
& $0.001(0.04)$ & $0.001(0.04)$ & $0.014(0.04)$ \\ 

& $\phi_2$ & {\color{blue}-$0.005(0.02)$} & $0.011(0.03)$ & {\color{red}-$0.001(0.02)$} 
& {\color{blue}-$0.001(0.03)$} & -$0.001(0.03)$ & $0.016(0.04)$ \\ 

& $I_{n}(f;f_\theta)$ & {\color{red}$0.049(0.07)$} & $0.053(0.07)$ & {\color{blue} $0.049(0.07)$}
& $0.053(0.07)$ & $0.054(0.08)$ & $0.058(0.08)$ \\ \hline

\multicolumn{8}{l}{Best fitting AR$(1)$ coefficients $\theta = (\phi_1, \phi_2)$ and spectral divergence:} \\
\multicolumn{8}{l}{~~$-$ $\theta_{20}=(1.367, -0.841)$, $\theta_{50}=(1.364,-0.803)$, $\theta_{300}=(1.365,-0.802)$. } \\
\multicolumn{8}{l}{~~$-$ $I_{20}(f; f_{\theta}) = 2.902$, $I_{50}(f; f_{\theta}) = 2.937$, $I_{300}(f; f_{\theta}) = 2.916$.} \\

\end{tabular} 
\caption{\textit{Same as in Table \ref{tab:arma11} but fitting AR(2) model.
}}
\label{tab:ar2}
\end{table}

We first discuss the parameter estimates. 
Comparing the asymptotic bias of the Gaussian likelihood with the boundary
corrected Whittle likelihood (see Appendix \ref{sec:multi}),
the Gaussian likelihood has an additional
bias term of form $\sum_{r=1}^{d} I^{(j,r)}\Ex[ \frac{\partial
  \mathcal{L}_n}{\partial \theta_{r}}]\rfloor_{\theta_{n}}$. But there
is no guarantee that the inclusion of this term
increases or decreases the bias. This is borne out in the simulations,
where we observe that overall the Gaussian likelihood or the new
likelihoods tend to have a smaller parameter bias (there is no
clear winner). The tapered
likelihood is a close contender, performing very well for the moderate
sample sizes $n=50$. 
Similarly, in terms of the RMSE, again there is no clear winner between
the Gaussian and the new likelihoods. Overall (in the simulations) the hybrid Whittle
likelihood tends to outperform the Gaussian likelihood. 

We next turn our attention to the estimated spectral divergence
$I_{n}(f,f_{\widehat{\theta}})$. For the fitted ARMA$(1,1)$ model, the
estimated spectral divergence of the new likelihood estimators tends to
be the smallest or second smallest in terms of the RMSE (its nearest competitor is the
tapered likelihood). On the other hand, for the AR$(2)$ model the
spectral divergence of Gaussian likelihood has the smallest RMSE for
all the sample sizes. The new likelihood comes in second for
sample sizes $n=20$ and $300$.

In the simulations above we select $p$ using the AIC. As mention at
the start of the section, this leads to an additional error of
$O(p^{3}n^{-3/2})$ in the new likelihoods. Thus, if a large $p$ is
selected the error $O(p^{3}n^{-3/2})$ will be large. In order to
understand the impact $p$ has on the estimator, in Appendix
\ref{sec:fixedP} we compare the the likelihoods constructed using the
predictive DFT based on the AIC with the likelihoods constructed using
the predictive DFT based on the best fitting estimated AR$(1)$ model. We simulate from the
ARMA$(3,2)$ model described above and fit an ARMA$(1,1)$ and AR$(2)$ model.
As is expected, the bias tends to be a little larger when the
order is fixed to $p=1$. But even when fixing $p=1$, we do observe an improvement over
the Whittle likelihood (in some cases an improvement over the Gaussian
likelihood).


%% file: appendix_main_proof.tex
\section*{Summary of results in the Supplementary material} 

To navigate the supplementary material, we briefly
summarize the contents of each section.

\begin{enumerate}
\item In Appendix \ref{sec:proofbio} we prove the results stated in
  Section  2 (which concern representing the Gaussian likelihood in the
  frequency domain and obtaining an explicit expression for the
  predictive DFT for AR models). Some of these proofs will use results
  from Appendix \ref{sec:sup}. 
\item In Appendix  \ref{sec:approxproofs} we prove both the first
  order approximation and higher order approximation 
  results stated in Sections 3. The proof of Theorem
  \ref{theorem:approx}, uses the extended Baxter's inequality 
(for completeness we prove this result in Appendix \ref{sec:baxter1},
though we believe the result is well known). The proof Theorem
\ref{thm:higherorder} uses an explicit expression for finite predictors.
\item In Appendix \ref{sec:approxproofs2} we prove Lemma \ref{lemma:12}.
The proof is similar to the proof of Theorem \ref{theorem:approx} but
with some subtle differences. 
\item Appendix \ref{sec:baxter} mainly deals with Baxter-type
  inequalities. In Appendix \ref{sec:baxter} we give a proof of the
  extended Baxter inequality. In Appendix  \ref{sec:baxterderivatives}
  we obtain Baxter type bounds for the derivatives of the finite
  predictors is stated (with respect to the unknown parameter). 
These results are used in Appendix \ref{sec:prelim} to bound the
difference between the derivatives of the Gaussian and Whittle
likelihood.
\item In Appendix \ref{sec:consistent} we show consistency of the new
  likelihood estimators (defined in Section
  \ref{sec:newlikelihood}). We also state and prove the necessary
  lemmas for proving the asymptotic equivalence between the feasible
  and infeasible estimators (proof of Theorem \ref{theorem:equivalence}).
\item In Appendix \ref{sec:bias} we obtain the bias of the  Gaussian, 
Whittle,  boundary corrected and hybrid Whittle likelihoods under quite general assumptions on the underlying
time series $\{X_{t}\}$. In particular, in Appendix
\ref{sec:oneparameter} we state the results in the one-parameter case
(the proof is given in Appendix \ref{sec:proofbias}). The results for
the special case of the AR$(1)$ is given in Appendix
\ref{sec:AR1bias}. The bias for multi-parameters is given in Appendix \ref{sec:multi}.
\item In Appendix \ref{sec:sim_appendix}, we supplement the
  simulations in Section \ref{sec:sim}. 
\item In Appendix \ref{sec:sim-long} we apply the new likelihood
  estimation methods to
  long memory time series. We consider both parametric methods and
  also variants on the local Whittle likelihood estimator for the long
  memory parameter. In Appendix \ref{sec:alternative} we construct
  alternative estimators of the predictive DFT, which are used to
  build the new likelihoods. Through simulations we compare the
  estimators based on these different new likelihoods. 
\end{enumerate}


\setcounter{section}{0}
\section{Proof of Theorems \ref{theorem:bio}, \ref{lemma:ARp1}
and \ref{theorem:longmemory}} \label{sec:proofbio}

We start this section by giving the proof of Theorems
\ref{theorem:bio}. This result is instrumental to the subsequent
results in the paper.

\vspace{2mm}
\noindent{\bf PROOF of Theorem \ref{theorem:bio}}
First, to prove Theorem \ref{theorem:bio}, we recall that it entails obtaining a transform 
$U_{n}\underline{X}_{n}$ where 
$\cov_{f}\left(U_{n}\underline{X}_{n}, F_n\underline{X}_n\right) =
\Delta_{n}(f_{})$. Pre and post multiplying this covariance with $F_{n}^{*}$ and $F_{n}$ gives 
\begin{equation*}
F_{n}^{*}\cov_{f}\left(U_{n}\underline{X}_{n},
F_n\underline{X}_n\right)F_{n} = \cov_{\theta}\left(F_{n}^{*}U_{n}\underline{X}_{n},
\underline{X}_n\right)=
F_{n}^{*}\Delta_{n}(f_{})F_{n} = C_{n}(f_{}).
\end{equation*}
Thus our
objective is to find the transform $\underline{Y}_{n} = F_{n}^{*}U_{n}\underline{X}_{n}$ such that
$\cov_{f}(\underline{Y}_{n},
\underline{X}_{n})=C_{n}(f)$.
Then, the vector $F_{n}\underline{Y}_n =U_{n}\underline{X}_{n}$ will be 
biorthogonal to $F_{n}\underline{X}_n$, as required.
We observe that the entries of the circulant matrix $C_{n}(f)$
are 
\begin{equation*}
(C_{n}(f))_{u,v} = n^{-1}\sum_{k=1}^{n}f(\omega_{k,n})\exp(-i(u-v)\omega_{k,n})
= \sum_{\ell \in \mathbb{Z}}c_{f}(u-v+\ell n),
\end{equation*}
where the second equality is due to the 
Poisson summation. The random vector 
$\underline{Y}_{n}=\{Y_{u,n}\}_{u=1}^{n}$ is such that $\cov_{f}(Y_{u,n}, X_{v}) =
\sum_{\ell \in
\mathbb{Z}}c_{f}(u-v+\ell n)$ and $Y_{u} \in \spa(\Xunder_{n})$.
Since $\cov_{f}(X_{u+\ell n},X_{v}) = c_{f}(u-v+\ell n)$, at least ``formally'' 
$\cov_{f}(\sum_{\ell \in \mathbb{Z}}X_{u+\ell n},X_{v}) =
\sum_{\ell \in
\mathbb{Z}}c_{f}(u-v+\ell n)$. However, $\sum_{\ell\in
\mathbb{Z}}X_{u+\ell n}$ is neither a well defined random
variable nor does not it belong to $\spa(\Xunder_{n})$. 
We replace each element in the sum $\sum_{\ell \in
\mathbb{Z}}X_{u+\ell n}$ with an element that belongs to $\spa(\Xunder_{n})$ and gives the
same covariance. 
To do this we use the following well known result. Let
$Z$ and $\underline{X}$ denote a random variable and vector
respectively. Let $P_{\underline{X}}(Z)$ denote the projection of $Z$ onto
$\spa(\underline{X})$, i.e., the best linear
predictor of $Z$ given $\underline{X}$, then $\cov_{f}(Z,\underline{X}) =
\cov_{f}(P_{\underline{X}}(Z), \underline{X})$. Let $\widehat{X}_{\tau,n}$
denote best linear predictor of $X_{\tau}$ given
$\Xunder_{n}=(X_{1},\ldots,X_{n})$ (as defined in 
(\ref{eq:Xtaun})). $\widehat{X}_{\tau,n}$ retains the pertinent properties of
$X_{\tau}$ in the sense that $\cov_{f}(\widehat{X}_{\tau,n},X_{t})=c_{f}(\tau-t)$ for
all $\tau\in \mathbb{Z}$ and $1\leq t\leq n$. Define
\begin{equation*}
Y_{u,n} =\sum_{\ell\in
\mathbb{Z}}\widehat{X}_{u+\ell n,n} =
\sum_{s=1}^{n}\left(\sum_{\ell \in \mathbb{Z}}\phi_{s,n}(u+\ell n;f)\right)X_{s} \in \spa(\Xunder_{n}),
\end{equation*}
where we note that $Y_{u,n}$ a well defined random variable, since by
using Lemma \ref{lemma:BaxterEx} it can be shown that
$\sup_{n}\sum_{s=1}^{n}\sum_{\ell =
-\infty}^{\infty}|\phi_{s,n}(u+\ell n;f)|<\infty$.
Thus by definition of $Y_{u,n}$ the following holds 
\begin{equation}
\label{eq:XtYtau}
\cov_{f}\left( Y_{u,n}, X_{v} \right)=
\sum_{\ell \in \mathbb{Z}}c_{f}(u-v+\ell n) =
\left( C_{n}(f)\right)_{u,v},
\end{equation}
and $\underline{Y}_{n}=F_{n}^{*}U_{n}\underline{X}_{n}$, gives the
desired transformation of the time series. 
Thus, based on this construction, $F_{n}\underline{Y}_{n} = U_{n}\underline{X}_{n}$ and
$F_{n}\underline{X}_{n}$ are biorthogonal transforms,
with entries
$(F_{n}\underline{X}_{n})_{k} = J_n(\omega_{k,n})$ and 
\begin{eqnarray}
(U_{n}\underline{X}_{n})_{k} = (F_{n}\underline{Y}_{n})_{k}
&=& n^{-1/2}\sum_{\ell \in \mathbb{Z}}\sum_{u=1}^{n}\widehat{X}_{u+\ell
n,n}e^{iu\omega_{k,n}} \nonumber \\
&=& n^{-1/2}\sum_{\tau \in
\mathbb{Z}}\widehat{X}_{\tau,n}e^{i\tau\omega_{k,n}} \nonumber\\
&=& n^{-1/2}\sum_{t=1}^{n}X_{t}\sum_{\tau\in \mathbb{Z}}\phi_{t,n}(\tau;f)e^{i\tau\omega_{k,n}}.
\label{eq:FYUX}
\end{eqnarray}
The entries of the matrix $U_{n}$ are $(U_{n})_{k,t} = n^{-1/2}\sum_{\tau\in \mathbb{Z}}\phi_{t,n}(\tau;f) e^{i\tau\omega_{k,n}}$.
To show that $U_{n}$ ``embeds'' the regular DFT, we observe that for $1\leq \tau \leq n$,
$\phi_{t,n}(\tau;f) = \delta_{\tau,t}$, furthermore, due to second order stationarity
the coefficients $\phi_{t,n}(\tau;f)$ are reflective i.e.
the predictors of $X_{m}$ (for $m > n$) and
$X_{n+1-m}$ share the same set of prediction coefficients (just reflected) such that 
\begin{equation}
\label{eq:reflective}
\phi_{t,n}(m;f) = \phi_{n+1-t,n}(n+1-m;f) \qquad \textrm{ for }m > n.
\end{equation}
Using these two observations we can decompose
$(U_{n})_{k,t}$ as 
\begin{eqnarray*}
(U_{n})_{k,t} &=& 
n^{-1/2} \left( e^{it\omega_{k,n}} + 
\sum_{\tau\leq 0}\phi_{t,n}(\tau;f) e^{i\tau\omega_{k,n}} + 
\sum_{\tau\geq n+1}\phi_{t,n}(\tau;f)e^{i\tau\omega_{k,n}} \right) \\
&=&n^{-1/2} e^{it\omega_{k,n}} + n^{-1/2}
\sum_{\tau\leq 0}\left( \phi_{t,n}(\tau;f)e^{i\tau\omega_{k,n}} + 
\phi_{n+1-t,n}(\tau;f)e^{-i(\tau-1-n)\omega_{k,n}}\right).
\end{eqnarray*}
It immediately follows from the above decomposition that $U_{n} =
F_{n} + D_{n}(f)$ where $D_{n}(f)$ is defined in
(\ref{eq:Dftheta}). Thus proving (\ref{eq:UX}). 

To prove (\ref{eq:DFTPred}), we first observe that (\ref{eq:UX})
implies
\begin{equation*}
\cov_{f}\left( ((F_{n}+D_{n}(f))\underline{X}_{n})_{k_1}, (F_{n}\underline{X}_{n})_{k_2}\right)
= f_{}(\omega_{k_1,n})\delta_{k_1,k_2}.
\end{equation*}
It is clear that $(F_{n}\underline{X}_{n})_{k}=J_{n}(\omega_{k,n })$
and from the representation of $F_{n}\underline{Y}_{n}$
given in (\ref{eq:FYUX}) we have 
\begin{eqnarray*}
(F_{n}\underline{Y}_{n})_{k} 
&=& n^{-1/2}\sum_{\tau =1}^{n}X_{\tau}e^{i\tau\omega_{k,n}} +
n^{-1/2}\sum_{\tau \notin \{1, ..., n\}}\widehat{X}_{\tau,n}e^{i\tau\omega_{k,n}} \nonumber\\
&=& J_{n}(\omega_{k,n}) + \widehat{J}_{n} (\omega_{k,n};f_{}).
\end{eqnarray*}
This immediately proves (\ref{eq:DFTPred}). 
\hfill $\Box$

\vspace{2mm} 
\noindent Note that equation (\ref{eq:DFTPred}) can be verified
directly by using the properties of linear predictors and covariances
discussed in the above proof. 

\vspace{1em}

\vspace{3mm}
\noindent To prove Theorem \ref{lemma:ARp1} we study the predictive DFT
for autoregressive processes. We start by obtaining an explicit expression for
$\widehat{J}_{n}(\omega;f_{\theta})$ where $f_{\theta}(\omega) =
\sigma^{2}|1-\sum_{u=1}^{p}\phi_{u}e^{-iu\omega}|^{-2}$ 
(the spectral density
corresponding to an AR$(p)$ process). It is straightforward to show 
that predictive DFT predictor based on the AR$(1)$ model is
\begin{eqnarray*}
\widehat{J}_{n}(\omega;f_{\theta})
& =&
n^{-1/2}\sum_{\tau=-\infty}^{0}\phi^{-\tau+1}X_{1}e^{i\tau\omega}+
n^{-1/2}\sum_{\tau=n+1}^{\infty}\phi^{\tau+1-n}X_{n}e^{i\tau\omega}\\
&=& \frac{n^{-1/2}\phi}{ \phi_{1}(\omega)}X_{1} + \frac{n^{-1/2} \phi}{ \overline{\phi_{1}(\omega)}}X_{n}e^{i(n+1)\omega},
\end{eqnarray*}
where $\phi_{1}(\omega) = 1- \phi e^{-i\omega}$. In order to prove Theorem
\ref{lemma:ARp1}, which generalizes the above expression to AR$(p)$
processes, we partition $\widehat{J}_{n}(\omega;f_{\theta})$ into the
predictions involving the past and 
future terms
\begin{equation*}
\widehat{J}_{n}(\omega;f_{\theta}) = \widehat{J}_{n,L}(\omega;f_{\theta})+\widehat{J}_{n,R}(\omega;f_{\theta})
\end{equation*}
where 
\begin{equation*}
\widehat{J}_{n,L}(\omega;f_{\theta})= n^{-1/2}\sum_{\tau=-\infty}^{0}
\widehat{X}_{\tau, n} e^{i \tau \omega} \quad\textrm{and}\quad
\widehat{J}_{n,R}(\omega;f_{\theta}) = 
n^{-1/2}\sum_{\tau=n+1}^{\infty} \widehat{X}_{\tau, n} e^{i \tau \omega}.
\end{equation*}
We now obtain expressions for $\widehat{J}_{n,L}(\omega;f_{\theta})$
and $\widehat{J}_{n,R}(\omega;f_{\theta})$ separately, in the case the
predictors are based on the AR$(p)$ parameters where 
$f_{\theta}(\omega)=\sigma^{2}|1-\sum_{j=1}^{p}\phi_{j}e^{ij\omega}|^{-2}$
and the $\{\phi_{j}\}_{j=1}^{p}$ correspond to the causal AR$(p)$ representation. To do so, 
we define the $p$-dimension vector $\underline{\phi}^{\prime} = (\phi_{1},\ldots,\phi_{p})$
and the matrix $A_{p}(\underline{\phi})$ as 
\begin{eqnarray}
\label{eq:ARP}
A_{p}(\underline{\phi}) = \left(
\begin{array}{ccccc}
\phi_{1} &\phi_{2} & \ldots & \phi_{p-1} & \phi_{p} \\
1 & 0 & \ldots & 0 & 0 \\
0 & 1 & \ldots & 0 & 0 \\
\vdots & \vdots & \ddots & 0 & 0 \\
0 & 0 & \ldots & 1 & 0 \\ 
\end{array}
\right).
\end{eqnarray}
Therefore, for $\tau \leq 0$, since $\widehat{X}_{\tau,n} =
\left[A_{p}(\underline{\phi})^{|\tau|+1}\underline{X}_{p}\right]_{(1)}$,
where $\underline{X}_{p} = (X_{1},\ldots,X_{p})$, we can write
\begin{equation}
\label{eq:JhatP}
\widehat{J}_{n,L}(\omega;f_\theta) = n^{-1/2} 
\sum_{\tau=-\infty}^{0}\left[A_{p}(\underline{\phi})^{|\tau|+1}\underline{X}_{p}\right]_{(1)}e^{i\tau\omega}. 
\end{equation}

\begin{lemma}\label{lemma:DFTinfty}
Let $\widehat{J}_{n,L}(\omega)$ be defined as in (\ref{eq:JhatP}),
where the parameters $\underline{\phi}$ are such that the roots of
$\phi(z) = 1 - \sum_{j=1}^{p}\phi_{j}z^{j}$ lie outside the unit
circle. Then an analytic
expression for $\widehat{J}_{n,L}(\omega;f_\theta)$ is 
\begin{equation}
\label{eq:JPPP}
\widehat{J}_{n,L}(\omega;f_\theta) =
\frac{n^{-1/2}}{\phi_{p}(\omega)}\sum_{\ell=1}^{p}X_{\ell}\sum_{s=0}^{p-\ell}\phi_{\ell+s}e^{-is\omega}.
\end{equation}
where $\phi_{p}(\omega) = 1-\sum_{s=1}^{p}\phi_{s}e^{-is\omega}$.
\end{lemma}
PROOF. By using (\ref{eq:APexpand}) we have
\begin{equation*}
[A_{p}(\underline{\phi})^{|\tau|+1}\underline{X}_{p}]_{(1)} = \sum_{\ell=1}^{p}X_{\ell}\sum_{s=0}^{p-\ell}\phi_{\ell+s}\psi_{|\tau|-s}. 
\end{equation*}
Therefore, using (\ref{eq:JhatP}) and the change of variables $\tau \leftarrow -\tau$
\begin{eqnarray*}
\widehat{J}_{n,L}(\omega;f_\theta) 
&=& n^{-1/2} \sum_{\ell=1}^{p}X_{\ell}\sum_{s=0}^{p-\ell}\phi_{\ell+s}\sum_{\tau=0}^{\infty}\psi_{\tau-s}e^{-i\tau
\omega} \\
&=&n^{-1/2} \sum_{\ell=1}^{p}X_{\ell}\sum_{s=0}^{p-\ell}\phi_{\ell+s}e^{-is\omega}\sum_{\tau=0}^{\infty}\psi_{\tau-s}e^{-i(\tau-s)
\omega} \\
&=& n^{-1/2} \sum_{\ell=1}^{p}X_{\ell}\sum_{s=0}^{p-\ell}\phi_{\ell+s}e^{-is\omega}\sum_{\tau=s}^{\infty}\psi_{\tau-s}e^{-i(\tau-s) \omega}. 
\end{eqnarray*}
Let $\sum_{s=0}^{\infty}\psi_{s}e^{-is\omega} = \psi(\omega) = \phi_{p}(\omega)^{-1}$, and
substitute this into the above to give 
\begin{equation}
\label{eq:hatJtilde2}
\widehat{J}_{n,L}(\omega;f_\theta)= \frac{n^{-1/2}}{\phi_{p}(\omega)}\sum_{\ell=1}^{p}X_{\ell}\sum_{s=0}^{p-\ell}\phi_{\ell+s}e^{-is\omega},
\end{equation}
Thus we obtain the desired result. \hfill $\Box$

\vspace{1em}

\noindent {\bf PROOF of Theorem \ref{lemma:ARp1}} 
To prove (\ref{eq:JAR}), we note that 
the same proof as that above can be
used to prove that the right hand side predictive DFT
$\widehat{J}_{n,R}(\omega;f_\theta)$ has the representation
\begin{equation*}
\widehat{J}_{n,R}(\omega;f_\theta) = e^{in\omega} \frac{n^{-1/2}}{
\overline{\phi_{p}(\omega)}} \sum_{\ell=1}^{p}X_{n+1-\ell}\sum_{s=0}^{p-\ell}
\phi_{\ell+s}e^{i(s+1)\omega}.
\end{equation*}
Since $\widehat{J}_{n}(\omega;f_{\theta}) = \widehat{J}_{n,L}(\omega;f_{\theta})+\widehat{J}_{n,R}(\omega;f_{\theta})$,
Lemma \ref{lemma:DFTinfty} and the above give an explicit
expression for $\widehat{J}_{n}(\omega;f_{\theta})$, thus proving
equation (\ref{eq:JAR}). 

To prove (\ref{eq:DDD}) we use that
\begin{equation*}
(\widehat{J}_{n}(\omega_{1,n};f_\theta),\ldots,\widehat{J}_{n}(\omega_{n,n};f_\theta))^{\prime}
= D_{n}(f_{\theta})\underline{X}_{n}.
\end{equation*}
Now by using (\ref{eq:JAR}) together with the above we immediately
obtain (\ref{eq:DDD}). 

Finally, we prove (\ref{eq:DDDmatrix}). We use the result
$n^{-1}\sum_{k=1}^{n}\phi_{p}(\omega_{k,n})\exp(is\omega_{k,n}) =
\widetilde{\phi}_{s\bmod n}$ where $\phi_{p}(\omega) =
\sum_{r=0}^{n-1}\widetilde{\phi}_{r}e^{-ir\omega}$ and
$\widetilde{\phi}_{r}=0$ for $p+1\leq r \leq n$. 
For $1\leq t \leq p$ we use have 
\begin{eqnarray*}
(F_{n}^{*}\Delta_{n}(f_{\theta}^{-1})D_{n}(f_\theta) )_{s,t} 
&=& \frac{1}{n}\sum_{k=1}^{n}\frac{\phi_{t,p}(\omega_{k,n})
}{f_{\theta}(\omega_{k,n})}\exp(-is\omega_{k,n})\\
&=& \frac{\sigma^{-2}}{n}\sum_{k=1}^{n}\overline{\phi_{p}(\omega_{k,n})}
\sum_{\ell = 0}^{p-t}\phi_{\ell+t}\exp(-i\ell\omega_{k,n})
\exp(-is\omega_{k,n})\\
&=& \sigma^{-2} \sum_{\ell = 0}^{p-t}\phi_{\ell+t}\frac{1}{n}\sum_{k=1}^{n}
\overline{\phi_{p}(\omega_{k,n})}
\exp(-i(\ell+s)\omega_{k,n}) \\
&=& \sigma^{-2} \sum_{\ell = 0}^{p-t}\phi_{\ell+t} \widetilde{\phi}_{(\ell + s)\bmod n}.
\end{eqnarray*}
Similarly, for $1\leq t \leq p$,
\begin{eqnarray*}
(F_{n}^{*}\Delta_{n}(f_{\theta}^{-1})D_{n}(f_\theta) )_{s,n-t+1} 
&=& \frac{1}{n}\sum_{k=1}^{n}\frac{\overline{\phi_{t,p}(\omega_{k,n})}
}{f_{\theta}(\omega_{k,n})}\exp(i(1-s)\omega_{k,n})\\
&=& \frac{\sigma^{-2}}{n}\sum_{k=1}^{n}\phi_{p}(\omega_{k,n})
\sum_{\ell = 0}^{p-t}\phi_{\ell+t}\exp(i\ell\omega_{k,n})
\exp(i(1-s)\omega_{k,n})\\
&=& \sigma^{-2} \sum_{\ell = 0}^{p-t}\phi_{\ell+t} \frac{1}{n}\sum_{k=1}^{n}
\phi_{p}(\omega_{k,n})
\exp(i(\ell+1-s)\omega_{k,n}) \\
&=& 
\sigma^{-2} \sum_{\ell = 0}^{p-t}\phi_{\ell+t} \widetilde{\phi}_{(\ell +1 -s)\bmod n}.
\end{eqnarray*}
\hfill $\Box$

\vspace{3mm}

\noindent {\bf PROOF of Equation (\ref{eq:LKdiff})} We use that 
$\frac{1}{\phi_{p}(\omega)}f_{\theta}(\omega)^{-1} =
\sigma^{-2} \overline{\phi_{p}(\omega)}$. This gives 
\begin{eqnarray*}
\mathcal{L}_{n}(\phi) - K_{n}(\phi) = I + II
\end{eqnarray*}
where 
\begin{eqnarray*}
I &=&
\frac{1}{n^{3/2}}\sum_{k=1}^{n}\frac{\overline{J_{n}(\omega_{k,n})}}{f_{\theta}(\omega_{k,n})}\left\{
\frac{1}{\phi_{p}(\omega_{k,n})}\sum_{\ell=1}^{p}X_{\ell}\sum_{s=0}^{p-\ell}\phi_{\ell+s} e^{-is\omega_{k,n}}
\right\}\\
&=& \frac{\sigma^{-2}}{n^{3/2}}\sum_{k=1}^{n}\overline{J_{n}(\omega_{k,n})}\overline{\phi_{p}(\omega_{k,n})}\sum_{\ell=1}^{n}\sum_{s=0}^{p-\ell}\phi_{s+k} e^{-is\omega_{k,n}}
\phi_{\ell+s} e^{-is\omega_{k,n}} \\
&=& \frac{\sigma^{-2}}{n}\sum_{\ell=1}^{p}X_{\ell}\sum_{s=0}^{p-\ell}\phi_{s+k} \frac{1}{n^{1/2}}
\sum_{k=1}^{n}\overline{J_{n}(\omega_{k,n})}\overline{\phi_{p}(\omega_{k,n})}
e^{-is\omega_{k,n}}
\end{eqnarray*}
and 
\begin{equation}
II = \frac{\sigma^{-2}}{n^{3/2}}\sum_{k=1}^{n}\frac{\overline{J_{n}(\omega_{k,n})}}{f_{\theta}(\omega_{k,n})}\frac{1}{\overline{\phi_{p}(\omega_{k,n})}} 
\sum_{\ell=1}^{p}X_{n+1-\ell}\sum_{s=0}^{p-\ell}
\phi_{\ell+s}e^{i(s+1)\omega_{k,n}}. 
\end{equation}
We first consider $I$.
Using that $\overline{\phi_{p}(\omega_{k,n})} = 1 -
\sum_{j=1}^{p}\phi_{j}e^{ij\omega_{k,n}}$ and
$n^{-1/2}\sum_{k=1}^{n}\overline{J_{n}(\omega_{k,n})}e^{is\omega_{k,n}}
= X_{s\bmod n}$, gives 
\begin{eqnarray*}
I &=& -\frac{\sigma^{-2}}{n}\sum_{\ell=1}^{p}X_{\ell}\sum_{s=0}^{p-\ell}\sum_{j=0}^{p}\phi_{j}\phi_{s+\ell} \frac{1}{n^{1/2}}
\sum_{k=1}^{n}\overline{J_{n}(\omega_{k,n})}e^{-i(s-r)\omega_{k,n}}~ (\textrm{set }\phi_{0}=-1) \\
&=&
-\frac{\sigma^{-2}}{n}\sum_{\ell=1}^{p}X_{\ell}\sum_{s=0}^{p-\ell}\sum_{j=0}^{p}\phi_{j}\phi_{s+\ell} X_{-(s-j)\bmod n}.
\end{eqnarray*}
The proof of $II$ is similar. Altogether this proves the
result. \hfill $\Box$

\vspace{3mm}

Finally, in this section we prove Theorem \ref{theorem:longmemory}
(from Section \ref{sec:conc}). This result is analogous to  Theorem \ref{theorem:bio}.

\noindent{\bf PROOF of Theorem \ref{theorem:longmemory}} Under the
stated condition that $\widetilde{D}_{n}$ is a finite
matrix, then $(\widetilde{F}_{n}+\widetilde{D}_{n})\underline{X}_{n}$
is a well defined random variable. Thus 
\begin{eqnarray} 
\label{eq:biortho2}
\cov((\widetilde{F}_{n}+\widetilde{D}_{n})\underline{X}_{n}, \widetilde{F}_{n}\underline{X}_{n}) 
=(\widetilde{F}_{n}+\widetilde{D}_{n}) \Gamma_{n} \widetilde{F}_{n}^{*} = \widetilde{\Delta}_{n}.
\end{eqnarray} 
Therefore, using (\ref{eq:biortho2}) and  $F_{n}^{*} = [\widetilde{F}_n^{*}, \overline{e_n}] = [\widetilde{F}_n^{*}, e_{n}]$ gives
\begin{eqnarray*}
(\widetilde{F}_{n}+\widetilde{D}_{n}) \Gamma_{n} F_{n}^{*}  
=(\widetilde{F}_{n}+\widetilde{D}_{n}) \Gamma_{n} 
[\widetilde{F}_{n}^{*}, e_n] 
&=& [(\widetilde{F}_{n}+\widetilde{D}_{n}) \Gamma_{n} \widetilde{F}_{n}^{*}   , (\widetilde{F}_{n}+\widetilde{D}_{n}) \Gamma_{n}  e_n] \\
&=& [\widetilde{\Delta}_{n}, (\widetilde{F}_{n}+\widetilde{D}_{n}) \Gamma_{n}  e_n].
\end{eqnarray*} 
Note that for $1\leq k \leq n-1$,
\begin{eqnarray*}
(e_k + d_k)^{\prime} \Gamma_n e_{n} &=& \cov( \widetilde{J}_n(\omega_{k,n};f), J_n(0) ) = 0.
\end{eqnarray*} 
Therefore, we have
\begin{eqnarray*}
(\widetilde{F}_{n}+\widetilde{D}_{n}) \Gamma_{n} F_{n}^{*}   = [\widetilde{\Delta}_{n}, \textbf{0}_{n-1}].
\end{eqnarray*}  
where $\textbf{0}_{n-1} = (0,...,0)^{\prime}$ is a $(n-1)$-dimension
zero vector. 
Right multipling the above with $F_{n} \Gamma_n^{-1}$ gives 
\begin{eqnarray*}
(\widetilde{F}_{n}+\widetilde{D}_{n}) = [\widetilde{\Delta}_{n}, \textbf{0}_{n-1}] F_{n} \Gamma_n^{-1}.
\end{eqnarray*} 
Left multiplying the above with $[\widetilde{\Delta}_{n}^{-1}, \textbf{0}_{n-1}]^{T}$ gives
\begin{eqnarray*}
\begin{bmatrix}
\widetilde{\Delta}_{n}^{-1} \\
\textbf{0}_{n-1}^{\prime} 
\end{bmatrix}
(\widetilde{F}_{n}+\widetilde{D}_{n}) = 
\begin{bmatrix}
\widetilde{\Delta}_{n}^{-1} \\
\textbf{0}_{n-1}^{\prime} 
\end{bmatrix}
[\widetilde{\Delta}_{n}, \textbf{0}_{n-1}] F_{n} \Gamma_n^{-1}
&=& 
\begin{bmatrix}
I_{n-1} & \textbf{0}_{n-1} \\
\textbf{0}_{n-1}^{\prime} & 0 
\end{bmatrix} F_{n} \Gamma_n^{-1} \\
&=& \begin{bmatrix}
I_{n-1} & \textbf{0}_{n-1} \\
\textbf{0}_{n-1}^{\prime} & 0 
\end{bmatrix} 
\begin{bmatrix}
\widetilde{F}_{n}  \\
e_n^{\prime} 
\end{bmatrix} \Gamma_n^{-1} \\
&=& 
\begin{bmatrix}
\widetilde{F}_{n}  \\
\textbf{0}_{n-1}^{\prime} 
\end{bmatrix} \Gamma_n^{-1}
\end{eqnarray*} 
Left multiplying the above with $[\widetilde{F}_n^{*},
\textbf{0}_{n-1}]$ gives 
\begin{eqnarray}\label{eq:matrixeq1} 
\widetilde{F}_n^{*} \widetilde{\Delta}_{n}^{-1}(\widetilde{F}_{n}+\widetilde{D}_{n})
=
[\widetilde{F}_n^{*}, \textbf{0}_{n-1}]
\begin{bmatrix}
\widetilde{\Delta}_{n}^{-1} \\
\textbf{0}_{n-1}^{\prime} 
\end{bmatrix}
(\widetilde{F}_{n}+\widetilde{D}_{n}) 
=
 [\widetilde{F}_n^{*}, \textbf{0}_{n-1}] 
\begin{bmatrix}
\widetilde{F}_{n}  \\
\textbf{0}_{n-1}^{\prime} 
\end{bmatrix} \Gamma_n^{-1}
= \widetilde{F}_n^{*} \widetilde{F}_n \Gamma_n^{-1}.
\end{eqnarray} Finally, using that $e_n = n^{-1/2} \textbf{1}_n$
\begin{eqnarray*}
I_{n} = F_{n}^{*} F_{n} = [\widetilde{F}_n^{*}, n^{-1/2} \textbf{1}_n] 
\begin{bmatrix}
\widetilde{F}_{n}  \\
n^{-1/2} \textbf{1}_n^{\prime}  
\end{bmatrix} = \widetilde{F}_n^{*} \widetilde{F}_n + n^{-1} \textbf{1}_n \textbf{1}_n^{\prime}.
\end{eqnarray*} 
Substituting the above into (\ref{eq:matrixeq1}) gives
\begin{eqnarray*}
\widetilde{F}_n^{*} \widetilde{\Delta}_{n}^{-1}(\widetilde{F}_{n}+\widetilde{D}_{n})
= (I_n - n^{-1} \textbf{1}_n \textbf{1}_n^{\prime}) \Gamma_n^{-1}.
\end{eqnarray*} 
This proves the identity (\ref{eq:longinverse}). 
\hfill $\Box$

\subsection{Some auxiliary lemmas}\label{sec:sup}

In the case that the spectral density $f$ corresponds to  an
AR$(p)$ model, $\phi_{j}(\tau;f) = \sum_{s=0}^{p-j}\phi_{j+s}\psi_{|\tau|-s}$ for
$\tau \leq 0$. This result is well known (see \cite{p:ino-06}, page
980). However we could not find the proof, thus for completeness we
give the proof below. 

\begin{lemma}\label{lemma:Amatrix}
Suppose
$f_{\theta}(\omega)=\sigma^{2}|1-\sum_{j=1}^{p}\phi_{j}e^{-ij\omega}|^{-2}
= \sigma^{2}|\sum_{j=0}^{\infty}\psi_{j}e^{-ij\omega}|^{2}$, where
$\{\phi_{j}\}_{j=1}^{p}$ correspond to the causal AR$(p)$ representation. 
Let $\phi_{j}(\tau;f)$ be defined as in (\ref{eq:Xtaun}). 
Then $\phi_{j}(\tau;f) = \sum_{s=0}^{p-j}\phi_{j+s}\psi_{|\tau|-s}$.
\begin{eqnarray}
\label{eq:APexpand}
[A_{p}(\underline{\phi})^{|\tau|+1}\underline{X}_{p}]_{(1)} = \sum_{\ell=1}^{p}X_{\ell}\sum_{s=0}^{p-\ell}\phi_{\ell+s}\psi_{|\tau|-s}. 
\end{eqnarray}
where we set
$\psi_{j}=0$ for $j<0$.
\end{lemma}
PROOF. To simplify notation let $A = A_{p}(\underline{\phi})$. 
The proof is based on the observation that the $j$th row of $A^{m}$ ($m\geq 1$) is the
$(j-1)$th row of $A^{m-1}$ (due to the structure of $A$).
Let $(a_{1,m},\ldots,a_{p,m})$ denote the first row of $A^{m}$. Using
this notation we have 
\begin{eqnarray*}
\begin{pmatrix}
a_{1,m} & a_{2,m} & \ldots & a_{p,m} \\
a_{1,m-1} & a_{2,m-1} & \ldots & a_{p,m-1} \\
\vdots & \vdots & \ddots & \vdots \\ 
a_{1,m-p+1} & a_{2,m-p+1} & \ldots & a_{p,m-p+1} \\
\end{pmatrix}
= 
\begin{pmatrix}
\phi_{1} & \phi_{2} & \ldots & \phi_{p-1} & \phi_{p} \\
1 & 0 & \ldots & 0 & 0\\
0 & 1 & \ldots & 0 & 0 \\
\vdots & \vdots & \ddots & \vdots & \vdots \\
0 & 0 & \ldots & 1 & 0 \\ 
\end{pmatrix}
\begin{pmatrix}
a_{1,m-1} & a_{2,m-1} & \ldots & a_{p,m-1} \\
a_{1,m-2} & a_{2,m-2} & \ldots & a_{p,m-2} \\
\vdots & \vdots & \ddots & \vdots \\
a_{1,m-p} & a_{2,m-p} & \ldots & a_{p,m-p} \\
\end{pmatrix}.
\end{eqnarray*}
From the above we observe that $a_{\ell,m}$ satisfies the system of equations
\begin{eqnarray}
a_{\ell,m} &=& \phi_{\ell}a_{1,m-1} + a_{\ell+1,m-1} \qquad 1\leq \ell \leq p-1\nonumber\\
a_{p,m} &=& \phi_{p}a_{1,m-1}.\label{eq:asystem}
\end{eqnarray}
Our aim is to obtain an expression for $a_{\ell,m}$ in terms of
$\{\phi_{j}\}_{j=1}^{p}$ and $\{\psi_{j}\}_{j=0}^{\infty}$ which we now define. Since the roots
of $\phi(\cdot)$ lies outside the unit circle the function $(1-\sum_{j=1}^{p}\phi_{j}z^{j})^{-1}$
is well defined for $|z|\leq 1$ and has the power series expansion
$(1-\sum_{i=1}^{p}\phi_{i}z)^{-1} = \sum_{i=0}^{\infty}\psi_{i}z^{i}$
for $|z|\leq 1$.  We use the well know result
$[A^{m}]_{1,1} = a_{1,m} = \psi_{m}$ (which can be proved by induction).
 Using this we obtain 
an expression for the coefficients $\{a_{\ell,m};2\leq \ell\leq p\}$ in terms of
$\{\phi_{i}\}$ and $\{\psi_{i}\}$. Solving the system of
equations in (\ref{eq:asystem}), starting with $a_{1,1}=\psi_{1}$
and recursively solving for $a_{p,m},\ldots,a_{2,m}$
we have 
\begin{eqnarray*}
a_{p,r} &=& \phi_{p}\psi_{r-1}  \qquad \qquad\qquad m-p \leq r \leq m \\
a_{\ell,r} &=& \phi_{\ell}a_{1,r-1} + a_{\ell+1,r-1} \qquad 1\leq \ell
               \leq p-1, \quad m-p \leq r \leq m\nonumber\\
\end{eqnarray*}
This gives $a_{p,m} = \phi_{p}\psi_{m-1}$, for $\ell=p-1$
\begin{eqnarray*}
a_{p-1,m} &=& \phi_{p-1}a_{1,m-1} + a_{p,m-1} \\
             &=& \phi_{p-1}\psi_{m-1} + \psi_{p}\psi_{m-2} 
\end{eqnarray*}
\begin{eqnarray*}
a_{p-2,m} &=& \phi_{p-2}a_{1,m-1} + a_{p-1,m-1} \\
             &=& \phi_{p-2}\psi_{m-1} + \phi_{p-1}\psi_{m-2} + \psi_{p}\psi_{m-3} 
\end{eqnarray*}
up to 
\begin{eqnarray*}
a_{1,m} &=& \phi_{1}a_{1,m-1} + a_{2,m-1} \\
             &=& \sum_{s=0}^{p-1}\phi_{1+s}\psi_{m-1-s} = (\psi_{m}).
\end{eqnarray*}
This gives the general expression
\begin{eqnarray*}
a_{p-r,m} &=& \sum_{s=0}^{r}\phi_{p-r+s}\psi_{m-1-s} \qquad 0\leq r\leq p-1. 
\end{eqnarray*}
In the last line of the above we change variables with $\ell = p-r$ to
give for $m\geq 1$
\begin{eqnarray*}
a_{\ell,m} &=& \sum_{s=0}^{p-\ell}\phi_{\ell+s}\psi_{m-1-s} \qquad 1\leq
\ell\leq p,
\end{eqnarray*}
where we set $\psi_{0}=1$ and for $t<0$, $\psi_{t}=0$.
Therefore 
\begin{eqnarray*}
[A_{}^{|\tau|+1}\underline{X}_{p}]_{(1)} = \sum_{\ell=1}^{p}X_{\ell}\sum_{s=0}^{p-\ell}\phi_{\ell+s}\psi_{|\tau|-s}. 
\end{eqnarray*}
Thus we obtain the desired result. \hfill $\Box$

\vspace{2mm}

In the proof of Lemma \ref{lemma:DFTinfty} we obtained an expression in terms for the best
linear predictor based on the parameters of an AR$(p)$
process. For completeness, we obtain an expression for the
left hand side predictive DFT for a general second order stationary
time series which is based on infinite future.
This expression can be used to prove the identity in equation
(\ref{eq:JARIn}).

 For $\tau \leq 0$, let  $\widehat{X}_{\tau}$ be the best linear predictor of $X_{\tau}$ given
the infinite future $\{X_{t}\}_{t=1}^{\infty}$ i.e.
\begin{eqnarray}
\label{eq:bestLin}
\widehat{X}_{\tau} = \sum_{s=1}^{\infty}\phi_{s}(\tau;f)X_{s}.
\end{eqnarray}
The left hand side predictive DFT given the infinite future is defined
as 
\begin{eqnarray}
\label{eq:JL}
\widehat{J}_{\infty,L}(\omega;f) = n^{-1/2}\sum_{\tau=-\infty}^{0}\widehat{X}_{\tau}e^{i\tau\omega}.
\end{eqnarray}

\begin{corollary}\label{cor:predictionTau}
Suppose that $f_{}$ satisfies Assumption \ref{assum:A}.
 Let $\{\phi_{j}(f)\}_{j=1}^{\infty}$ denote the
AR$(\infty)$ coefficients associated with $f$.
Let $\widehat{J}_{\infty,L}(\omega;f)$ be defined as in (\ref{eq:JL}).
Then, the best linear predictor of $X_{\tau}$ given
$\{X_{t}\}_{t=1}^{\infty}$ where $\tau\leq 0$ (defined in
(\ref{eq:bestLin})) can be evaluated using the recursion
$\widehat{X}_{\tau} =
\sum_{s=1}^{\infty}\phi_{s}(f)\widehat{X}_{\tau+s}$, where we set
$\widehat{X}_{t} = X_{t}$ for $t\geq
1$. Further, $\phi_{\ell}(\tau;f) = \sum_{s=0}^{\infty}\phi_{\ell+s}(f)\psi_{|\tau|-s}(f)$
and 
\begin{eqnarray}
\label{eq:Jinfty}
\widehat{J}_{\infty,L}(\omega;f) =
\frac{n^{-1/2}}{\phi(\omega;f)}\sum_{\ell=1}^{\infty}X_{\ell}\sum_{s=0}^{\infty}\phi_{\ell+s}(f)e^{-is\omega}.
\end{eqnarray}
\end{corollary}
PROOF. We recall that for general processes $X_{t}$ with $f_{}$ bounded away from 0 has the
AR$(\infty)$ representation $X_{t}=\sum_{j=1}^{\infty}\phi_{j}X_{j+t}+\varepsilon_{t}$ where
$\{\varepsilon_{t}\}$ are uncorrelated random variables. This
immediately implies that the best linear prediction of $X_{\tau}$
given $\{X_{t}\}_{t=1}^{\infty}$ can be evaluated using the recursion 
$\widehat{X}_{\tau} = \sum_{s=1}^{\infty}\phi_{s}(f)\widehat{X}_{\tau+s}$.
By using (\ref{eq:APexpand}) where we let $p\rightarrow \infty$ 
we have $\phi_{\ell}(\tau;f) =
\sum_{s=0}^{\infty}\phi_{\ell+s}(f)\psi_{|\tau|-s}(f)$. This gives the
first part of the result. 

To obtain an expression for $\widehat{J}_{L}(\cdot;f)$ we use 
(\ref{eq:JPPP}) where we let $p\rightarrow \infty$ to obtain the desired result.
\hfill $\Box$

\vspace{1em}

By a similar argument we can show that 
\begin{eqnarray}
\label{eq:JinftyR}
\widehat{J}_{\infty,R}(\omega;f) =
                                     \sum_{\tau=n+1}^{\infty}\widehat{X}_{\tau}e^{i\tau\omega}
= e^{i(n+1)\omega}
\frac{n^{-1/2}}{ \overline{\phi(\omega;f_\theta)}} \sum_{\ell=1}^{n}X_{n+1-t}\overline{\phi_{t}^{\infty}(\omega;f_\theta)}.
\end{eqnarray}
Since $\widehat{J}_{\infty,n}(\omega;f) =
\widehat{J}_{\infty,L}(\omega;f)+\widehat{J}_{\infty,R}(\omega;f)$, by
using (\ref{eq:Jinfty}) and (\ref{eq:JinftyR}) we immediately obtain
the identity (\ref{eq:JARIn}).

\subsubsection{A fast algorithm for computing the predictive DFT}
To end this section, for $p \leq n$, we provide an $O(\min(n \log n,np))$ algorithm to compute 
the predictive DFT of an AR$(p)$ spectral density at frequencies $\omega = \omega_{k,n}$, $1\leq k \leq n$.
 Recall from (\ref{eq:JAR}), 
\begin{eqnarray*}
&&\widehat{J}_{n}(\omega_{k,n} ;f_p) \\
&&\quad =\frac{n^{-1/2}}{\phi_{p}(\omega_{k,n})} \sum_{\ell=1}^{p}X_{\ell}\sum_{s=0}^{p-\ell}\phi_{\ell+s}e^{-is\omega_{k,n}}+
\frac{n^{-1/2}}{ \overline{\phi_{p}(\omega_{k,n})}} \sum_{\ell=1}^{p}X_{n+1-\ell}\sum_{s=0}^{p-\ell}
\phi_{\ell+s}e^{i(s+1)\omega_{k,n}},
\end{eqnarray*} where $f_p(\cdot) = |\phi_p(\cdot)|^2$ and $\phi_p(\omega_{k,n}) = 1-\sum_{j=1}^{p} \phi_{j} e^{-ij \omega_{k,n}}$.
We focus on the first term of $\widehat{J}_{n}(\omega_{k,n} ;f_p)$ since the second term is almost identical. Interchange the summation, the 
first term is
\begin{eqnarray*}
&&\frac{n^{-1/2}}{\phi_{p}(\omega_{k,n})} \sum_{\ell=1}^{p}X_{\ell}\sum_{s=0}^{p-\ell}\phi_{\ell+s}e^{-is\omega_{k,n}} \\
&&= \frac{n^{-1/2}}{\phi_{p}(\omega_{k,n})} \sum_{s=0}^{p-1} \left( \sum_{\ell=1}^{p-s} X_{\ell} \phi_{\ell+s} \right) e^{-is\omega_{k,n}} 
= \frac{n^{-1/2}}{\phi_{p}(\omega_{k,n})} \sum_{s=0}^{p-1} Y_{s} e^{-is\omega_{k,n}}
\end{eqnarray*} where $Y_{s} = \sum_{\ell=1}^{p-s} X_{\ell} \phi_{\ell+s}$ for $0 \leq s \leq p-1$. 
Note that $Y_s$ can be viewed as a convolution between
$(X_{1}, ..., X_{p})$ and $(0,0,\ldots,0,\phi_1, ..., \phi_p,0,\ldots,0)$. 
Based on  this  observation, the FFT can be utilized to evaluate $\{Y_{s}: 0\leq s \leq p-1\}$ in $O(p \log p)$ operations.

By direct calculation $\{\phi_p(\omega_{k,n}): 0\leq k \leq n-1\}$ and $\{ \sum_{s=0}^{p-1}Y_{s} e^{-is\omega_{k,n}}: 0\leq k \leq n-1\}$
has $O(np)$ complexity. An alternative method of calculation is based
on the observation that both $\phi_p(\omega_{k,n})$ and $\sum_{s=0}^{p-1}Y_{s} e^{-is\omega_{k,n}}$
can be viewed as the $k$th component of the DFT of length $n$ sequences $(1,-\phi_{1}, ..., -\phi_p, 0,...,0)$ and
$(Y_0, ..., Y_{p-1}, 0,..,0)$ respectively. Thus the FFT can be used
to evaluate both $\{\phi_p(\omega_{k,n}): 0\leq k \leq n-1\}$ and $\{
\sum_{s=0}^{p-1}Y_{s} e^{-is\omega_{k,n}}: 0\leq k \leq n-1\}$ in
$O(n\log n)$ operations. Therefore, since either method can be used to
evaluate these terms the total number of operations for evaluation of 
$\{\phi_p(\omega_{k,n}): 0\leq k \leq n-1\}$ and $\{
\sum_{s=0}^{p-1}Y_{s} e^{-is\omega_{k,n}}: 0\leq k \leq n-1\}$ is $O(\min(n \log n,np))$. 

Therefore, the overall computational complexity is $O(p \log p +n \log n
\wedge np)=O(\min(n \log n,np))$.

\section{Proof of results in Sections \ref{sec:approx} and 
\ref{sec:boundary}}\label{sec:proofs2}

\subsection{Proof of Theorems \ref{theorem:approx} and 
\ref{theorem:Bound} }\label{sec:approxproofs}

Many of the results below hinge on a small generalisation of Baxter's
inequality which we summarize below. 

\begin{lemma}[Extended Baxter's inequality]\label{lemma:BaxterEx}
Suppose $f(\cdot)$ is a spectral density function which satisfies
Assumption \ref{assum:A}. Let $\psi(\cdot)$ and $\phi(\cdot)$ be
defined as in (\ref{eq:thetaphi}) (for the simplicity, we omit the notation $f$ inside the $\psi(\cdot)$ and $\phi(\cdot)$). 
Let $\phi_{p+1}^{\infty}(\omega) =
\sum_{s=p+1}^{\infty}\phi_{s}e^{-is\omega}$. Further, let $\{\phi_{s,n}(\tau)\}$ denote the 
coefficients in the best linear predictor of $X_{\tau}$ given
$\Xunder_{n} = \{X_{t}\}_{t=1}^{n}$ and $\{\phi_{s}(\tau)\}$ the corresponding the
coefficients in the best linear predictor of $X_{\tau}$ given
$\Xunder_{\infty} = \{X_{t}\}_{t=1}^{\infty}$, where $\tau \leq 0$.
Suppose $p$ is large enough such that 
$\left\|\phi_{p}^{\infty}\right\|_{K}
\left\| \psi\right\|_{K}\leq 
\varepsilon<1$. Then for all $n>p$ we have 
\begin{equation}
\label{eq:baxterExt}
\sum_{s=1}^{n}(2^{K}+s^{K})\left|\phi_{s,n}(\tau)-\phi_{s}(\tau)\right| 
\leq C_{f,K}
\sum_{s=n+1}^{\infty}(2^{K}+s^{K})\left| \phi_{s}(\tau)\right|,
\end{equation}
where $C_{f,K} = \frac{3-\varepsilon}{1-\varepsilon}\left\|\phi\right\|_{K}^{2}
\left\|\psi\right\|_{K}^{2}$ and 
$\phi_{s}(\tau) =
\sum_{j=0}^{\infty}\phi_{s+j}\psi_{|\tau|-j}$ (we set $\psi_{0}=1$ and $\psi_{j}=0$
for $j < 0$).
\end{lemma}
\noindent PROOF. For completeness we give the proof in Appendix \ref{sec:baxter}. \hfill $\Box$

\vspace{1em}
\noindent {\bf PROOF of Equation (\ref{eq:widetildeDalt})} 
Since 
\begin{equation}
\label{eq:DF}
(D_{\infty,n}(f_{\theta}))_{k,t} =
n^{-1/2}\sum_{\tau\leq0}\left(\phi_{t}(\tau)e^{i\tau\omega_{k,n}}
+\phi_{n+1-t}(\tau)e^{-i(\tau-1)\omega_{k,n}}\right)
\end{equation}
we replace $\phi_{t}(\tau)$ in the above with the coefficients of the
MA and AR infinity expansions;
$\phi_{t}(\tau)=\sum_{s=0}^{\infty}\phi_{t+s}\psi_{|\tau|-s}$. Substituting
this into the first term in (\ref{eq:DF}) gives 
\begin{eqnarray*}
n^{-1/2}\sum_{\tau\leq0}\phi_{t}(\tau)e^{i\tau\omega_{k,n}}
&=&
n^{-1/2}\sum_{\tau\leq0}\sum_{s=0}^{\infty} \phi_{t+s}\psi_{-\tau-s}
e^{i\tau\omega_{k,n}} \\
&=& n^{-1/2}\sum_{s=0}^{\infty}\phi_{t+s}e^{-is\omega_{k,n}} \sum_{\tau\leq0}\psi_{-\tau-s} e^{-i(-\tau-s)\omega_{k,n}} \\
&=& n^{-1/2} \psi(\omega_{k,n}) \sum_{s=0}^{\infty} \phi_{t+s}e^{-is\omega_{k,n}} \\
&=& n^{-1/2}\phi(\omega_{k,n})^{-1} \phi_{t}^{\infty}(\omega_{k,n}),
\end{eqnarray*}
which gives the first term in (\ref{eq:widetildeDalt}). The second term
follows similarly. Thus giving the identity in equation
(\ref{eq:widetildeDalt}). \hfill $\Box$

\noindent Next we prove Theorem \ref{theorem:approx}. To do this we note that the entries of
$F_{n}^{*}\Delta_{n}(f_{\theta}^{-1})D_{\infty,n}(f)$ are 
\begin{eqnarray}
&& (F_{n}^{*}\Delta_{n}(f_{\theta}^{-1})D_{\infty,n}(f))_{s,t} \nonumber \\
&& \qquad =
\sum_{\tau \leq 0}\left[\phi_{t}(\tau;f)
G_{1,n}(s,\tau;f_{\theta})+\phi_{n+1-t}(\tau;f) G_{2,n}(s,\tau;f_{\theta})\right],
\label{eq:DDeltaDtilde}
\end{eqnarray} 
where $G_{1,n}$ and $G_{2,n}$ are defined as in (\ref{eq:DDeltaD}).
Thus 
\begin{eqnarray}
&&\left(F_{n}^{*}\Delta_{n}(f_{\theta}^{-1})\left[D_{n}(f) -
D_{\infty,n}(f)\right]\right)_{s,t}\nonumber \\
&& \quad = \sum_{\tau \leq 0}[ \left\{\phi_{t,n}(\tau;f)-\phi_{t}(\tau;f)\right\}
G_{1,n}(s,\tau;f_{\theta}) \nonumber \\
&& \qquad \quad \quad + \left\{\phi_{n+1-t,n}(\tau;f)-\phi_{n+1-t}(\tau;f)\right\} G_{2,n}(s,\tau;f_{\theta}) ].
\label{eq:DDeltaDtildeD}
\end{eqnarray} 
To prove Theorem \ref{theorem:approx} we bound the above terms. 
\vspace{1em}

\noindent {\bf PROOF of Theorem \ref{theorem:approx}}. 
To simplify
notation we only emphasis the coefficients associated with $f_\theta$
and not the coefficients associated with $f$. I.e. we set
$\phi_{s,n}(\tau;f) = \phi_{s,n}(\tau)$, $\phi_{s}(\tau;f) =
\phi_{s}(\tau)$, $\phi_{f} = \phi$ and $\psi_{f} = \psi$.

The proof of (\ref{eq:approx1AA}) simply follows from the definitions
of $D_{n}(f)$ and $D_{\infty,n}(f)$.

Next we prove (\ref{eq:approx1A}).
By using (\ref{eq:DDeltaDtildeD}) we have
\begin{equation*}
\left\|F_{n}^{*}\Delta_{n}(f_{\theta}^{-1})D_{n}(f)-
F_{n}^{*}\Delta_{n}(f_{\theta}^{-1})D_{\infty,n}(f)
\right\|_{1} \leq T_{1,n} + T_{2,n},
\end{equation*}
where 
\begin{eqnarray*}
T_{1,n} &=& 
\sum_{s,t=1}^{n}\sum_{\tau=-\infty}^{0}|\phi_{s,n}(\tau) -
\phi_{s}(\tau)||G_{1,n}(t,\tau;f_\theta)| \\
T_{2,n} &=& 
\sum_{s,t=1}^{n}\sum_{\tau=-\infty}^{0}|\phi_{n+1-s,n}(\tau) -
\phi_{n+1-s}(\tau)| |G_{2,n}(t,\tau;f_\theta)|.
\end{eqnarray*}
We focus on $T_{1,n}$, noting that the method for bounding $T_{2,n}$
is similar. Exchanging the summands we have 
\begin{equation*}
T_{1,n} \leq \sum_{\tau=-\infty}^{0}\sum_{t=1}^{n}|G_{1,n}(t,\tau;f_\theta)|
\sum_{s=1}^{n}|\phi_{s,n}(\tau) -
\phi_{s}(\tau)|.
\end{equation*}
To bound $\sum_{s=1}^{n}|\phi_{s,n}(\tau) -
\phi_{s}(\tau)|$ we require the generalized Baxter's inequality stated in 
Lemma \ref{lemma:BaxterEx}.
Substituting the bound in Lemma \ref{lemma:BaxterEx} into the above
(and for a sufficiently large $n$) we have 
\begin{equation*}
T_{1,n} \leq
C_{f,0}\sum_{\tau=-\infty}^{0}\sum_{t=1}^{n}|G_{1,n}(t,\tau;f_\theta)|\sum_{s=n+1}^{\infty}|\phi_{s}(\tau)|.
\end{equation*}
Using that $G_{1,n}(t,\tau) = \sum_{a\in
\mathbb{Z}}K_{f_\theta^{-1}}(\tau - t +an)$ we have the bound
\begin{eqnarray*}
T_{1,n} &\leq & C_{f,0}\sum_{\tau=-\infty}^{0}\sum_{t=1}^{n}\sum_{a\in
\mathbb{Z}}|K_{f_{\theta}^{-1}}(t-\tau+an)|\sum_{s=n+1}^{\infty}|\phi_{s}(\tau)| \\
&=& C_{f,0} \sum_{r\in
\mathbb{Z}}|K_{f_{\theta}^{-1}}(r)|\sum_{\tau=-\infty}^{0}\sum_{s=n+1}^{\infty}|\phi_{s}(\tau)|.
\end{eqnarray*}
Therefore,
\begin{eqnarray*}
T_{1,n}&\leq& C_{f,0}\sum_{r\in
\mathbb{Z}}|K_{f_{\theta}^{-1}}(r)|\sum_{\tau=-\infty}^{0}\sum_{s=n+1}^{\infty}|\phi_{s}(\tau)|\\
&\leq& C_{f,0}\sum_{r\in
\mathbb{Z}}|K_{f_{\theta}^{-1}}(r)|\sum_{\tau=-\infty}^{0}\sum_{s=n+1}^{\infty}\sum_{j=0}^{\infty}|\phi_{s+j}||\psi_{-\tau-j}|
~~ \textrm{ (use }\phi_{s}(\tau)=\sum_{j=0}^{\infty}\phi_{s+j}\psi_{|\tau|-j})\\
&=& C_{f,0}\sum_{r\in
\mathbb{Z}}|K_{f_{\theta}^{-1}}(r)|\sum_{\tau=0}^{\infty}|\psi_{\tau-j}|\sum_{s=n+1}^{\infty}\sum_{j=0}^{\infty}|\phi_{s+j}|
\qquad (\textrm{change limits of }\sum_{\tau}) 
\\
&\leq&C_{f,0}\sum_{r\in
\mathbb{Z}}|K_{f_{\theta}^{-1}}(r)|\sum_{\ell}|\psi_{\ell}|\sum_{u=n+1}^{\infty}|u\phi_{u}| 
\quad (\textrm{change of variables } u = s+j).
\end{eqnarray*}
Next we use Assumption \ref{assum:A}(i) to give
\begin{eqnarray*}
T_{1,n} 
&\leq& C_{f,0}\sum_{r\in
\mathbb{Z}}|K_{f_{\theta}^{-1}}(r)|\sum_{\ell}|\psi_{\ell}|\sum_{s=n+1}^{\infty}\frac{s^{K}}{s^{K-1}}|\phi_{s}|
\\
&\leq& \frac{C_{f,0}}{n^{K-1}}\sum_{r\in
\mathbb{Z}}|K_{f_{\theta}^{-1}}(r)|\sum_{\ell}|\psi_{\ell}|\sum_{s=n+1}^{\infty}|s^{K}\phi_{s}|
\\
&\leq&\frac{C_{f,0}}{n^{K-1}}\rho_{n,K}(f)\|\psi\|_0\|\phi\|_{K}\sum_{r\in
\mathbb{Z}}|K_{f_{\theta}^{-1}}(r)|.
\end{eqnarray*}
We note that the inverse covariance  $K_{f_\theta^{-1}}(r) =
\int_{0}^{2\pi}f_{\theta}^{-1}(\omega)e^{ir\omega}d\omega =
\sigma_{f_\theta}^{-2}\int_{0}^{2\pi}|\phi_{f_\theta}(\omega)|^{2}e^{ir\omega}d\omega
= 
\sigma_{f_\theta}^{-2}\sum_{j}\phi_{j}(f_\theta)\phi_{j+r}(f_\theta)$. Therefore
\begin{equation}
\label{eq:BoundK}
\sum_{r=-\infty}^{\infty}|K_{f_\theta}(r)| \leq \sigma_{f_\theta}^{-2}\|\phi_{f_\theta}\|_{0}^{2}.
\end{equation}
Substituting this into the above yields the bound 
\begin{equation*}
T_{1,n} 
\leq \frac{C_{f,0}}{\sigma_{f_\theta}^{2}n^{K-1}}\rho_{n,K}(f)\|\psi\|_0 \|\phi_{f_\theta}\|_{0}^{2}\|\phi\|_{K}.
\end{equation*}
The same bound holds for $T_{2,n}$. Together the bounds for $T_{1,n}$
and $T_{2,n}$ give 
\begin{equation*}
\left\|F_{n}^{*}\Delta_{n}(f_{\theta}^{-1})D_{n}(f_{\theta})-
F_{n}^{*}\Delta_{n}(f_{\theta}^{-1})D_{\infty,n}(f_{\theta})
\right\|_{1} \leq \frac{2C_{f,0}}{\sigma_{f_\theta}^{2}n^{K-1}}\rho_{n,K}(f)\|\psi\|_0 \|\phi_{f_\theta}\|_{0}^{2}\|\phi\|_{K}.
\end{equation*}
Replacing $\|\psi_f\|_0 = \|\psi\|_0$ and $\|\phi_{f}\|_{K} = \|\phi\|_{K}$,
this proves (\ref{eq:approx1A}).

To prove (\ref{eq:approx1B}) we recall
\begin{equation*}
\|X_{t}X_{s}\|_{\Ex,q} = \left( \Ex |X_{t}X_{s}|^{q}\right)^{1/q} \leq
\left( \Ex |X_{t}|^{2q}\right)^{1/2q} \left( \Ex |X_{s}|^{2q}\right)^{1/2q}
\leq \|X\|_{\Ex,2q}^{2}.
\end{equation*} 
Therefore,
\begin{eqnarray*}
&& n^{-1}\left\|\underline{X}_{n}^{\prime}F_{n}^{*}\Delta_{n}(f_{\theta}^{-1})\left(D_{n}(f)-
D_{\infty,n}(f)\right)\underline{X}_{n}
\right\|_{\Ex,q} \\
&& \qquad \leq 
n^{-1}\sum_{s,t=1}^{n}
\left|\left( F_{n}^{*}\Delta_{n}(f_{\theta}^{-1})\left(D_{n}(f)-
D_{\infty,n}(f)\right)\right)_{s,t}\right| \|X_{t}X_{s}\|_{\Ex,q} \\
&& \qquad \leq n^{-1}\left\| F_{n}^{*}\Delta_{n}(f_{\theta}^{-1})\left(D_{n}(f)-
D_{\infty,n}(f)\right) \right\|_{1} \|X\|_{\Ex,2q}^{2} \\
&& \qquad \leq
\frac{2C_{f,0}}{\sigma_{f_\theta}^{2}n^{K}}\rho_{n,K}(f)\|\psi_{f}\|_0 \|\phi_{f_\theta}\|_{0}^{2}\|\phi_{f}\|_{K} \|X\|_{\Ex,2q}^{2},
\end{eqnarray*} 
where the last line follows from the inequality in (\ref{eq:approx1A}).
This proves (\ref{eq:approx1B}).
\hfill $\Box$

\vspace{1em}

\noindent {\bf PROOF of Theorem \ref{theorem:Bound}}
For notational simplicity, we omit the parameter dependence on
$f_{\theta}$. We first prove (\ref{eq:Bound1}). We observe that
\begin{eqnarray*}
\left\| F_{n}^{*}\Delta_{n}(f_{\theta}^{-1})D_{\infty,n}(f_{\theta})
\right\|_{1} &\leq& 
\sum_{s,t=1}^{n}\sum_{\tau=-\infty}^{0}\left( |
\phi_{s}(\tau)||G_{1,n}(t,\tau)| + |\phi_{n+1-s}(\tau)||G_{2,n}(t,\tau)| \right) \\
&=& S_{1,n} + S_{2,n}.
\end{eqnarray*}
As in the proof of Theorem \ref{theorem:approx}, we bound each term
separately. Using a similar set of bounds to those used in the proof
of Theorem \ref{theorem:approx} we have 
\begin{eqnarray*}
S_{1,n} &\leq&
\sum_{r\in
\mathbb{Z}}|K_{f_{\theta}^{-1}}(r)|\sum_{\ell}|\psi_{\ell}|\sum_{s=1}^{n}\sum_{j=0}^{\infty}|\phi_{s+j}| \\
&\leq& \sum_{r\in
\mathbb{Z}}|K_{f_{\theta}^{-1}}(r)|\sum_{\ell}|\psi_{\ell}|\sum_{s=1}^{\infty}|s\phi_{s}|
\leq \frac{1}{\sigma_{f_\theta}^{2}}\|\psi_{f_\theta}\|_0 \|\phi_{f_\theta}\|_{0}^{2}\|\phi_{f_\theta}\|_{1},
\end{eqnarray*}
where the bound $\sum_{r\in
\mathbb{Z}}|K_{f_{\theta}^{-1}}(r)| \leq \sigma_{f_\theta}^{-2}\|\phi_{f_\theta}\|_{0}^{2}$ follows from
(\ref{eq:BoundK}).
Using a similar method we obtain the bound $S_{2,n}\leq
\sigma_{f_\theta}^{-2} \|\psi_{f_\theta}\|_0
\|\phi_{f_\theta}\|_{0}^{2}\|\phi_{f_\theta}\|_{1}$. Altogether the bounds
for $S_{1,n}$ and $S_{2,n}$ give
\begin{equation*}
\left\| F_{n}^{*}\Delta_{n}(f_{\theta}^{-1})D_{\infty,n}(f_{\theta})
\right\|_{1} \leq \frac{2}{\sigma_{f_\theta}^{2}}\|\psi_{f_\theta}\|_0 \|\phi_{f_\theta}\|_{0}^{2}\|\phi_{f_\theta}\|_{1},
\end{equation*}
this proves (\ref{eq:Bound1}).

The proof of (\ref{eq:Bound2}) uses 
the triangle inequality 
\begin{eqnarray*}
\left\| \Gamma_{n}(f_{\theta})^{-1} -
C_{n}(f_{\theta}^{-1})\right\|_{1}&=&\left\|F_{n}^{*}\Delta_{n}(f_\theta^{-1})D_{n}(f_\theta)
\right\|_{1} \\
&\leq &\left\|F_{n}^{*}\Delta_{n}(f_\theta^{-1})\left(D_{n}(f_\theta) - D_{\infty,n}(f_{\theta})\right)
\right\|_{1} \\
&& + \left\|F_{n}^{*}\Delta_{n}(f_\theta^{-1})D_{\infty,n}(f_\theta)
\right\|_{1}.
\end{eqnarray*}
Substituting the bound Theorem \ref{theorem:approx} (equation (\ref{eq:approx1A})) and
(\ref{eq:Bound1}) into the above gives (\ref{eq:Bound2}). 

The proof of (\ref{eq:BD1}) uses the bound in
(\ref{eq:Bound2}) together with
similar arguments to those in the proof of Theorem
\ref{theorem:approx}, we omit the details. 
\hfill $\Box$

\vspace{3mm}

\subsection{Proof of results in Section \ref{sec:higherexpan}}\label{sec:higherproof}


To prove Theorem \ref{thm:higherorder} 
we recall some notation. Let $P_{[1,\infty)} X$ and $P_{(-\infty,n]}
X$ denote the projection of $X$ onto $\spa (X_{t}; t\geq 1)$ or $\spa (X_{t}; t \leq n)$. 
Then using the definition of $\{\phi_{j}(\tau)\}$ we have 
\begin{eqnarray*}
P_{[1,\infty)}X_{\tau} = \sum_{j=1}^{\infty}\phi_{j}(\tau)X_{j}
  \textrm{ for }\tau\leq 1
\quad \text{and} \quad
P_{(-\infty,n]}X_{\tau}  =  \sum_{j=1}^{\infty}h_{j}(\tau)X_{j}
  \textrm{ for }\tau>n.
\end{eqnarray*} 
By stationarity we have $h_{j}(\tau) = \phi_{n+1-j}(n+1-\tau)$. We use
this notation below. 

\vspace{2mm}
\noindent {\bf PROOF of Theorem \ref{thm:higherorder}} We first focus on the
prediction coefficient associated with $X_{j}$ where $1\leq j \leq n$ and show that
\begin{eqnarray*}
\sum_{\tau=-\infty}^{0}\phi_{j,n}(\tau)e^{i\tau\omega} = \phi(\omega)^{-1}\sum_{s=1}^{\infty}\zeta_{j,n}^{(s)}(\omega;f),
\end{eqnarray*}
where $\zeta_{j,n}^{(s)}(\omega;f)$ is defined in
(\ref{eq:Dinftys}).  This will lead to the series expansion of
$\widehat{J}_{n}(\omega;f)$. 

Let us assume $\tau\leq 0$. By using
\cite{p:ino-06}, Theorem 2.5, we obtain an expression for 
the difference between the 
coefficients of $X_{j}$ for the
finite predictor (based on the space $\spa(X_{1},\ldots,X_{n})$) and the
infinite predictor (based on the space
$\spa(X_{1},X_{2},\ldots)$):
\begin{eqnarray*}
&& \phi_{j,n}(\tau) - \phi_{j}(\tau) \\
&& =
                                      \sum_{u_{1}=n+1}^{\infty}\phi_{u_1}(\tau)h_{j}(u_{1})
 +
\sum_{u_{1}=n+1}^{\infty}\sum_{v_{1}=-\infty}^{0}\phi_{u_1}(\tau)h_{v_{1}}(u_{1})\phi_{j}(v_{1})+\\
&&\sum_{u_{1}=n+1}^{\infty}\sum_{v_{1}=-\infty}^{0}\sum_{u_{2}=n+1}^{\infty}\phi_{u_1}(\tau)h_{v_{1}}(u_{1})\phi_{u_2}(v_{1})
   h_{j}(u_2)+\ldots\\
\end{eqnarray*}
Replacing $h_{j}(\tau) = \phi_{n+1-j}(n+1-\tau)$ we have 
\begin{eqnarray*}
&& \phi_{j,n}(\tau) - \phi_{j}(\tau) \\
&&=   \sum_{u_{1}=n+1}^{\infty}\phi_{u_1}(\tau)\phi_{n+1-j}(n+1-u_{1}) +
\sum_{u_{1}=n+1}^{\infty}\sum_{v_{1}=-\infty}^{0}\phi_{u_1}(\tau)\phi_{n+1-v_{1}}(n+1-u_{1})\phi_{j}(v_{1})+\\
&&\sum_{u_{1}=n+1}^{\infty}\sum_{v_{1}=-\infty}^{0}\sum_{u_{2}=n+1}^{\infty}\phi_{u_1}(\tau)\phi_{n+1-v_{1}}(n+1-u_{1})\phi_{u_2}(v_{1})
   \phi_{n+1-j}(n+1-u_2)+\ldots\\
\end{eqnarray*}
Changing variables with $u_{1}\rightarrow u_{1}-n-1$, $v_{1}\rightarrow
-v_{1},\ldots$ gives
\begin{eqnarray*}
&& \phi_{j,n}(\tau) - \phi_{j}(\tau) \\
&&=   \sum_{u_{1}=0}^{\infty}\phi_{n+1+u_1}(\tau)\phi_{n+1-j}(-u_{1}) +
\sum_{u_{1}=0}^{\infty}\sum_{v_{1}=0}^{\infty}\phi_{n+1+u_1}(\tau)\phi_{n+1+v_{1}}(-u_{1})\phi_{j}(-v_{1})+\\
&&\sum_{u_{1}=0}^{\infty}\sum_{v_{1}=0}^{\infty}\sum_{u_{2}=0}^{\infty}\phi_{n+1+u_1}(\tau)\phi_{n+1+v_{1}}(-u_{1})\phi_{n+1+u_2}(-v_{1})
   \phi_{n+1-j}(-u_2)+\ldots\\
&& =
    \sum_{s=2}^{\infty}\sum_{u_{1},\ldots,u_{s}=0}^{\infty}\phi_{n+1+u_1}(\tau)\left[\prod_{a=1}^{s-1}\phi_{n+1+u_{a+1}}(-u_{a})\right]
\left[\phi_{n+1-j}(-u_{s})\delta_{s\bmod 2=1}+\phi_{j}(-u_{s})\delta_{s\bmod 2=0}\right].
\end{eqnarray*}
Next our focus will be on the term inside the sum
$\sum_{s=2}^{\infty}$, where we define 
\begin{eqnarray*}
\phi_{j,n}^{(s)}(\tau) = \sum_{u_{1},\ldots,u_{s}=0}^{\infty}\phi_{n+1+u_1}(\tau)\left[\prod_{a=1}^{s-1}\phi_{n+1+u_{a+1}}(-u_{a})\right]
\left[\phi_{n+1-j}(-u_{s})\delta_{s\bmod 2=1}+\phi_{j}(-u_{s})\delta_{s\bmod 2=0}\right].
\end{eqnarray*}
Using this notation we have
\begin{eqnarray*}
\phi_{j,n}(\tau) - \phi_{j}(\tau) =  \sum_{s=2}^{\infty}\phi_{j,n}^{(s)}(\tau).
\end{eqnarray*}
We will rewrite $\phi_{j,n}^{(s)}(\tau)$ as a convolution. To do so, we first note
that from Lemma \ref{lemma:DFTinfty}
\begin{eqnarray*}
\phi_{j}(\tau) = \sum_{s=0}^{\infty}\phi_{j+s}\psi_{|\tau|-s} = 
\sum_{\ell=0}^{|\tau|}\psi_{\ell}\phi_{j+|\tau|-\ell} \quad\textrm{ for
  }\tau\leq 0.
\end{eqnarray*}
This can be written as an integral
\begin{eqnarray*}
\phi_{j}(\tau) = 
  \frac{1}{2\pi}\int_{0}^{2\pi}\phi(\lambda)^{-1}\phi_{j}^{\infty}(\lambda)e^{-i\tau\lambda}d\lambda, \quad \tau \leq 0
\end{eqnarray*}
where $\phi(\lambda)^{-1}=\sum_{s=0}^{\infty}\psi_{s}e^{-is\lambda}$
and
$\phi_{j}^{\infty}(\lambda)=\sum_{s=0}^{\infty}\phi_{j+s}e^{-is\lambda}$. Using
this representation we observe that for $\tau\geq 1$ 
\begin{eqnarray*}
\frac{1}{2\pi}\int_{0}^{2\pi}\phi(\lambda)^{-1}\phi_{j}^{\infty}(\lambda)e^{-i\tau\lambda}d\lambda = 0.
\end{eqnarray*}
Based on the above we define the ``extended'' coefficient
\begin{eqnarray}
\widetilde{\phi}_{j}(\tau) =  \frac{1}{2\pi}\int_{0}^{2\pi}\phi(\lambda)^{-1}\phi_{j}^{\infty}(\lambda)e^{-i\tau\lambda}d\lambda=
\left\{
\begin{array}{cc}
\phi_{j}(\tau) & \tau \leq 0 \\
0   & \tau \geq 1.
\end{array}
\right. \label{eq:phitilde}
\end{eqnarray}
$\phi_{j}(\tau)$ can be treated as the first term in the
expansion $\phi_{j,n}(\tau) = \phi_{j}(\tau) +
\sum_{s=2}^{\infty}\phi_{j,n}^{(s)}(\tau)$. 
In the same way we have extended the definition of $\phi_{j,n}(\tau)$
over $\tau \in \mathbb{Z}$ we do the same for the higher order terms
and define 
\begin{eqnarray*}
\widetilde{\phi}_{j,n}^{(s)}(\tau) = \sum_{u_{1},\ldots,u_{s}=0}^{\infty}\widetilde{\phi}_{n+1+u_1}(\tau)\left[\prod_{a=1}^{s-1}\phi_{n+1+u_{a+1}}(-u_{a})\right]
\left[\phi_{n+1-j}(-u_{s})\delta_{s\bmod 2=1}+\phi_{j}(-u_{s})\delta_{s\bmod 2=0}\right].
\end{eqnarray*}
Observe from the above definition that $\widetilde{\phi}_{j,n}^{(s)}(\tau) = \phi_{j,n}^{(s)}(\tau)$ for
$\tau \leq 0$ and $\widetilde{\phi}_{j,n}^{(s)}(\tau)=0$ for $\tau\geq
1$. Substituting (\ref{eq:phitilde}) into $\widetilde{\phi}_{j,n}^{(s)}(\tau)$
for $\tau \in \mathbb{Z}$ gives rise to a convolution for $\widetilde{\phi}_{j,n}^{(s)}(\tau)$
\begin{eqnarray}
\widetilde{\phi}_{j,n}^{(s)}(\tau) 
&=& \frac{1}{(2\pi)^{s}}\int_{[0,2\pi]^{s}}e^{-i\tau\lambda_{1}}
\left[\prod_{a=1}^{s}\phi(\lambda_{a})^{-1}\right]
\left[\prod_{a=1}^{s-1}\sum_{u_{a}=0}^{\infty}\phi_{n+1+u_a}^{\infty}(\lambda_{a})e^{iu_{a}\lambda_{a+1}}\right] \nonumber \\
&&\left[\phi_{n+1-j}^{\infty}(\lambda_{s})\delta_{s\bmod
   2=0}+\phi_{j}^{\infty}(\lambda_{s})\delta_{s\bmod
   2=1}\right]d\underline{\lambda}_{s} \nonumber  \\
&=&  \frac{1}{(2\pi)^{s}}\int_{[0,2\pi]^{s}}e^{-i\tau\lambda_{1}}
\left[\prod_{a=1}^{s}\phi(\lambda_{a})^{-1}\right]\prod_{a=1}^{s-1}\Phi_{n}(\lambda_{a},\lambda_{a+1})  \nonumber \\
&&\times \left[\phi_{n+1-j}^{\infty}(\lambda_{s})\delta_{s\bmod
   2=0}+\phi_{j}^{\infty}(\lambda_{s})\delta_{s\bmod
   2=1}\right]d\underline{\lambda}_{s}.
\label{eq:leftpredict}
\end{eqnarray}
Evaluating the Fourier transform of
$\{\widetilde{\phi}_{j,n}^{(s)}(\tau)\}_{\tau\in \mathbb{Z}}$ gives 
\begin{eqnarray*}
\sum_{\tau \leq 0}\phi_{j,n}^{(s)}(\tau)e^{i\tau\omega} &=& 
\sum_{\tau \in
                                                            \mathbb{Z}}\widetilde{\phi}_{j,n}^{(s)}(\tau)e^{i\tau\omega} \\
&=& \sum_{\tau \in \mathbb{Z}} \frac{1}{(2\pi)^{s}}\int_{[0,2\pi]^{s}}e^{-i\tau(\lambda_{1}-\omega)}
\left[\prod_{a=1}^{s}\phi(\lambda_{a})^{-1}\right]\prod_{a=1}^{s-1}\Phi_{n}(\lambda_{a},\lambda_{a+1})\\
&&\quad \times \left[\phi_{n+1-j}^{\infty}(\lambda_{s})\delta_{s\bmod
   2=0}+\phi_{j}^{\infty}(\lambda_{s})\delta_{s\bmod
   2=1}\right]d\underline{\lambda}_{s} \\
&=&  \frac{1}{(2\pi)^{s}}\int_{[0,2\pi]^{s}} \left( \sum_{\tau \in \mathbb{Z}} e^{-i\tau(\lambda_{1}-\omega)} \right)
\left[\prod_{a=1}^{s}\phi(\lambda_{a})^{-1}\right]\prod_{a=1}^{s-1}\Phi_{n}(\lambda_{a},\lambda_{a+1})\\
&&\quad \times \left[\phi_{n+1-j}^{\infty}(\lambda_{s})\delta_{s\bmod
   2=0}+\phi_{j}^{\infty}(\lambda_{s})\delta_{s\bmod
   2=1}\right]d\underline{\lambda}_{s} \\
&=&  \frac{1}{(2\pi)^{s-1}}\int_{[0,2\pi]^{s}}
\left[\prod_{a=1}^{s}\phi(\lambda_{a})^{-1}\right]\prod_{a=1}^{s-1}\Phi_{n}(\lambda_{a},\lambda_{a+1})\\
&&\quad  \times \left[\phi_{n+1-j}^{\infty}(\lambda_{s})\delta_{s\bmod
   2=0}+\phi_{j}^{\infty}(\lambda_{s})\delta_{s\bmod
   2=1}\right] \delta_{\lambda_1 = \omega}   d\underline{\lambda}_{s} \\
&=&   \phi(\omega)^{-1} \frac{1}{(2\pi)^{s-1}}\int_{[0,2\pi]^{s}}
\prod_{a=1}^{s-1}\phi(\lambda_{a+1})^{-1}\Phi_{n}(\lambda_{a},\lambda_{a+1})\\
&&\quad  \times \left[\phi_{n+1-j}^{\infty}(\lambda_{s})\delta_{s\bmod
   2=0}+\phi_{j}^{\infty}(\lambda_{s})\delta_{s\bmod
   2=1}\right] \delta_{\lambda_1 = \omega}   d\underline{\lambda}_{s}
  \\
&=& \phi(\omega)^{-1}\zeta_{j,n}^{(s)}(\omega;f), 
\end{eqnarray*} 
where 
\begin{eqnarray*}
\zeta_{j,n}^{(s)}(\omega;f)  &=& \frac{1}{(2\pi)^{s-1}}\int_{[0,2\pi]^{s}}
\prod_{a=1}^{s-1}\phi(\lambda_{a+1})^{-1}\Phi_{n}(\lambda_{a},\lambda_{a+1})\\
&&\quad  \times \left[\phi_{n+1-j}^{\infty}(\lambda_{s})\delta_{s\bmod
   2=0}+\phi_{j}^{\infty}(\lambda_{s})\delta_{s\bmod
   2=1}\right] \delta_{\lambda_1 = \omega}   d\underline{\lambda}_{s}.
\end{eqnarray*}
The above holds for $s\geq 2$. But a similar representation also holds
for the first term, $\phi_{j}(\tau)$, in the expansion of
$\phi_{j,n}(\tau)$. Using the same argument as above we have 
\begin{eqnarray*}
\sum_{\tau \leq 0}\phi_{j}(\tau) e^{i\tau\omega} =
 \phi(\lambda_1)^{-1} \zeta_{j,n}^{(1)}(\omega; f) = \int_{0}^{2\pi}
  \phi(\lambda_1)^{-1} \phi_{j}^{\infty}(\lambda_1) 
\delta_{\lambda_1 = \omega} d\lambda_1.
\end{eqnarray*} 
Altogether this gives an expression for the Fourier transform of the
predictor coefficients corresponding to $X_{j}$ ($1\leq j \leq n$) at all
lags $\tau\leq 0$:
\begin{eqnarray}
\label{eq:phileft}
\sum_{\tau \leq 0}\phi_{j,n}(\tau) e^{i\tau\omega} =
\sum_{\tau \leq 0}\phi_{j}(\tau) e^{i\tau\omega} +
\sum_{s=2}^{\infty}\sum_{\tau \leq 0}\phi_{j,n}^{(s)}(\tau)e^{i\tau\omega}
=  \phi(\omega)^{-1}\sum_{s=1}^{\infty} \zeta_{j,n}^{(s)}(\omega; f).
\end{eqnarray}
Using a similar set of arguments we have 
\begin{eqnarray}
\label{eq:phiright}
\sum_{\tau \geq n+1}\phi_{j,n}(n+1-\tau) e^{i\tau\omega} =
  \overline{\phi(\omega)}^{-1}
e^{i(n+1)\omega}\sum_{s=1}^{\infty}\overline{\zeta_{n+1-j,n}^{(s)}(\omega; f)}.
\end{eqnarray}
Therefore by using the above and setting $j=t$ we have the series expansion
\begin{eqnarray*}
\widehat{J}_{n}(\omega;f) &=&
  \frac{1}{n}\sum_{t=1}^{n}X_{t}\left(\sum_{\tau\leq
  0}\phi_{t,n}(\tau)e^{i\tau\omega} + \sum_{\tau\geq
                            n+1}\phi_{n+1-t,n}(n+1-\tau)e^{i\tau\omega}\right) \\
&=& \sum_{s=1}^{\infty}\widehat{J}_{n}^{(s)}(\omega;f),
\end{eqnarray*}
where $\widehat{J}_{n}^{(s)}(\omega;f)$ is defined in Theorem \ref{thm:higherorder}.
Thus proving (\ref{eq:Japprox1}).  To prove (\ref{eq:Dapprox1}) we use that 
\begin{eqnarray*}
(D_{n}(f))_{k,t} =  n^{-1/2} \sum_{\tau \leq 0}\phi_{t,n}(\tau)
  e^{i\tau\omega_{k,n}} +  n^{-1/2}e^{i\omega_{k,n}}\sum_{\tau \leq 0}\phi_{n+1-t,n}(\tau) e^{-i\tau\omega_{k,n}}.
\end{eqnarray*}
This together with (\ref{eq:phileft}) and (\ref{eq:phiright}) proves
the (\ref{eq:Dapprox1}). 

Finally, we will obtain a bound for $\zeta_{t,n}^{(s)}(\omega; f)$,
which results in bound for $\widehat{J}_{n}^{(s)}(\omega;f)$ (which we
use to prove (\ref{eq:Jhatapprox2})). 
By using  (\ref{eq:Dinftys}) we have
\begin{eqnarray}
|\zeta_{t,n}^{(s)}(\omega;f)| &\leq &
\left(\sup_{\lambda}|\phi(\lambda;f)|^{-1}\sup_{\lambda_{1},\lambda_{2}}|\Phi_{n}(\lambda_{1},\lambda_{2})|\right)^{s-1}\times \nonumber\\
&& \bigg(\|\phi_{t}^{\infty}(\lambda_s;f)\|_{0}\delta_{ s\equiv 1 (\bmod 2)} 
+\|\phi_{n+1-t}^{\infty}(\lambda_s;f)\|_{0}\delta_{ s\equiv 0(\bmod
   2)}\bigg). \label{eq:phiii3}
\end{eqnarray}
To bound the terms above we note that
\begin{eqnarray*}
\sup_{\lambda_{1},\lambda_{2}}|\Phi_{n}(\lambda_{1},\lambda_{2})|\leq
  \sum_{j=1}^{\infty}|j\phi_{n+j}(f)| = O(n^{-K+1}),
\end{eqnarray*}
where the above follows from Assumption \ref{assum:A}. Further 
\begin{eqnarray*}
|\phi(\lambda;f)|^{-1} \leq \sum_{j=0}^{\infty}|\psi_{j}(f)|,
\end{eqnarray*}
where $\{\psi_{j}(j)\}$ are the MA$(\infty)$ coefficients
corresponding to the spectral density $f$. Substituting these bounds
into (\ref{eq:phiii3}) gives 
\begin{eqnarray*}
|\zeta_{t,n}^{(s)}(\omega;f)| &\leq & \left(\sum_{j=0}^{\infty}|\psi_{j}(f)|\sum_{j=1}^{\infty}|j\phi_{n+j}(f)|\right)^{s-1}
\bigg(\sum_{j=t}^{\infty}|\phi_{j}(f)|+\sum_{j=n+1-t}^{\infty}|\phi_{j}(f)|\bigg).
\end{eqnarray*}
Substituting the above into $\widehat{J}_{n}^{(s)}(\omega;f)$ gives
the bound
\begin{eqnarray*}
|\widehat{J}_{n}^{(s)}(\omega;f)| &\leq& \sum_{j=0}^{\infty}|\psi_{j}(f)|
  \left(\sum_{j=0}^{\infty}|\psi_{j}(f)|\sum_{j=1}^{\infty}|j\phi_{n+j}(f)|\right)^{s-1}\\
&&\times \frac{2}{\sqrt{n}}\sum_{t=1}^{n}|X_{t}|\left(\sum_{j=t}^{\infty}|\phi_{j}(f)|+\sum_{j=n+1-t}^{\infty}|\phi_{j}(f)|\right).
\end{eqnarray*}
Therefore taking expectation of $|\widehat{J}_{n}^{(s)}(\omega;f)|$
gives the bound 
\begin{eqnarray*}
\Ex|\widehat{J}_{n}^{(s)}(\omega;f)| \leq \Ex|X_{0}| \sum_{j=0}^{\infty}|\psi_{j}(f)|
  \left(\sum_{j=0}^{\infty}|\psi_{j}(f)|\sum_{j=1}^{\infty}|j\phi_{n+j}(f)|\right)^{s-1}
\frac{4}{\sqrt{n}}\sum_{j=1}^{\infty}|j\phi_{j}(f)|.
\end{eqnarray*}
For large enough $n$,
$\sum_{j=0}^{\infty}|\psi_{j}(f)| \cdot \sum_{j=1}^{\infty}|j\phi_{n+j}(f)| \leq
Cn^{-K+1}<1$.  Therefore
\begin{eqnarray*}
\sum_{s=m+1}^{\infty}\Ex|\widehat{J}_{n}^{(s)}(\omega;f)| &\leq&
                                                                 n^{-1/2}\sum_{s=m+1}^{\infty}\left(\frac{C}{n^{K-1}}\right)^{s-1}
= O(n^{-m(K-1)-1/2})
\end{eqnarray*}
thus 
\begin{eqnarray*}
\sum_{s=m+1}^{\infty}\widehat{J}_{n}^{(s)}(\omega;f) 
= O_{p}(n^{-m(K-1)-1/2}).
\end{eqnarray*}
This proves the approximation in (\ref{eq:Jhatapprox2}). Thus we have
proved the result. \hfill $\Box$

\vspace{1em}

\noindent {\bf Proof of equation (\ref{eq:recursionZeta})} Our aim is to show for $s \geq 3$,
\begin{eqnarray}
\zeta_{t,n}^{(s)}(\omega;f) =
  \frac{1}{(2\pi)^{2}}\int_{[0,2\pi]^{2}}\phi(y_{1};f)^{-1}\phi(y_{2};f)^{-1}\Phi_{n}(\omega,y_1)
\Phi_{n}(y_1,y_2) \zeta_{t,n}^{(s-2)}(y_2;f)dy_{1}dy_{2}. 
\label{eq:recursionZeta2}
\end{eqnarray}
We recall the definition of $\zeta_{t,n}^{(s)}(\omega;f)$ in equation (\ref{eq:Dinftys}) 
\begin{eqnarray*}
\zeta_{t,n}^{(s)}(\omega;f) &=&
  \frac{1}{(2\pi)^{s-1}}\int_{[0,2\pi]^{s-1}}\Phi_{n}(\omega,\lambda_{2}) \phi(\lambda_{2};f)^{-1}
\left(\prod_{a=2}^{s-1}\phi(\lambda_{a+1};f)^{-1}\Phi_{n}(\lambda_{a},\lambda_{a+1})\right)\times \nonumber\\
&& \bigg(\phi_{t}^{\infty}(\lambda_s;f)\delta_{ s\equiv 1 (\bmod 2)} 
+\phi_{n+1-t}^{\infty}(\lambda_s;f)\delta_{ s\equiv 0(\bmod
   2)}\bigg)d\lambda_{2}d\lambda_{3}\ldots\lambda_{s}.
\end{eqnarray*}
We change notation and set
$u_{s}=\omega,u_{s-1}=\lambda_{2}\ldots,u_{2}=\lambda_{s-1}$
This gives
\begin{eqnarray*}
\zeta_{t,n}^{(s)}(u_{s};f) &=&
  \frac{1}{(2\pi)^{s-1}}\int_{[0,2\pi]^{s-1}}\Phi_{n}(u_{s},u_{s-1}) \phi(u_{s-1};f)^{-1}
\left(\prod_{a=s-1}^{2}\phi(u_{s-a-1};f)^{-1}\Phi_{n}(u_{s-a},u_{s-a-1})\right) \nonumber\\
&&\times \bigg(\phi_{t}^{\infty}(u_1;f)\delta_{ s\equiv 1 (\bmod 2)} 
+\phi_{n+1-t}^{\infty}(u_1;f)\delta_{ s\equiv 0(\bmod
   2)}\bigg)du_{1}du_{2}\ldots du_{s-1}.
\end{eqnarray*}
Again we write the above as (by rearranging the product term)
\begin{eqnarray*}
\zeta_{t,n}^{(s)}(u_{s};f) &=&
  \frac{1}{(2\pi)^{s-1}}\int_{[0,2\pi]^{s-1}}\Phi_{n}(u_{s},u_{s-1}) \phi(u_{s-1};f)^{-1}
\left(\prod_{a=2}^{s-1}\phi(u_{a};f)^{-1}\Phi_{n}(u_{a+1},u_{a})\right)\times \nonumber\\
&& \bigg(\phi_{t}^{\infty}(u_1;f)\delta_{ s\equiv 1 (\bmod 2)} 
+\phi_{n+1-t}^{\infty}(u_1;f)\delta_{ s\equiv 0(\bmod
   2)}\bigg)du_{1}du_{2}\ldots du_{s-2} du_{s-1} \\
&=& \frac{1}{(2\pi)}\int_{[0,2\pi]}\Phi_{n}(u_{s},u_{s-1})
    \phi(u_{s-1};f)^{-1}\frac{1}{(2\pi)^{s-2}}
\bigg\{\int_{[0,2\pi]^{s-2}}
\left(\prod_{a=2}^{s-1}\phi(u_{a};f)^{-1}\Phi_{n}(u_{a+1},u_{a})\right)\nonumber\\
&&\times  \bigg(\phi_{t}^{\infty}(u_1;f)\delta_{ s\equiv 1 (\bmod 2)} 
+\phi_{n+1-t}^{\infty}(u_1;f)\delta_{ s\equiv 0(\bmod
   2)}\bigg)du_{1}du_{2}\ldots du_{s-2}\bigg\} du_{s-1}. \\
\end{eqnarray*}
The term inside the integral is analogus to
$\zeta_{t,n}^{(s-1)}(u_{s-1};f)$ (though it is not this). To obtain
the exactly expression we apply the same procedure to that described above to the inner integral of the
above. This proves  (\ref{eq:recursionZeta2}). 
\hfill $\Box$

\vspace{1em}

It is worth mentioning that analogous to the recursion for
$\zeta_{t,n}^{(s)}(\omega;f)$ a recursion can also be obtained for 
\begin{eqnarray*}
\widehat{J}_{n}^{(s)}(\omega;f) &=&
  \frac{n^{-1/2}}{\phi(\omega;f)}\sum_{t=1}^{n}X_{t}\zeta_{t,n}^{(s)}(\omega;f)
  + e^{i(n+1)\omega}\frac{n^{-1/2}}{\overline{\phi(\omega;f)}}\sum_{t=1}^{n}X_{n+1-t}\overline{\zeta_{t,n}^{(s)}(\omega;f)} \\
&=& \widehat{J}_{L,n}^{(s)}(\omega;f) +e^{i(n+1)\omega}\overline{\widehat{J}_{R,n}^{(s)}(\omega;f)},
\end{eqnarray*}
where 
\begin{eqnarray*}
\widehat{J}_{L,n}^{(s)}(\omega;f) &=&
  \frac{n^{-1/2}}{\phi(\omega;f)}\sum_{t=1}^{n}X_{t}\zeta_{t,n}^{(s)}(\omega;f) \\
\textrm{and} \qquad
\widehat{J}_{R,n}^{(s)}(\omega;f) &=&
\frac{n^{-1/2}}{\phi(\omega;f)}\sum_{t=1}^{n}X_{n+1-t}\zeta_{t,n}^{(s)}(\omega;f).
\end{eqnarray*}
By using the recursion for $\zeta_{t,n}^{(s)}(\omega;f)$ we observe
that for $s\geq 3$ we can write $\widehat{J}_{L,n}^{(s)}(\omega;f) $
and $\widehat{J}_{R,n}^{(s)}(\omega;f) $ as 
\begin{eqnarray*}
\widehat{J}_{L,n}^{(s)}(\omega;f) =
 \frac{1}{(2\pi)^{2}}\int_{[0,2\pi]^{2}}\frac{\Phi_{n}(\omega,y_1)}{\phi(\omega;f)}
\frac{\Phi_{n}(y_1,y_2)}{\phi(y_{1};f)} \widehat{J}_{L,n}^{(s-2)}(y_2;f)dy_{1}dy_{2}
\end{eqnarray*}
and 
\begin{eqnarray*}
\widehat{J}_{R,n}^{(s)}(\omega;f) =
 \frac{1}{(2\pi)^{2}}\int_{[0,2\pi]^{2}}\frac{\Phi_{n}(\omega,y_1)}{\phi(\omega;f)}
\frac{\Phi_{n}(y_1,y_2)}{\phi(y_{1};f)} \widehat{J}_{R,n}^{(s-2)}(y_2;f)dy_{1}dy_{2}.
\end{eqnarray*}

\subsection{Proof of Lemma \ref{lemma:12}}\label{sec:approxproofs2}

\vspace{1em}

We now prove Lemma \ref{lemma:12}. The proof is similar to the proof
of Theorem \ref{theorem:approx}, but with some subtle differences. 
Rather than bounding the best finite
predictors with the best infinite predictors, we bound the best
infinite predictors with the plug-in estimators based on the best
fitting AR$(p)$ parameters. For example, the bounds use the regular Baxter's
inequality rather than the generalized Baxter's inequality. 

\vspace{1em}

\noindent {\bf PROOF of Lemma \ref{lemma:12}} 
We first prove (\ref{eq:BD12}). By using 
the triangular inequality we have 
\begin{eqnarray}
&& \left\|F_{n}^{*}\Delta_{n}(f_{\theta}^{-1})\left(D_{n}(f_{})
- D_{n}(f_p)\right)\right\|_{1} \nonumber \\
&\leq& \left\|F_{n}^{*}\Delta_{n}(f_{\theta}^{-1})\left(D_{n}(f_{})
- D_{\infty,n}(f)\right)\right\|_{1} +
\left\|F_{n}^{*}\Delta_{n}(f_{\theta}^{-1})\left(D_{\infty,n}(f_{})
- D_{n}(f_p)\right)\right\|_{1} \nonumber \\
&& \qquad \qquad \leq \frac{C_{f,0}\rho_{n,K}(f)}{n^{K-1}}A_{K}(f,f_\theta) +
\left\|F_{n}^{*}\Delta_{n}(f_{\theta}^{-1})\left(D_{\infty,n}(f_{})
- D_{n}(f_p)\right)\right\|_{1}, 
\label{eq:lemma4.1}
\end{eqnarray}
where the first term of the right hand side of the above follows from (\ref{eq:approx1A}).
Now we bound the second term on the right hand side of the above. We
observe that since the AR($p$) process only uses the first and last
$p$ observations for the predictions that
$D_{n}(f_{p}) =D_{\infty, n}(f_{p})$, thus we can 
write the second term as 
\begin{equation*}
F_{n}^{*}\Delta_{n}(f_{\theta}^{-1})\left(D_{\infty,n}(f_{})
- D_{n}(f_p)\right) = F_{n}^{*}\Delta_{n}(f_{\theta}^{-1})\left(D_{\infty,n}(f_{})
- D_{\infty,n}(f_p)\right). 
\end{equation*}
Recall that $\{a_{j}(p)\}_{j=1}^{p}$ are the best fitting AR$(p)$ parameters based
on the autocovariance function associated with the spectral density $f$.
Let $a_{p}(\omega) = 1 - \sum_{s=1}^{p}a_{s}(p)e^{-is\omega}$,
$a_{j,p}^{\infty}(\omega) = 1 -
\sum_{s=1}^{p-j}a_{s+j}(p)e^{-is\omega}$ and $a_{p}(\omega)^{-1} =
\psi_{p}(\omega) = \sum_{j=0}^{\infty}\psi_{j,p}e^{-ij\omega}$. By
using the expression for $D_{\infty,n}(f)$ given in
(\ref{eq:widetildeDalt}) we have 
\begin{eqnarray*}
\left[F_{n}^{*}\Delta_{n}(f_\theta^{-1})\left(D_{\infty,n}(f_{})
- D_{\infty,n}(f_p)\right)\right]_{t,j} &=&U^{j,t}_{1,n} + U_{2,n}^{j,t}
\end{eqnarray*}
where 
\begin{eqnarray*}
U^{j,t}_{1,n} &=& \frac{1}{n}\sum_{k=1}^{n}\frac{e^{-it\omega_{k,n}}}{f_{\theta}(\omega_{k,n})}
\left(\frac{\phi_{j}^{\infty}(\omega_{k,n})}{\phi(\omega_{k,n})}
-
\frac{a_{j,p}^{\infty}(\omega_{k,n})}{a_{p}(\omega_{k,n})}\right)\\
U^{j,t}_{2,n}&=&\frac{1}{n}\sum_{k=1}^{n}\frac{e^{-i(t-1)\omega_{k,n}}}{f_{\theta}(\omega_{k,n})}\left(
\frac{\overline{\phi_{n+1-j}^{\infty}(\omega_{k,n})}}{\overline{\phi(\omega_{k,n})}}
-
\frac{\overline{a_{n+1-j,p}^{\infty}(\omega_{k,n})}}{\overline{a_{p}(\omega_{k,n})}}\right).
\end{eqnarray*}
We focus on $U_{1,n}^{j,t}$, and partition it into two terms
$U_{1,n}^{j,t} = U_{1,n,1}^{j,t} + U_{1,n,2}^{j,t}$,
where 
\begin{eqnarray*}
U_{1,n,1}^{j,t} = \frac{1}{n}\sum_{k=1}^{n}\frac{e^{-it\omega_{k,n}}}{\phi(\omega_{k,n})f_{\theta}(\omega_{k,n})}
\left(\phi_{j}^{\infty}(\omega_{k,n})
- a_{j,p}^{\infty}(\omega_{k,n})\right)
\end{eqnarray*}
and 
\begin{eqnarray*}
U_{1,n,2}^{j,t} &=&\frac{1}{n}\sum_{k=1}^{n}\frac{e^{-it\omega_{k,n}}a_{j,p}^{\infty}(\omega_{k,n})}{f_{\theta}(\omega_{k,n})}
\left(\phi(\omega_{k,n})^{-1} - a_{p}(\omega_{k,n})^{-1}\right)\\
&=&\frac{1}{n}\sum_{k=1}^{n}\frac{e^{-it\omega_{k,n}}a_{j,p}^{\infty}(\omega_{k,n})}{f_{\theta}(\omega_{k,n})}
\left(\psi(\omega_{k,n}) - \psi_{p}(\omega_{k,n})\right).
\end{eqnarray*}
We first consider $U_{1,n,1}^{j,t}$. We observe
$\phi(\omega_{k,n})^{-1} =
\psi(\omega_{k,n}) = \sum_{\ell =
0}^{\infty}\psi_{\ell}e^{-i\ell\omega_{k,n}}$. Substituting this
into $U_{1,n,1}^{j,t}$ gives 
\begin{eqnarray*}
U_{1,n,1}^{j,t} &=& \sum_{s=0}^{\infty}\left(\phi_{j+s} - a_{j+s}(p) \right)
\frac{1}{n}\sum_{k=1}^{n}\frac{e^{-i(t+s)\omega_{k,n}}}{\phi(\omega_{k,n})f_{\theta}(\omega_{k,n})}\\
&=&\sum_{s=0}^{\infty}\left(\phi_{j+s} - a_{j+s}(p) \right)\sum_{\ell=0}^{\infty}\psi_{\ell}
\frac{1}{n}\sum_{k=1}^{n}f_{\theta}(\omega_{k,n})^{-1}e^{-i(t+\ell+s)\omega_{k,n}}\\
&=& \sum_{s=0}^{\infty}\left(\phi_{j+s} - a_{j+s}(p) \right)\sum_{\ell=0}^{\infty}\psi_{\ell}
\sum_{r\in \mathbb{Z}}^{}K_{f_\theta^{-1}}(t+\ell+s+rn),
\end{eqnarray*} where 
$K_{f_\theta^{-1}}(r) = \int_{0}^{2\pi}
f_{\theta}(\omega)^{-1}e^{ir\omega}d\omega$.
Therefore, the absolute sum of the above gives 
\begin{eqnarray*}
\sum_{j,t=1}^{n}|U_{1,n,1}^{j,t}| &\leq& \sum_{j,t=1}^{n}
\sum_{s=0}^{\infty}|\phi_{j+s} - a_{j+s}(p)|\sum_{\ell=0}^{\infty}|\psi_{\ell}|
\sum_{r\in \mathbb{Z}}^{}|K_{f_\theta^{-1}}(t+\ell+s+rn)| \\
&=&
\sum_{j=1}^{n}\sum_{s=0}^{\infty}|\phi_{j+s} - a_{j+s}(p) |
\sum_{\ell=0}^{\infty}|\psi_{\ell}|
\sum_{t=1}^{n}\sum_{r\in \mathbb{Z}}|K_{f_{\theta}^{-1}}(t+\ell+s+rn)|\\
&\leq& \left( \sum_{j=1}^{n}\sum_{s=0}^{\infty}|\phi_{j+s} - a_{j+s}(p) | \right)
\|\psi_{f}\|_{0} \sum_{\tau \in \mathbb{Z}}|K_{f_{\theta}^{-1}}(\tau)|\\
&\leq& \left( \sum_{s=1}^{\infty} s|\phi_{s} - a_{s}(p) | \right)
\|\psi_{f}\|_{0} \sum_{\tau \in \mathbb{Z}}|K_{f_{\theta}^{-1}}(\tau)|.
\end{eqnarray*} 
By using (\ref{eq:BoundK}) we have
$\sum_{\tau \in \mathbb{Z}}|K_{f_{\theta}^{-1}}(\tau)| \leq
\sigma_{f_{\theta}}^{-2} \|\phi_{f_{\theta}}\|_{0}^2$. Further, by
using the regular Baxter inequality we have
\begin{equation*}
\sum_{s=1}^{\infty} s|\phi_{s} - a_{s}(p) | \leq
(1+C_{f,1})\sum_{s=p+1}^{\infty} s|\phi_{s}| \leq (1+C_{f,1}) 
p^{-K+1} \rho_{p,K}(f) \|\phi_{f}\|_{K}.
\end{equation*} 
Substituting these two bounds into $\sum_{j,t=1}^{n}|U_{1,n,1}^{j,t}|
$ yields
\begin{equation*}
\sum_{j,t=1}^{n}|U_{1,n,1}^{j,t}| \leq \frac{(1+C_{f,1})}{\sigma_{f_\theta}^2 p^{K-1}} \rho_{p,K}(f) \|\phi_{f}\|_{K} \|\psi_{f}\|_{0} \|\phi_{f_{\theta}}\|_{0}^2.
\end{equation*}
Next we consider the second term $U_{1,n,2}^{j,t}$. Using that 
$\psi(\omega_{k,n}) = \sum_{s=0}^{\infty}\psi_{s}e^{-is\omega}$ and 
$\psi_{p}(\omega_{k,n}) = \sum_{s=0}^{\infty}\psi_{s,p}e^{-is\omega}$
we have 
\begin{eqnarray*}
U_{1,n,2}^{j,t} &=&
\frac{1}{n}\sum_{k=1}^{n}\frac{e^{-it\omega_{k,n}}a_{j,p}^{\infty}(\omega_{k,n})}{f_{\theta}(\omega_{k,n})}
\left(\psi(\omega_{k,n}) - \psi_{p}(\omega_{k,n})\right) \\
&=& \sum_{s=0}^{\infty}\left(\psi_{s} - \psi_{s,p}\right) 
\frac{1}{n}\sum_{k=1}^{n}\frac{e^{-i(t+s)\omega_{k,n}}a_{j,p}^{\infty}(\omega_{k,n})}{f_{\theta}(\omega_{k,n})}\\
&=& \sum_{s=0}^{\infty}\left(\psi_{s} - \psi_{s,p}\right) \sum_{\ell=0}^{\infty}a_{j+\ell}(p)
\frac{1}{n}\sum_{k=1}^{n}\frac{e^{-i(t+s+\ell)\omega_{k,n}}}{f_{\theta}(\omega_{k,n})}
\\
&=& \sum_{s=0}^{\infty}\left(\psi_{s} - \psi_{s,p}\right)
\sum_{\ell=0}^{\infty}a_{j+\ell}(p)\sum_{r\in
\mathbb{Z}}K_{f^{-1}_\theta}(t+s+\ell+rn).
\end{eqnarray*}
Taking the absolute sum of the above gives 
\begin{eqnarray*}
\sum_{j,t=1}^{n}|U_{1,n,2}^{j,t}| 
&\leq& \sum_{j,t=1}^{n}\sum_{s=0}^{\infty}|\psi_{s} - \psi_{s,p}|
\sum_{\ell=0}^{\infty}|a_{j+\ell}(p)|\sum_{r\in
\mathbb{Z}}|K_{f_{\theta}^{-1}}(t+s+\ell+rn)| \\
&=&\sum_{s=0}^{\infty}|\psi_{s} - \psi_{s,p}|
\sum_{j=1}^{n}\sum_{\ell=0}^{\infty}|a_{j+\ell}(p)|\sum_{r\in
\mathbb{Z}}|K_{f_{\theta}^{-1}}(r)| \quad \textrm{(apply the bound (\ref{eq:BoundK})})\\
&\leq& \sigma_{f_\theta}^{-2} \|\phi_{f_\theta}\|_0^{2}\left( \sum_{s=0}^{\infty}|\psi_{s}
- \psi_{s,p}|\right)
\sum_{u=0}^{\infty}|ua_{u}(p)| \\
&\leq& \sigma_{f_\theta}^{-2} \|\phi_{f_\theta}\|_0^{2}\|a_{p}\|_{1}\sum_{s=0}^{\infty}|\psi_{s}
- \psi_{s,p}|. 
\end{eqnarray*}
Next we bound $\|a_{p}\|_{1}$ and $\sum_{s=0}^{\infty}|\psi_{s}
- \psi_{s,p}|$. Let $\phi_{p}(\omega) =
1-\sum_{j=1}^{p}\phi_{j}e^{ij\omega}$ (the truncated AR$(\infty)$
process). Then by applying Baxter's inequality, it is straightforward
to show that
\begin{equation}
\label{eq:BoundA}
\|a_p\|_{1} \leq \|\phi_p\|_{1} + \|a_p - \phi_p\|_{1} \leq (C_{f,1}+1)\|\phi_{f}\|_{1}.
\end{equation}
To bound $\sum_{s=0}^{\infty}|\psi_{s}
- \psi_{s,p}|$ we use the inequality in \cite{p:kre-11}, page 2126
\begin{equation*}
\sum_{s=0}^{\infty}|\psi_{s} - \psi_{s,p}| \leq
\frac{\|\psi_{f}\|_{0}^{2}\cdot
\sum_{j=1}^{\infty}|\phi_{j}-a_{j}(p)|}{ 1 - \|\psi_{f}\|\cdot \|a_{p}
- \phi\|_0}.
\end{equation*}
Applying Baxter's inequality to the numerator of the above gives 
\begin{equation}
\label{eq:BoundPsi}
\sum_{s=0}^{\infty}|\psi_{s} - \psi_{s,p}|
\leq \frac{\|\psi_{f}\|_{0}^{2} (C_{f,0}+1) \rho_{p,K}(f) \|\phi_{f}\|_{K} }{p^{K}(1 - \|\psi_{f}\|_0\cdot \|a_{p} - \phi\|_0)}
\end{equation}
Substituting the bound in (\ref{eq:BoundA}) and (\ref{eq:BoundPsi})
into $\sum_{j,t=1}^{n}|U_{1,n,2}^{j,t}| $ gives 
\begin{equation*}
\sum_{j,t=1}^{n}|U_{1,n,2}^{j,t}| 
\leq \frac{(C_{f,1}+1)^{2}}{\sigma_{f_\theta}^{2}p^{K}} \cdot \frac{
\|\psi_f\|_{0}^{2}\|\phi_{f}\|_{1}
\|\phi_{f}\|_{K}\|\phi_{f_\theta}\|_0^{2} \rho_{p,K}(f)}{
1 - \|\psi_f\|_0 \|a_{p} - \phi\|_0}
\end{equation*}
Altogether, for sufficiently large $p$, where $\|\psi_{f}\|_0\cdot \|a_{p} -
\phi\|_0\leq 1/2$ we have
\begin{eqnarray*}
\sum_{t,j=1}^{n}|U_{1,n}^{j,t}| &\leq&
\frac{(1+C_{f,1})}{\sigma_{f_\theta}^2 p^{K-1}} \rho_{p,K}(f) \|\phi_{f}\|_{K} \|\psi_{f}\|_{0} \|\phi_{f_{\theta}}\|_{0}^2 \\
&&+
\frac{2(C_{f,1}+1)^{2}}{\sigma_{f_\theta}^{2}p^{K}} \|\psi_f\|_{0}^{2}\|\phi_{f}\|_{1}
\|\phi_{f}\|_{K}\|\phi_{f_\theta}\|_0^{2} \rho_{p,K}(f) \\
&\leq& \frac{(C_{f,1}+1)}{\sigma_{f_\theta}^{2}p^{K-1}} 
\rho_{p,K}(f) \|\phi_f\|_{K}\|\phi_{f_\theta}\|_{0}^2 \|\psi_f\|_{0} \left(1+
\frac{2(1+C_{f,1})}{p} \|\psi_f\|_{0} \|\phi_f\|_{1} \right)
\end{eqnarray*}
The same bound holds for $\sum_{t,j=1}^{n}|U_{2,n}^{j,t}|$, thus using (\ref{eq:lemma4.1}) and 
$\rho_{n,K}(f) \leq \rho_{p,K}(f)$ gives 
\begin{eqnarray*}
&& \left\|F_{n}^{*}\Delta_{n}(f_{\theta}^{-1})\left(D_{\infty,n}(f_{})
- D_{n}(f_p)\right)\right\|_{1} \\
&& \qquad \leq \rho_{p,K}(f) A_{K}(f,f_\theta) \left( \frac{(C_{f,1}+1)}{p^{K-1}} +
\frac{2(C_{f,1}+1)^2 }{p^{K}}\|\psi_f\|_{0} \|\phi_f\|_{1}\right). 
\end{eqnarray*}
Substituting the above into (\ref{eq:lemma4.1}) gives (\ref{eq:BD12}).

The proof of (\ref{eq:BD2}) is similar to the proof of Theorem
\ref{theorem:approx}, we omit the details.
\hfill $\Box$



%% file: appendix_baxter.tex
\section{An extension of Baxter's inequalities }\label{sec:baxter}

Let $\{X_{t}\}$ be a second order stationary time series with
absolutely summable autocovariance and spectral density $f$.
We can represent $f$ as 
$f(\omega) = \psi(\omega)\overline{\psi(\omega)} = 1/\left( \phi(\omega)\overline{\phi(\omega)}\right)$
where 
\begin{equation*}
\phi(\omega) = 1 - \sum_{s=1}^{\infty}\phi_{s}e^{-is\omega} \textrm{
and }
\psi(\omega) = 1 + \sum_{s=1}^{\infty}\psi_{s}e^{-is\omega}.
\end{equation*}
Note that $\{\phi_{s}\}$ and $\{\psi_{s}\}$ are the corresponding 
AR$(\infty)$ and MA$(\infty)$ coefficients respectively and 
$\psi(\omega) = \phi(\omega)^{-1}$. To simplify notation we have ignored the variance of the
innovation. 

\subsection{Proof of the extended Baxter inequality}\label{sec:baxter1}

Let $\{\phi_{s,p}(\tau)\}_{s=1}^{p}$ denote the the coefficients of the
best linear predictor of $X_{t+\tau}$ (for $\tau \geq 0$) given $\{X_{s}\}_{t-p}^{t-1}$
\begin{equation}
\label{eq:TimeNormal}
\Ex\left[\left(X_{t+\tau} - \sum_{s=1}^{p}\phi_{s,p}(\tau)X_{t-s}\right)X_{t-k}
\right] = 0\textrm{ for } k = 1,\ldots,p. 
\end{equation}
and $\{\phi_{s}(\tau)\}$ denote the coefficients of the best linear
predictor of $X_{t+\tau}$ given the infinite past $\{X_{s}\}_{s=-\infty}^{t-1}$
\begin{equation}
\label{eq:TimeNormalI}
\Ex\left[\left(X_{t+\tau} - \sum_{s=1}^{\infty}\phi_{s}(\tau)X_{t-s}\right)X_{t-k}
\right] = 0\textrm{ for } k = 1,2,\ldots 
\end{equation}

\noindent Before we begin, we define an appropriate norm on the subspace of $L_{2} {[0,2\pi]}$.

\begin{defin}[Norm on the subspace of $L_{2} {[0,2\pi]}$] \label{def:norm}
Suppose the sequence of positive weights $\{v(k)\}_{k\in \mathbb{Z}}$ satisfies 2 conditions:
(1) $v(n)$ is even, i.e., $v(-n) = v(n)$ for all $n\geq 0$;
(2) $v(n+m) \leq v(n) v(m)$ for all $n,m \in \mathbb{Z}$.

\noindent Given $\{v(k)\}$ satisfies 2 conditions above, define a subspace $A_{v}$ of $L_{2} {[0,2\pi]}$ by
\begin{equation*}
A_{v} = \{ f \in L_{2}[0,2\pi]: \sum_{k \in \mathbb{Z}} v(k)|f_k| < \infty \}.
\end{equation*} where, $f(\omega) = \sum_{k\in \mathbb{Z}}f_{k}e^{ik\omega}$.
We define a norm $\|f\|$ on $A_{v}$ by $\|f\| = \sum_{k \in \mathbb{Z}}v(k)|f_k|$, then
it is easy to check this is a valid norm.
\end{defin}

\begin{remark}[Properties of $\|\cdot\|$] \label{rmk:norm}
Suppose the sequence $\{v(k)\}_{k\in \mathbb{Z}}$ satisfies 2 conditions in Definition \ref{def:norm}, and
define the norm $\|\cdot\|$ with respect to $\{v(k)\}$. Then,
beside the triangle inequality, this norm also
satisfies $\|1\| = v(0)\leq1$, $\|f\| = \| \overline{f} \|$, and $\|fg\|\leq \|f\|\|g\|$ (which
does not hold for all norms but is an important component of the
(extended) Baxter's proof), i.e., $(A_{v}, \|\cdot\|)$ is a Banach algebra with involution operator. The proof for the multiplicative inequality
follows from the fact that 
$(fg)_{k}= \sum_{r}f_{r}g_{k-r}$, where $f_{k}$ and $g_{k}$ are $k$th Fourier coefficient of $f$ and $g$. Thus
\begin{eqnarray*}
\left\|fg\right\| &\leq& 
\sum_{k \in \mathbb{Z}}v(k)\left|\sum_{r\in 
\mathbb{Z}}f_{r}g_{k-r}\right| \\
&\leq& \sum_{k\in \mathbb{Z}}v(r)v(k-r)\left|\sum_{r\in 
\mathbb{Z}}f_{r}g_{k-r}\right| \leq 
\sum_{k,r\in \mathbb{Z}}v(r)v(k-r)|f_{r}||g_{k-r}|= \|f\|\|g\|.
\end{eqnarray*}

\noindent Examples of weights include $v(r)=(2^{q}+|r|^{q})$ 
or $v(r)=(1+|r|)^{q}$ for some $q \geq 0$. In these two examples, when $q=K$, under Assumption \ref{assum:A},
$\psi(\omega),~ \phi(\omega) \in A_{v}$ where $\psi(\omega)=
1+\sum_{j=1}^{\infty}\psi_{j}e^{-ij\omega}$ and $\phi(\omega) = 1 -
\sum_{j=1}^{\infty}\phi_{j}e^{-ij\omega}$ (see
\cite{p:kre-11}). 
\end{remark}
We believe that  Lemma \ref{lemma:BaxterEx} is well known. But as we
could not find a prove we give a proof. 
The proof below follows closely the proof of Baxter (1962, 1963)\nocite{p:bax-62,p:bax-63}.


\noindent {\bf PROOF of Lemma \ref{lemma:BaxterEx}} \hspace{2mm} We use the same proof as Baxter, which is based on rewriting the
normal equations in (\ref{eq:TimeNormal}) within the frequency domain
to yield
\begin{equation*}
\frac{1}{2\pi}\int_{0}^{2\pi}\left(e^{i\tau\omega}-\sum_{s=1}^{p}\phi_{s,p}(\tau)e^{-is\omega}\right)f(\omega)e^{-ik\omega}d\omega
= 0,\textrm{ for } k=1,\ldots,p 
\end{equation*}
Similarly, using the infinite past to do prediction yields the normal equations
\begin{equation*}
\frac{1}{2\pi}\int_{0}^{2\pi}\left(e^{i\tau\omega}-\sum_{s=1}^{\infty}\phi_{s}(\tau)e^{-is\omega}\right)f(\omega)e^{-ik\omega}d\omega
= 0,\textrm{ for } k\geq 1.
\end{equation*}
Thus taking differences of the above two equations for 
$k=1,\ldots,p$ gives
\begin{eqnarray}
&&\frac{1}{2\pi}\int_{0}^{2\pi}\left(\sum_{s=1}^{p}\left[\phi_{s,p}(\tau)-\phi_{s}(\tau)\right]e^{-is\omega}\right)f(\omega)e^{-ik\omega}d\omega 
\nonumber\\
&& \qquad =
\frac{1}{2\pi}\int_{0}^{2\pi}\left(\sum_{s=p+1}^{\infty}\phi_{s}(\tau)e^{-is\omega}\right)f(\omega)e^{-ik\omega}d\omega
\quad 1\leq k \leq p. 
\label{eq:alpha-diff}
\end{eqnarray}
These $p$-equations give rise to Baxter's Weiner-Hopf equations
and allow one to find a bound for
$\sum_{s=1}^{p}\left|\phi_{s,p}(\tau)-\phi_{s}(\tau)\right|$ in terms of
$\sum_{s=p+1}^{\infty}|\phi_{s}(\tau)|$. Interpreting the above, we have two
different functions 
$\left(\sum_{s=1}^{p}\left[\phi_{s,p}(\tau)-\phi_{s}(\tau)\right]e^{-is\omega}\right)f(\omega)$
and $\left(\sum_{s=p+1}^{\infty}\phi_{s}(\tau)e^{-is\omega}\right)f(\omega)$
whose first $p$ Fourier coefficients are the same. 

\noindent Define the polynomials 
\begin{equation}
\label{eq:hpgp}
h_{p}(\omega)
=\sum_{s=1}^{p}\left[\phi_{s,p}(\tau)-\phi_{s}(\tau)\right]e^{-is\omega}
\quad \textrm{and} \quad
g_{p}(\omega) = \sum_{k=1}^{p}g_{k,p}e^{ik\omega}
\end{equation}
where 
\begin{equation}
\label{eq:gsp}
g_{k,p}=
(2\pi)^{-1}\int_{0}^{2\pi}\left(\sum_{s=p+1}^{\infty}\phi_{s}(\tau) e^{-is\omega}\right)f(\omega)e^{-ik\omega}d\omega.
\end{equation}
For the general norm $\|\cdot\|$ defined in Definition \ref{def:norm}, will show that for a sufficiently large $p$, $\|h_p\|\leq
C_{f}\|g_{p}\|$, where the constant $C_{f}$ is a function of the
spectral density (that we will derive).

The Fourier expansion of $h_pf$ is 
\begin{equation*}
h_{p}(\omega)f(\omega) = \sum_{k=-\infty}^{\infty}\widetilde{g}_{k,p}e^{ik\omega},
\end{equation*}
where $\widetilde{g}_{k,p} = (2\pi)^{-1}\int_{0}^{2\pi}h_{p}(\omega)f(\omega)e^{-ik\omega}d\omega$.
Then, by (\ref{eq:alpha-diff}) for $1\leq k\leq p$, $\widetilde{g}_{k,p} =
g_{k,p}$ (where $g_{k,p}$ is defined in (\ref{eq:gsp})). Thus 
\begin{equation}
\label{eq:hp}
h_{p}(\omega)f(\omega) = G_{-\infty}^{0}(\omega) + g_{p}(\omega) + G_{p+1}^{\infty}(\omega)
\end{equation}
where 
\begin{equation*}
G_{-\infty}^{0}(\omega) =
\sum_{k=-\infty}^{0}\widetilde{g}_{k,p}e^{ik\omega} \quad \textrm{
and }\quad
G_{p+1}^{\infty}(\omega) =
\sum_{s=p+1}^{\infty}\widetilde{g}_{k,p}e^{ik\omega}.
\end{equation*}
Dividing by $f^{-1} = \phi\overline{\phi}$ and 
taking the $\|\cdot\|$-norm we have 
\begin{eqnarray}
\left\|h_{p}\right\| &\leq& 
\left\|f^{-1}G_{-\infty}^{0}\right\| + 
\left\|f^{-1} g_{p}\right\| +
\left\|f^{-1}G_{p+1}^{\infty}\right\| \nonumber\\
&\leq&\left\|f^{-1}G_{-\infty}^{0}\right\| + 
\left\|f^{-1} \right\|\left\|g_{p}\right\| +
\left\|f^{-1}G_{p+1}^{\infty}\right\| \nonumber\\
&\leq& \left\|\overline{\phi}\right\| \left\|\phi G_{-\infty}^{0}\right\| + 
\left\|f^{-1} \right\|\left\|g_{p}\right\| +
\left\|\phi
\right\|\left\|\overline{\phi} G_{p+1}^{\infty}\right\|. \label{eq:hp1}
\end{eqnarray}
First we obtain bounds for
$\left\|\phi G_{-\infty}^{0}\right\| $ and 
$\left\|\overline{\phi}
G_{p+1}^{\infty}\right\|$ in terms of $\|g_p\|$.
We will show that for a sufficiently large $p$
\begin{eqnarray*}
\left\|\phi G_{-\infty}^{0} \right\|
&\leq& \left\|\phi \right\|
\left\|g_{p}\right\| + 
\varepsilon
\left\|\overline{\phi} G_{p+1}^{\infty}\right\| \\
\left\|\overline{\phi} G_{p+1}^{\infty}\right\| 
&\leq& 
\left\|\overline{\phi}\right\|\left\|g_{p}\right\|+
\varepsilon
\left\|\phi
G_{-\infty}^{0}\right\|.
\end{eqnarray*}
The bound for these terms hinges on the Fourier coefficients of a function
being unique, which allows us to compare coefficients across functions. 
Some comments are in order that will help in the bounding of the above. We recall
that $f(\omega)^{-1} =\phi(\omega)
\overline{\phi(\omega)}$, where
\begin{equation*}
\phi(\omega) = 1-\sum_{s=1}^{\infty}\phi_{s}e^{-is\omega} \qquad 
\overline{\phi(\omega)} = 1-\sum_{s=1}^{\infty}\phi_{s}e^{is\omega}.
\end{equation*}
Thus $\phi(\omega)G_{-\infty}^{0}(\omega)$
and $\overline{\phi(\omega)}G_{p+1}^{\infty}(\omega)$ have Fourier expansions
with only less than the first and greater than the $p$th frequencies respectively. This
observation gives the important insight into the proof. 
Suppose $b(\omega) = \sum_{j=-\infty}^{\infty}b_{j}e^{ij\omega}$, 
we will make
the use of the notation $\{b(\omega)\}_{+} =
\sum_{j=1}^{\infty}b_{j}e^{ij\omega}$ and $\{b(\omega)\}_{-} =
\sum_{j=-\infty}^{0}b_{j}e^{ij\omega}$, thus $b(\omega) = \{b(\omega)\}_{-}+\{b(\omega)\}_{+}$.

We now return to (\ref{eq:hp}) using that $f =
\psi(\omega)\overline{\psi(\omega)}$ we multiply (\ref{eq:hp}) by 
$\psi(\omega)^{-1}=\phi(\omega)$ to give 
\begin{equation}
\label{eq:hp2}
h_{p}(\omega)\overline{\psi(\omega)} =
\phi(\omega) G_{-\infty}^{0}(\omega) + \phi(\omega) g_{p}(\omega) + 
\phi(\omega) G_{p+1}^{\infty}(\omega).
\end{equation}
Rearranging the above gives 
\begin{equation*}
-\phi(\omega)G_{-\infty}^{0}(\omega) = 
-h_{p}(\omega)\overline{\psi(\omega)} +\phi(\omega)g_{p}(\omega) +
\phi(\omega) G_{p+1}^{\infty}(\omega). 
\end{equation*}
We recall that $h_{p}(\omega)\overline{\psi(\omega)}$ only contain positive
frequencies, whereas 
$\phi(\omega)G_{-\infty}^{0}(\omega)$ only
contains non-positive frequencies. Based on these observations we have
\begin{eqnarray}
&& -\phi(\omega)G_{-\infty}^{0}(\omega) \nonumber \\
&& \qquad = \left\{ -\phi(\omega)G_{-\infty}^{0}(\omega) \right\}_{-}
= \left\{\phi(\omega)g_{p}(\omega)\right\}_{-} +
\left\{\phi(\omega) G_{p+1}^{\infty}(\omega)\right\}_{-}.
\label{eq:hp3}
\end{eqnarray}
We further observe that $G_{p+1}^{\infty}$ only contains non-zero coefficients for
positive frequencies of $p$+1 and greater, thus only the coefficients
of $\phi(\omega)$ with frequencies less or equal to $-(p+1)$ will
give non-positive frequencies when multiplied with
$G_{p+1}^{\infty}$. Therefore 
\begin{equation*}
-\phi(\omega)G_{-\infty}^{-1}(\omega) 
= \left\{\phi(\omega)g_{p}(\omega)\right\}_{-} + 
\left\{\phi_{p+1}^{\infty}(\omega) G_{p+1}^{\infty}(\omega)\right\}_{-},
\end{equation*}
where $\phi_{p+1}^{\infty}(\omega) =
\sum_{s=p+1}^{\infty}\phi_{s}e^{-is\omega}$.
Evaluating the norm of the above (using both the triangle and the
multiplicative inequality) we have 
\begin{eqnarray*}
\left\|\phi G_{-\infty}^{0} \right\|
&\leq& \left\|\phi\right\|\left\|g_{p}\right\| + 
\left\|\phi_{p+1}^{\infty}G_{p+1}^{\infty}\right\| \\
&\leq&\left\|\phi\right\|\left\|g_{p}\right\| + 
\left\|\phi_{p+1}^{\infty}\right\|\left\|\overline{\psi}\right\|
\left\|\overline{\phi} G_{p+1}^{\infty}\right\| ~~\textrm{ since
}\overline{\psi(\omega)} \overline{\phi(\psi)}= 1.
\end{eqnarray*}
This gives a bound 
for $\left\|\phi G_{-\infty}^{0}
\right\|$ in terms of $\|g_{p}\|$ and $\left\|\overline{\phi}
G_{p+1}^{\infty}\right\|$. 
Next we obtain a similar bound for
$\left\|\overline{\phi} G_{p+1}^{\infty}\right\|$ in terms of 
$\left\|g_{p}\right\|$ and 
$\left\|\phi G_{-\infty}^{0}
\right\|$. 

Again using (\ref{eq:hp}), $f(\omega) =
\psi(\omega)\overline{\psi(\omega)}$, but this time multiplying (\ref{eq:hp}) by 
$\overline{\psi(\omega)}^{-1} = \overline{\phi(\omega)}$, we have 
\begin{equation*}
h_{p}(\omega)\psi(\omega) =
\overline{\phi(\omega)}G_{-\infty}^{0}(\omega) +\overline{\phi(\omega)}g_{p}(\omega) + 
\overline{\phi(\omega)}G_{p+1}^{\infty}(\omega).
\end{equation*}
Rearranging the above gives
\begin{equation*}
\overline{\phi(\omega)} G_{p+1}^{\infty}(\omega) =
h_{p}(\omega) \psi(\omega) - 
\overline{\phi(\omega)} G_{-\infty}^{0}(\omega) - \overline{\phi(\omega)} g_{p}(\omega). 
\end{equation*}
We observe that $\overline{\phi(\omega)} G_{p+1}^{\infty}(\omega)$ contains frequencies greater than
$p$ whereas 
$h_{p}(\omega)\psi(\omega)$ only contains frequencies
less or equal to the order $p$ (since $h_p$ is a polynomial up to order
$p$). Therefore multiply $e^{-ip\omega}$ on both side and take $\{\}_{+}$ gives
\begin{eqnarray}
&& e^{-ip\omega}\overline{\phi(\omega)} G_{p+1}^{\infty}(\omega) \nonumber \\
&& \qquad -\left\{e^{-ip\omega} \overline{\phi(\omega)}G_{-\infty}^{0}(\omega)\right\}_{+} -
\left\{e^{-ip\omega} \overline{\phi(\omega)} g_{p}(\omega)\right\}_{+},
\label{eq:hp4}
\end{eqnarray}
By the similar technique from the previous, it is easy to show 
\begin{equation}
\left\{e^{-ip\omega} \overline{\phi(\omega)}G_{-\infty}^{0}(\omega)\right\}_{+}
= \left\{e^{-ip\omega} \overline{\phi_{p+1}^{\infty}(\omega)}G_{-\infty}^{0}(\omega)\right\}_{+}.
\end{equation}
Multiplying $e^{ip\omega}$ and evaluating the $\|\cdot\|$-norm of the above yields the inequality
\begin{eqnarray*}
\left\|\overline{\phi} G_{p+1}^{\infty}\right\| &\leq& 
\left\|\overline{\phi} g_{p}\right\| +
\left\|\overline{\phi_{p+1}^{\infty}}G_{-\infty}^{0}\right\| \\
&\leq& 
\left\|\overline{\phi}\right\| \left\|g_{p}\right\|+
\left\|\overline{\phi_{p+1}^{\infty}}\right\|
\left\|\psi \right\| \left\|\phi G_{-\infty}^{0}\right\|.
\end{eqnarray*}
We note that $\|\phi_{p+1}^{\infty}\| = 
\|\overline{\phi_{p+1}^{\infty}}\|$. For $\phi \in A_{v}$ (see Definition \ref{def:norm} and Remark \ref{rmk:norm}),
$\|\overline{\phi_{p+1}^{\infty}}\| = \sum_{s=p+1}^{\infty}v(s)|\phi_{s}|\rightarrow 0$ as 
$p\rightarrow \infty$, for a large enough $p$, 
$\| \psi(\omega)\|\cdot\|\phi_{p+1}^{\infty}\|<1$. Suppose that 
$p$ is such that $\left\|\phi_{p+1}^{\infty}(\omega)\right\|
\left\| \psi(\omega)\right\|\leq 
\varepsilon<1$, then we have the desired bounds
\begin{eqnarray*}
\left\|\phi G_{-\infty}^{0} \right\|
&\leq& \left\|\phi \right\|
\left\|g_{p}\right\| + 
\varepsilon
\left\|\overline{\phi} G_{p+1}^{\infty}\right\| \\
\left\|\overline{\phi} G_{p+1}^{\infty}\right\| 
&\leq& 
\left\|\overline{\phi}\right\|\left\|g_{p}\right\|+
\varepsilon
\left\|\phi
G_{-\infty}^{0}\right\|.
\end{eqnarray*}
The above implies that
$\left\|\phi G_{-\infty}^{0} \right\| + 
\left\|\overline{\phi} G_{p+1}^{\infty}\right\| \leq
2(1-\varepsilon)^{-1} \left\|\phi \right\|
\left\|g_{p}\right\|$.
Substituting the above in (\ref{eq:hp1}), and using 
that $\|\phi \|\geq 1$ (since $\phi 
= 1-\sum_{s=1}^{\infty}\phi_{s}e^{-is\omega}$, $\|\phi\| \geq \|1\| = v(0) \geq 1$) we have 
\begin{eqnarray*}
\left\|h_{p} \right\| &\leq&
\frac{2\left\|\phi \right\| \left\| g_{p} \right\|
}{1-\varepsilon}+
\left\|f ^{-1} \right\|\left\|g_{p} \right\|\\
&\leq&(1-\varepsilon)^{-1}\left( 
2\left\|\phi \right\| + 
(1-\varepsilon)\left\|\phi \right\|^{2}
\right) \left\|g_{p} \right\| \leq
\frac{3-\varepsilon}{1-\varepsilon}
\left\|\phi \right\|^{2}
\left\|g_{p}\right\|. 
\end{eqnarray*} 
Thus based on the above we have 
\begin{equation}
\label{eq:gp1}
\left\|h_{p} \right\|
\leq \frac{3-\varepsilon}{1-\varepsilon}
\left\|\phi\right\|^{2}\left\|g_{p}\right\|.
\end{equation}
Finally, we obtain a bound for $\|g_{p}\|$ in terms of 
$\sum_{s=p+1}^{\infty}|\phi_{s}(\tau)|$. 
We define an extended version of the function 
$g_{p}(\omega)$. Let $\widetilde{g}_{p}(\omega)= 
\sum_{k \in \mathbb{Z}}g_{k,p}e^{ik\omega}$ 
where $g_{k,p}$ is as in (\ref{eq:hpgp}). By definition,
$\widetilde{g}_{p}(\omega) = \left(\sum_{s=p+1}^{\infty}\phi_{s}(\tau) e^{-is\omega}\right)f(\omega)$
and the Fourier coefficients of $g_{p}(\omega)$ are contained within
$\widetilde{g}_{p}(\omega)$, which implies 
\begin{equation}
\label{eq:gp2}
\left\|g_{p}\right\|\leq \left\|\widetilde{g}_{p}\right\| = 
\left\|\phi_{p+1}^{\infty}(\tau) f \right\|\leq 
\left\|\phi_{p+1}^{\infty}(\tau)\right\| \left\|f \right\| \leq
\left\|\phi_{p+1}^{\infty}\right\| 
\|\psi\|^{2}.
\end{equation} where $\phi_{p+1}^{\infty}(\tau)(\omega) = \sum_{s=p+1}^{\infty}\phi_{s}(\tau) e^{-is\omega}$.
Finally, substituting (\ref{eq:gp2}) into (\ref{eq:gp1}), implies that 
if $p$ is large enough such that 
$\left\|\phi_{p+1}^{\infty}\right\|
\left\| \psi\right\|\leq 
\varepsilon<1$, then 
\begin{equation*}
\left\|h_{p}\right\| \leq \frac{3-\varepsilon}{1-\varepsilon}
\left\|\phi\right\|^{2}\left\|\psi\right\|^{2} 
\left\|\phi_{p+1}^{\infty}(\tau)\right\|. 
\end{equation*}
Thus, if the weights in the norm are $v(m) = (2^K + m^K)$ (it is well-defined weights, see Remark \ref{rmk:norm}) we have
\begin{eqnarray}
&& \sum_{s=1}^{p}(2^{K}+s^{K})\left|\phi_{s,p}(\tau)-\phi_{s}(\tau)\right| \nonumber \\
&& \qquad \leq \frac{3-\varepsilon}{1-\varepsilon}\left\|\phi\right\|_{K}^{2}
\left\|\psi \right\|_{K}^{2}
\sum_{s=p+1}^{\infty}(2^{K}+s^{K})\left| \phi_{s}(\tau)\right|.
\label{eq:baxter2}
\end{eqnarray}
Using Corollary \ref{cor:predictionTau} we have for $\tau\leq 0$
$\phi_{s}(\tau) =
\sum_{j=0}^{\infty}\phi_{s+j}\psi_{|\tau|-j}$ (noting that
$\psi_{j}=0$ for $j<0$), and the desired result. \hfill $\Box$


\subsection{Baxter's inequality on the derivatives of the coefficients}\label{sec:baxterderivatives}

Our aim is to obtain a Baxter-type inequality for the derivatives
of the linear predictors. These bounds will be used when obtaining
expression for the bias of the Gaussian and Whittlelikelihoods.
However, they may also be of independent
interest. It is interesting to note that the
  following result can be used to show that the Gaussian and Whittle
  likelihood estimators are asymptotically equivalent in the sense
  that $\sqrt{n}|\widehat{\theta}_{n}^{(G)} -
  \widehat{\theta}_{n}^{(K)}|_{1}\Pcon 0$ as $n\rightarrow \infty$.

The proof of the result is based on the novel proof strategy developed
in Theorem 3.2 of \cite{p:kre-17}  (for spatial processes).
We require the following definitions. 
Define the two $n$-dimension vectors  
\begin{eqnarray}
\underline{\varphi}_{n}(\tau;f_{\theta}) &=& 
\left(\phi_{1,n}(\tau;f_{\theta}),\ldots, \phi_{n,n}(\tau;f_{\theta})\right)^{\prime} \quad (\textrm{best linear
                                        finite future predictor})\label{eq:varphidef}\\
\underline{\phi}_{n}(\tau;f_{\theta}) &=& 
\left(\phi_{1}(\tau;f_{\theta}),\ldots, \phi_{n}(\tau;f_{\theta})\right)^{\prime} \quad
                                        (\textrm{truncated best linear
                                        infinite future predictor})\nonumber.
\end{eqnarray}

\begin{lemma}\label{lemma:PHI}
Let $\theta$ be a $d$-dimension vector. Let $\{c_{\theta}(r)\}$, $\{\phi_{j}(f_{\theta})\}$ and $\{\psi_{j}(f_{\theta})\}$
  denote the autocovariances, AR$(\infty)$, and  MA$(\infty)$ coefficients
  corresponding to the spectral density  $f_\theta$.
  For all $\theta \in \Theta$ and for $0\leq i \leq \kappa$ we assume
\begin{equation}
\label{eq:armaAssumption}
\sum_{j=1}^{\infty}\|j^{K}\nabla_{\theta}^{i}\phi_{j}(f_{\theta})\|_{1}<\infty\quad
\sum_{j=1}^{\infty}\|j^{K}\nabla_{\theta}^{i}\psi_{j}(f_{\theta})\|_{1}<\infty,\quad 
\end{equation} 
where $K>1$. Let $\underline{\varphi}_{n}(\tau;f_{\theta})$ and $\underline{\phi}_{n}(\tau;f_{\theta})$,
be defined as in (\ref{eq:varphidef}). 
We assume that $\tau \leq 0$. Then for all $0\leq i \leq \kappa$, we have 
\begin{eqnarray*}
\left\|\frac{\partial^{i}}{\partial\theta_{r_1}\ldots\partial
  \theta_{r_i}}\left[ \underline{\varphi}_{n}(\tau;f_{\theta}) - \underline{\phi}_{n}(\tau;f_{\theta})\right]
\right\|_2 &\leq&
       f_0\bigg(\sum_{\stackrel{a_1+a_2=i}{a_2\neq  i}}
C_{a_1} \binom{i}{a_1}
\left\|\nabla_{\theta}^{a_2}[\underline{\varphi}_{n}(\tau;f_{\theta})-
  \underline{\phi}_{n}(\tau;f_{\theta})]\right\|_{2} \nonumber\\
&& +\sum_{b_1+b_2=i}
C_{b_1} \binom{i}{b_1} \sum_{j=n+1}^{\infty}
\left\|\nabla_{\theta}^{b_2}\phi_{j}(\tau;f_{\theta})\right\|_1
\bigg),
\end{eqnarray*}
where 
$f_0 = (\inf_{\omega} f_{\theta}(\omega))^{-1}$ and  $C_{a} =
\sum_{r}\|\nabla_{\theta}^{a} c_\theta(r)\|_{1}$,
$\nabla_{\theta}^{a}g(f_{\theta})$ is the $a$th order partial derivative of
$g$ with respect to $\theta=(\theta_{1},\ldots,\theta_{d})$ and
$\|\nabla_{\theta}^{a}g(f_{\theta})\|_{p}$ denotes the $\ell_p-$norm of the matrix with elements containing all the
partial derivatives in $\nabla_{\theta}^{a}g(f_{\theta})$. 
\end{lemma}
\noindent PROOF. To prove the result, we
define the  $n$-dimension vector
\begin{eqnarray}
\underline{c}_{n,\tau} =
                                    \left(c(\tau-1),c(\tau-2),\ldots,c(\tau-n)\right)^{\prime} \quad
                                    (\textrm{covariances from lag
                                    $\tau-1$ to lag $\tau-n$}). \nonumber
\end{eqnarray}
To simplify notation we drop the $f_{\theta}$ notation from the prediction
coefficients $\phi_{j,n}(\tau;f_{\theta})$ and $\phi_{j}(\tau;f_{\theta})$.

\vspace{2mm}
\noindent \underline{Proof for the case $i=0$} This
  is the regular  Baxter inequality but with the $\ell_2$-norm rather than
$\ell_1$-norm.  We recall that for $\tau\leq 0$ we have the best linear predictors 
\begin{eqnarray*}
\widehat{X}_{\tau,n} = \sum_{j=1}^{n}\phi_{j,n}(\tau)X_{j}
\quad \text{and} \quad
 X_{\tau} = \sum_{j=1}^{\infty}\phi_{j}(\tau)X_{j}. 
\end{eqnarray*}
Thus by evaluating the covariance of the above with $X_{r}$ for all $1\leq r
\leq n$ gives the sequence of $r$ normal equation, which can be written
in matrix form
\begin{equation*}
\Gamma_{n}(f_\theta)\underline{\varphi}_{n}(\tau) =
                                                   \underline{c}_{n,\tau} \quad \textrm{and}\quad
\Gamma_{n}(f_\theta)\underline{\phi}_{n}(\tau) + \sum_{j=n+1}^{\infty} \phi_{j}(\tau)\underline{c}_{n,-j+\tau} =
                                                   \underline{c}_{n,\tau}.
\end{equation*}
Taking differences of the above gives 
 \begin{eqnarray}
\Gamma_{n}(f_\theta)\left[\underline{\varphi}_{n}(\tau)
  -\underline{\phi}_{n}(\tau) \right] &=&
\sum_{j=n+1}^{\infty}\phi_{j}(\tau)\underline{c}_{n,-j+\tau} \nonumber \\
\Rightarrow \quad \left[\underline{\varphi}_{n}(\tau)
  -\underline{\phi}_{n}(\tau) \right] &=&\Gamma_{n}(f_\theta)^{-1}
\sum_{j=n+1}^{\infty}\phi_{j}(\tau)\underline{c}_{n,-j+\tau} 
\label{eq:GammaExpand}
\end{eqnarray}
The $\ell_2$-norm of the above gives 
\begin{equation*}
\left\|\underline{\varphi}_{n}(\tau)
  -\underline{\phi}_{n}(\tau) \right\|_{2} \leq  \|\Gamma_{n}(f_\theta)^{-1}\|_{spec}
\sum_{j=n+1}^{\infty}|\phi_{j}(\tau)|\cdot\|\underline{c}_{n,-j+\tau}\|_{2}.
\end{equation*}
To bound the above we use the 
well known result $\|\Gamma_{n}(f_\theta)^{-1}\|_{spec} \leq
1/\inf_{\omega}f(\omega)=f_0$ and $\|\underline{c}_{n,-j+\tau}
\|_{2} \leq \sum_{r\in \mathbb{Z}}|c_\theta(r)|=C_{0}$. This gives the bound 
\begin{equation*}
\left\|\underline{\varphi}_{n}(\tau)
  -\underline{\phi}_{n}(\tau) \right\|_{2} \leq  f_{0}C_{0}\sum_{j=n+1}^{\infty}|\phi_{j}(\tau)|.
\end{equation*}

\noindent \underline{Proof for the case $i=1$.}
As our aim is to bound the derivative of the difference 
$\underline{\varphi}_{n}(\tau)
  -\underline{\phi}_{n}(\tau)$, we evaluate the partial derivative of
(\ref{eq:GammaExpand}) with respective to $\theta_{r}$ and isolate  $\partial [\underline{\varphi}_{n}(\tau)
  -\underline{\phi}_{n}(\tau)]/\partial \theta_{r}$. 
Differentiating both sides of (\ref{eq:GammaExpand}) with respect to
$\theta_r$ gives 
 \begin{eqnarray}
\label{eq:Gammaderiv}
&& \frac{\partial \Gamma_{n}(f_\theta)}{\partial \theta_{r}}
\left[\underline{\varphi}_{n}(\tau)
   -\underline{\phi}_{n}(\tau) \right] +
\Gamma_{n}(f_\theta)\frac{\partial }{\partial \theta_r}
\left[\underline{\varphi}_{n}(\tau)
   -\underline{\phi}_{n}(\tau) \right] \nonumber\\
&& \quad =
    \sum_{j=n+1}^{\infty}\left[\frac{\partial \phi_{j}(\tau)}{\partial \theta_r}\underline{c}_{n,-j+\tau}+
\phi_{j}(\tau)\frac{\partial \underline{c}_{n,-j+\tau}}{\partial \theta_r}\right]. 
\end{eqnarray}
Isolating $\partial [\underline{\varphi}_{n}(\tau)
  -\underline{\phi}_{n}(\tau)]/\partial \theta_{r}$ gives 
\begin{eqnarray}
\frac{\partial }{\partial \theta_r}
\left[\underline{\varphi}_{n}(\tau)
   -\underline{\phi}_{n}(\tau) \right] &=&
                                           -\Gamma_{n}(f_\theta)^{-1}\frac{\partial
                                           \Gamma_{n}(f_\theta)}{\partial
                                           \theta_r}
\left[\underline{\varphi}_{n}(\tau)
   -\underline{\phi}_{n}(\tau) \right] \nonumber\\
&& +
    \Gamma_{n}(f_\theta)^{-1}\sum_{j=n+1}^{\infty}\left[\frac{\partial
   \phi_{j}(\tau)}{\partial \theta_r}\underline{c}_{n,-j+\tau}+
\phi_{j}(\tau)\frac{\partial \underline{c}_{n,-j+\tau}}{\partial
   \theta_r}\right]. \label{eq:GammaInverse}
\end{eqnarray}
Evaluating
the $\ell_{2}$ norm of the above and using
$\|AB\underline{x}\|_{2}\leq
\|A\|_{spec}\|B\|_{spec}\|\underline{x}\|_{2}$ gives the bound
\begin{eqnarray}
&&\left\|\frac{\partial }{\partial \theta_r}
\left[\underline{\varphi}_{n}(\tau)
   -\underline{\phi}_{n}(\tau) \right]\right\|_{2}
\nonumber\\
&& \leq \|\Gamma_{n}(f_\theta)^{-1}\|_{spec}
   \bigg(\left\|\frac{\partial \Gamma_{n}(f_\theta)}{\partial \theta_r}\right\|_{spec}
\left\|\underline{\varphi}_{n}(\tau)
   -\underline{\phi}_{n}(\tau) \right\|_{2}\nonumber\\
&&\qquad\qquad\qquad\qquad +
\sum_{j=n+1}^{\infty}\left|\frac{\partial \phi_{j}(\tau)}{\partial\theta_r}\right|\|\underline{c}_{n,-j+\tau}\|_{2}+ 
\sum_{j=n+1}^{\infty}|\phi_{j}(\tau)|\left\|\frac{\partial
   \underline{c}_{n,-j+\tau}}{\partial \theta_r}\right\|_{2}\bigg) \nonumber\\
&&\leq f_0
   \bigg(\left\|\frac{\partial \Gamma_{n}(f_\theta)}{\partial \theta_r}\right\|_{spec}
\left\|\underline{\varphi}_{n}(\tau)
   -\underline{\phi}_{n}(\tau) \right\|_{2}
\nonumber\\
&&\qquad\qquad\quad +C_0
\sum_{j=n+1}^{\infty}\left|\frac{\partial
   \phi_{j}(\tau)}{\partial\theta_r}\right|+ 
\left(\sum_{r\in \mathbb{Z}}\|\nabla_{\theta}c_{\theta}(r)\|_2\right)
\sum_{j=n+1}^{\infty}|\phi_{j}(\tau)|\bigg) \nonumber\\
&&\leq  f_0
   \bigg(\left\|\frac{\partial \Gamma_{n}(f_\theta)}{\partial \theta_r}\right\|_{spec}
\left\|\underline{\varphi}_{n}(\tau)
   -\underline{\phi}_{n}(\tau) \right\|_{2}
 +C_{0}
\sum_{j=n+1}^{\infty}\left|\frac{\partial
   \phi_{j}(\tau)}{\partial\theta_r}\right|+ 
C_{1}
\sum_{j=n+1}^{\infty}|\phi_{j}(\tau)|\bigg) \qquad \label{eq:B26}
\end{eqnarray}
where the last line in the above uses the bound
$\sum_{r\in
    \mathbb{Z}}\|\nabla_{\theta}^{a}c_{\theta}(r)\|_2 \leq \sum_{r\in
    \mathbb{Z}}\|\nabla_{\theta}^{a}c_{\theta}(r)\|_1=C_{a}$ (for
  $a=0$ and $1$). We
  require a bound for  $\|\partial \Gamma_{n}(f_\theta)/\partial
  \theta_r\|_{spec}$. Since $\Gamma_{n}(f_\theta)$ is a symmetric Toeplitz
  matrix, then $\partial \Gamma_{n}(f_\theta)/\partial \theta_r$ is
  also a symmetric Toeplitz matrix (though not necessarily positive definite) with entries 
\begin{eqnarray*}
\left[\frac{\partial \Gamma_{n}(f_\theta)}{\partial
  \theta_{r}}\right]_{s,t} =
  \frac{1}{2\pi}\int_{0}^{2\pi}\frac{\partial
  f_{\theta}(\omega)}{\partial \theta_{r}}\exp(i(s-t)\omega)d\omega.
\end{eqnarray*}
We mention that the symmetry is clear, since $\frac{\partial
  c(s-t;f_{\theta})}{\partial \theta_{r}} = \frac{\partial
  c(t-s;f_{\theta})}{\partial \theta_{r}}$. Since the matrix is
symmetric the spectral norm is the spectral radius. This gives 
\begin{eqnarray*}
\left\|\frac{\partial \Gamma_{n}(f_\theta)}{\partial
  \theta_{r}}\right\|_{spec} &=& 
\sup_{\|x\|_{2}=1}\left|\sum_{s,t=1}^{n}x_{s}x_{t}\frac{\partial
  c(s-t;f_{\theta})}{\partial \theta_{r}} \right| 
=
    \sup_{\|x\|_{2}=1}\left|\frac{1}{2\pi}\int_{0}^{2\pi}|\sum_{s=1}^{n}x_{s}e^{is\omega}|^{2}
\frac{\partial f_{\theta}(\omega)}{\partial \theta_{r}}d\omega \right| \\
    &\leq& \sup_{\omega}\left|\frac{\partial
    f_{\theta}(\omega)}{\partial \theta_{r}}\right|\sup_{\|x\|_{2}=1}
\frac{1}{2\pi}\int_{0}^{2\pi}\left|\sum_{s=1}^{n}x_{s}e^{is\omega}\right|^{2}d\omega
           = \sup_{\omega}\left|\frac{\partial
    f_{\theta}(\omega)}{\partial \theta_{r}}\right|.
\end{eqnarray*}
By using the same argument one can show that the $a$th derivative is
\begin{eqnarray}
\label{eq:specB}
\left\|\frac{\partial^{a} \Gamma_{n}(f_\theta)}{\partial
    \theta_{r_1}\ldots\partial \theta_{r_a}}\right\|_{spec}\leq \sup_{\omega}\left|\frac{\partial^{a}
  f_{\theta}(\omega)}{\partial \theta_{r_1}\ldots\partial
  \theta_{r_a}}\right| = 
\sum_{r\in \mathbb{Z}}\left|\frac{\partial^{a} c(r;f_{\theta})}{\partial 
    \theta_{r_1}\ldots\partial \theta_{r_a}}\right| \leq C_{a}.
\end{eqnarray} 
This general bound will be useful when evaluating the higher order
derivatives below.
Substituting (\ref{eq:specB}) into (\ref{eq:B26})
gives 
\begin{eqnarray*}
\left\|\frac{\partial }{\partial \theta_r}
\left[\underline{\varphi}_{n}(\tau)
   -\underline{\phi}_{n}(\tau) \right]\right\|_{2}
 \leq f_0\bigg(C_{1}\left\|\underline{\varphi}_{n}(\tau)-\underline{\phi}_{n}(\tau)\right\|_{2}
+C_{1}\sum_{j=n+1}^{\infty}|\phi_{j}(\tau)| +
C_{0}\sum_{j=n+1}^{\infty}\left\|\nabla_{\theta}\phi_{j}(\tau)\right\|_1
\bigg).
\end{eqnarray*}
This proves the result for $i=1$. 

\noindent \underline{Proof for the case $i=2$} We
differentiate both sides of (\ref{eq:Gammaderiv}) with respect to $\theta_{r_2}$
 to give the second derivative 
\begin{eqnarray}
&&\left(\frac{\partial^{2}\Gamma_{n}(f_\theta)}{\partial \theta_{r_1}\partial \theta_{r_2}}\right)
\left[\underline{\varphi}_{n}(\tau)
   -\underline{\phi}_{n}(\tau) \right] +
\frac{\partial\Gamma_{n}(f_\theta)}{\partial
   \theta_{r_2}}\frac{\partial }{\partial \theta_{r_1}}
\left[\underline{\varphi}_{n}(\tau)
   -\underline{\phi}_{n}(\tau) \right] +\nonumber \\
&& \frac{\partial\Gamma_{n}(f_\theta)}{\partial
   \theta_{r_1}}\frac{\partial }{\partial \theta_{r_2}}
\left[\underline{\varphi}_{n}(\tau)
   -\underline{\phi}_{n}(\tau) \right] +
\Gamma_{n}(f_\theta)
\left(\frac{\partial^{2}}{\partial\theta_{r_2}\partial \theta_{r_1}}\left[\underline{\varphi}_{n}(\tau)
   -\underline{\phi}_{n}(\tau) \right] \right)\nonumber\\
&&  =
    \sum_{j=n+1}^{\infty}\bigg[
\frac{\partial^{2}\phi_{j}(\tau)}{\partial \theta_{r_2}\partial \theta_{r_1}}\underline{c}_{n,-j+\tau}+
\frac{\partial \phi_{j}(\tau)}{\partial \theta_{r_1}}\frac{\partial
   \underline{c}_{n,-j+\tau}}{\partial \theta_{r_2}}+ 
 \frac{\partial \phi_{j}(\tau)}{\partial \theta_{r_2}}\frac{\partial
   \underline{c}_{n,-j+\tau}}{\partial \theta_{r_1}}+
\phi_{j}(\tau)\frac{\partial^{2}\underline{c}_{n,-j+\tau}}{\partial
   \theta_{r_2}\partial \theta_{r_1}}\bigg]. \nonumber\\
&& 
\end{eqnarray}
Rearranging the above to isolate
$\frac{\partial^{2}}{\partial\theta_{r_2}
\partial \theta_{r_1}}\left[\underline{\varphi}_{n}(\tau)
   -\underline{\phi}_{n}(\tau) \right]$ gives
\begin{eqnarray*}
&&\frac{\partial^{2}}{\partial \theta_{r_1}\partial \theta_{r_2}}\left[\underline{\varphi}_{n}(\tau)
   -\underline{\phi}_{n}(\tau) \right] \nonumber\\
&&= -\Gamma_{n}(f_\theta)^{-1}
\left(\frac{\partial^{2}\Gamma_{n}(f_\theta)}{\partial \theta_{r_1}\partial \theta_{r_2}}\right)
\left[\underline{\varphi}_{n}(\tau)
   -\underline{\phi}_{n}(\tau) \right] -
\Gamma_{n}(f_\theta)^{-1}\frac{\partial\Gamma_{n}(f_\theta)}{\partial
   \theta_{r_2}}\frac{\partial }{\partial \theta_{r_1}}
\left[\underline{\varphi}_{n}(\tau)
   -\underline{\phi}_{n}(\tau) \right] \nonumber \\
&&\quad - 	\Gamma_{n}(f_\theta)^{-1}\frac{\partial\Gamma_{n}(f_\theta)}{\partial
   \theta_{r_1}}\frac{\partial }{\partial \theta_{r_2}}
\left[\underline{\varphi}_{n}(\tau)
   -\underline{\phi}_{n}(\tau) \right] \nonumber\\
&&\quad   +\Gamma_{n}(f_\theta)^{-1}
    \sum_{j=n+1}^{\infty}\bigg[
\frac{\partial^{2}\phi_{j}(\tau)}{\partial \theta_{r_2}\partial \theta_{r_1}}\underline{c}_{n,-j+\tau}+
\frac{\partial \phi_{j}(\tau)}{\partial \theta_{r_1}}\frac{\partial
   \underline{c}_{n,-j+\tau}}{\partial \theta_{r_2}} \nonumber\\
&&\quad + \frac{\partial \phi_{j}(\tau)}{\partial \theta_{r_2}}\frac{\partial
   \underline{c}_{n,-j+\tau}}{\partial \theta_{r_1}}+
\phi_{j}(\tau)\frac{\partial^{2}\underline{c}_{n,-j+\tau}}{\partial \theta_{r_2}\partial \theta_{r_1}}\bigg]. 
\end{eqnarray*}
Taking the $\ell_{2}$-norm of $\frac{\partial^{2}}{\partial \theta_{r_1}\partial \theta_{r_2}}\left[\underline{\varphi}_{n}(\tau)
   -\underline{\phi}_{n}(\tau) \right]$ and using (\ref{eq:specB})
 gives 
\begin{eqnarray*}
&&\left\|\frac{\partial^{2}}{\partial \theta_{r_1}\partial \theta_{r_2}}
\left[\underline{\varphi}_{n}(\tau)   -\underline{\phi}_{n}(\tau) \right]\right\|_{2} \\
&& \leq f_{0}\bigg(
C_2\left\| \underline{\varphi}_{n}(\tau)
   -\underline{\phi}_{n}(\tau) \right\|_{2}+
2C_{1}\left\|\nabla_{\theta} [\underline{\varphi}_{n}(\tau)
   -\underline{\phi}_{n}(\tau)]\right\|_{2} \nonumber\\
&&   \quad +C_{0}\sum_{j=n+1}^{\infty}\left\|\nabla^{2}\phi_{j}(\tau)\right\|_1
+ 2C_{1}\sum_{j=n+1}^{\infty}\left\|\nabla \phi_{j}(\tau)\right\|_1 +
C_{2}\sum_{j=n+1}^{\infty}|\phi_{j}(\tau)|\bigg). \nonumber\\
\end{eqnarray*}
This proves the result for $i=2$. The proof for $i > 2$ follows using
a similar argument (we omit the details). \hfill $\Box$

\vspace{3mm}
\noindent The above result gives an $\ell_2$-bound between the
derivatives of the finite and infinite predictors. However, for our
purposes an $\ell_1$-bound is more useful. Thus we use the Cauchy-Schwarz
inequality and norm inequality $\|\cdot\|_{2}\leq \|\cdot\|_{1}$
to give the $\ell_{1}$-bound
\begin{eqnarray}
\label{eq:BaxterL0}
&& \left\|\frac{\partial^{i}}{\partial\theta_{r_1}\ldots\partial
  \theta_{r_i}}\left[ \underline{\varphi}_{n}(\tau;f_{\theta}) - \underline{\phi}_{n}(\tau;f_{\theta})\right]
\right\|_1 \nonumber \\ 
&& \qquad \leq
       n^{1/2}f_0\bigg(\sum_{\stackrel{a_1+a_2=i}{a_2\neq  i}}
\binom{i}{a_1} C_{a_1}
\left\|\nabla_{\theta}^{a_2}[\underline{\varphi}_{n}(\tau;f_{\theta})-
  \underline{\phi}_{n}(\tau;f_{\theta})]\right\|_{1} + \nonumber\\
&&\qquad \qquad \qquad  \sum_{b_1+b_2=i}
\binom{i}{b_1}C_{b_1}\sum_{j=n+1}^{\infty}
\left\|\nabla_{\theta}^{b_2}\phi_{j}(\tau;f_{\theta})\right\|_1
\bigg),
\end{eqnarray}
this incurs an additional $n^{1/2}$ term. Next, considering all the
partial derivatives with respect to $\theta$ of order $i$ and using (\ref{eq:varphidef})
we have
\begin{eqnarray}
\label{eq:BaxterL1}
&& \sum_{t=1}^{n}\|\nabla_{\theta}^{i}[\phi_{t,n}(\tau;f_{\theta})-\phi_{t}(\tau;f_{\theta})]\|_{1}\nonumber\\
&& \qquad \leq
      d^{i} n^{1/2}f_0\bigg(\sum_{\stackrel{a_1+a_2=i}{a_2\neq  i}}
\binom{i}{a_1} C_{a_1}
\left\|\nabla_{\theta}^{a_2}[\underline{\varphi}_{n}(\tau;f_{\theta})-
  \underline{\phi}_{n}(\tau;f_{\theta})]\right\|_{1} + \nonumber\\
&& \qquad \qquad \qquad \sum_{b_1+b_2=i}
\binom{i}{b_1}C_{b_1}\sum_{j=n+1}^{\infty}
\left\|\nabla_{\theta}^{b_2}\phi_{j}(\tau;f_{\theta})\right\|_1
\bigg),
\end{eqnarray}
where $d$ is the dimension of the vector $\theta$.
The above gives a bound in terms of the infinite predictors. We now
obtain a bound in terms of the corresponding AR$(\infty)$ and
MA$(\infty)$ coefficients. To do this, 
we recall that for $\tau \leq0$,
$\phi_{j}(\tau;f_{\theta}) =
\sum_{s=0}^{\infty}\phi_{s+j}(f_\theta)\psi_{|\tau|-j}(f_\theta)$. Thus
the partial derivatives of $\phi_{j}(\tau;f_{\theta})$ give the bound
\begin{equation*}
\sum_{j=n+1}^{\infty}\left\|\nabla_{\theta}\phi_{j}(\tau;f_\theta)\right\|_{1} \leq 
\sum_{s=0}^{\infty}\sum_{j=n+1}^{\infty}\left(|\psi_{|\tau|-j}(f_\theta)|\cdot
\|\nabla_{\theta}\phi_{s+j}(f_\theta)\|_{1} + |\phi_{s+j}(f_\theta)|\cdot\|\nabla_{\theta}\psi_{|\tau|-j}(f_\theta)\|_{1}
\right).
\end{equation*}
Substituting the above bound into (\ref{eq:BaxterL1}) and using Lemma \ref{lemma:BaxterEx} for the case $i=1$ gives
\begin{eqnarray}
&&\sum_{s=1}^{n}\left
\|\nabla_{\theta}[\phi_{s,n}(\tau;f_\theta) - \phi_{s}(\tau;f_\theta)]
\right\|_{1} \nonumber\\
&\leq&  n^{1/2}df_{0} \bigg\{C_1(C_{0} + 1)
\sum_{j=n+1}^{\infty}|\phi_{j}(\tau;f_\theta)| + \nonumber\\
&& 
C_{0}\sum_{s=0}^{\infty}\sum_{j=n+1}^{\infty}\left(|\psi_{|\tau|-j}(f_\theta)|\cdot
\|\nabla_{\theta}\phi_{s+j}(f_\theta)\|_{1} + |\phi_{s+j}(f_\theta)|\cdot\|\nabla_{\theta}\psi_{|\tau|-j}(f_\theta)\|_{1}
\right)\bigg\}.
\label{eq:LLLB}
\end{eqnarray}
The above results are used to obtain bounds between the derivatives of the Whittle
and Gaussian likelihood in Appendix \ref{sec:prelim}. 
Similar bounds can also be obtained for the higher order derivatives  $\sum_{s=1}^{n}\left
\|\nabla_{\theta}^{i}[\phi_{s,n}(\tau;f_\theta) - \phi_{s}(\tau;f_\theta)]
\right\|_{1}$ in terms of the derivatives of the MA$(\infty)$ and
AR$(\infty)$ coefficients. 

\begin{remark}
It is interesting to note that it is  possible to prove
analgous bounds to those in
Lemma \ref{lemma:PHI} 
using the derivatives of the series expansion
in the proof of Theorem \ref{thm:higherorder}, that is
\begin{eqnarray*}
\frac{\partial }{\partial \theta_{r}}[\phi_{j,n}(\tau;f_{\theta}) -
  \phi_{j}(\tau;f_{\theta})] =  
\sum_{s=2}^{\infty}\frac{\partial }{\partial
  \theta_{r}}\phi_{j,n}^{(s)}(\tau;f_{\theta}). 
\end{eqnarray*}
This would give bounds on
\begin{eqnarray*}
\sum_{\tau \leq 0}\frac{\partial }{\partial \theta_{r}}[\phi_{j,n}(\tau;f_{\theta}) -\phi_{j}(\tau;f_{\theta})] e^{i\tau\omega} 
=  \phi(\omega)^{-1}\sum_{s=2}^{\infty} \frac{\partial }{\partial
  \theta_{r}}
\left[\phi(\omega;f_{\theta})^{-1}\zeta_{j,n}^{(s)}(\omega; f_{\theta})\right].
\end{eqnarray*}
This approach would immediately give bounds based on the derivatives
of the AR$(\infty)$ coefficients. 
\end{remark}

We now state and prove a lemma which will be useful in a later
section (it is not directly related to Baxter's inequality).

\begin{lemma}\label{lemma:logcoeff}
Let $\{\phi_{j}(f_{\theta})\}$ and $\{\psi_{j}(f_{\theta})\}$
  denote the AR$(\infty)$, and  MA$(\infty)$ coefficients
  corresponding to the spectral density  $f_\theta$.
Suppose the same set of Assumptions in Lemma \ref{lemma:PHI} holds. 
Let  
$\{\alpha_{j}(f_\theta)\}$ denote the Fourier coefficients in the one-sided expansion
\begin{eqnarray}
\label{eq:log_MA}
\log (\sum_{j=0}^{\infty} \psi_j(f_\theta) z^{j})=
 \sum_{j=1}^{\infty} \alpha_j(f_\theta) z^j  
\quad \text{for } |z|<1,
\end{eqnarray}
Then for all $\theta \in \Theta$ and for  $0\leq s \leq \kappa$ we have 
\begin{equation*}
\sum_{j=1}^{\infty}\|j^{K}\nabla_{\theta}^{s}\alpha_{j}(f_\theta)\|_1<\infty,
\end{equation*} 
 \end{lemma}
PROOF. 
We first consider the case  $s=0$. The derivative of (\ref{eq:log_MA})
with respect to $z$ together with
$\psi(z;f_\theta)^{-1} = 
\phi(z;f_\theta) = 1-\sum_{j=1}^{\infty} \phi_j(f_\theta)z^{j}$ gives
\begin{eqnarray} 
\sum_{j=1}^{\infty} j\alpha_j(f_\theta) z^{j-1} &=& \left( \sum_{j=1}^{\infty}j \psi_j(f_\theta) z^{j-1} \right) \left(\sum_{j=0}^{\infty} \psi_{j} (f_\theta) z^{j} \right)^{-1} \nonumber \\
&=&
\left( \sum_{j=1}^{\infty}j \psi_j(f_\theta) z^{j-1} \right) \left(\sum_{j=0}^{\infty} \widetilde{\phi}_{j} (f_\theta) z^{j} \right),
\label{eq:log_identity}
\end{eqnarray}
where $\sum_{j=0}^{\infty} \widetilde{\phi}_{j} (f_\theta) z^{j} = 1- \sum_{j=1}^{\infty} \phi_j(f_\theta) z^j$.
 Comparing the coefficients of $z^{j-1}$ from both side of above
 yields the identity
\begin{eqnarray} 
&& j\alpha_j(f_\theta) = \sum_{\ell=0}^{j-1} (j-\ell)\psi_{j-\ell}
  (f_\theta) \widetilde{\phi}_{\ell}(f_\theta) \nonumber \\
\Rightarrow && \alpha_j(f_\theta) = j^{-1}\sum_{\ell=0}^{j-1} (j-\ell)\psi_{j-\ell}
  (f_\theta) \widetilde{\phi}_{\ell}(f_\theta) \quad \textrm{for} \quad j\geq 1. 
\label{eq:log_coeff_theta}
\end{eqnarray} 
Therefore, using the above and taking the absolute into the summand we have
\begin{eqnarray*}
\sum_{j=1}^{\infty} j^{K} |\alpha_j(f_\theta)| 
&\leq& \sum_{j=1}^{\infty} \sum_{\ell=0}^{j-1} j^{K-1} (j-\ell) |\psi_{j-\ell}
  (f_\theta)| |\widetilde{\phi}_{\ell}(f_\theta)| \\
&=& \sum_{\ell=1}^{\infty} |\widetilde{\phi}_{\ell}(f_\theta)| \sum_{j=\ell+1}^{\infty} j^{K-1} (j-\ell) |\psi_{j-\ell}
  (f_\theta)|  \qquad \text{(exchange summation)} \\
&=& \sum_{\ell=1}^{\infty} |\widetilde{\phi}_{\ell}(f_\theta)| \sum_{s=1}^{\infty} (s+\ell)^{K-1} s |\psi_{s}
  (f_\theta)|  \qquad \text{(change of variable $s=j+\ell$)} \\
&\leq& \sum_{\ell=1}^{\infty} |\widetilde{\phi}_{\ell}(f_\theta)| \sum_{s=1}^{\infty} (s+\ell)^{K} |\psi_{s} (f_\theta)|.
\qquad ((s+\ell)^{-1} \leq s^{-1})
\end{eqnarray*} Since $K \geq 1$, using inequality $(a+b)^{K} \leq 2^{K-1} (a^K + b^K)$ for $a,b>0$, we have
\begin{eqnarray*}
\sum_{j=1}^{\infty} j^{K} |\alpha_j(f_\theta)| 
&\leq& \sum_{\ell=1}^{\infty} |\widetilde{\phi}_{\ell}(f_\theta)| \sum_{s=1}^{\infty} (s+\ell)^{K} |\psi_{s} (f_\theta)| \\
&\leq& 2^{K-1} \sum_{\ell=1}^{\infty} |\widetilde{\phi}_{\ell}(f_\theta)| \sum_{s=1}^{\infty} (s^K+\ell^K) |\psi_{s} (f_\theta)| \\
&=& 2^{K-1} \bigg( \sum_{\ell=1}^{\infty} \ell^K |\widetilde{\phi}_{\ell}(f_\theta)| \cdot \sum_{s=1}^{\infty} |\psi_{s} (f_\theta)|
+ \sum_{\ell=1}^{\infty} |\widetilde{\phi}_{\ell}(f_\theta)| \cdot \sum_{s=1}^{\infty} s^{K}|\psi_{s} (f_\theta)| \bigg) \leq \infty
\end{eqnarray*}
and this proves the lemma when $s=0$.

To prove lemma for $s=1$, we differentiate (\ref{eq:log_coeff_theta}) with $\theta$, then, by Assumption \ref{assum:B} (iii),
\begin{eqnarray*}
&& \sum_{j=1}^{\infty} j^{K}\|\nabla_{\theta}\alpha_j(f_\theta)\|_1 \\
&& \quad \leq \sum_{j=1}^{\infty} \sum_{\ell=0}^{j-1} \| j^{K-1}(j-\ell) \nabla_{\theta}\psi_{j-\ell} (f_\theta) \widetilde{\phi}_{\ell}(f_\theta)\|_1 
+ \sum_{j=1}^{\infty} \sum_{\ell=0}^{j-1} \| j^{K-1}(j-\ell) \psi_{j-\ell} (f_\theta) \nabla_{\theta}\widetilde{\phi}_{\ell}(f_\theta)\|_1 \\
\end{eqnarray*} 
Using similar technique to prove $s=0$, we show $\sum_{j=1}^{\infty} \|j^{K} \nabla_{\theta} \alpha_j(f_\theta)\|_1 <\infty$ and the
proof for $s\geq 2$ is similar (we omit the detail).
\hfill $\Box$

\vspace{1em}

\subsection{The difference between the derivatives of the Gaussian and
  Whittle likelihoods}\label{sec:prelim}

We now obtain an expression for the difference between the derivatives of the Gaussian
likelihood and the Whittle likelihood. These expression will be used
later for obtaining the bias of the Gaussian likelihood (as compared
with the Whittle likelihood).

For the Gaussian likelihood, we have shown in Theorem \ref{theorem:gauss} that 
\begin{eqnarray*}
\underline{X}_{n}^{\prime}\Gamma_{n}(\theta)^{-1}\underline{X}_{n}
  =
  \underline{X}_{n}^{\prime}F_{n}^{*}\Delta_{n}(f_\theta^{-1})F_{n}\underline{X}_{n}
  + \underline{X}_{n}^{\prime}F_{n}^{*}\Delta_{n}(f_\theta^{-1})D_{n}(f_\theta)\underline{X}_{n},
\end{eqnarray*}
where the first term is the Whittle likelihood and the second term the
additional term due to the Gaussian likelihood. Clearly the
derivative with respect to $\theta^{\prime} =
(\theta_{1},\ldots,\theta_{d})$ is 
\begin{eqnarray*}
\underline{X}_{n}^{\prime}\nabla_{\theta}^{i}\Gamma_{n}(\theta)^{-1}\underline{X}_{n}
  =
  \underline{X}_{n}^{\prime}F_{n}^{*}\nabla_{\theta}^{i}\Delta_{n}(f_\theta^{-1})F_{n}\underline{X}_{n}
  + \underline{X}_{n}^{\prime}F_{n}^{*}\nabla_{\theta}^{i}[\Delta_{n}(f_\theta^{-1})D_{n}(f_\theta)]\underline{X}_{n}.
\end{eqnarray*}
The first term on the right hand side is the derivative of the Whittle
likelihood with respect
to $\theta$, the second term is the additional term due to the
Gaussian likelihood. 

For the simplicity, assume $\theta$ is univariate. Our objective in the next few lemmas is to show that 
\begin{eqnarray*}
\left\|\underline{X}_{n}^{\prime}F_{n}^{*}\frac{d^{i}}{d\theta^{i}}\Delta_{n}(f_\theta^{-1})D_{n}(f_\theta)\underline{X}_{n}\right\|_{1}
  = O_{}(1),
\end{eqnarray*}
which is a result analogous to Theorem \ref{theorem:Bound}, but for
the derivatives. We will use this result to prove Theorem
\ref{theorem:bias}, in particular to show the derivatives of the Whittle likelihood and
the Gaussian likelihood (after normalization by $n^{-1}$) differ by $O(n^{-1})$. 

Just as in the proof of Theorem \ref{theorem:Bound}, the derivative of 
this term with respect to $\theta$ does not (usually) have a simple
analytic form. Therefore, analogous to Theorem \ref{theorem:approx} it is
easier to replace the derivatives of $D_{n}(f_\theta)$ with the derivatives
of $D_{\infty,n}(f_\theta)$, and show that the replacement error is
``small''. 

\begin{lemma}\label{lemma:DDeriv}
Suppose Assumption \ref{assum:B}(i),(iii)
holds and $g$ is a bounded function. Then for $1\leq i \leq 3$ we have
\begin{equation}
\label{eq:DDerivS}
\left\|F_{n}^{*}\Delta_{n}(g) \frac{d^{i}}{d\theta^{i}}\left(D_{n}(f_{\theta})-
  D_{\infty,n}(f_{\theta})\right)\right\|_{1}
  = O(n^{-K+3/2}),
\end{equation}
and 
\begin{eqnarray}
\label{eq:DDerivSS}
\left\|
F_{n}^{*}\sum_{k=0}^{i}\binom{i}{k}
\frac{d^{k}\Delta_{n}\left(f_{\theta}^{-1}\right)}{d\theta^k}\frac{d^{k - i}D_{\infty,n}(f_\theta)}{d\theta^{k-i}}\right\|_{1} = O(1).
\end{eqnarray}
\end{lemma}
\noindent PROOF.
To bound (\ref{eq:DDerivS}), we  use the expression for $F_{n}^{*}\Delta_{n}(g^{-1}) \left(D_{n}(f_{\theta})-
  D_{\infty,n}(f_{\theta})\right)$ given in (\ref{eq:DDeltaDtildeD}) 
\begin{eqnarray*}
\left(F_{n}^{*}\Delta_{n}(g^{-1})\left[D_{n}(f_\theta) -
D_{\infty,n}(f_\theta)\right]\right)_{s,t}
&=&  \sum_{\tau \leq 0}[ \left\{\phi_{t,n}(\tau;f_\theta)-\phi_{t}(\tau;f_\theta)\right\}
G_{1,n}(s,\tau;g) \nonumber \\
&& \quad \quad + \left\{\phi_{n+1-t,n}(\tau;f_\theta)-\phi_{n+1-t}(\tau;f_\theta)\right\} G_{2,n}(s,\tau;g) ].
\end{eqnarray*}
Differentiating the above with respect to $\theta$ gives
\begin{eqnarray*}
&&\left[F_{n}^{*}\Delta_{n}(g) \frac{d^{}}{d\theta^{}}\left(D_{n}(f_{\theta})-
  D_{\infty,n}(f_{\theta})\right)\right]_{s,t}\\ 
&& \quad= 
\sum_{\tau \leq 0}\left[G_{1,n}(s,\tau) \frac{d}{d\theta}[\phi_{t,n}(\tau)-\phi_{t}(\tau)]
    +G_{2,n}(s,\tau)
\frac{d}{d\theta}\left[\phi_{n+1-t}(\tau)-\phi_{n+1-t}(\tau)\right]
   \right] \\
&&\quad = T_{s,t,1} + T_{s,t,2}.
\end{eqnarray*}
We recall that equation (\ref{eq:LLLB}) gives the bound
\begin{eqnarray*}
&&\sum_{s=1}^{n}
\left|\frac{d}{d\theta}[\phi_{s,n}(\tau;f_\theta) - \phi_{s}(\tau;f_\theta)]
\right| \nonumber\\
&\leq&  n^{1/2}f_{0} \bigg\{C_1(C_{f,0} + 1)
\sum_{j=n+1}^{\infty}|\phi_{j}(\tau;f_\theta)| + \nonumber\\
&& 
C_{0}\sum_{s=0}^{\infty}\sum_{j=n+1}^{\infty}\left(|\psi_{|\tau|-j}(f_\theta)|\cdot
\left|\frac{d}{d\theta}\phi_{s+j}(f_\theta)\right| + |\phi_{s+j}(f_\theta)|\cdot\left|\frac{d}{d\theta}\nabla_{\theta}\psi_{|\tau|-j}(f_\theta)\right|
\right)\bigg\}.
\end{eqnarray*}
Substituting this into $T_{s,t,1}$ gives the bound
\begin{eqnarray*}
|T_{s,t,1}| &\leq& 
Cn^{1/2}\sum_{\tau \leq 0}G_{1,n}(s,\tau) 
\bigg(\sum_{j=n+1}^{\infty}\sum_{s=0}^{\infty}|\phi_{s+j}||\psi_{|\tau|-j}|
              +
\sum_{s=0}^{\infty}\sum_{j=n+1}^{\infty}\left|\frac{d\phi_{s+j}}{d\theta}\psi_{|\tau|-j}
+\phi_{s+j}\frac{d\psi_{|\tau|-j}}{d\theta}
\right| \bigg).
\end{eqnarray*}
Using the same techniques used to prove Theorem
\ref{theorem:approx} yields
\begin{eqnarray*}
\sum_{s,t=1}^{n}|T_{s,t,1}| = O\left(n^{1/2}n^{-K+1}\right) = O(n^{-K+3/2}).
\end{eqnarray*}
Similarly, we can show that
$\sum_{s,t=1}^{n}|T_{s,t,2}| = O\left(n^{1/2}n^{-K+1}\right) = O(n^{-K+3/2})$.
Altogether this gives 
\begin{equation*}
\left\|F_{n}^{*}\Delta_{n}(g) \left(D_{n}(f_{\theta})-
  D_{\infty,n}(f_{\theta})\right)\right\|_{1} = O(n^{-K+3/2}).
\end{equation*}
This proves (\ref{eq:DDerivS}) for the case $i=1$. The
proof for  the cases $i=2, 3$ is similar.

To prove (\ref{eq:DDerivSS}) we use the same method used to prove
Theorem \ref{theorem:approx}, equation (\ref{eq:approx1A}). But with 
$\frac{d^{k}f_{\theta}^{-1}}{d\theta^{k}}$ replacing
$f_{\theta}$ in $\Delta_{n}(\cdot)$ and $\frac{d^{i-k}}{d\theta^k}\phi_{j}(\tau;f_\theta) =
\frac{d^{i-k}}{d\theta^k}\sum_{s=0}^{\infty}\phi_{s+j}\psi_{|\tau|-j}$
replacing $\phi_{j}(\tau;f_\theta) =
\sum_{s=0}^{\infty}\phi_{s+j}\psi_{|\tau|-j}$ in
$D_{n}(f_{\theta})$. We omit the details. \hfill $\Box$

\vspace{2mm}

\noindent 
We now apply the above results to quadratic forms of random variables.
\begin{corollary}
\label{corollary:DDeriv}
Suppose Assumptions \ref{assum:B} (i),(iii) hold 
 and $g$
is a bounded function. Further, if $\{X_{t}\}$ is a time
series where $\sup_{t} \|X_{t}\|_{\Ex,2q} = \|X\|_{\Ex,2q}<\infty$ (for some
$q>1$), then
\begin{equation}
\label{eq:DDerivSC}
\left\|\frac{1}{n}\underline{X}_n^{\prime}\left[F_{n}^{*}\Delta_{n}(g) \frac{d^{i}}{d\theta^{i}}\left(D_{n}(f_{\theta})-
  D_{\infty,n}(f_{\theta})\right)\right]\underline{X}_n\right\|_{\Ex,q}
  = O(n^{-K+1/2}),
\end{equation}
and 
\begin{eqnarray}
\label{eq:DDerivSSC}
\left\|\frac{1}{n}\underline{X}_{n}^{\prime}
F_{n}^{*}\sum_{\ell=0}^{i}\binom{i}{\ell}
\frac{d^{\ell}\Delta_{n}\left(f_{\theta}^{-1}\right)}{d\theta^{\ell}}\frac{d^{\ell-i}D_{\infty,n}(f_\theta)}{d\theta^{\ell-i}}
\underline{X}_n\right\|_{\Ex,q} = O(n^{-1})
\end{eqnarray}
for $i=1,2$ and $3$. 
\end{corollary}
PROOF. To prove (\ref{eq:DDerivSC}), we observe that
\begin{eqnarray*}
&&\left\|\frac{1}{n}\underline{X}_n^{\prime}\left[F_{n}^{*}\Delta_{n}(g) \frac{d^{i}}{d\theta^{i}}\left(D_{n}(f_{\theta})-
  D_{\infty,n}(f_{\theta})\right)\right]\underline{X}_n\right\|_{\Ex,q}
  \\
&&\quad \leq \frac{1}{n}\sum_{s,t=1}^{n}\left|
\left[F_{n}^{*}\Delta_{n}(g) \frac{d^{i}}{d\theta^{i}}\left(D_{n}(f_{\theta})-
  D_{\infty,n}(f_{\theta})\right)\right]_{s,t}\right|\|X_{s}X_{t}\|_{\Ex,q} \\
&&\quad =\frac{1}{n}\sup_{t}\|X_{t}\|_{\Ex,2q}^{2}\sum_{s,t=1}^{n}\left|
\left[F_{n}^{*}\Delta_{n}(g) \frac{d^{i}}{d\theta^{i}}\left(D_{n}(f_{\theta})-
  D_{\infty,n}(f_{\theta})\right)\right]_{s,t}\right| \\
&&\quad = \frac{1}{n} \|X\|_{\Ex,2q}^{2}\left\|
F_{n}^{*}\Delta_{n}(g) \frac{d^{i}}{d\theta^{i}}\left(D_{n}(f_{\theta})-
  D_{\infty,n}(f_{\theta})\right)\right\|_{1} = O(n^{-K+1/2})
\end{eqnarray*}
where the above follows from Lemma \ref{lemma:DDeriv}, equation
(\ref{eq:DDerivS}). This proves (\ref{eq:DDerivSC}).

To prove (\ref{eq:DDerivSSC}) we use the the bound in
(\ref{eq:DDerivSS}) together with a similar proof to that described above. This
immediately proves (\ref{eq:DDerivSSC}). \hfill $\Box$

\vspace{1em}

We now apply the above result to the difference in the derivatives of
the Gaussian and Whittle likelihood. It is straightforward to show that
\begin{eqnarray}
\label{eq:DDDexpand}
&&\underline{X}_{n}^{\prime}F_{n}^{*}\frac{d^{i}}{d\theta^{i}}\Delta_{n}(f_\theta^{-1})D_{n}(f_\theta) \underline{X}_{n} \\
&&= \underline{X}_{n}^{\prime}F_{n}^{*}\left[\sum_{\ell=0}^{i}\binom{i}{\ell}
\frac{d^{\ell}\Delta_{n}\left(f_{\theta}^{-1}\right)}{d\theta^{\ell}}\frac{d^{\ell-i}D_{n}(f_\theta)}{d\theta^{\ell-i}}\right]
\underline{X}_{n}
  \nonumber\\
&&= \underline{X}_{n}^{\prime}F_{n}^{*}\left[\sum_{\ell=0}^{i}\binom{i}{\ell}
\frac{d^{\ell}\Delta_{n}\left(f_{\theta}^{-1}\right)}{d\theta^{\ell}}\frac{d^{\ell-i}D_{\infty,n}(f_\theta)}{d\theta^{\ell-i}}\right]
\underline{X}_{n} \nonumber\\
&& +\underline{X}_{n}^{\prime}F_{n}^{*}\left(
\sum_{\ell=0}^{i}\binom{i}{\ell}
\frac{d^{\ell}\Delta_{n}\left(f_{\theta}^{-1}\right)}{d\theta^{\ell}}
\left[\frac{d^{\ell-i}D_{n}(f_\theta)}{d\theta^{\ell-i}}-\frac{d^{\ell-i}D_{\infty,n}(f_\theta)}{d\theta^{\ell-i}}\right]\right)\underline{X}_{n}.
\end{eqnarray}
First we study the second term on the right hand side of the above. 
By applying Corollary \ref{corollary:DDeriv} (and under
Assumption \ref{assum:B}) for  
$1\leq i \leq 3$ we have
\begin{eqnarray}
\label{eq:DiffDerivatives0}
&&\left\|n^{-1}
\underline{X}_{n}^{\prime}F_{n}^{*}\left(\frac{d^{i}}{d\theta^{i}}\Delta_{n}(f_\theta^{-1})\left[D_{n}(f_\theta)
  - D_{\infty,n}(f_\theta)\right]\right)
\underline{X}_{n}\right\|_{\Ex,1} \nonumber\\
&&=\left\|n^{-1}\underline{X}_{n}^{\prime}F_{n}^{*}
\sum_{\ell=0}^{i}\binom{i}{\ell}
\frac{d^{\ell}\Delta_{n}\left(f_{\theta}^{-1}\right)}{d\theta^{\ell}}
\left[\frac{d^{\ell-i}D_{n}(f_\theta)}{d\theta^{\ell-i}}-\frac{d^{\ell-i}D_{\infty,n}(f_\theta)}{d\theta^{\ell-i}}\right]\underline{X}_{n}
\right\|_{\Ex,1} \nonumber \\
&&= O(n^{-K+1/2}). 
\end{eqnarray}
On the other hand, the first term on the right hand side of
(\ref{eq:DDDexpand}) has the bound
\begin{eqnarray}
\label{eq:DiffDerivatives}
\left\|n^{-1}
\underline{X}_{n}^{\prime}F_{n}^{*}\frac{d^{i}}{d\theta^{i}}\left[\Delta_{n}(f_\theta^{-1})D_{n}(f_\theta)\right]\underline{X}_{n}\right\|_{\Ex,1}
  = O(n^{-1}).
\end{eqnarray}

%% file: appendix_consistency.tex
\section{Rates of convergence of the new likelihood estimators}\label{sec:consistent}

In this section we study the sampling properties of the new criteria. 

\subsection{The criteria}

To begin with, we state the assumptions required to  obtain rates of
convergence of the new criteria and asymptotic equivalence to the
infeasible criteria. These results will be used to derive the
asymptotic sampling properties of the new likelihood estimators,
including their asymptotic bias (in a later section). To do this, we start
by defining the criteria we will be considering. 

We assume that $\{X_{t}\}$ is a stationary time series with spectral density $f$, where $f$ is
bounded away from zero (and bounded above). We fit the model with
spectral density $f_{\theta}$ to the observed time series. We do not
necessarily assume that there exists a $\theta_{0} \in \Theta$
where $f = f_{\theta_{0}}$. 
Since we allow the misspecified case, for a given $n$, it seems
natural that the ``ideal'' best fitting parameter is 
\begin{equation}
\label{eq:thetan}
\theta_{n} =
  \arg\min_{\theta}I_{n}(f,f_\theta).
\end{equation}
where $I_{n}(f,f_\theta)$ is defined in (\ref{eq:Ifthetaf}).
Note that in the case the spectral density is correctly specified, 
then $\theta_{n}=\theta_{0}$ for all $n$ where $f = f_{\theta_{0}}$. 

We now show that Assumption \ref{assum:B}(ii,iii) 
allows us to ignore the 
$n^{-1}\sum_{k=1}^{n}\log f_{\theta}(\omega_{k,n})$ in the Whittle,
boundary corrected Whittle and hybrid Whittle likelihoods. 
To show why this is true, we obtain the  
Fourier expansion of $\log f_\theta(\omega) 
= \sum_{r \in \mathbb{Z}} \alpha_r(f_\theta) e^{ir \omega}$, where $
\alpha_0(f_\theta)=\log \sigma^{2}$, in terms of
the corresponding MA$(\infty)$ coefficients. We use the well known 
Szeg{\"o}'s  identity 
\begin{eqnarray*}
\log f_\theta(\cdot) = \log \sigma^2
|\psi(\cdot;f_\theta)|^2 = \log \sigma^2 + \log \psi(\cdot;f_\theta) + 
\log \overline{\psi(\cdot;f_\theta)}
\end{eqnarray*}
where $\psi(\omega;f_\theta) =
\sum_{j=0}^{\infty} \psi_j(f_\theta)e^{-ij\omega}$ with
$\psi_0(f_\theta)=1$ and the roots of the MA transfer function
$\sum_{j=0}^{\infty} \psi_j(f_\theta)z^{j}$ lie outside the unit
circle (minimum phased). Comparing 
\\$\log f_\theta(\omega) 
= \sum_{r \in \mathbb{Z}} \alpha_r(f_\theta) e^{ir \omega}$ with 
the positive half of the above expansion gives 
\begin{eqnarray*}
\log (\sum_{j=0}^{\infty} \psi_j(f_\theta) z^{j})=
 \sum_{j=1}^{\infty} \alpha_j(f_\theta) z^j  
\quad \text{for } |z|<1,
\end{eqnarray*}
and since $\log f_\theta$ is real and symmetric about $\pi$, $\alpha_{-j}(f_\theta) = \alpha_{j}(f_\theta) \in \mathbb{R}$.
This allows us to obtain coefficients $\{\alpha_{j}(f_\theta)\}$
in terms of the MA$(\infty)$ coefficients (it is interesting to note that \cite{b:pou-01}
gives a recursion for $\alpha_{j}(f_\theta)$ in terms of the
MA$(\infty)$ coefficient). The result is given in Lemma
\ref{lemma:logcoeff}, but we summarize it below.  Under Assumption \ref{assum:B}(iii)
we have for $0 \leq s \leq \kappa$ (for some $\kappa \geq 4$)
\begin{equation*}
\sum_{j=1}^{\infty}\|j^{K}\nabla_{\theta}^{s}\alpha_{j}(f_\theta)\|_1<\infty.
\end{equation*}
Using this result, we bound $n^{-1}\sum_{k=1}^{n}\log
f_{\theta}(\omega_{k,n})$. Applying the Poisson summation formula to
this sum we have
\begin{eqnarray}
\frac{1}{n} \sum_{k=1}^{n} \log f_{\theta}(\omega_{k,n}) =
\sum_{r \in \mathbb{Z}} \alpha_{rn} (f_\theta) &=& \alpha_{0}(f_\theta) +
\sum_{r \in\mathbb{Z}\setminus \{0\} }\alpha_{rn}(f_\theta) \nonumber\\
&=& \log(\sigma^2) + \sum_{r \in \mathbb{Z} \setminus \{0\} }
    \alpha_{rn} (f_\theta). \label{eq:logf}
\end{eqnarray}
The $s$th-order derivative ($s\geq 1$) with respect to $\theta$ (and using 
Assumption \ref{assum:B}(ii) that $\sigma^{2}$ does not depend on $\theta$) we have 
\begin{eqnarray*}
\frac{1}{n} \sum_{k=1}^{n} \nabla_{\theta}^{s}\log f_{\theta}(\omega_{k,n}) 
= \sum_{r \in \mathbb{Z} \setminus \{0\} } \nabla_{\theta}^{s}\alpha_{rn} (f_\theta).
\end{eqnarray*}
By using Lemma \ref{lemma:logcoeff} 
for $0\leq s\leq \kappa$ we have 
\begin{eqnarray}
\label{eq:alphaBound}
\| \sum_{r \in \mathbb{Z} \setminus \{0\} } \nabla_{\theta}^{s}\alpha_{rn} (f_\theta)\|_1 \leq
2 \sum_{j\geq n} \|\nabla_{\theta}^{s} \alpha_{j} (f_\theta)\|_1 = O(n^{-K}).
\end{eqnarray}
Substituting the bound in (\ref{eq:alphaBound})  (for $s=0$) into
(\ref{eq:logf}) gives 
\begin{eqnarray*}
\left|\frac{1}{n} \sum_{k=1}^{n} \log f_{\theta}(\omega_{k,n})
  -\log \sigma^2 \right| = O(n^{-K}).
\end{eqnarray*}
Using (\ref{eq:alphaBound})  for $1\leq s \leq \kappa$ we have 
\begin{eqnarray*}
 \left\| \frac{1}{n} \sum_{k=1}^{n}\nabla_{\theta}^{s} \log f_{\theta}(\omega_{k,n}) \right\|_1 = O(n^{-K}).
\end{eqnarray*} 
Therefore if $K>1$, the log deterministic 
term in the Whittle, boundary corrected, and hybrid Whittle likelihood is negligible
as compared with $O(n^{-1})$ (which we show is the leading order in the bias).

However, for the Gaussian likelihood, the log determinant cannot be ignored.
Specifically, by applying  the strong Szeg{\"o}'s theorem (see e.g., Theorem 10.29 of
\cite{b:bot-13}) to $\Gamma_n(f_\theta)$ we have 
\begin{eqnarray*}
\frac{1}{n} \log |\Gamma_{n}(f_\theta)| = \log \sigma^2 + \frac{1}{n}E(\theta) + o(n^{-1}) 
\end{eqnarray*}where 
$E(\theta) = \sum_{k=1}^{\infty} \alpha_{k}(f_\theta)^2$. Therefore,
unlike the other three
quasi-likelihoods, the error in
\\$\log |\Gamma_{n}(f_\theta)|$  is of order $O(n^{-1})$, which is of
the same order as the bias. In Section \ref{sec:AR1bias}, we show that
 the inclusion and exclusion of $n^{-1} \log
 |\Gamma_{n}(f_\theta)|$ leads to Gaussian likelihood estimators with substantial
 differences in their bias.  Further, there is no clear rule whether the
 inclusion of the $n^{-1}\log |\Gamma_{n}(f_\theta)|$ in the Gaussian
 likelihood improves the bias or makes it worse. In the case that
 $n^{-1}\log |\Gamma_{n}(f_\theta)|$ is included in the Gaussian
 likelihood, then the expression for the bias  will include 
the derivatives of $E(\theta)$. Except for a few simple models (such
as the AR$(1)$ model) the expression for the derivatives of
$E(\theta)$ will be extremely unwieldy. 

Based on the above, to make the derivations cleaner, we define all the
quasi-likelihoods without the log term and let 
\begin{eqnarray}
\mathcal{L}_{n}(\theta) &=&
                            n^{-1}\underline{X}_{n}^{\prime}\Gamma_{n}(f_\theta)^{-1}\underline{X}_{n}
  =  \frac{1}{n}\sum_{k=1}^{n}\frac{\widetilde{J}_{n}(\omega_{k,n};f_\theta) \overline{J_{n}(\omega_{k,n})}}{f_\theta(\omega_{k,n})}\nonumber\\
K_{n}(\theta) &=&\frac{1}{n}\sum_{k=1}^{n}\frac{|J_{n}(\omega_{k,n})|^{2}}{f_\theta(\omega_{k,n})}
                            \nonumber\\
\widehat{W}_{p,n}(\theta) &=& 
                    \frac{1}{n}\sum_{k=1}^{n}\frac{\widetilde{J}_{n}(\omega_{k,n};\widehat{f}_p)\overline{J_{n}(\omega_{k,n})}}{f_\theta(\omega_{k,n})}
  \nonumber\\
\widehat{H}_{p,n}(\theta) &=&
  \frac{1}{n}\sum_{k=1}^{n}\frac{\widetilde{J}_{n}(\omega_{k,n};
 \widehat{f}_{p})\overline{J_{n,\underline{h}_{n}}(\omega_{k,n})}
                              }{f_{\theta}(\omega_{k,n})}.
\label{eq:LIKE}
\end{eqnarray} 
In the case of the hybrid Whittle likelihood, 
we make the assumption the data taper 
$\{h_{t,n}\}$ is such that $h_{t,n} = c_{n}h_{n}(t/n)$ where
$c_{n} = n/H_{1,n}$ and $h_{n}:[0,1] \rightarrow \mathbb{R}$ is a sequence of taper functions
 which satisfy the taper conditions in Section 5, \cite{p:dah-88}.

We define the parameter estimators as 
\begin{eqnarray}
\widehat{\theta}_{n}^{(G)} &=& \arg\min \mathcal{L}_{n}(\theta), \quad 
\widehat{\theta}_{n}^{(K)} = \arg\min K_{n}(\theta), \nonumber \\
\widehat{\theta}_{n}^{(W)} &=& \arg\min \widehat{W}_{p,n}(\theta),
                             \textrm{ and  } 
\widehat{\theta}_{n}^{(H)} = \arg\min \widehat{H}_{p,n}(\theta) 
\label{eq:thetaBiased}
\end{eqnarray}

\subsection{Asymptotic equivalence to the infeasible criteria}

In this section we  analyze the feasible estimators 
$\widehat{\theta}_{n}^{(W)}=\arg\min \widehat{W}_{p,n}(\theta)$ and 
$\widehat{\theta}_{n}^{(H)}=\arg\min \widehat{H}_{p,n}(\theta)$. We
show that they are asymptotically equivalent to the corresponding infeasible estimators 
$\widetilde{\theta}_{n}^{(W)}=\arg\min W_{n}(\theta)$ and
$\widetilde{\theta}_{n}^{(H)}=\arg\min H_{n}(\theta)$, in the sense that
\begin{eqnarray*}
|\widehat{\theta}_{n}^{(W)}-\widetilde{\theta}_{n}^{(W)}|_1 = 
O_{p}\left(\frac{p^{3}}{n^{3/2}}+\frac{1}{np^{K-1}}\right) 
\quad \textrm{and} \quad
|\widehat{\theta}_{n}^{(H)}-\widetilde{\theta}_{n}^{(H)}|_1 = 
O_{p}\left(\frac{p^{3}}{n^{3/2}}+\frac{1}{np^{K-1}}\right)
\end{eqnarray*} 
where
$|\underline{a}|_{1}=\sum_{j=1}^{d}|a_{j}|$, where $\underline{a} =
(a_{1},\ldots,a_{d})$. We will make the assumption that that the ratio
of tapers satisfy the condition $H_{2,n}/H_{1,n}^2 \sim n^{-1}$. This
has some benefits. The first is that the rates for the hybrid Whittle
and the boundary corrected Whittle are the same. In particular, 
by using Corollary 2.1 and Theorem 3.1
in \cite{p:dsy-20} (under Assumption \ref{assum:TS}) we have 
\begin{eqnarray}
\label{eq:Precise1}
[\widehat{J}_{n}(\omega_{};\widehat{f}_{p})-\widehat{J}_{n}(\omega_{};f)]
\overline{J_{n,\underline{h}_n}(\omega_{})}
 =  O_{p}\left(\frac{p^{2}}{n}+\frac{p^{3}}{n^{3/2}} \right)
\end{eqnarray}
and 
\begin{eqnarray}
\label{eq:Precise2}
\widehat{H}_{p,n}(\theta) = H_{n}(\theta)+ O_{p}\left(\frac{p^{3}}{n^{3/2}} + \frac{1}{np^{K-1}}\right).
\end{eqnarray}
Using this, we show below that the Hybrid
Whittle  estimator has the classical $n^{1/2}$--rate. If we
were to relax the rate on $H_{2,n}/H_{1,n}^2 \sim n^{-1}$,
then the $n^{1/2}$--rate and the rates in (\ref{eq:Precise1}) and (\ref{eq:Precise2})
would change. This will make the proofs more
technical. Thus for ease of notation and presentation we will assume
that $H_{2,n}/H_{1,n}^2 \sim n^{-1}$. 

We start by obtaining  a ``crude'' bound for $\nabla_{\theta}^{s}\widehat{W}_{p,n}(\theta)-\nabla_{\theta}^{s} W_{n}(\theta)$. 
\begin{lemma}\label{lemma:WWW}
Suppose that Assumptions \ref{assum:B}(i,iii) and \ref{assum:TS}(i,ii)
hold. Then for $0\leq s \leq \kappa$ (for some $\kappa \geq 4$) we have
\begin{eqnarray*}
\sup_{\theta\in \Theta}\left\|\nabla_{\theta}^{s}\widehat{W}_{p,n}(\theta)-\nabla_{\theta}^{s} W_{n}(\theta)\right\|_{1}
  = O_{p}\left(\frac{p^{2}}{n}\right)
\end{eqnarray*}
and
\begin{eqnarray*}
\sup_{\theta\in \Theta}\left\|\nabla_{\theta}^{s}\widehat{H}_{p,n}(\theta)-\nabla_{\theta}^{s} H_{n}(\theta)\right\|_{1}
  = O_{p}\left(\frac{p^{2}}{n}\right).
\end{eqnarray*}
\end{lemma}
PROOF. We first prove the result in the case that $s=0$ and for
$\widehat{W}_{p,n}(\cdot)$. In this case
\begin{eqnarray*}
\widehat{W}_{p,n}(\theta)-W_{n}(\theta) =
  \frac{1}{n}\sum_{k=1}^{n}f_{\theta}(\omega_{k,n})^{-1} \left[\widehat{J}_{n}(\omega_{k,n};\widehat{f}_{p})-
\widehat{J}_{n}(\omega_{k,n};f)\right] \overline{J_{n}(\omega_{k,n})}.
\end{eqnarray*}
Thus 
\begin{eqnarray*}
\sup_{\theta\in \Theta}|\widehat{W}_{p,n}(\theta)-W_{n}(\theta)| \leq
 \sup_{\theta,\omega}f_{\theta}(\omega_{})^{-1}\times
 \frac{1}{n}\sum_{k=1}^{n}\left|\left[\widehat{J}_{n}(\omega_{k,n};\widehat{f}_{p})-
\widehat{J}_{n}(\omega_{k,n};f)\right]\overline{J_{n}(\omega_{k,n})}\right|.
\end{eqnarray*}
By using Corollary 2.1 in \cite{p:dsy-20} (under Assumption \ref{assum:TS}) we have 
\begin{eqnarray*}
\left[\widehat{J}_{n}(\omega_{};\widehat{f}_{p})-\widehat{J}_{n}(\omega_{};f)\right]\overline{J_{n}(\omega_{})}
 =  \Delta(\omega) + O_{p}\left(\frac{p^{3}}{n^{3/2}} \right)
\end{eqnarray*}
where  $\Delta(\omega)$ is defined in Corollary 2.1 in
\cite{p:dsy-20}, $O_{p}\left(p^{3}n^{-3/2}\right)$ bound is uniform all
frequencies, $\sup_{\omega}\Ex[\Delta(\omega)] = O((np^{K-1})^{-1}+p^{3}/n^{2})$ and 
$\sup_{\omega}\var[\Delta(\omega)] = O(p^{4}/n^{2})$. Thus using this
we have 
\begin{eqnarray*}
\sup_{\theta\in \Theta}|\widehat{W}_{p,n}(\theta)-W_{n}(\theta)| 
&=&  \sup_{\theta,\omega}f_{\theta}(\omega_{})^{-1}\times
 \frac{1}{n}\sum_{k=1}^{n}|\Delta(\omega_{k,n})| +
    O_{p}\left(\frac{p^{3}}{n^{3/2}} \right) \\
&=& O_{p}\left(\frac{p^{2}}{n}+\frac{p^{3}}{n^{3/2}}\right) = O_{p}\left(\frac{p^{2}}{n}\right).
\end{eqnarray*}
This proves the result for $s=0$. 
A similar argument applies for the derivatives of
$\widehat{W}_{p,n}(\theta)$ (together with Assumption \ref{assum:B}(iii)) and  $\widehat{H}_{p,n}(\theta)$,  we
omit the details. 
\hfill $\Box$

\vspace{1em}

\begin{lemma}\label{lemma:WWWcon}
Suppose that Assumptions \ref{assum:B}(i,iii) and \ref{assum:TS}(i,ii) hold.
Then 
\begin{eqnarray*}
| \widehat{\theta}_n^{(W)} - \theta_{n}|_1 \Pcon 0
\qquad \text{and} \qquad
|\widehat{\theta}_n^{(H)} - \theta_{n}|_1 \Pcon 0
\end{eqnarray*}
with $p^{2}/n\rightarrow 0$ as $p,n\rightarrow \infty$. 
\end{lemma}
PROOF. We start with the infeasible criterion $W_{n}(\theta)$. 
Let $\Ex[W_{n}(\theta)] = \mathcal{W}_{n}(\theta)$.
We first show the uniformly convergence of $W_n(\theta)$, i.e.,
\begin{eqnarray}
\label{eq:Wtheta}
\sup_{\theta\in \Theta}|W_{n}(\theta)-\mathcal{W}_{n}(\theta)|\Pcon 0.
\end{eqnarray} Using \cite{p:dsy-20}, Theorem A.1 and the classical
result $\var[K_n(\theta)] = O(n^{-1})$ we have 
\begin{eqnarray*}
\var[W_n(\theta)] &=& \var\left( K_n(\theta) + n^{-1}\sum_{k=1}^{n} \frac{\widehat{J}_n(\omega_{k,n};f) \overline{J_n(\omega_{k,n})}}{f_\theta(\omega_{k,n})}\right)  \\
&\leq& 2 \var[K_n(\theta)] + \frac{2}{n}\sum_{k=1}^{n}  \var[\widehat{J}_n(\omega_{k,n};f) \overline{J_n(\omega_{k,n})}] / f_{\theta}(\omega_{k,n})^2 \\
&\leq& 2 \var[K_n(\theta)] + O(n^{-2}) = O(n^{-1}).
\end{eqnarray*} Therefore, by Markov's inequality, $W_n(\theta) \Pcon \mathcal{W}_n(\theta)$ for each $\theta \in \Theta$. To show 
a uniform convergence, since $\Theta$ is compact, it is enough to show
that $\{W_n(\theta); \theta \in \Theta\}$ is equicontinuous in 
probability. For arbitrary $\theta_1, \theta_2 \in \Theta$,
\begin{eqnarray*}
W_n(\theta_1) - W_n(\theta_2) &=& 
n^{-1} \sum_{k=1}^{n} (f_{\theta_1}^{-1}(\omega_{k,n}) - f_{\theta_2}^{-1}(\omega_{k,n})) \widetilde{J}_n(\omega_{k,n};f) \overline{J_n(\omega_{k,n})}  \\
&+& n^{-1} \sum_{k1}^{n} (\log f_{\theta_1} (\omega_{k,n}) - \log f_{\theta_2} (\omega_{k,n}) )= I_1(\theta_1,\theta_2) + I_2(\theta_1,\theta_2).
\end{eqnarray*} 
To (uniformly) bound $I_1(\theta_1,\theta_2)$, we use the mean value theorem
\begin{eqnarray*}
I_1(\theta_1,\theta_2) 
&=& n^{-1} \sum_{k=1}^{n} \left( f_{\theta_1}^{-1}(\omega_{k,n}) - f_{\theta_2}^{-1}(\omega_{k,n})\right) \widetilde{J}_n(\omega_{k,n};f) \overline{J_n(\omega_{k,n})} \\
&=& n^{-1} \sum_{k=1}^{n} \nabla_{\theta}f_{\theta}^{-1}(\omega_{k,n})\rfloor_{\theta = \overline{\theta}_k}^{\prime}  (\theta_1 -\theta_2) \widetilde{J}_n(\omega_{k,n};f) \overline{J_n(\omega_{k,n})} \\
&=& \mathcal{K}_n(\overline{\theta})^{\prime} (\theta_1 -\theta_2), 
\end{eqnarray*} 
where $\mathcal{K}_n(\overline{\theta}) = n^{-1} \sum_{k=1}^{n} \widetilde{J}_n(\omega_{k,n};f) \overline{J_n(\omega_{k,n})}
\nabla_{\theta}f_{\theta}^{-1}(\omega_{k,n}) \rfloor_{\theta = \overline{\theta}_k}$ and
$\overline{\theta}_1, ..., \overline{\theta}_n$ are convex
combinations of $\theta_1$ and $\theta_2$. It is clear that 
\begin{eqnarray*}
\|\mathcal{K}_n(\overline{\theta})\|_{1} \leq 
  \sup_{\theta,\omega}\|\nabla_{\theta}f_{\theta}^{-1}(\omega)\|_{1}\frac{1}{n}\sum_{k=1}^{n} \left|\widetilde{J}_n(\omega_{k,n};f) \overline{J_n(\omega_{k,n})}
\right| = K_{n}.
\end{eqnarray*}
Thus 
\begin{eqnarray} \label{eq:I1bound}
|I_1(\theta_1,\theta_2)|\leq K_{n}|\theta_{1}-\theta_{2}|_{1}.
\end{eqnarray}
We need to show that $K_{n} = O_{p}(1)$ (it is enough to show that
$\sup_n \Ex[K_{n}]<\infty$). To show this, we use the classical results on DFT
\begin{eqnarray*}
\Ex \left|\widetilde{J}_n(\omega_{k,n};f) \overline{J_n(\omega_{k,n})}
\right| &\leq& \Ex |J_n(\omega_{k,n})|^2 + \Ex \left|\widehat{J}_n(\omega_{k,n};f) \overline{J_n(\omega_{k,n})}
\right| \\
&\leq& f(\omega_{k,n}) + \var (\widehat{J}_n(\omega_{k,n};f))^{1/2} \var(J_n(\omega_{k,n}))^{1/2} \\
&\leq& f(\omega_{k,n}) (1+O(n^{-1})).
\end{eqnarray*} Using above and Assumption \ref{assum:B}(iii-a) gives
\begin{eqnarray*}
\sup_n \Ex [K_n] \leq \sup_{\theta,\omega}\|\nabla_{\theta}f_{\theta}^{-1}(\omega)\|_{1}
\cdot \sup_n
\frac{1}{n} \sum_{k=1}^{n}  f(\omega_{k,n}) (1+O(n^{-1})) < \infty.
\end{eqnarray*}
Therefore, $K_n = O_p(1)$ and from (\ref{eq:I1bound}), $I_{1}(\theta_1, \theta_2)$ is 
equicontinuous in probability. Using similar argument, we can show that $I_2(\theta_1,\theta_2)$ is equicontinous in probability and thus,
$\{W_n(\theta); \theta \in \Theta\}$ is equicontinous in
probability. This imples $\sup_{\theta\in
  \Theta}|W_{n}(\theta)-\mathcal{W}_{n}(\theta)|\Pcon 0$, thus we have shown (\ref{eq:Wtheta}).

\vspace{1em}

\noindent Next, let $\widetilde{\theta}_{n}^{(W)} = \arg\min_{\theta\in
  \Theta}W_{n}(\theta)$. Since $\theta_{n}=\arg\min_{\theta\in \Theta}
\mathcal{W}_{n}(\theta)$ we have 
\begin{eqnarray*}
W_{n}(\widetilde{\theta}_{n}^{(W)}) - \mathcal{W}_{n}(\widetilde{\theta}_n^{(W)})   \leq
W_{n}(\widetilde{\theta}_{n}^{(W)}) - \mathcal{W}_{n}(\theta_n) \leq
W_{n}(\theta_{n}) - \mathcal{W}_{n}(\theta_{n}). 
\end{eqnarray*}
Thus 
\begin{eqnarray*}
|W_{n}(\widetilde{\theta}_{n}^{(W)}) - \mathcal{W}_{n}(\theta_n)| \leq
\sup_{\theta}|W_{n}(\theta)-\mathcal{W}_{n}(\theta)|\Pcon 0.
\end{eqnarray*}
If $\theta_{n}$ uniquely minimises $I_{n}(f,f_{\theta})$, then by
using the above we have that
$|\widetilde{\theta}_{n}^{(W)}-\theta_{n}|_1 \Pcon 0$. However,
$W_{n}(\theta)$ is an infeasible criterion. To show consistency we
need to obtain a uniform bound on the feasible criterion
$\widehat{W}_{p,n}(\theta)$. That is 
\begin{eqnarray}
\label{eq:Wtheta1}
\sup_{\theta}|\widehat{W}_{p,n}(\theta)- \mathcal{W}_{n}(\theta)|\Pcon 0.
\end{eqnarray}
Now by using the triangular inequality, together with 
 (\ref{eq:Wtheta}) and Lemma \ref{lemma:WWW},
 (\ref{eq:Wtheta1}) immediately follows. Therefore, by
using the same arguments those given above we have 
$|\widehat{\theta}_n^{(W)} - \theta_{n}|_1 \Pcon 0$, which is the
desired result.  

By the same set of arguments we have $|\widehat{\theta}_n^{(H)} - \theta_{n}|_1 \Pcon 0$.
\hfill $\Box$

\vspace{1em}

For the simplicity, we assume $\theta$ is univariate and state the following lemma. It can
be easily generalized to the multivariate case.

\begin{lemma}\label{lemma:secondorderexpansion}
Suppose Assumptions \ref{assum:B}(i,iii) and \ref{assum:TS} hold. Then for $i=1,2$ we have 
\begin{eqnarray}
\label{eq:W1}
\frac{d^{i}\widehat{W}_{p,n}(\theta)}{d\theta^{i}}\rfloor_{\theta =
  \theta_{n}} = \frac{d^{i}W_{n}(\theta)}{d\theta^{i}}\rfloor_{\theta =
  \theta_{n}} + O_{p}\left(\frac{p^{3}}{n^{3/2}}+ \frac{1}{np^{K-1}}\right)
\end{eqnarray}
and 
\begin{eqnarray}
\label{eq:W2}
\left|\frac{d^{3}\widehat{W}_{p,n}(\theta)}{d\theta^{3}}\rfloor_{\theta =
  \overline{\theta}_{n}} - \frac{d^{3}W_{n}(\theta)}{d\theta^{3}}\rfloor_{\theta =
 \theta_{n}}\right| =
   O_{p}\left(\frac{p^{2}}{n}\right) + |\widehat{\theta}_n^{(W)} - \theta_n|O_{p}(1),
\end{eqnarray}
where $\overline{\theta}_n$ is a convex combination of $\widehat{\theta}_n^{(W)}$ and $\theta_n$.
This gives rise to the first order and second expansions
\begin{eqnarray}
\label{eq:W3}
(\widehat{\theta}_{n}^{(W)}  - \theta_{n}) = -\left[\Ex\left[\frac{d^{2}
                                        W_{n}(\theta_{n})}{d\theta_{n}^{2}} \right]\right]^{-1}
\frac{d W_{n}(\theta_{n})}{d\theta_{n}} +
O_{p}\left( \frac{1}{n} + \frac{p^{3}}{n^{3/2}}+\frac{1}{np^{K-1}}\right)
\end{eqnarray}
and 
\begin{eqnarray}
&& \frac{dW_{n}(\theta)}{d\theta}\rfloor_{\theta =
  \theta_{n}} + (\widehat{\theta}_{n}^{(W)}-\theta_{n}) 
\frac{d^{2}W_{n}(\theta)}{d\theta^{2}}\rfloor_{\theta =
  \theta_{n}} + \frac{1}{2}(\widehat{\theta}_{n}^{(W)}-\theta_{n})^{2}
\frac{d^{3}W_{n}(\theta)}{d\theta^{3}}\rfloor_{\theta =\theta_{n}}\nonumber \\
&& \qquad  \qquad =
  O_{p}\left(\frac{p^{3}}{n^{3/2}}+\frac{1}{np^{K-1}}\right). 
\label{eq:W4}
\end{eqnarray}
\end{lemma}
PROOF.  By using Theorem 3.1, \cite{p:dsy-20} we have for $i=1$ and $2$
\begin{eqnarray*}
\frac{d^{i}\widehat{W}_{p,n}(\theta)}{d\theta^{i}}\rfloor_{\theta =
  \theta_{n}} = \frac{d^{i}W_{n}(\theta)}{d\theta^{i}}\rfloor_{\theta =
  \theta_{n}} + O_{p}\left(\frac{p^{3}}{n^{3/2}} + \frac{1}{np^{K-1}}\right),
\end{eqnarray*}
this immediately gives (\ref{eq:W1}). 
Let $\overline{\theta}_n$ denote a convex combination of $\theta_n$ and
$\widehat{\theta}_n^{(W)}$ (note that $\widehat{\theta}_{n}^{(W)}$ is a consistent
estimator of $\theta_{n}$). To evaluate $\frac{d^{3}\widehat{W}_{n}(\theta)}{d\theta^{3}}$ 
at the (consistent) estimator $\overline{\theta}_{n}$, a slightly
different approach is required (due to the additional 
random parameter $\overline{\theta}_{n}$). By using triangular inequality and   
Lemma \ref{lemma:WWW} we have
\begin{eqnarray*}
&&\left|\frac{d^{3}\widehat{W}_{p,n}(\theta)}{d\theta^{3}}\rfloor_{\theta =
  \overline{\theta}_{n}} - \frac{d^{3}W_{n}(\theta)}{d\theta^{3}}\rfloor_{\theta =
 \theta_{n}}\right| \\
&& \leq \left|\frac{d^{3}\widehat{W}_{p,n}(\theta)}{d\theta^{3}}\rfloor_{\theta =
  \overline{\theta}_{n}} - \frac{d^{3}W_{n}(\theta)}{d\theta^{3}}\rfloor_{\theta =
 \overline{\theta}_{n}}\right| \\
&& \quad + \frac{1}{n}\sum_{k=1}^{n}
 \left|\frac{d^{3}}{d\theta^{3}}\left[f_{\overline{\theta}_{n}}(\omega_{k,n})^{-1}-f_{\theta_{n}}(\omega_{k,n})^{-1}\right]\right|
    \left|\widetilde{J}_{n}(\omega_{k,n};f) \overline{J_{n}(\omega_{k,n})} \right| \\
&& =   O_{p}\left(\frac{p^{2}}{n} \right)
+ \frac{1}{n}\sum_{k=1}^{n}
 \left|\frac{d^{3}}{d\theta^{3}}\left[f_{\overline{\theta}_{n}}(\omega_{k,n})^{-1}-f_{\theta_{n}}(\omega_{k,n})^{-1}\right]\right|
    \left|\widetilde{J}_{n}(\omega_{k,n};f) \overline{J_{n}(\omega_{k,n})} \right|.
\end{eqnarray*}
For the second term in the above, we apply the mean value theorem to
$\frac{d^{3}}{d\theta^{3}} f_{\theta}^{-1}$ to give 
\begin{eqnarray*}
\left|\frac{d^{3}}{d\theta^{3}}(f_{\bar{\theta}_{n}}^{-1} -
  f_{\theta_{n}}^{-1}) \right|\leq
  \sup_{\theta}\left|\frac{d^{4}}{d\theta^{4}}f_{\theta}^{-1}\right|\cdot
  |\bar{\theta}_{n}-\theta_{n}| \leq 
  \sup_{\theta}\left|\frac{d^{4}}{d\theta^{4}}f_{\theta}^{-1}\right|\cdot
  |\widehat{\theta}_{n}^{(W)}-\theta_{n}|,
\end{eqnarray*}
note that to bound the 
fourth derivation we require Assumption \ref{assum:B}(iii) for $\kappa=4$.
Substituting this into the previous inequality gives 
\begin{eqnarray*}
&&\left|\frac{d^{3}\widehat{W}_{p,n}(\theta)}{d\theta^{3}}\rfloor_{\theta =
  \overline{\theta}_{n}} - \frac{d^{3}W_{n}(\theta)}{d\theta^{3}}\rfloor_{\theta =
 \theta_{n}}\right| \\
&& \leq   O_{p}\left(\frac{p^{2}}{n} \right)
+ \frac{1}{n}\sum_{k=1}^{n}
 \left|\frac{d^{3}}{d\theta^{3}}\left[f_{\overline{\theta}_{n}}(\omega_{k,n})^{-1}-f_{\theta_{n}}(\omega_{k,n})^{-1}\right]\right|
    \left|\widetilde{J}_{n}(\omega_{k,n};f) \overline{J_{n}(\omega_{k,n})} \right| \\
&& =   O_{p}\left(\frac{p^{2}}{n}\right) + |\widehat{\theta}_n^{(W)} - \theta_n|O_{p}(1).
\end{eqnarray*}
The above proves  (\ref{eq:W2}). 

Using (\ref{eq:W1}) and (\ref{eq:W2}) we now obtain the first and second order expansions
in (\ref{eq:W3}) and (\ref{eq:W4}).
In order to prove (\ref{eq:W3}), we will show that 
\begin{eqnarray*}
(\widehat{\theta}_{n}^{(W)}  - \theta_{n})&=& O_{p}\left(\frac{1}{n^{1/2}}+\frac{p^{3}}{n^{3/2}}\right).
\end{eqnarray*}
if $p^{2}/n\rightarrow 0$ we make 
a second order expansion of $\frac{d
  \widehat{W}_{p,n}(\widehat{\theta}_n^{(W)})}{d\theta}$
about $\theta_{n}$ and assuming that $\widehat{\theta}_{n}^{(W)}$ lies
inside the parameter space we have 
\begin{eqnarray*}
0=\frac{d \widehat{W}_{p,n}(\widehat{\theta}_{n}^{(W)})}{d\widehat{\theta}_{n}^{(W)}} = 
\frac{d \widehat{W}_{p,n}(\theta_{n})}{d\theta_{n}} +
  (\widehat{\theta}_{n}^{(W)}  - \theta_{n})
\frac{d^{2} \widehat{W}_{p,n}(\theta_{n})}{d\theta_{n}^{2}} + 
\frac{1}{2}(\widehat{\theta}_{n}^{(W)} - \theta_{n})^{2}
\frac{d^{3} \widehat{W}_{p,n}(\overline{\theta}_{n})}{d\theta_{n}^{3}}
\end{eqnarray*}
where $\overline{\theta}_{n}$ is a convex combination of $\theta_{n}$
and $\widehat{\theta}_{n}^{(W)}$. Now by using (\ref{eq:W1}) and
(\ref{eq:W2}) we can replace in the above
$\widehat{W}_{p,n}(\theta_{n})$ and its derivatives with
$W_{n}(\theta_n)$ and its derivatives. 
Therefore,  
\begin{eqnarray}
&&\frac{d W_{n}(\theta_{n})}{d\theta_{n}} +  (\widehat{\theta}_{n}^{(W)}  - \theta_{n})
\frac{d^{2} W_{n}(\theta_{n})}{d\theta_{n}^{2}} + \frac{1}{2}(\widehat{\theta}_{n}^{(W)} - \theta_{n})^{2}
\frac{d^{3} W_{n}(\theta_{n})}{d\theta_{n}^{3}} \nonumber\\
&=&O_{p}\left(\frac{p^{3}}{n^{3/2}}+\frac{1}{np^{K-1}}\right) + (\widehat{\theta}_{n}^{(W)} - \theta_{n})^{2}O_{p}\left(\frac{p^{2}}{n}\right)
+|\widehat{\theta}_{n}^{(W)} - \theta_{n}|^{3}O_{p}(1). \label{eq:secondExpan}
\end{eqnarray}
Rearranging the above gives 
\begin{eqnarray}
\label{eq:thetahatT}
(\widehat{\theta}_{n}^{(W)}  - \theta_{n})
&=& -\left[\frac{d^{2} W_{n}(\theta_{n})}{d\theta_{n}^{2}} \right]^{-1}
\frac{d W_{n}(\theta_{n})}{d\theta_{n}} -
\frac{1}{2}\left[\frac{d^{2} W_{n}(\theta_{n})}{d\theta_{n}^{2}} \right]^{-1}
\frac{d^{3} W_{n}(\theta_{n})}{d\theta_{n}^{3}} (\widehat{\theta}_{n}^{(W)} - \theta_{n})^{2}
\nonumber\\
&&+O_{p}\left(\frac{p^{3}}{n^{3/2}}+\frac{1}{np^{K-1}}\right) + (\widehat{\theta}_{n}^{(W)} - \theta_{n})^{2}O_{p}\left(\frac{p^{2}}{n}\right)
+|\widehat{\theta}_{n}^{(W)} - \theta_{n}|^{3}O_{p}(1). 
\end{eqnarray}
Next we obtain a bound for $\frac{dW_{n}(\theta_{n})}{d\theta_{n}}$
(to substitute into the above).
Since
$\Ex[\frac{dW_{n}(\theta_{n})}{d\theta_{n}}]=O(n^{-K})$ (from equation (\ref{eq:alphaBound})) and
$\var[\frac{d W_{n}(\theta_{n})}{d\theta_{n}}]=O_{p}(n^{-1})$ 
we have $\frac{d W_{n}(\theta_{n})}{d\theta_{n}} =
O_{p}(n^{-1/2})$. Substituting this into (\ref{eq:thetahatT}) gives 
\begin{eqnarray*}
(\widehat{\theta}_{n}^{(W)}  - \theta_{n})&=& 
\frac{1}{2}\left[\frac{d^{2} W_{n}(\theta_{n})}{d\theta_{n}^{2}} \right]^{-1}
\frac{d^{3} W_{n}(\theta_{n})}{d\theta_{n}^{3}} (\widehat{\theta}_{n}^{(W)} - \theta_{n})^{2}
\nonumber\\
&&+O_{p}(n^{-1/2})+ O_{p}\left(\frac{p^{3}}{n^{3/2}}+\frac{1}{np^{K-1}}\right) + (\widehat{\theta}_{n}^{(W)} - \theta_{n})^{2}O_{p}\left(\frac{p^{2}}{n}\right)
+|\widehat{\theta}_{n}^{(W)} - \theta_{n}|^{3}O_{p}(1).
\end{eqnarray*} 
Using that $\left[\frac{d^{2} W_{n}(\theta_{n})}{d\theta_{n}^{2}} \right]^{-1}
\frac{d^{3} W_{n}(\theta_{n})}{d\theta_{n}^{3}}=O_{p}(1)$ and
substituting this into the above gives 
\begin{eqnarray}
\label{eq:widehattheta}
(\widehat{\theta}_{n}^{(W)}  - \theta_{n}) &=& 
  O_{p}\left(\frac{1}{n^{1/2}}+\frac{p^{3}}{n^{3/2}}+\frac{1}{np^{K-1}}\right) + 
(\widehat{\theta}_{n}^{(W)} - \theta_{n})^{2}O_{p}\left(\frac{p^{2}}{n}+1\right)
+|\widehat{\theta}_{n}^{(W)} - \theta_{n}|^{3}O_{p}(1). \nonumber\\
&& 
\end{eqnarray} 
Thus, from the above and the consistency result in Lemma
\ref{lemma:WWWcon} 
($|\widehat{\theta}_{n}^{(W)}  - \theta_{n}|=o_{p}(1)$) we have\footnote{The precise proof for (\ref{eq:thetarate}): By using
  (\ref{eq:widehattheta}) we have 
\begin{eqnarray*}
(\widehat{\theta}_{n}^{(W)}  - \theta_{n}) = 
A_n + (\widehat{\theta}_{n}^{(W)}  - \theta_{n})^2 B_n + 	(\widehat{\theta}_{n}^{(W)}  - \theta_{n})^3 C_n
\end{eqnarray*} where 
the random variables $A_{n}, B_{n}$ and $C_{n}$ are such that
$A_n = O_{p}\left(\frac{1}{n^{1/2}}+\frac{p^{3}}{n^{3/2}}+\frac{1}{np^{K-1}}\right)$,
$B_n = O_{p}\left(\frac{p^{2}}{n}+1\right)$, and $C_n = O_p(1)$.
Then, moving the second and third term in the RHS to LHS
\begin{eqnarray*}
(\widehat{\theta}_{n}^{(W)}  - \theta_{n}) [1 + (\widehat{\theta}_{n}^{(W)}  - \theta_{n}) B_n + (\widehat{\theta}_{n}^{(W)}  - \theta_{n})^2 C_n]
=O_{p}\left(\frac{1}{n^{1/2}}+\frac{p^{3}}{n^{3/2}}+\frac{1}{np^{K-1}}\right).
\end{eqnarray*}
Using consistency, $(\widehat{\theta}_{n}^{(W)}  - \theta_{n}) B_n = o_p(1) O_{p}\left(\frac{p^{2}}{n}+1\right) = o_p(\frac{p^{2}}{n}+1)$
and $(\widehat{\theta}_{n}^{(W)}  - \theta_{n})^2 C_n = o_p(1)$.
Finally, we use that $(1+o_p(1))^{-1} = O_p(1)$, then,
\begin{eqnarray*}
(\widehat{\theta}_{n}^{(W)}  - \theta_{n}) 
=O_p(1) O_{p}\left(\frac{1}{n^{1/2}}+\frac{p^{3}}{n^{3/2}}+\frac{1}{np^{K-1}}\right).
\end{eqnarray*}
Thus giving the required probabilistic rate.} 
\begin{eqnarray}
\label{eq:thetarate}
(\widehat{\theta}_{n}^{(W)}  - \theta_{n})
= O_{p}\left(\frac{1}{n^{1/2}}+\frac{p^{3}}{n^{3/2}}\right).
\end{eqnarray}
We use the above bound to obtain an exact expression for the
dominating rate $O_{p}(n^{-1/2})$. 
Returning to equation (\ref{eq:thetahatT}) and substituting this bound
into the quadratic term in
(\ref{eq:thetahatT}) gives 
\begin{eqnarray*}
(\widehat{\theta}_{n}^{(W)}  - \theta_{n})
= -\left[\frac{d^{2} W_{n}(\theta_{n})}{d\theta_{n}^{2}} \right]^{-1}
\frac{d W_{n}(\theta_{n})}{d\theta_{n}} +O_{p}\left(\frac{1}{n} + \frac{p^{3}}{n^{3/2}}+ \frac{1}{np^{K-1}}\right). 
\end{eqnarray*}
Using that $\frac{d^{2} W_{n}(\theta_{n})}{d\theta_{n}^{2}}
= \Ex\left(\frac{d^{2}    W_{n}(\theta_{n})}{d\theta_{n}^{2}}\right)+O_{p}(n^{-1/2})$ and
under Assumption \ref{assum:TS}(iii) we have 
\begin{eqnarray*}
(\widehat{\theta}_{n}^{(W)}  - \theta_{n})
= -\left[\Ex\left[\frac{d^{2}W_{n}(\theta_{n})}{d\theta_{n}^{2}} \right]\right]^{-1}
\frac{d W_{n}(\theta_{n})}{d\theta_{n}} +
O_{p}\left( \frac{1}{n} + \frac{p^{3}}{n^{3/2}}+ \frac{1}{np^{K-1}}\right).
\end{eqnarray*}
This proves (\ref{eq:W3}). 

To prove (\ref{eq:W4}) we return to (\ref{eq:secondExpan}). By
substituting (\ref{eq:thetarate}) into (\ref{eq:secondExpan}) we have 
\begin{eqnarray*}
\frac{dW_{n}(\theta)}{d\theta}\rfloor_{\theta =
  \theta_{n}} + (\widehat{\theta}_{n}^{(W)}-\theta_{n}) 
\frac{d^{2}W_{n}(\theta)}{d\theta^{2}}\rfloor_{\theta =
  \theta_{n}} + \frac{1}{2}(\widehat{\theta}_{n}^{(W)}-\theta_{n})^{2}
\frac{d^{3}W_{n}(\theta)}{d\theta^{3}}\rfloor_{\theta =\theta_{n}} =
  O_{p}\left(\frac{p^{3}}{n^{3/2}}+
\frac{1}{np^{K-1}} \right). 
\end{eqnarray*}
This proves (\ref{eq:W4}).  \hfill $\Box$

\vspace{1em}
The second order expansion (\ref{eq:W4}) is instrumental in proving
the equivalence result Theorem \ref{theorem:equivalence}. By following
a similar set of arguments to those in Lemma
\ref{lemma:secondorderexpansion} for the multivariate parameter
$\theta = (\theta_{1},\ldots,\theta_{d})$, the feasible estimator
satisfies the expansion
\begin{eqnarray}
&&\frac{\partial W_{n}(\theta)}{\partial \theta_{r}} + 
\sum_{s=1}^{d}(\widehat{\theta}_{s,n}^{(W)} - \theta_{s,n})
  \frac{\partial^{2}W_{n}(\theta)}{\partial \theta_{s}\partial
    \theta_{r}}\rfloor_{\theta  =\theta_{n}} +\nonumber\\
&& \frac{1}{2}\sum_{s_1,s_2=1}^{d}(\widehat{\theta}_{s_1,n}^{(W)} - \theta_{s_1,n}) (\widehat{\theta}_{s_2,n}^{(W)} - \theta_{s_2,n})
  \frac{\partial^{3}W_{n}(\theta)}{\partial    \theta_{s_1}\partial\theta_{s_2}\partial \theta_r}\rfloor_{\theta
    =\theta_{n}} = \left(\frac{p^{3}}{n^{3/2}}+
\frac{1}{np^{K-1}} \right). 
\label{eq:W4M}
\end{eqnarray}
By using the same set of arguments we can obtain a first and second
order expansion for the hybrid Whittle estimator
\begin{eqnarray}
\label{eq:H3}
(\widehat{\theta}_{n}^{(H)}  - \theta_{n})= -\left[\Ex\left[\frac{d^{2}
                                        H_{n}(\theta_{n})}{d\theta_{n}^{2}} \right]\right]^{-1}
\frac{d H_{n}(\theta_{n})}{d\theta_{n}} +
O_{p}\left( \frac{1}{n} + \frac{p^{3}}{n^{3/2}}+ \frac{1}{np^{K-1}}\right)
\end{eqnarray}
and 
\begin{eqnarray}
&& \frac{dH_{n}(\theta)}{d\theta}\rfloor_{\theta =
  \theta_{n}} + (\widehat{\theta}_{n}^{(W)}-\theta_{n}) 
\frac{d^{2}H_{n}(\theta)}{d\theta^{2}}\rfloor_{\theta =
  \theta_{n}} + \frac{1}{2}(\widehat{\theta}_{n}^{(H)}-\theta_{n})^{2}
\frac{d^{3}H_{n}(\theta)}{d\theta^{3}}\rfloor_{\theta =\theta_{n}} \nonumber \\
&& \qquad \qquad =
  \left(\frac{p^{3}}{n^{3/2}}+
\frac{1}{np^{K-1}} \right). 
\label{eq:H4}
\end{eqnarray}

\vspace{1em}

\noindent \textbf{PROOF of Theorem \ref{theorem:equivalence}}. 
We first prove the result for the one parameter case when $p\geq 1$. 
By using (\ref{eq:W4}) for the feasible estimator
$\widehat{\theta}^{(W)}_{n} = \arg\min \widehat{W}_{p,n}(\theta)$ we
have 
\begin{eqnarray*}
\frac{dW_{n}(\theta)}{d\theta}\rfloor_{\theta =
  \theta_{n}} + (\widehat{\theta}_{n}^{(W)}-\theta_{n}) 
\frac{d^{2}W_{n}(\theta)}{d\theta^{2}}\rfloor_{\theta =
  \theta_{n}} + \frac{1}{2}(\widehat{\theta}_{n}^{(W)}-\theta_{n})^{2}
\frac{d^{3}W_{n}(\theta)}{d\theta^{3}}\rfloor_{\theta =\theta_{n}} =
  O_{p}\left(\frac{p^{3}}{n^{3/2}}+\frac{1}{np^{K-1}}\right). 
\end{eqnarray*}
Whereas for the infeasible estimator
$\widetilde{\theta}^{(W)}_{n} = \arg\min \widehat{W}_{p,n}(\theta)$ we
have 
\begin{eqnarray*}
\frac{dW_{n}(\theta)}{d\theta}\rfloor_{\theta =
  \theta_{n}} + (\widetilde{\theta}_{n}^{(W)}-\theta_{n}) 
\frac{d^{2}W_{n}(\theta)}{d\theta^{2}}\rfloor_{\theta =
  \theta_{n}} + \frac{1}{2}(\widetilde{\theta}_{n}^{(W)}-\theta_{n})^{2}
\frac{d^{3}W_{n}(\theta)}{d\theta^{3}}\rfloor_{\theta =\theta_{n}} =
  O_{p}\left(\frac{1}{n^{3/2}}\right). 
\end{eqnarray*}
Taking differences for the two expansions above we have
\begin{eqnarray*}
&&(\widetilde{\theta}_{n}^{(W)}-\widehat{\theta}_{n}^{(W)}) 
\frac{d^{2}W_{n}(\theta)}{d\theta^{2}}\rfloor_{\theta =
  \theta_{n}} + \frac{1}{2}(\widetilde{\theta}_{n}^{(W)}-\widehat{\theta}_{n}^{(W)})\left[
(\widetilde{\theta}_{n}^{(W)} -
   \theta_{n})+(\widehat{\theta}_{n}^{(W)} - \theta_{n})\right] 
\frac{d^{3}W_{n}(\theta)}{d\theta^{3}}\rfloor_{\theta =\theta_{n}} \\
&&\qquad =
  O_{p}\left(\frac{p^3}{n^{3/2}}+\frac{1}{np^{K-1}}\right). 
\end{eqnarray*}
Now replacing $\frac{d^{2}W_{n}(\theta)}{d\theta^{2}}\rfloor_{\theta =
  \theta_{n}} $  with its expectation and using that $|\widetilde{\theta}_{n}^{(W)} -
   \theta_{n}| = o_{p}(1)$ and $|\widehat{\theta}_{n}^{(W)} -
   \theta_{n}| = o_{p}(1)$ we have 
\begin{eqnarray*}
(\widetilde{\theta}_{n}^{(W)}-\widehat{\theta}_{n}^{(W)}) 
\Ex\left(\frac{d^{2}W_{n}(\theta)}{d\theta^{2}}\rfloor_{\theta =
  \theta_{n}}\right) + o_p(1)
(\widetilde{\theta}_{n}^{(W)}-\widehat{\theta}_{n}^{(W)})
=  O_{p}\left(\frac{p^3}{n^{3/2}}+\frac{1}{np^{K-1}}\right). 
\end{eqnarray*}
Since $\Ex\left(\frac{d^{2}W_{n}(\theta)}{d\theta^{2}}\rfloor_{\theta =
  \theta_{n}}\right)$ is greater than 0, the above implies
\begin{eqnarray*}
(\widetilde{\theta}_{n}^{(W)}-\widehat{\theta}_{n}^{(W)}) 
=  O_{p}\left(\frac{p^3}{n^{3/2}}+\frac{1}{np^{K-1}}\right). 
\end{eqnarray*}
Now we prove the result for the case $p=0$. If $p=0$, then
$\widehat{W}_{p,n}(\theta) = K_{n}(\theta)$ (the Whittle likelihood). Let 
\begin{eqnarray*}
\widehat{\theta}_{n}^{(K)} = \arg\min K_{n}(\theta) 
\quad \textrm{and} \quad
\widetilde{\theta}_{n}^{(W)} = \arg\min W_{n}(\theta).
\end{eqnarray*}
Our aim is to show that $|\widehat{\theta}_{n}^{(K)}  -
\widetilde{\theta}_{n}^{(W)}| = O_{p}(n^{-1})$. 
Note that $W_{n}(\theta) = K_{n}(\theta) + C_{n}(\theta)$, 
where 
\begin{eqnarray*}
C_{n}(\theta) = \frac{1}{n}\sum_{k=1}^{n}\frac{\widehat{J}_{n}(\omega;f)\overline{J_{n}(\omega_{k,n})}}{f_{\theta}(\omega_{k,n})}.
\end{eqnarray*}
Using a Taylor expansion, similar to the above, we have 
\begin{eqnarray*}
\frac{dK_{n}(\theta)}{d\theta}\rfloor_{\theta =
  \theta_{n}} + (\widehat{\theta}_{n}^{(K)}-\theta_{n}) 
\frac{d^{2}K_{n}(\theta)}{d\theta^{2}}\rfloor_{\theta =
  \theta_{n}} + \frac{1}{2}(\widehat{\theta}_{n}^{(K)}-\theta_{n})^{2}
\frac{d^{3}K_{n}(\theta)}{d\theta^{3}}\rfloor_{\theta =\theta_{n}} =
  O_{p}\left(\frac{1}{n^{3/2}}\right) 
\end{eqnarray*}
and
\begin{eqnarray*}
\frac{dW_{n}(\theta)}{d\theta}\rfloor_{\theta =
  \theta_{n}} + (\widetilde{\theta}_{n}^{(W)}-\theta_{n}) 
\frac{d^{2}W_{n}(\theta)}{d\theta^{2}}\rfloor_{\theta =
  \theta_{n}} + \frac{1}{2}(\widetilde{\theta}_{n}^{(W)}-\theta_{n})^{2}
\frac{d^{3}W_{n}(\theta)}{d\theta^{3}}\rfloor_{\theta =\theta_{n}} =
  O_{p}\left(\frac{1}{n^{3/2}}\right). 
\end{eqnarray*}
Taking differences of the two expansions
\begin{eqnarray}
&&\frac{dC_{n}(\theta)}{d\theta}\rfloor_{\theta =
  \theta_{n}} + (\widehat{\theta}_{n}^{(K)}-\widetilde{\theta}_{n}^{(W)}) 
\frac{d^{2}K_{n}(\theta)}{d\theta^{2}}\rfloor_{\theta =
  \theta_{n}} -(\widehat{\theta}_{n}^{(W)}-\theta_{n}) \frac{d^{2}C_{n}(\theta)}{d\theta^{2}}\rfloor_{\theta =
  \theta_{n}}  \nonumber\\
&&
  + \frac{1}{2}(\widehat{\theta}_{n}^{(K)}-\widetilde{\theta}_{n}^{(W)})
\left[  (\widehat{\theta}_{n}^{(K)}-\theta_{n})+(\widetilde{\theta}_{n}^{(W)}-\theta_{n})\right]
\frac{d^{3}K_{n}(\theta)}{d\theta^{3}}\rfloor_{\theta =\theta_{n}} \nonumber\\
&& -\frac{1}{2}(\widetilde{\theta}_{n}^{(W)} - \theta_{n})^2 \frac{d^{3}C_{n}(\theta)}{d\theta^{3}}\rfloor_{\theta =\theta_{n}} =  O_{p}\left(\frac{1}{n^{3/2}}\right). \label{eq:thetannn}
\end{eqnarray}
To bound the above we use that
\begin{eqnarray*}
|\widehat{\theta}_{n}^{(K)}-\theta_{n}| = O_{p}(n^{-1/2}) \textrm{ and
  }|\widetilde{\theta}_{n}^{(W)}-\theta_{n}| = O_{p}(n^{-1/2}). 
\end{eqnarray*}
In addition by using a proof analogous to the proves of Theorem
\ref{theorem:Bound}, equation (\ref{eq:BD1})  
we have
\begin{eqnarray*}
\frac{d^{s}C_{n}(\theta)}{d\theta^{s}} = 
\frac{1}{n}\sum_{k=1}^{n}\widehat{J}_{n}(\omega_{k,n};f)\overline{J_{n}(\omega_{k,n})}\frac{d^{s}}{d\theta^{s}}f_{\theta}(\omega_{k,n})^{-1} =
O_{p}(n^{-1}) \quad \textrm{for} \quad 0\leq s \leq3.
\end{eqnarray*}
Substituting the above bounds into (\ref{eq:thetannn}) gives 
\begin{eqnarray*}
(\widehat{\theta}_{n}^{(K)}-\widetilde{\theta}_{n}^{(W)}) 
\frac{d^{2}K_{n}(\theta)}{d\theta^{2}}\rfloor_{\theta =
  \theta_{n}}  + \frac{1}{2}(\widehat{\theta}_{n}^{(K)}-\widetilde{\theta}_{n}^{(W)})O_{p}(n^{-1/2})
 = -
\frac{dC_{n}(\theta)}{d\theta}\rfloor_{\theta =
  \theta_{n}} + O_{p}\left(\frac{1}{n^{3/2}}\right). 
\end{eqnarray*}
Since $[\frac{d^{2}K_{n}(\theta)}{d\theta^{2}}\rfloor_{\theta =
  \theta_{n}}]^{-1}=O_{p}(1)$ we have
\begin{eqnarray*}
|\widehat{\theta}_{n}^{(K)}-\widetilde{\theta}_{n}^{(W)}|=O_{p}(n^{-1}),
\end{eqnarray*}
thus giving the desired rate. 

For the multiparameter case we use (\ref{eq:W4M}) and the same
argument to give 
\begin{eqnarray*}
&&\sum_{s=1}^{d}(\widehat{\theta}_{s,n}^{(W)} - \widetilde{\theta}_{s,n}^{(W)})
  \frac{\partial^{2}W_{n}(\theta)}{\partial \theta_{s}\partial
    \theta_{r}}\rfloor_{\theta  =\theta_{n}} +\nonumber\\
&& \frac{1}{2}\sum_{s_1,s_2=1}^{d}\left[
(\widehat{\theta}_{s_1,n}^{(W)} - \widetilde{\theta}_{s_1,n}^{(W)})
   (\widehat{\theta}_{s_2,n}^{(W)} - \theta_{s_2,n}) +
(\widehat{\theta}_{s_2,n}^{(W)} - \widetilde{\theta}_{s_2,n}^{(W)})
   (\widetilde{\theta}_{s_1,n}^{(W)} - \theta_{s_1,n}) \right]\times
  \frac{\partial^{3}W_{n}(\theta)}{\partial    \theta_{s_1}\partial\theta_{s_2}\partial \theta_r}\rfloor_{\theta
    =\theta_{n}} \\
&=& O_{p}\left(\frac{p^3}{n^{3/2}}+\frac{1}{np^{K-1}}\right). 
\end{eqnarray*}
Replacing $\frac{\partial^{2}W_{n}(\theta)}{\partial \theta_{s}\partial
    \theta_{r}}\rfloor_{\theta  =\theta_{n}}$ with its expectation gives 
\begin{eqnarray*}
(\widehat{\theta}_{n}^{(W)} - \widetilde{\theta}_{n}^{(W)})^{\prime}
  \Ex\left[\nabla^{2}_{\theta}W_{n}(\theta)\rfloor_{\theta  =\theta_{n}}\right] +
 |\widehat{\theta}_{n}^{(W)} - \widetilde{\theta}_{n}^{(W)}|_{1}o_{p}(1)
&=& O_{p}\left(\frac{p^3}{n^{3/2}}+\frac{1}{np^{K-1}}\right). 
\end{eqnarray*}
Thus under the assumption that
$\Ex\left[\nabla^{2}_{\theta}W_{n}(\theta)\rfloor_{\theta
    =\theta_{n}}\right]$ is invertible we have 
\begin{eqnarray*}
|\widetilde{\theta}_{n}^{(W)}-\widehat{\theta}_{n}^{(W)}|_{1}
=  O_{p}\left(\frac{p^3}{n^{3/2}}+\frac{1}{np^{K-1}}\right). 
\end{eqnarray*}
By a similar argument we have 
\begin{eqnarray*}
|\widetilde{\theta}_{n}^{(H)}-\widehat{\theta}_{n}^{(H)}|_{1} 
=  O_{p}\left(\frac{p^3}{n^{3/2}}+\frac{1}{np^{K-1}}\right). 
\end{eqnarray*}
The case when  $p=0$ is analogous to the uniparameter case and we omit
the details. This concludes the proof. \hfill $\Box$

%% file: appendix_bias.tex
\section{The bias of the different criteria}\label{sec:bias}

In this section we derive the approximate bias of the Gaussian, 
Whittle,  boundary corrected and hybrid Whittle likelihoods under
quite general assumptions on the underlying time series $\{X_{t}\}$.  The bias 
we evaluate will be in the sense of \cite{p:bar-53} and will
be based on the second order expansion of the loss function. We
mention that for certain specific models (such as the specified AR, or certain MA or
ARMA) the bias of the least
squares, Whittle likelihood or maximum likeihood estimators are given
in \cite{p:tan-83,p:tan-84,p:sha-sti-88}.

\subsection{Bias for the estimator of one unknown parameter}\label{sec:oneparameter}

In order to derive the limiting bias, we require the  following
definitions
\begin{equation*}
I(\theta) = -\frac{1}{2\pi}\int_{0}^{2\pi}\left(\frac{d^{2}f_{\theta}(\omega)^{-1}}{d\theta^{2}}\right)f(\omega)d\omega
\quad \textrm{and} \quad
J(g) = \frac{1}{2\pi}\int_{0}^{2\pi}g(\omega)f_{}(\omega)d\omega.
\end{equation*}
For real functions $g,h\in L^{2}[0,2\pi]$ we define
\begin{eqnarray}
V(g,h) &=& \frac{2}{2\pi}\int_{0}^{2\pi}g(\omega)h(\omega)f(\omega)^{2}d\omega \nonumber \\
&&+
  \frac{1}{(2\pi)^2}\int_{0}^{2\pi}\int_{0}^{2\pi}g(\omega_1)h(\omega_2)f_{4}(\omega_1,-\omega_1,\omega_{2})
d\omega_{1}d\omega_{2}, \label{eq:Vgh}
\end{eqnarray}
where $f_{4}$ denotes the fourth order cumulant density of the time
series $\{X_{t}\}$. Further, we define 
\begin{eqnarray*}
B_{G,n}(\theta) &=&  
\R \frac{2}{n}\sum_{s,t=1}^{n}c(s-t)\frac{1}{n}\sum_{k=1}^{n}
  e^{-is\omega_{k,n}}\frac{d}{d\theta}\left[\overline{\phi(\omega_{k,n};f_\theta)}\phi_{t}^{\infty}(\omega_{k,n};f_\theta)\right] \\
B_{K,n}(\theta) &=& \frac{1}{n}\sum_{k=1}^{n}
f_{n}(\omega_{k,n})\frac{df_\theta(\omega_{k,n})^{-1}}{d\theta} 
\end{eqnarray*}
where 
$\phi(\omega;f_\theta)$ and $\phi_{t}^{\infty}(\omega;f_\theta)$ are
defined in Section \ref{sec:approx},
$c(r) = \cov(X_{0},X_{r})$ and $f_{n}(\omega_k) = \int
F_{n}(\omega - \lambda) f(\lambda)d\lambda$ and $F_{n}(\cdot)$ is the Fej\'{e}r
kernel of order $n$. 


\begin{theorem}\label{theorem:bias}
Suppose that the parametric spectral densities $\{f_{\theta};\theta \in \Theta\}$ satisfy
Assumptions \ref{assum:B}. Suppose the underlying time series $\{X_{t}\}$ is a stationary time series with
spectral density $f$ and satisfies Assumption \ref{assum:TS}. Let 
$\widehat{\theta}_{n}^{(G)}$, $\widehat{\theta}_{n}^{(K)}$ and
$\widehat{\theta}_{n}^{(W)}$, and $\widehat{\theta}_{n}^{(H)}$ be defined as
in (\ref{eq:thetaBiased}). Then the asymptotic bias is
\begin{eqnarray*}
\Ex_{\theta}[\widehat{\theta}_{n}^{(G)} - \theta_n] &=&
 I(\theta)^{-1}\left(B_{K,n}(\theta_{n}) + B_{G,n}(\theta_{n})\right) +n^{-1}G(\theta_{n})+ O(n^{-3/2})\\
\Ex_{\theta}[\widehat{\theta}_{n}^{(K)} - \theta_n] &=&
I(\theta)^{-1}B_{K,n}(\theta_{n})+ n^{-1}G(\theta_{n}) + O(n^{-3/2})\\
\Ex_{\theta}[\widehat{\theta}_{n}^{(W)} - \theta_n] &=&
n^{-1}G(\theta_{n})+O\left(p^{3}n^{-3/2}+ n^{-1}p^{-K+1}\right) \\
\textrm{and} \qquad \Ex_{\theta}[\widehat{\theta}_{n}^{(H)} - \theta_n] &=&
\frac{H_{2,n}}{H_{1,n}^2} G(\theta_{n})+O\left(p^{3}n^{-3/2}+ n^{-1}p^{-K+1}\right) \\
\end{eqnarray*}
where $H_{q,n} = \sum_{t=1}^{n} h_{n}(t/n)^{q}$,
\begin{equation*}
G(\theta) = I(\theta)^{-2}
V\left(\frac{d f^{-1}_{\theta}}{d\theta},\frac{d^{2}f^{-1}_{\theta}}{d\theta^{2}}\right)  +
  2^{-1}I(\theta)^{-3}
V\left(\frac{d f^{-1}_{\theta}}{d\theta}, \frac{d f^{-1}_{\theta}}{d\theta}\right)
J\left(\frac{d^{3}f_{\theta}^{-1}}{d\theta^{3}}\right),
\end{equation*}
and $V(g,h)$ is defined in (\ref{eq:Vgh}).
\end{theorem}
PROOF. See Supplementary \ref{sec:proofbias}. \hfill $\Box$

\begin{remark}
In the case that the model is 
  linear,  then $f_{4}(\omega_{1},-\omega_{1},\omega_{2}) =
  (\kappa_{4}/\sigma^{4})f(\omega_1)f(\omega_{2})$ where $\sigma^{2}$
  and $\kappa_{4}$ is the $2$nd and $4$th order
cumulant of the innovation in the model. 

Furthermore, in the case the model is correct specification
  and linear, we can show that Assumption \ref{assum:B}(ii) implies
  that fourth order cumulant term in $V\left(\frac{d
      f^{-1}_{\theta}}{d\theta},\frac{d^{2}f^{-1}_{\theta}}{d\theta^{2}}\right)$
  and $V\left(\frac{d f^{-1}_{\theta}}{d\theta}, \frac{d
      f^{-1}_{\theta}}{d\theta}\right)$ is zero. 
This results in the   fourth order cumulant term in $G(\cdot)$ being zero. 
\end{remark}

\subsection{The bias for the AR(1) model}\label{sec:AR1bias}

In general, it is difficult to obtain a simple expression for
the bias defined in Theorem \ref{theorem:bias}, but
in the special case a model $AR(1)$ is fitted to the data the bias can
be found. In the calculation below let $\theta$ denote the AR$(1)$
coefficient for the best fitting AR$(1)$ parameter. 
We assume Gaussianity, which avoids dealing with
the fourth order spectral density. 

If the true model
is a Gaussian AR$(1)$ the bias for the various criteria is
\begin{itemize}
\item The Gaussian likelihood 
\begin{eqnarray*}
\Ex[\widehat{\theta}_{n}^{G}-\theta] =  -\frac{1}{n}\theta+O(n^{-3/2})
\end{eqnarray*}
\item The Whittle likelihood 
\begin{eqnarray*}
\Ex[\widehat{\theta}_{n}^{K}-\theta] = 
-\frac{3}{n}\theta
  +\frac{1}{n}\theta^{n-1} + O(n^{-3/2})
\end{eqnarray*}
\item The boundary corrected Whittle likelihood 
\begin{eqnarray*}
\Ex[\widehat{\theta}_{n}^{W}-\theta] = 
-\frac{2}{n}\theta
+ O(p^3n^{-3/2} +(np^{K-1})^{-1} )
\end{eqnarray*}
\item The hybrid Whittle likelihood
\begin{eqnarray*}
\Ex[\widehat{\theta}_{n}^{H}-\theta] =
-2\frac{H_{2,n}}{H_{1,n}^{2}}\theta+ O(p^3n^{-3/2} +(np^{K-1})^{-1}) .
\end{eqnarray*}
\end{itemize}
Moreover, if the Gaussian likelihood included the determinant term in
the Gaussian likelihood,
i.e. $\widetilde{\theta}_{n}^{G} = \arg\min_{\theta}[
\mathcal{L}_{n}(\theta) + n^{-1}\log|\Gamma_{n}(f_\theta)|]$, then 
\begin{eqnarray*}
\Ex[\widetilde{\theta}_{n}^{G}-\theta] =  -\frac{2}{n}\theta+O(n^{-3/2}). 
\end{eqnarray*}
We observe for the AR$(1)$ model (when the true time series is Gaussian with
an AR$(1)$ representation) that the ``true'' Gaussian likelihood with
the log-determinant term has a larger bias than the Gaussian
likelihood without the Gaussian determinant term.

The above bounds show that the Gaussian likelihood with the
log-determinant term and the boundary corrected Whittle likelihood have the
same asymptotic bias. This is substantiated in the
simulations. However, in the simulations in Section \ref{sec:specified}, we do
observe that the bias of the Gaussian likelihood is a little less than the boundary corrected Whittle. The
difference between two likelihoods is likely due to differences in the higher order
terms which are of order $O(n^{-3/2})$ (for the Gaussian likelihood)
and $O(p^{3}n^{-3/2})$ (for the boundary corrected Whittle likelihood, due to additional
estimation of the predictive DFT). 
	
\vspace{1em}

\noindent PROOF.
The inverse of the spectral density function and autocovariance function is
\begin{eqnarray*}
f_{\theta}(\omega)^{-1}= \sigma^{-2}\left(1+\theta^{2}-2\theta\cos(\omega) \right)
\quad \text{and} \quad
c(r) = \frac{\sigma^2 \theta^{|r|}}{1-\theta^2}.
\end{eqnarray*}
Thus
\begin{eqnarray*}
\frac{d}{d\theta} f_{\theta}(\omega)^{-1} = 2\sigma^{-2} (\theta-\cos \omega)
\quad \text{and} \quad
\frac{d^{2}}{d\theta^{2}} f_{\theta}(\omega)^{-1} = 2 \sigma^{-2}. 
\end{eqnarray*}
This gives 
\begin{eqnarray}
\label{eq:Itheta}
I(\theta) = -\frac{1}{2\pi}\int_{0}^{2\pi}\frac{d^{2}f_{\theta}(\omega)^{-1}}{d\theta^{2}} f(\omega)d\omega 
=-\frac{1}{\pi \sigma^2}\int_{0}^{2\pi} f(\omega) d\omega = -\frac{2}{\sigma^2}c(0).
\end{eqnarray}
Next we calculate $B_{G,n}$, since 
$\phi(\omega) = 1-\theta e^{-i\omega}$ it is is easy to show 
\begin{equation*}
\phi_{1}^{\infty}(\omega) = \theta \quad \text{and} \quad
\phi_{j}^{\infty}(\omega) = 0\textrm{ for }j\geq2.
\end{equation*} 
Therefore,
\begin{eqnarray*}
B_{G,n}(\theta) &=& \R \frac{2}{n}\sum_{t,j=1}^{n}c(t-j)\frac{1}{n}\sum_{k=1}^{n}
  e^{-it\omega_{k,n}}\frac{d}{d\theta}\left[\overline{\phi(\omega_{k,n};f_\theta)}\phi_{j}^{\infty}(\omega_{k,n};f_\theta)\right] \\
&=& \frac{2\sigma^{-2}}{n}\sum_{t=1}^{n}c(t-1)\frac{1}{n}\sum_{k=1}^{n}e^{-it\omega_{k,n}}
\frac{d}{d\theta}\left[(1-\theta e^{i\omega_{k,n}}) \theta\right] \\
&=& \frac{2\sigma^{-2}}{n}\sum_{t=1}^{n}c(t-1)
\frac{1}{n}\sum_{k=1}^{n} \left(e^{-it\omega_{k,n}} - 2\theta
    e^{-i(t-1)\omega_{k,n}} \right) \\
&=& \frac{2\sigma^{-2}}{n}\sum_{t=1}^{n}c(t-1)\left[
\frac{1}{n}\sum_{k=1}^{n}e^{-it\omega_{k,n}} - \frac{2\theta}{n}\sum_{k=1}^{n}
    e^{-i(t-1)\omega_{k,n}} \right] \\
\end{eqnarray*} 
The second summation (over $k$) is 0 unless $t\in \{1,n\}$. Therefore, 
\begin{equation}
\label{eq:BGn}
B_{G,n}(\theta) =  \frac{2\sigma^{-2}}{n}c(n-1) - \frac{4\sigma^{-2}\theta}{n}c(0). 
\end{equation} 
To calculate $B_{K,n}$ we have 
\begin{eqnarray}
\label{eq:BKn}
B_{K,n}(\theta) &=& \frac{1}{n}\sum_{k=1}^{n}
f_{n}(\omega_{k,n})\frac{df_\theta(\omega_{k,n})^{-1}}{d\theta} \nonumber\\
 &=& \frac{2\sigma^{-2}}{n}\sum_{k=1}^{n}
f_{n}(\omega_{k,n})\left(\theta - \cos(\omega_{k,n})\right) \nonumber\\
&=& 2\sigma^{-2}\left[\theta c(0) - \left(\frac{n-1}{n}\right)c(1) -
    \frac{1}{n}c(1-n) \right] \nonumber\\
&=& \frac{2\sigma^{-2}}{n}\left[ c(1) - c(n-1)\right].
\end{eqnarray}
Altogether this gives 
\begin{eqnarray}
\label{eq:Btheta}
I(\theta)^{-1}\left(B_{K,n}(\theta) +B_{G,n}(\theta)\right) &=& 
-\frac{\sigma^{2}}{2nc(0)}\left(2\sigma^{-2}\left[c(1) -
  c(n-1)\right]+2\sigma^{-2}c(n-1) - 4\sigma^{-2}\theta c(0)\right)
                                                                \nonumber\\
&=& -\frac{1}{nc(0)}\left(c(1)  - 2\theta c(0)\right)  = \frac{\theta}{n}
\end{eqnarray}
and 
\begin{eqnarray}
\label{eq:BthetaK}
I(\theta)^{-1}B_{K,n}(\theta)  =
-\frac{1}{nc(0)}\left[c(1) -
  c(n-1)\right]=-\frac{1}{n}(\theta - \theta^{n-1}). 
\end{eqnarray}

Next, we calculate $G(\theta)$. Since the third derivative of $f_{\theta}^{-1}$ with respect to
$\theta$ is zero we have 
\begin{eqnarray*}
G(\theta)  = I(\theta)^{-2}V\left(\frac{d}{d\theta} f_{\theta}^{-1},
\frac{d^{2}}{d\theta^{2}} f_{\theta}^{-1}\right)
\end{eqnarray*}
where 
\begin{eqnarray*}
V\left(\frac{d}{d\theta} f_{\theta}^{-1},
\frac{d^{2}}{d\theta^{2}} f_{\theta}^{-1}\right) &=& \frac{1}{\pi}\int_{0}^{2\pi}\left(\frac{2\theta
  - 2\cos(\omega)}{\sigma^{2}}\right)\left(\frac{2}{\sigma^{2}}\right)
  f(\omega)^{2}d\omega \\
 &=&  \frac{4}{\sigma^{4}}\frac{1}{\pi}\int_{0}^{2\pi}\left[\theta
  - \cos(\omega)\right]f(\omega)^{2}d\omega \\
 &=& \frac{8}{\sigma^{4}}\left(\theta c_2(0)- c_2(1)\right)
\end{eqnarray*}
where $\{c_{2}(r)\}$ is the autocovariance function associated with
$f(\omega)^{2}$, it is the convolution of $c(r)$ with itself;
\begin{eqnarray*}
c_{2}(r) = \sum_{\ell\in \mathbb{Z}}c(\ell)c(\ell+r)
\end{eqnarray*}
Using this expansion we have 
\begin{eqnarray*}
  V\left(\frac{d}{d\theta} f_{\theta}^{-1},
\frac{d^{2}}{d\theta^{2}} f_{\theta}^{-1}\right) 
 &=& \frac{8}{\sigma^{4}}  \sum_{\ell\in \mathbb{Z}}c(\ell)\left[\theta c(\ell) - c(\ell+1)\right] 
\end{eqnarray*}
and 
\begin{eqnarray}
\label{eq:Gtheta}
G(\theta)  = \frac{\sigma^{4}}{4c(0)^{2}}
\frac{8}{\sigma^{4}}\left(\theta c_2(0)- c_2(1)\right)
= \frac{2}{c(0)^{2}}\left(\theta c_2(0)- c_2(1)\right).
\end{eqnarray}
Putting (\ref{eq:Gtheta}) with (\ref{eq:Btheta}) gives 
\begin{eqnarray*}
\Ex[\widehat{\theta}_{n}^{G}-\theta] &\approx& 
I(\theta)^{-1}\left(B_{K,n}(\theta) +B_{G,n}(\theta)\right) +  n^{-1}G(\theta) \\
&=& 
\frac{\theta}{n}
+\frac{2}{nc(0)^{2}}\left(\theta c_2(0)- c_2(1)\right), \\
\Ex[\widehat{\theta}_{n}^{K}-\theta] &\approx& I(\theta)^{-1}B_{K,n}(\theta) +  n^{-1}G(\theta)\\
&=& 
-\frac{1}{n}(\theta -\theta^{n-1}) 
+\frac{2}{nc(0)^{2}}\left(\theta c_2(0)- c_2(1)\right), \\
\Ex[\widehat{\theta}_{n}^{W}-\theta] &\approx& \frac{2}{nc(0)^{2}}\left(\theta c_2(0)- c_2(1)\right), \\
\Ex[\widehat{\theta}_{n}^{H}-\theta] &\approx& \frac{2}{c(0)^{2}} \frac{H_{2,n}}{H_{1,n}^2}\left(\theta c_2(0)- c_2(1)\right).
\end{eqnarray*}
It is not entirely clear how to access the above. So now we consider
the case that the model is fully specified. 
Under correct specification we have 
\begin{eqnarray*}
V\left(\frac{d}{d\theta} f_{\theta}^{-1},
\frac{d^{2}}{d\theta^{2}} f_{\theta}^{-1}\right)  &=& 
\frac{2\sigma^2}{\pi}\int_{0}^{2\pi}f(\omega)^{2}\frac{df_{\theta}(\omega)^{-1}}{d\theta}d\omega \\
&=&
    -\frac{2\sigma^2}{\pi}\int_{0}^{2\pi}f(\omega)^{2}\left(\frac{1}{f_\theta(\omega)^{2}}\right)
\frac{df_{\theta}(\omega)}{d\theta}d\omega \\
&=&
    -\frac{2\sigma^2}{\pi}\int_{0}^{2\pi}\frac{df_{\theta}(\omega)}{d\theta}d\omega 
    =
 -\frac{2\sigma^2}{\pi}\frac{d}{d\theta}\int_{0}^{2\pi}f_{\theta}(\omega)d\omega \\
 &=& -4\sigma^2 \frac{d}{d\theta} c(0)
= -\frac{8\sigma^4 \theta}{(1-\theta^{2})^{2}}.  
\end{eqnarray*}
Thus $I(\theta)^{-2}V\left(\frac{d}{d\theta} f_{\theta}^{-1},
\frac{d^{2}}{d\theta^{2}} f_{\theta}^{-1}\right) = -2\theta/n$. 
Substituting this into the above we have 
\begin{eqnarray*}
\Ex[\widehat{\theta}_{n}^{(G)}-\theta] &=& \frac{1}{n}\theta
                                           -\frac{2}{n}\theta + O(n^{-3/2}) =
                                           -\frac{1}{n}\theta + O(n^{-3/2}), \\
\Ex[\widehat{\theta}_{n}^{(K)}-\theta] &=& 
-\frac{1}{n}\left[\theta -
  \theta^{n-1}\right]  -\frac{2}{n}\theta + o(n^{-1})\approx -\frac{3}{n}\theta +\frac{1}{n}\theta^{n-1}+ O(n^{-3/2}), \\
\Ex[\widehat{\theta}_{n}^{(W)}-\theta] &=& -\frac{2}{n}\theta + O(n^{-3/2}), \\
\Ex[\widehat{\theta}_{n}^{(H)}-\theta] &=&
                                                 -2\frac{H_{2,n}}{H_{1,n}^2}\theta + O(n^{-3/2}).
\end{eqnarray*}
This proves the main part of the assertion. 
To compare the above bias with the ``true'' Gaussian likelihood,
we consider the Gaussian likelihood with the log determinant term. First, consider the correlation
of AR(1) matrix 
$(A_n)_{s,t} = \theta^{|s-t|}$. Then,
\begin{eqnarray*}
A_{n+1} = 
\begin{pmatrix}
A_n & B_n \\
B_n^{\prime}  & 1
\end{pmatrix}, \qquad B_n = (\theta^{n}, ..., \theta)^{\prime}.
\end{eqnarray*}
Therefore, using block matrix determinant identity, 
$|A_{n+1}| = |A_n| (1-B_{n}^\prime A_n^{-1} B_{n})$. Moreover, it is easy to show
$A_n R_n = B_n$, where $R_n = (0,...,0,\theta)^\prime$. Thus
\begin{eqnarray*}
|A_{n+1}| = |A_n| ( 1- B_n^\prime R_n) = |A_n| (1-\theta^2).
\end{eqnarray*} Using iteration, $|A_n| = (1-\theta^2)^{n-2} |A_2| = (1-\theta^2)^{n-1}$ and 
thus,
\begin{eqnarray*}
|\Gamma_n(f_\theta)| = \left|\frac{\sigma^2}{1-\theta^2} A_n\right| = \left( \frac{\sigma^2}{1-\theta^2}\right)^{n} (1-\theta^2)^{n-1}
= \frac{(\sigma^2)^n}{1-\theta^2}.
\end{eqnarray*} Then, by simple calculus,
\begin{eqnarray*}
\frac{d}{d\theta} n^{-1}\log |\Gamma_n(f_\theta)| = \frac{2\theta}{n(1-\theta^2)} = \frac{2\sigma^{-2}}{n}c(1).
\end{eqnarray*} and thus,
\begin{eqnarray*}
\Ex[\widetilde{\theta}^{G}-\theta] &\approx&
I(\theta)^{-1} \left( B_{K,n}(\theta) + B_{G,n}(\theta) + \frac{d}{d\theta} n^{-1}\log |\Gamma_n(f_\theta)|\right)
+ n^{-1} G(\theta) \\
&=& I(\theta)^{-1} \frac{1}{n} (2\sigma^{-2} c(1) - 4\sigma^{-2} \theta c(0) +2\sigma^{-2} c(1))
- \frac{2\theta}{n} = - \frac{2\theta}{n},
\end{eqnarray*}
which proves the results.
\hfill $\Box$



\subsection{Proof of Theorem \ref{theorem:bias}}\label{sec:proofbias}

In Theorem \ref{theorem:equivalence} we showed that 
\begin{eqnarray*}
|\widehat{\theta}_{n}^{(W)}-\widetilde{\theta}_{n}^{(W)}| = 
O_{p}\left(\frac{p^{3}}{n^{3/2}}+\frac{1}{np^{K-1}}\right) 
\quad \textrm{and} \quad
|\widehat{\theta}_{n}^{(H)}-\widetilde{\theta}_{n}^{(H)}| = 
O_{p}\left(\frac{p^{3}}{n^{3/2}}+\frac{1}{np^{K-1}}\right), 
\end{eqnarray*}
where 
$\widehat{\theta}_{n}^{(W)}=\arg\min \widehat{W}_{p,n}(\theta)$,
$\widetilde{\theta}_{n}^{(W)}=\arg\min W_{n}(\theta)$, 
$\widehat{\theta}_{n}^{(H)}=\arg\min \widehat{H}_{p,n}(\theta)$,
and 
\\ $\widetilde{\theta}_{n}^{(H)}=\arg\min H_{n}(\theta)$.
We will show that the asymptotic bias of $\widetilde{\theta}_{n}^{(W)}$ and
$\widetilde{\theta}_{n}^{(H)}$ (under certain conditions on the taper)
are of order $O(n^{-1})$, thus if $p^{3}n^{-1/2}\rightarrow 0$
as $n, p \rightarrow \infty$, then the infeasible
estimators and feasible estimators share the same asymptotic
bias. Therefore in the proof we obtain the bias of the infeasible
estimators.

Now we obtain a general expansion (analogous to the Bartlett
correction). Let 
$L_{n}(\cdot)$ denote the general minimization criterion
(it can be $\mathcal{L}_{n}(\theta)$, $K_{n}(\theta)$, $W_{n}(\theta)$,  or $H_{n}(\theta)$) and
$\widehat{\theta} = \arg\min L_{n}(\theta)$. For all the criteria,
it is easily shown that
\begin{equation*}
(\widehat{\theta} - \theta) = U(\theta)^{-1}\frac{d L_{n}(\theta)}{d\theta}+O_{p}(n^{-1})
\end{equation*}
where $U(\theta) = -\Ex[\frac{d^{2} L_{n}}{d\theta^{2}}]$ and 
\begin{equation*}
\frac{dL_{n}(\theta)}{d\theta} + (\widehat{\theta} - \theta)
  \frac{d^{2}L_{n}(\theta)}{d\theta^{2}} +
 \frac{1}{2}(\widehat{\theta} - \theta)^{2}
  \frac{d^{3}L_{n}(\theta)}{d\theta^{3}} = O_{p}(n^{-3/2}).
\end{equation*}
Ignoring the
probabilistic error, the first and second order expansions are
\begin{equation}
\label{eq:expan1}
(\widehat{\theta} - \theta) \approx U(\theta)^{-1}\frac{d L_{n}(\theta)}{d\theta},
\end{equation}
and 
\begin{equation*}
\frac{dL_{n}(\theta)}{d\theta} + (\widehat{\theta} - \theta)
  \frac{d^{2}L_{n}(\theta)}{d\theta^{2}} +
 \frac{1}{2}(\widehat{\theta} - \theta)^{2}
  \frac{d^{3}L_{n}(\theta)}{d\theta^{3}}\approx 0.
\end{equation*}
The method described below follows the 
Bartlett correction described in \cite{p:bar-53} and \cite{p:cox-68}.
Taking expectation of the above we have 
\begin{eqnarray*}
&&\Ex\left[\frac{dL_{n}(\theta)}{d\theta}\right] + 
\Ex\left[(\widehat{\theta} - \theta)\frac{d^{2}L_{n}(\theta)}{d\theta^{2}}\right] +
 \frac{1}{2} \Ex\left[(\widehat{\theta} - \theta)^{2}
  \frac{d^{3}L_{n}(\theta)}{d\theta^{3}}\right] \\
&&= \Ex\left[\frac{dL_{n}(\theta)}{d\theta}\right] +
\Ex\left[(\widehat{\theta} -
    \theta)\right]\Ex\left[\frac{d^{2}L_{n}(\theta)}{d\theta^{2}}\right]
    +\cov\left[(\widehat{\theta} -
    \theta),\frac{d^{2}L_{n}(\theta)}{d\theta^{2}}\right] \\
&&\quad + \frac{1}{2} \Ex\left[(\widehat{\theta} - \theta)^{2}\right]
  \Ex\left[\frac{d^{3}L_{n}(\theta)}{d\theta^{3}}\right] + 
\frac{1}{2} \cov\left[(\widehat{\theta} - \theta)^{2},
\frac{d^{3}L_{n}(\theta)}{d\theta^{3}}\right].
\end{eqnarray*}
Substituting $(\widehat{\theta} - \theta) \approx
U(\theta)^{-1}\frac{d L_{n}(\theta)}{d\theta}$ into the last three
terms on the right hand side of the above gives
\begin{eqnarray*}
&&\Ex\left(\frac{d L_{n}(\theta)}{d\theta} \right) -
  U(\theta) \Ex(\widehat{\theta} - \theta)+
U(\theta)^{-1}\cov\left(\frac{d
  L_{n}(\theta)}{d\theta}, \frac{d^{2}
  L_{n}(\theta)}{d\theta^{2}}\right) \\
&& \quad +2^{-1}U(\theta)^{-2} 
\Ex\left(\frac{d
  L_{n}(\theta)}{d\theta}\right)^{2}\Ex\left(\frac{d^{3}
   L_{n}(\theta)}{d\theta^{3}}\right) \\
&& \quad + 2^{-1}U(\theta)^{-2} 
\cov\left(\left(\frac{d
  L_{n}(\theta)}{d\theta}\right)^{2},\frac{d^{3}
   L_{n}(\theta)}{d\theta^{3}}\right)\approx 0.
\end{eqnarray*}
Using the above to solve for $\Ex(\widehat{\theta} - \theta)$ gives 
\begin{eqnarray*}
\Ex(\widehat{\theta} - \theta) 
&=& U(\theta)^{-1}\Ex\left(\frac{d L_{n}(\theta)}{d\theta} \right) +
U(\theta)^{-2}\cov\left(\frac{d
  L_{n}(\theta)}{d\theta}, \frac{d^{2}
  L_{n}(\theta)}{d\theta^{2}}\right) \nonumber\\
&& +2^{-1}U(\theta)^{-3} 
\Ex\left(\frac{d
  L_{n}(\theta)}{d\theta}\right)^{2}\Ex\left(\frac{d^{3}
   L_{n}(\theta)}{d\theta^{3}}\right)+ 2^{-1}U(\theta)^{-3} 
\cov\left(\left(\frac{d
  L_{n}(\theta)}{d\theta}\right)^{2},\frac{d^{3}
   L_{n}(\theta)}{d\theta^{3}}\right) \nonumber\\
&=&  U(\theta)^{-1}\Ex\left(\frac{d L_{n}(\theta)}{d\theta} \right) +
U(\theta)^{-2}\cov\left(\frac{d
  L_{n}(\theta)}{d\theta}, \frac{d^{2}
  L_{n}(\theta)}{d\theta^{2}}\right) \nonumber\\
&& +2^{-1}U(\theta)^{-3} \left[\var\left(\frac{d
  L_{n}(\theta)}{d\theta}\right)+
\left\{\Ex\left[\frac{d
  L_{n}(\theta)}{d\theta}\right]\right\}^{2}\right]\Ex\left(\frac{d^{3}
   L_{n}(\theta)}{d\theta^{3}}\right) \nonumber\\
&&+ 2^{-1}U(\theta)^{-3} 
\cov\left(\left(\frac{d
  L_{n}(\theta)}{d\theta}\right)^{2},\frac{d^{3}
   L_{n}(\theta)}{d\theta^{3}}\right). 
\end{eqnarray*}
Thus 
\begin{eqnarray}
\label{eq:BIASmle}
\Ex(\widehat{\theta} - \theta) &=& I_{0}+I_{1} + I_{2} + I_{3} + I_{4}
\end{eqnarray}
where
\begin{eqnarray*}
I_{0}&=&  U(\theta)^{-1}\Ex\left(\frac{d L_{n}(\theta)}{d\theta}
         \right) \\
I_{1}&=&U(\theta)^{-2}\cov\left(\frac{d L_{n}(\theta)}{d\theta}, \frac{d^{2}
  L_{n}(\theta)}{d\theta^{2}}\right)  \nonumber \\
I_{2}&=&2^{-1}U(\theta)^{-3} \var\left(\frac{d
  L_{n}(\theta)}{d\theta}\right)\Ex\left(\frac{d^{3}
   L_{n}(\theta)}{d\theta^{3}}\right) \nonumber\\
I_{3} &=& 2^{-1}U(\theta)^{-3} \left\{\Ex\left(\frac{d
  L_{n}(\theta)}{d\theta}\right)\right\}^{2}\Ex\left(\frac{d^{3}
   L_{n}(\theta)}{d\theta^{3}}\right) \\
I_{4}&=& 2^{-1}U(\theta)^{-3} 
\cov\left(\left(\frac{d
  L_{n}(\theta)}{d\theta}\right)^{2},\frac{d^{3}
   L_{n}(\theta)}{d\theta^{3}}\right). 
\end{eqnarray*}
Note that the term $\Ex\left(\frac{d
    L_{n}(\theta)}{d\theta} \right)$ will be 
different for the four quasi-likelihoods (and will be of order $O(n^{-1})$). However
the remaining terms are asymptotically
the same for three quasi-likelihoods and will be slightly different for
the hybrid Whittle likelihood. 

\vspace{2mm}
\noindent \underline{The first derivative}
We first obtain expressions for $\Ex \left( \frac{d L_{n}(\theta)}{d\theta}\right)$ for the four quasi-likelihoods:
\begin{eqnarray*}
\Ex\left(\frac{d
    K_{n}(\theta)}{d\theta} \right) = \frac{1}{n}\sum_{k=1}^{n}
\Ex[|J_{n}(\omega_{k,n})|^{2}]\frac{d}{d\theta}f_{\theta}(\omega_{k,n})^{-1}
                                        = \frac{1}{n}\sum_{k=1}^{n}f_{n}(\omega_{k,n})
\frac{d}{d\theta}f_{\theta}(\omega_{k,n})^{-1} = B_{K,n}(\theta),
\end{eqnarray*}
where $f_{n}(\omega) = \int
F_{n}(\omega - \lambda) f(\lambda)d\lambda$ and $F_{n}$ is the Fej\'{e}r
kernel of order $n$. 

To obtain the expected derivative of $\mathcal{L}_{n}(\theta)$ we
recall that 
\begin{eqnarray*}
\Ex\left[\frac{d}{d\theta}\mathcal{L}_{n}(\theta)\right] &=&
                                                  \Ex\left[\frac{d}{d\theta} K_{n}(\theta)\right] 
+\Ex\left[n^{-1}\underline{X}_{n}^{\prime}F_{n}^{*}\frac{d}{d\theta}\Delta_{n}(f_\theta^{-1})D_{n}(f_\theta) \underline{X}_{n}\right].
\end{eqnarray*}
Now by replacing $D_{n}(f_\theta) $ with $D_{\infty,n}(f_\theta) $ and
using (\ref{eq:DDDexpand}) we have 
\begin{eqnarray*}
\Ex\left[\frac{d}{d\theta}\mathcal{L}_{n}(\theta)\right] 
&=&   \Ex\left[\frac{d}{d\theta} K_{n}(\theta)\right] 
+\Ex\left[n^{-1}\underline{X}_{n}^{\prime}F_{n}^{*}\frac{d}{d\theta}\Delta_{n}(f_\theta^{-1})D_{\infty,n}(f_\theta) 
\underline{X}_{n}\right] \\
&& +\Ex\left[n^{-1}\underline{X}_{n}^{\prime}F_{n}^{*}\frac{d}{d\theta} \Delta_n\left(f_{\theta}^{-1}\right) \left(D_{n}(f_\theta)-D_{\infty,n}(f_\theta) \right)\underline{X}_{n} \right]\\
&=&       \Ex\left[\frac{d}{d\theta} K_{n}(\theta)\right] +                            
n^{-1}\sum_{s,t=1}^{n}c(s-t)
\frac{1}{n}\sum_{k=1}^{n}
  e^{-is\omega_{k,n}}\frac{d}{d\theta} \varphi_{t,n}(\omega_{k,n};f_\theta) \\
&&+n^{-1}\Ex\left[\underline{X}_{n}^{\prime}F_{n}^{*} \frac{d}{d\theta} \Delta_n\left(f_{\theta}^{-1}\right) \left(D_{n}(f_\theta)-D_{\infty,n}(f_\theta) \right)\underline{X}_{n} \right]
\end{eqnarray*}
where 
$\varphi_{t,n}(\omega;f_\theta) = \sigma^{-2}
\left[ \overline{\phi(\omega;f_\theta)}\phi_{t}^{\infty}(\omega;f_\theta)    +
  e^{i\omega}\phi(\omega;f_\theta)\overline{\phi_{n+1-t}^{\infty}(\omega;f_\theta)}\right]$. The
first term on the RHS of the above is $B_{K,n}(\theta)$. Using the
change of variables $t^{\prime}=n+1-t$, 
the second term in RHS above can be written as 
\begin{eqnarray*}
&& n^{-1}\sum_{s,t=1}^{n}c(s-t)
\frac{1}{n}\sum_{k=1}^{n}
  e^{-is\omega_{k,n}}\frac{d}{d\theta} \varphi_{t,n}(\omega_{k,n};f_\theta) \\
&&= n^{-1}\sum_{s,t=1}^{n}c(s-t)
\frac{1}{n}\sum_{k=1}^{n}\frac{d}{d\theta} \left[e^{-is\omega_{k,n}}\overline{\phi(\omega_{k,n};f_\theta)}\phi_{t}^{\infty}(\omega_{k,n};f_\theta)+
  e^{-i(s-1)\omega_{k,n}}\phi(\omega_{k,n};f_\theta)\overline{\phi_{n+1-t}^{\infty}(\omega_{k,n};f_\theta)}\right]
  \\
&&= n^{-1}\sum_{s,t=1}^{n}c(s-t)
\frac{1}{n}\sum_{k=1}^{n}e^{-is\omega_{k,n}}\frac{d}{d\theta} \overline{\phi(\omega_{k,n};f_\theta)}\phi_{t}^{\infty}(\omega_{k,n};f_\theta)\\
&&\quad + n^{-1}\sum_{s,t^{\prime}=1}^{n}c(s-n-1+t^{\prime})
\frac{1}{n}\sum_{k=1}^{n}e^{-i(s-1)\omega_{k,n}}\frac{d}{d\theta} \phi(\omega_{k,n};f_\theta)\overline{\phi_{t^{\prime}}^{\infty}(\omega_{k,n};f_\theta)}\quad 
\quad (\textrm{let} \quad t^{\prime} = n+1-t) \\
&&= \R \frac{2}{n}\sum_{s,t=1}^{n}c(s-t)
\frac{1}{n}\sum_{k=1}^{n}e^{-is\omega_{k,n}}\frac{d}{d\theta} \overline{\phi(\omega_{k,n};f_\theta)}\phi_{t}^{\infty}(\omega_{k,n};f_\theta) = B_{G,n}(\theta).
\end{eqnarray*}
Finally, by using Corollary \ref{corollary:DDeriv} 
we have 
\begin{eqnarray*}
n^{-1}\left\|
\underline{X}_{n}^{\prime}F_{n}^{*} \frac{d}{d\theta^{}}\Delta_{n}(f_\theta^{-1})\left[D_{n}(f_\theta)
  - D_{\infty,n}(f_\theta)\right]
\underline{X}_{n}\right\|_{\Ex,1} = O(n^{-K+1/2}).
\end{eqnarray*}
Thus the derivative of the Gaussian likelihood is 
\begin{eqnarray*}
\Ex\left(\frac{d\mathcal{L}_{n}(\theta)}{d\theta}\right)
= B_{K,n}(\theta) + B_{G,n}(\theta)+ O(n^{-K+1/2}).
\end{eqnarray*} 
Next we consider the boundary corrected Whittle likelihood. By using
that
\begin{eqnarray*}
\Ex[\widetilde{J}_{n}(\omega_{k,n};f)\overline{J_{n}(\omega_{k,n})}]
= f(\omega_{k,n}) 
\end{eqnarray*} we have 
\begin{eqnarray*}
\Ex\left(\frac{dW_{n}(\theta)}{d\theta}\right) &=&
\frac{1}{n}\sum_{k=1}^{n} 
\Ex[ \widetilde{J}_{n}(\omega_{k,n};f_{})\overline{J_{n}(\omega_{k,n})}]
\frac{d}{d\theta}f_{\theta}(\omega_{k,n})^{-1} \\
 &=& \frac{1}{n}\sum_{k=1}^{n} 
f(\omega_{k,n}) \frac{d}{d\theta}f_{\theta}(\omega_{k,n})^{-1}.
\end{eqnarray*}
Finally, the analysis of $H_{n}(\theta)$ is identical to the
analysis of $W_{n}(\theta)$ and we obtain
\begin{eqnarray*}
\Ex\left(\frac{dH_{n}(\theta)}{d\theta}\right) &=&
\frac{1}{n}\sum_{k=1}^{n} 
\Ex[ \widetilde{J}_{n}(\omega_{k,n};f) \overline{J_{\underline{h}_n,n}(\omega_{k,n})} ] 
\frac{d}{d\theta}f_{\theta}(\omega_{k,n})^{-1} \\
 &=& \frac{1}{n}\sum_{k=1}^{n} 
f(\omega_{k,n}) \frac{d}{d\theta}f_{\theta}(\omega_{k,n})^{-1}.
\end{eqnarray*}
In summary, evaluating the above at the best fitting parameter $\theta_{n}$
and by Assumption \ref{assum:B}(ii) gives 
\begin{eqnarray}
\Ex\left(\frac{d K_{n}(\theta)}{d\theta} \right)\rfloor_{\theta = \theta_n} &=& B_{K,n}(\theta_n)\nonumber\\
\Ex\left(\frac{d\mathcal{L}_{n}(\theta)}{d\theta}\right) \rfloor_{\theta = \theta_n}
&=& B_{K,n}(\theta_n) + B_{G,n}(\theta_n) +  O(n^{-K+1/2}) \nonumber\\
\text{and} \qquad  \Ex\left(\frac{dW_{n}(\theta)}{d\theta}\right) \rfloor_{\theta = \theta_n}
 &=& \Ex\left(\frac{dH_{n}(\theta)}{d\theta}\right) \rfloor_{\theta = \theta_n}
 = 0.
\label{eq:firstterm}
\end{eqnarray}
It can be shown that $ B_{K,n}(\theta_n)  = O(n^{-1})$ and
$B_{G,n}(\theta_n) = O(n^{-1})$. These terms could be negative
or positive so there is no clear cut answer as to whether $B_{K,n}(\theta_n) $ or
$B_{K,n}(\theta_n) +B_{G,n}(\theta_n)$ is larger (our
simulations results suggest that often $B_{K,n}(\theta_n)$ tends to be
larger).


\vspace{1em}

\noindent \underline{The second and third order derivatives} 
The analysis of all the higher order terms will require comparisons between the
derivatives of $\mathcal{L}_{n}(\theta), K_{n}(\theta),
W_{n}(\theta)$ and $H_{n}(\theta)$. We first represent the
derivatives of the Gaussian likelihood in terms of the Whittle likelihood
\begin{eqnarray*}
\frac{d^{i}\mathcal{L}_{n}(\theta)}{d\theta^{i}} = 
\frac{d^{i}K_{n}(\theta)}{d\theta^{i}} + 
\Ex\left[n^{-1}\underline{X}_{n}^{\prime}F_{n}^{*} \frac{d^{i}}{d\theta^{i}}\Delta_{n}(f_\theta^{-1})D_{n}(f_\theta) \underline{X}_{n}\right].
\end{eqnarray*}
By using (\ref{eq:DiffDerivatives}), for $1\leq i \leq 3$ we have
\begin{eqnarray}
\label{eq:Bounds11}
\left\|n^{-1}
\underline{X}_{n}^{\prime}F_{n}^{*}\frac{d^{i}}{d\theta^{i}}\Delta_{n}(f_\theta^{-1})D_{n}(f_\theta)\underline{X}_{n}\right\|_{\Ex,1}
  = O(n^{-1}).
\end{eqnarray}
Similarly, we represent the derivatives of
$W_{n}(\theta)$ and $H_{p,n}(\theta)$ in terms of the derivatives of $K_{n}(\theta)$ 
\begin{eqnarray*}
\frac{d^{i}W_{n}(\theta)}{d\theta^{i}} &=&
                                             \frac{d^{i}K_{n}(\theta)}{d\theta^{i}}
                                             + C_{i,n} \\
\frac{d^{i}H_{n}(\theta)}{d\theta^{i}} &=&
                                             \frac{d^{i}K_{n,\underline{h}_n}(\theta)}{d\theta^{i}}
                                             + D_{i,n} 
\end{eqnarray*}
where  $K_{n,\underline{h}_n}(\theta) = n^{-1}\sum_{k=1}^{n}\frac{
                            J_{n}(\omega_{k,n}) \overline{J_{n,\underline{h}_n}(\omega_{k,n})} }{f_\theta(\omega_{k,n})}$
and 
\begin{eqnarray*}
C_{i,n} &=& \frac{1}{n}\sum_{k=1}^{n}\frac{d^{i}}{d\theta^{i}}\frac{
                            \widehat{J}_{n}(\omega_{k,n};f)\overline{J_{n}(\omega_{k,n})} }{f_\theta(\omega_{k,n})}
=
    \frac{1}{n}\underline{X}_{n}^{\prime}F_{n}^{*}\Delta_n(\frac{d^{i}}{d\theta^{i}}f_\theta^{-1})
D_{n}(f)\underline{X}_{n} \\
D_{i,n} &=& \frac{1}{n}\sum_{k=1}^{n}\frac{d^{i}}{d\theta^{i}}\frac{
                            \widehat{J}_{n}(\omega_{k,n};f)\overline{J_{n,\underline{h}_n}(\omega_{k,n})}}{f_\theta(\omega_{k,n})}
= \frac{1}{n}\underline{X}_{n}^{\prime}H_{n}F_{n}^{*}\Delta_n(\frac{d^{i}}{d\theta^{i}}f_\theta^{-1})
D_{n}(f)\underline{X}_{n},
\end{eqnarray*}
where $H_{n}=\diag(h_{1,n},\ldots,h_{n,n})$.
In the analysis of the first order derivative obtaining an
exact bound between each ``likelihood'' and the Whittle likelihood was
important. However, for the higher order derivatives we simply require a moment bound on
the difference. To bound $C_{i,n}$, we use that
\begin{eqnarray*}
&& \left\|F_{n}^{*}\Delta_n(\frac{d^{i}}{d\theta^{i}}f_\theta^{-1})
D_{n}(f)\right\|_{1}\\
&& \quad \leq \left\|F_{n}^{*}\Delta_n(\frac{d^{i}}{d\theta^{i}}f_\theta^{-1})
[D_{n}(f)-D_{\infty,n}(f)]\right\|_{1} 
+ \left\|F_{n}^{*}\Delta_n(\frac{d^{i}}{d\theta^{i}}f_\theta^{-1})
D_{\infty,n}(f)\right\|_{1}.
\end{eqnarray*}
We use a similar method to the proof of Theorem \ref{theorem:approx}, equation (\ref{eq:approx1A})
and Theorem \ref{theorem:Bound}, equation (\ref{eq:Bound1}) with
$\Delta_{n}(\frac{d^{i}}{d\theta^{i}} f_\theta^{-1})$ and
$D_{\infty,n}(f)$ replacing $\Delta_{n}(f_\theta^{-1})$ and
$D_{\infty,n}(f_\theta)$ respectively together with Assumption
\ref{assum:B} and \ref{assum:TS}.
By using the proof of Theorem \ref{theorem:approx}, equation
(\ref{eq:approx1A}), we have $\|F_{n}^{*}\Delta_n(\frac{d^{i}}{d\theta^{i}}f_\theta^{-1})
[D_{n}(f)-D_{\infty,n}(f)]\|_{1} = O(n^{-K+1})$. Similarly, by
using the proof of Theorem \ref{theorem:Bound}, equation
(\ref{eq:Bound1}) we have $\|F_{n}^{*}\Delta_n(\frac{d^{i}}{d\theta^{i}}f_\theta^{-1})
D_{\infty,n}(f)\|_{1}=O(1)$. Altogether this gives 
\begin{eqnarray}
\label{eq:Bounds12}
(\Ex |C_{i,n}|^2)^{1/2} = 
n^{-1} \left\|\underline{X}_{n}^\prime F_{n}^{*}\Delta_n(\frac{d^{i}}{d\theta^{i}}f_\theta^{-1})
D_{n}(f) \underline{X}_{n}\right\|_{\Ex, 2} =
O(n^{-1}).
\end{eqnarray}
For the hybrid likelihood, we use that $\sup_{t,n}|h_{t,n}| < \infty$,
this gives 
\begin{eqnarray*}
\|H_{n}F_{n}^{*}\Delta_n(\frac{d^{i}}{d\theta^{i}}f_\theta^{-1})
D_{n}(f)\|_1 \leq ( \sup_{t} h_{t,n} ) \times \|F_{n}^{*}\Delta_n(\frac{d^{i}}{d\theta^{i}}f_\theta^{-1})
D_{n}(f)\|_1 = O(1).
\end{eqnarray*} 
Therefore, under the condition that $\{h_{t,n}\}$ is a bounded sequence 
\begin{eqnarray}
\label{eq:Bounds121}
(\Ex |D_{i,n}|^2)^{1/2} = 
n^{-1} \left\|H_n\underline{X}_{n}^\prime F_{n}^{*}\Delta_n(\frac{d^{i}}{d\theta^{i}}f_\theta^{-1})
D_{n}(f) \underline{X}_{n}\right\|_{\Ex, 2} =
O\left(n^{-1}\right).
\end{eqnarray}
Thus the expectations of the derivatives are 
\begin{eqnarray*}
\Ex\left(\frac{d^{i}\mathcal{L}_{n}(\theta)}{d\theta^{i}}\right) &=&  
\Ex\left(\frac{d^{i}K_{n}(\theta)}{d\theta^{i}}\right) + O(n^{-1}) \\
\Ex\left(\frac{d^{i}W_{n}(\theta)}{d\theta^{i}}\right) &=&  
\Ex\left(\frac{d^{i}K_{n}(\theta)}{d\theta^{i}}\right) +  O(n^{-1}) \\
\Ex\left(\frac{d^{i}H_{n}(\theta)}{d\theta^{i}}\right) &=&  
\Ex\left(\frac{d^{i}K_{n}(\theta)}{d\theta^{i}}\right) +  O\left(n^{-1}\right). 
\end{eqnarray*}
This gives the expectation of the second and third derivatives
of all likelihoods in
terms of $I(\theta)$ and $J(\frac{d^{3}f_\theta^{-1}}{d\theta^{3}})$:
\begin{eqnarray*}
\Ex\left(\frac{d^{2}L_{n}(\theta)}{d\theta^{2}}\right) = -I(\theta) + O(n^{-1}),
\quad \text{and} \quad
\Ex\left(\frac{d^{3}L_{n}(\theta)}{d\theta^{3}}\right) = J(\frac{d^{3}f_\theta^{-1}}{d\theta^{3}}) + O(n^{-1}).
\end{eqnarray*}

\vspace{1em}

\noindent \underline{Bounds for the covariances between the derivatives}
The terms $I_{1},I_{2}$ and $I_{4}$ all contain the covariance between
various likelihoods and its derivatives.  Thus to obtain expression and bounds for
these terms we use that 
\begin{eqnarray}
\label{eq:KKKLLL}
\var\left(\frac{d^{i}}{d\theta^{i}}K_{n}(\theta) \right) = O(n^{-1}),
\end{eqnarray}
where the above can be proved using \cite{b:bri-01}, Theorem
4.3.2. Further, if the 
data taper $\{h_{t,n}\}$ is such that $h_{t,n} = c_{n}h_{n}(t/n)$ where
$c_{n} = n/H_{1,n}$ and $h_{n}:[0,1] \rightarrow \mathbb{R}$ is a sequence of taper functions
 which satisfy the taper conditions
in Section 5, \cite{p:dah-88}, then 
\begin{eqnarray} \label{eq:Ktaperbound}
\var\left(\frac{d^{i}}{d\theta^{i}}K_{n,\underline{h}_n}(\theta) \right) = O\left(\frac{H_{2,n}}{H_{1,n}^{2}}\right).
\end{eqnarray}
By using (\ref{eq:Bounds11}), (\ref{eq:Bounds12}), and (\ref{eq:KKKLLL}) we have 
\begin{eqnarray*}
\cov\left(\frac{d
  \mathcal{L}_{n}(\theta)}{d\theta}, \frac{d^{2}
  \mathcal{L}_{n}(\theta)}{d\theta^{2}}\right) &=& \cov\left(\frac{d
  K_{n}(\theta)}{d\theta}, \frac{d^{2}
  K_{n}(\theta)}{d\theta^{2}}\right) + O(n^{-3/2}) \\
\cov\left(\frac{d
  W_{n}(\theta)}{d\theta}, \frac{d^{2}
  W_{n}(\theta)}{d\theta^{2}}\right)  &=& \cov\left(\frac{d
  K_{n}(\theta)}{d\theta}, \frac{d^{2}
  K_{n}(\theta)}{d\theta^{2}}\right) + O(n^{-3/2}) \\
\var\left(\frac{d
  \mathcal{L}_{n}(\theta)}{d\theta}\right) &=& \var\left(\frac{d
  K_{n}(\theta)}{d\theta}\right) + O(n^{-3/2}) \\
 \var\left(\frac{d
  W_{n}(\theta)}{d\theta}\right)  &=& \var\left(\frac{d
  K_{n}(\theta)}{d\theta}\right) + O(n^{-3/2}).
\end{eqnarray*} 
For the hybrid Whittle likelihood, by using (\ref{eq:Bounds121}) and (\ref{eq:Ktaperbound})
\begin{eqnarray*}
\cov\left(\frac{d
  H_{n}(\theta)}{d\theta}, \frac{d^{2}
  H_{n}(\theta)}{d\theta^{2}}\right) &=& \cov\left(\frac{d
  K_{n,\underline{h}_n}(\theta)}{d\theta}, \frac{d^{2}
  K_{n,\underline{h}_n}(\theta)}{d\theta^{2}}\right) +
O\left(\frac{H_{2,n}^{1/2}}{nH_{1,n}}\right)
\\
 \var\left(\frac{d
  H_{n}(\theta)}{d\theta}\right)  &=& \var\left(\frac{d
  K_{n,\underline{h}_n}(\theta)}{d\theta}\right) + O\left(\frac{H_{2,n}^{1/2}}{nH_{1,n}}\right).
\end{eqnarray*}
Using that $H_{2,n}/H_{1,n}^2 \sim n^{-1}$, we show that the above
error terms $O(H_{2,n}^{1/2}/(nH_{1,n}))$ 
(for the hybrid Whittle likelihood) is the same as the other likelihoods.
Next, having reduced the above covariances to those of the derivatives of
$K_{n}(\theta)$  and $K_{n,\underline{h}_n}(\theta)$. 
We first focus on $K_{n}(\theta)$. By using the expressions for
cumulants of DFTs given in \cite{b:bri-01}, Theorem
4.3.2 and well-known cumulant arguments we can show that
\begin{eqnarray*}
 \cov\left(\frac{d
  K_{n}(\theta)}{d\theta}, \frac{d^{2}
  K_{n}(\theta)}{d\theta^{2}}\right) &=&
                                           n^{-1}V\left(\frac{df_{\theta}^{-1}}{d\theta}, 
\frac{d^{2}f_{\theta}^{-1}}{d\theta^{2}}\right) + O(n^{-2})\\
\textrm{ and }\var\left(\frac{d
  K_{n}(\theta)}{d\theta}\right)&=&  n^{-1}V\left(\frac{df_{\theta}^{-1}}{d\theta}, 
\frac{df_{\theta}^{-1}}{d\theta}\right) + O(n^{-2}).
\end{eqnarray*}
To obtain expressions for the covariance involving
$K_{n,\underline{h}_n}(\theta)$, we
 apply similar techniques as those developed in \cite{p:dah-83}, Lemma
 6 together with cumulant arguments. This gives 
\begin{eqnarray*}
 \cov\left(\frac{d
  K_{n,\underline{h}_n}(\theta)}{d\theta}, \frac{d^{2}
  K_{n,\underline{h}_n}(\theta)}{d\theta^{2}}\right) &=&
                                           \frac{H_{2,n}}{H_{1,n}^{2}} V\left(\frac{df_{\theta}^{-1}}{d\theta}, 
\frac{d^{2}f_{\theta}^{-1}}{d\theta^{2}}\right) +O\left(\frac{H_{2,n}}{nH_{1,n}^2} \right)\\
\textrm{ and }\var\left(\frac{d
  K_{n,\underline{h}_n}(\theta)}{d\theta}\right)&=&  \frac{H_{2,n}}{H_{1,n}^{2}} V\left(\frac{df_{\theta}^{-1}}{d\theta}, 
\frac{df_{\theta}^{-1}}{d\theta}\right) + O\left(\frac{H_{2,n}}{nH_{1,n}^2}\right).
\end{eqnarray*}
These results yield expressions for $I_{1}$ and $I_{2}$ (we obtain
these below).

\vspace{1em}

\noindent \underline{Expression for $I_0$ and a bound for $I_{3}$}.
Using the results above we have 
\begin{itemize}
\item The Gaussian likelihood 
\begin{eqnarray}
I_{0} = I(\theta_n)^{-1}\left[B_{K,n}(\theta_n)+B_{G,n}(\theta_n)\right] + O(n^{-2})\label{eq:I0G}
\end{eqnarray}
\item The Whittle likelihood 
\begin{eqnarray*}
I_{0} = I(\theta_n)^{-1}B_{K,n}(\theta_n)+ O(n^{-2})
\end{eqnarray*}
\item The boundary corrected Whittle and hybrid Whittle likelihood
\begin{eqnarray*}
I_{0} = 0
\end{eqnarray*}
\end{itemize}
However, since for all the likelihoods
$\Ex[\frac{dL_{n}(\theta)}{d\theta}\rfloor_{\theta
  =\theta_n}]=O(n^{-1})$, this implies that for all the likelihoods
the term $I_{3}$ is 
\begin{eqnarray*}
I_{3} &=& 2^{-1}U(\theta)^{-3} \left\{\Ex\left(\frac{d
  L_{n}(\theta)}{d\theta}\right)\right\}^{2}\Ex\left(\frac{d^{3}
   L_{n}(\theta)}{d\theta^{3}}\right) = O(n^{-2}).
\end{eqnarray*}

\vspace{1em}

\noindent \underline{Expression for $I_1$ and $I_{2}$}. For the
Gaussian, Whittle, and boundary corrected Whittle likelihoods we
have 
\begin{eqnarray*}
I_{1}&=& n^{-1}I(\theta_{n})^{-2}V\left(\frac{df_{\theta}^{-1}}{d\theta}, 
\frac{d^{2}f_{\theta}^{-1}}{d\theta^{2}}\right) + O(n^{-3/2})  \nonumber \\
I_{2}&=&n^{-1}2^{-1}I(\theta_n)^{-3} V\left(\frac{df_{\theta}^{-1}}{d\theta}, 
\frac{df_{\theta}^{-1}}{d\theta}\right)J\left(\frac{d^{3}f_{\theta}^{-1}}{d\theta^{3}}\right)
         + O(n^{-3/2}).
\end{eqnarray*}
For the hybrid Whittle likelihood we obtain a similar expression 
\begin{eqnarray*}
I_{1}&=& \frac{H_{2,n}}{H_{1,n}^{2}} I(\theta_{n})^{-2}V\left(\frac{df_{\theta}^{-1}}{d\theta}, 
\frac{d^{2}f_{\theta}^{-1}}{d\theta^{2}}\right) + O\left(n^{-3/2}\right)  \nonumber \\
I_{2}&=& \frac{H_{2,n}}{H_{1,n}^{2}} 2^{-1}I(\theta_n)^{-3} V	\left(\frac{df_{\theta}^{-1}}{d\theta}, 
\frac{df_{\theta}^{-1}}{d\theta}\right)J\left(\frac{d^{3}f_{\theta}^{-1}}{d\theta^{3}}\right)
         +O\left(n^{-3/2}\right).
\end{eqnarray*}

\vspace{2mm}
\noindent \underline{A bound for $I_{4}$} We now show that $I_{4}$ has
a lower order term than the dominating terms $I_{0},I_{1}$ and $I_{2}$. We
recall that 
\begin{eqnarray*}
I_{4}&=& 2^{-1}U(\theta)^{-3} 
\cov\left(\left(\frac{d
  L_{n}(\theta)}{d\theta}\right)^{2},\frac{d^{3}
   L_{n}(\theta)}{d\theta^{3}}\right). 
\end{eqnarray*}
To bound the above we focus on $\cov\left(\left(\frac{d
  L_{n}(\theta)}{d\theta}\right)^{2},\frac{d^{3}
   L_{n}(\theta)}{d\theta^{3}}\right)$. By using indecomposable 
partitions we have 
\begin{eqnarray*}
\cov\left(\left(\frac{d
  L_{n}(\theta)}{d\theta}\right)^{2},\frac{d^{3}
   L_{n}(\theta)}{d\theta^{3}}\right) &=& 2\cov\left(\frac{d
  L_{n}(\theta)}{d\theta},\frac{d^{3}
   L_{n}(\theta)}{d\theta^{3}}\right)\Ex\left(\frac{d
  L_{n}(\theta)}{d\theta}\right)   \\
&&+ \cum\left(\frac{d
  L_{n}(\theta)}{d\theta}, \frac{d
  L_{n}(\theta)}{d\theta}\frac{d^{3}
   L_{n}(\theta)}{d\theta^{3}}\right) \\
&& +
\left[\Ex\left(\frac{d
  L_{n}(\theta)}{d\theta}\right)\right]^{2}\Ex\left(\frac{d^{3}
   L_{n}(\theta)}{d\theta^{3}}\right).
\end{eqnarray*}
We use (\ref{eq:KKKLLL}), (\ref{eq:Bounds11}) and
(\ref{eq:Bounds12}) to replace $L_{n}(\theta)$ with 
$K_{n}(\theta)$ or $K_{\underline{h},n}(\theta)$. Finally by using 
the expressions for cumulants of DFTs given in \cite{b:bri-01}, Theorem
4.3.2 we have that for the non-hybrid likelihoods 
\begin{eqnarray*}
I_{4} = O(n^{-2})
\end{eqnarray*}
and for the hybrid Whittle likelihood
\begin{eqnarray*}
I_{4} = O\left( \frac{H_{2,n}}{nH_{1,n}^{2}}\right).
\end{eqnarray*}
Thus, altogether for all the estimators we have that
\begin{eqnarray*}
(\widehat{\theta}_{n}-\theta_{n}) = I_{0} + I_{1}+I_{2} + O(n^{-2}),
\end{eqnarray*}
where for the Gaussian, Whittle and boundary corrected Whittle likelihoods 
\begin{eqnarray*}
I_{1} + I_{2} &=& n^{-1}\left[I(\theta_{n})^{-2}V\left(\frac{df_{\theta}^{-1}}{d\theta}, 
\frac{d^{2}f_{\theta}^{-1}}{d\theta^{2}}\right) +2^{-1}I(\theta_n)^{-3} V\left(\frac{df_{\theta}^{-1}}{d\theta}, 
\frac{df_{\theta}^{-1}}{d\theta}\right)J\left(\frac{d^{3}f_{\theta}^{-1}}{d\theta^{3}}\right)
                  \right] +O(n^{-3/2}) \\
 &=& n^{-1}G(\theta_n) + O(n^{-3/2})
\end{eqnarray*}
and for the hybrid Whittle likelihood
\begin{eqnarray*}
I_{1} + I_{2}  &=& \frac{H_{2,n}}{H_{1,n}^{2}} G(\theta_n)+O\left(n^{-3/2}\right).
\end{eqnarray*}
The terms for $I_{0}$ are given in (\ref{eq:I0G}).  This proves the
result. \hfill $\Box$

\subsection{Bias for estimators of multiple parameters}\label{sec:multi}

We now generalize the ideas above to multiple unknown
parameters. Suppose we fit the spectral density $f_{\theta}(\omega)$
to the time series $\{X_{t}\}$ where $\theta =
(\theta_{1},\ldots,\theta_{d})$ are the unknown parameters in $\Theta \subset \mathbb{R}^{d}$. 
$\mathcal{L}_{n}(\theta)$, $K_{n}(\theta)$,
$\widehat{W}_{p,n}(\theta)$ and $\widehat{H}_{p,n}(\theta)$ denote the
Gaussian likelihood, Whittle likelihood, boundary corrected Whittle and hybrid
Whittle likelihood defined in (\ref{eq:LIKE}).
Let $\widehat{\theta}_{n}^{(G)}$, $\widehat{\theta}_{n}^{(W)}$,
$\widehat{\theta}_{n}^{(W)}$ and $\widehat{\theta}_{n}^{(H)} $ be the
corresponding estimators defined in 
 (\ref{eq:thetaBiased}) and $\theta_n = (\theta_{1,n}, ..., \theta_{d,n})$ is the best fitting parameter defined as in (\ref{eq:thetan}). 
Then under Assumption \ref{assum:B} and
 \ref{assum:TS} we have the following asymptotic bias:
\begin{itemize}
\item The Gaussian likelihood (excluding the term $n^{-1} \log |\Gamma_{n}(\theta)|$)
\begin{eqnarray*}
\Ex[\widehat{\theta}_{j,n}^{(G)} - \theta_{j,n}] = \sum_{r=1}^{d}I^{(j,r)}\left[
B_{r,K,n}(\theta) + B_{r,G,n}(\theta)+ n^{-1}G_{r}(\theta)\right] +O\left(n^{-3/2}\right)
\end{eqnarray*}
\item The Whittle likelihood has bias
\begin{eqnarray*}
\Ex[\widehat{\theta}_{j,n}^{(K)} - \theta_{j,n}] = \sum_{r=1}^{d}I^{(j,r)}\left[
B_{r,K,n}(\theta) + n^{-1}G_{r}(\theta)\right]+O\left(n^{-3/2}\right).
\end{eqnarray*}
\item The boundary corrected Whittle likelihood has bias
\begin{eqnarray*}
\Ex[\widehat{\theta}_{j,n}^{(W)} - \theta_{j,n}] =
                                       n^{-1}\sum_{r=1}^{d}I^{(j,r)}G_{r}(\theta) +
O\left(p^{3}n^{-3/2}+(np^{K-1})^{-1}\right).
\end{eqnarray*}
\item The hybrid Whittle likelihood has bias 
\begin{eqnarray}
\Ex[\widehat{\theta}_{j,n}^{(H)} - \theta_{j,n}] &=& 
\frac{H_{2,n}}{H_{1,n}^{2}}\sum_{r=1}^{d}I^{(j,r)}G_{r}(\theta)+
O\left(p^{3}n^{-3/2}+(np^{K-1})^{-1}\right).
\label{eq:biasMultivar}
\end{eqnarray}
\end{itemize} Where $I^{(j,r)}$, $B_{r,G,n}(\cdot)$, $B_{r,K,n}(\cdot)$, and $G_r(\cdot)$ are
defined as in Section \ref{sec:sample}.

\vspace{1em}

\noindent PROOF. Let $L_{n}(\theta)$ be the criterion and $\widehat{\theta}_n = \arg\min
L_{n}(\theta)$ and $\theta_{n}$ the best fitting parameter. We use a 
similar technique used to prove Theorem \ref{theorem:bias}. The first
order expansion is 
\begin{eqnarray*}
\widehat{\theta}_n-\theta_n = U(\theta_n)^{-1}\nabla_{\theta}L_{n}(\theta_{n})
\end{eqnarray*}
where $U(\theta)$ is the $d\times d$ matrix
\begin{eqnarray*}
U(\theta) = -\Ex\left[\nabla_{\theta}^{2}L_{n}(\theta) \right].
\end{eqnarray*}
Thus entrywise we have
\begin{eqnarray*}
\widehat{\theta}_{r,n}-\theta_{r,n}=
  \sum_{s=1}^{d}U^{r,s}\frac{\partial L_{n}(\theta)}{\partial \theta_{s}}
\end{eqnarray*}
where $U^{(r,s)}$ denotes the $(r,s)$-entry of the $d\times d$ matrix
$U(\theta_n)^{-1}$. To obtain the ``bias'' we make a second order
expansion. For the simplicity, we omit the subscript $n$ from $\widehat{\theta}_{r,n}$ and $\theta_{r,n}$.
For $1\leq r \leq d$ we evaluate the partial derivative 
\begin{equation*}
\frac{\partial L_{n}(\theta)}{\partial \theta_{r}} + 
\sum_{s=1}^{d}(\widehat{\theta}_{s} - \theta_{s})
  \frac{\partial^{2}L_{n}(\theta)}{\partial \theta_{s}\partial \theta_{r}} +
 \frac{1}{2}\sum_{s_1,s_2=1}^{d}(\widehat{\theta}_{s_1} - \theta_{s_1}) (\widehat{\theta}_{s_2} - \theta_{s_2})
  \frac{\partial^{3}L_{n}(\theta)}{\partial
    \theta_{s_1}\partial\theta_{s_2}\partial \theta_r}\approx 0.
\end{equation*}
Taking expectation of the above gives 
\begin{eqnarray*}
&&\Ex\left[\frac{\partial L_{n}(\theta)}{\partial \theta_{r}}\right] + 
\sum_{s=1}^{d}
\Ex\left[(\widehat{\theta}_{s} -
   \theta_{s})\frac{\partial^{2}L_{n}(\theta)}{\partial
   \theta_{s}\partial \theta_{r}}\right] +
 \frac{1}{2}\sum_{s_1,s_2=1}^{d}\Ex\left[(\widehat{\theta}_{s_1} - \theta_{s_1}) (\widehat{\theta}_{s_2} - \theta_{s_2})
  \frac{\partial^{3}L_{n}(\theta)}{\partial
   \theta_{s_1}\partial\theta_{s_2}\partial \theta_r}\right]\approx 0.
\end{eqnarray*}
We now replace the product of random variables with their covariances 
\begin{eqnarray*}
&&\Ex\left[\frac{\partial L_{n}(\theta)}{\partial \theta_{r}}\right] + 
\sum_{s=1}^{d}
\Ex [\widehat{\theta}_{s} -
   \theta_{s} ] \Ex\left[\frac{\partial^{2}L_{n}(\theta)}{\partial
   \theta_{s}\partial \theta_{r}}\right] +
\sum_{s=1}^{d}
\cov\left[\widehat{\theta}_{s} -
   \theta_{s}, \frac{\partial^{2}L_{n}(\theta)}{\partial
   \theta_{s}\partial \theta_{r}}\right] \\
&&\quad+
 \frac{1}{2}\sum_{s_1,s_2=1}^{d}\cov\left(\widehat{\theta}_{s_1} - \theta_{s_1},\widehat{\theta}_{s_2} - \theta_{s_2}\right)
\Ex\left[\frac{\partial^{3}L_{n}(\theta)}{\partial
   \theta_{s_1}\partial\theta_{s_2}\partial \theta_r}\right]\\
&&\quad +\frac{1}{2}\sum_{s_1,s_2=1}^{d}\Ex[\widehat{\theta}_{s_1} - \theta_{s_1}] \Ex [\widehat{\theta}_{s_2} - \theta_{s_2}]
\Ex\left[\frac{\partial^{3}L_{n}(\theta)}{\partial
   \theta_{s_1}\partial\theta_{s_2}\partial \theta_r}\right]\\
&&\quad+ \frac{1}{2}\sum_{s_1,s_2=1}^{d}\cov\left[(\widehat{\theta}_{s_1} - \theta_{s_1})(\widehat{\theta}_{s_2} - \theta_{s_2}),
  \frac{\partial^{3}L_{n}(\theta)}{\partial
   \theta_{s_1}\partial\theta_{s_2}\partial \theta_r}\right] \approx 0.
\end{eqnarray*}
With the exception of $\Ex[\widehat{\theta}_{s} - \theta_{s}]$, we 
replace $\widehat{\theta}_{s} - \theta_{s}$ in
the above with their first order expansions $\sum_{j=1}^{d}U^{(s,j)}\frac{\partial L_{n}(\theta)}{\partial \theta_{j}}$. This gives 
\begin{eqnarray*}
&&\Ex\left[\frac{\partial L_{n}(\theta)}{\partial \theta_{r}}\right] - 
\sum_{s=1}^{d}
\Ex [\widehat{\theta}_{s} -
   \theta_{s}] U_{s,r} +
\sum_{s_1,s_2=1}^{d}U^{(s_1,s_2)}
\cov\left[\frac{\partial L_{n}(\theta)}{
   \partial \theta_{s_2}},\frac{\partial^{2}L_{n}(\theta)}{\partial
   \theta_{s_1}\partial \theta_{r}}\right] \\
&&\quad + \frac{1}{2}\sum_{s_1,s_2,s_3,s_4=1}^{d}U^{(s_1,s_3)}U^{(s_2,s_4)}
\cov\left(\frac{\partial L_{n}(\theta)}{\partial \theta_{s_3}}, \frac{\partial L_{n}(\theta)}{\partial \theta_{s_4}}\right)
\Ex\left[\frac{\partial^{3}L_{n}(\theta)}{\partial
   \theta_{s_1}\partial\theta_{s_2}\partial \theta_r}\right]\\
&&\quad +\frac{1}{2}\sum_{s_1,s_2,s_3,s_4=1}^{d}U^{(s_1,s_3)}U^{(s_2,s_4)}
\Ex\left[\frac{\partial L_{n}(\theta)}{\partial \theta_{s_3}}\right]\Ex\left[\frac{\partial L_{n}(\theta)}{\partial \theta_{s_4}}\right]
\Ex\left[\frac{\partial^{3}L_{n}(\theta)}{\partial
   \theta_{s_1}\partial\theta_{s_2}\partial \theta_r}\right]\\
&&\quad+ \frac{1}{2}\sum_{s_1,s_2,s_3,s_4=1}^{d}U^{(s_1,s_3)}U^{(s_2,s_4)}
\cov\left[\frac{\partial L_{n}(\theta)}{\partial \theta_{s_3}}\frac{\partial L_{n}(\theta)}{\partial \theta_{s_4}},
  \frac{\partial^{3}L_{n}(\theta)}{\partial
   \theta_{s_1}\partial\theta_{s_2}\partial \theta_r}\right] \approx 0,
\end{eqnarray*} where $U_{s,r}$ denotes the $(s,r)$-entry of the $d\times d$ matrix
$U(\theta_n)$

Now we consider concrete examples of likelihoods. Using the same
arguments as those used in the proof of Theorem  \ref{theorem:bias} we
have the last two terms of the above are of order $O(n^{-2})$ or
$O(H_{2,n}/(nH_{1,n}^{2}))$ depending on the likelihood used. This
implies that 
\begin{eqnarray*}
&&\Ex\left[\frac{\partial L_{n}(\theta)}{\partial \theta_{r}}\right] - 
\sum_{s=1}^{d}
\Ex[\widehat{\theta}_{s} -
   \theta_{s}] U_{s,r} +
\sum_{s_1,s_2=1}^{d}U^{(s_1,s_2)}
\cov\left[\frac{\partial L_{n}(\theta)}{
   \partial \theta_{s_2}},\frac{\partial^{2}L_{n}(\theta)}{\partial
   \theta_{s_1}\partial \theta_{r}}\right] \\
&& \quad + \frac{1}{2}\sum_{s_1,s_2,s_3,s_4=1}^{d}U^{(s_1,s_3)}U^{(s_2,s_4)}
\cov\left(\frac{\partial L_{n}(\theta)}{\partial \theta_{s_3}}, \frac{\partial L_{n}(\theta)}{\partial \theta_{s_4}}\right)
\Ex\left[\frac{\partial^{3}L_{n}(\theta)}{\partial
   \theta_{s_1}\partial\theta_{s_2}\partial \theta_r}\right]\approx 0.
\end{eqnarray*}
Let 
\begin{eqnarray*}
J\left(g \right) &=& 
\frac{1}{2\pi}\int_{0}^{2\pi} g(\omega)f(\omega)d\omega \\
I(\theta) &=& -\frac{1}{2\pi}\int_{0}^{2\pi}[\nabla_{\theta}^{2}f_{\theta}(\omega)^{-1}]f(\omega)d\omega
\end{eqnarray*}
and  $I_{s,r}$ (and $I^{(s,r)}$) corresponds to the $(s,r)$-th element
of $I(\theta_n)$ (and $I^{-1}(\theta_n)$).
So far, we have no specified the likelihood $L_{n}(\theta)$. But to
write a second order expansion for all four likelihoods we set 
$H_{2,n}/H_{1,n}^{2} = n^{-1}$ for the Gaussian, Whittle, and boundary
corrected Whittle likelihood and using the notation a similar proof
to Theorem  \ref{theorem:bias} we have 
\begin{eqnarray*}
&&\Ex\left[\frac{\partial L_{n}(\theta)}{\partial \theta_{r}}\right] -
\sum_{s=1}^{d}I_{s,r}
\Ex[\widehat{\theta}_{s} -
   \theta_{s}] +
 \frac{H_{2,n}}{H_{1,n}^{2}}\sum_{s_1,s_2=1}^{d}I^{(s_1,s_2)}V\left(\frac{\partial f_\theta^{-1}}{
   \partial \theta_{s_2}}, \frac{\partial^{2}f_\theta^{-1}}{\partial
   \theta_{s_1}\partial \theta_{r}}\right)\\
&& + \frac{H_{2,n}}{2H_{1,n}^{2}}\sum_{s_1,s_2,s_3,s_4=1}^{d}I^{(s_1,s_3)}I^{(s_2,s_4)}
V\left(\frac{\partial f_{\theta}^{-1}}{\partial \theta_{s_3}}, \frac{\partial f_{\theta}^{-1}}{\partial \theta_{s_4}}\right)
J\left(\frac{\partial^{3}f_{\theta}^{-1}}{\partial  \theta_{s_1}\partial\theta_{s_2}\partial \theta_r} \right)
\approx 0.
\end{eqnarray*}
Thus 
\begin{eqnarray*}
\sum_{s=1}^{d}I_{s,r}
\Ex[\widehat{\theta}_{s} -
   \theta_{s}] 
&\approx &\Ex\left[\frac{\partial L_{n}(\theta)}{\partial \theta_{r}}\right] +
 \frac{H_{2,n}}{H_{1,n}^{2}}\sum_{s_1,s_2=1}^{d}I^{(s_1,s_2)}V\left(\frac{\partial f_\theta^{-1}}{
   \partial \theta_{s_2}}, \frac{\partial^{2}f_\theta^{-1}}{\partial
   \theta_{s_1}\partial \theta_{r}}\right)\\
&& + \frac{H_{2,n}}{2H_{1,n}^{2}}\sum_{s_1,s_2,s_3,s_4=1}^{d}I^{(s_1,s_3)}I^{(s_2,s_4)}
V\left(\frac{\partial f_{\theta}^{-1}}{\partial \theta_{s_3}}, \frac{\partial f_{\theta}^{-1}}{\partial \theta_{s_4}}\right)
J\left(\frac{\partial^{3}f_{\theta}^{-1}}{\partial  \theta_{s_1}\partial\theta_{s_2}\partial \theta_r} \right).
\end{eqnarray*}
In the final stage, to extract $\Ex[\widehat{\theta}_{s} -
   \theta_{s}]$ from the above we define the $d$-dimensional column vector
   $\underline{D}^{\prime} = (D_{1},\ldots,D_{d})$, where 
$D_{r} = \sum_{s=1}^{d}I_{s,r} \Ex[\widehat{\theta}_{s} -
   \theta_{s}] = [I(\theta_n)(\widehat{\theta}_n -
   \theta_n)]_{r}$. Substituting this in the above gives
\begin{eqnarray*}
D_{r} &\approx& \Ex\left[\frac{\partial L_{n}(\theta)}{\partial \theta_{r}}\right] +
 \frac{H_{2,n}}{H_{1,n}^{2}}\sum_{s_1,s_2=1}^{d}I^{(s_1,s_2)}V\left(\frac{\partial f_\theta(\omega)^{-1}}{
   \partial \theta_{s_2}}, \frac{\partial^{2}f_\theta(\omega)^{-1}}{\partial
   \theta_{s_1}\partial \theta_{r}}\right)\\
&& + \frac{H_{2,n}}{2H_{1,n}^{2}}\sum_{s_1,s_2,s_3,s_4=1}^{d}I^{(s_1,s_3)}I^{(s_2,s_4)}
V\left(\frac{\partial f_{\theta}^{-1}}{\partial \theta_{s_3}}, \frac{\partial f_{\theta}^{-1}}{\partial \theta_{s_4}}\right)
J\left(\frac{\partial^{3}f_{\theta}^{-1}}{\partial  \theta_{s_1}\partial\theta_{s_2}\partial \theta_r} \right).
\end{eqnarray*}
Using that $\Ex[\widehat{\theta}_n - \theta_n] \approx
I(\theta_n)^{-1}\underline{D}$ and substituting this into the above
gives the bias for $\widehat{\theta}_j$
\begin{eqnarray}
\Ex[\widehat{\theta}_j - \theta_j] &\approx& \sum_{r=1}^{d}I^{(j,r)}\bigg[
\Ex\left[\frac{\partial L_{n}(\theta)}{\partial \theta_{r}}\right] +
 \frac{H_{2,n}}{H_{1,n}^{2}}\sum_{s_1,s_2=1}^{d}I^{(s_1,s_2)}V\left(\frac{\partial f_\theta(\omega)^{-1}}{
   \partial \theta_{s_2}}, \frac{\partial^{2}f_\theta(\omega)^{-1}}{\partial
   \theta_{s_1}\partial \theta_{r}}\right) \nonumber \\
&& + \frac{H_{2,n}}{2H_{1,n}^{2}}\sum_{s_1,s_2,s_3,s_4=1}^{d}I^{(s_1,s_3)}I^{(s_2,s_4)}
V\left(\frac{\partial f_{\theta}^{-1}}{\partial \theta_{s_3}}, \frac{\partial f_{\theta}^{-1}}{\partial \theta_{s_4}}\right)
J\left(\frac{\partial^{3}f_{\theta}^{-1}}{\partial  \theta_{s_1}\partial\theta_{s_2}\partial \theta_r} \right)
\bigg]. \label{eq:bias_expansion}
\end{eqnarray}
The above is a general result. We now obtain the bias for the different
criteria. Let 
\begin{eqnarray*}
B_{r,G,n}(\theta) &=&  
\R \frac{2}{n}\sum_{s,t=1}^{n}c(s-t)\frac{1}{n}\sum_{k=1}^{n}
  e^{-is\omega_{k,n}}\frac{\partial }{\partial \theta_r}\left[\overline{\phi(\omega_{k,n};f_\theta)}\phi_{t}^{\infty}(\omega_{k,n};f_\theta)\right] \\
B_{r,K,n}(\theta) &=& \frac{1}{n}\sum_{k=1}^{n}
f_{n}(\omega_{k,n})\frac{\partial f_\theta(\omega_{k,n})^{-1}}{\partial \theta_{r}} \\
\text{and} \qquad G_{r}(\theta) &=& \sum_{s_1,s_2=1}^{d}I^{(s_1,s_2)}V\left(\frac{\partial f_\theta(\omega)^{-1}}{
   \partial \theta_{s_2}}, \frac{\partial^{2}f_\theta(\omega)^{-1}}{\partial
   \theta_{s_1}\partial \theta_{r}}\right)\\
&& + \frac{1}{2}\sum_{s_1,s_2,s_3,s_4=1}^{d}I^{(s_1,s_3)}I^{(s_2,s_4)}
V\left(\frac{\partial f_{\theta}^{-1}}{\partial \theta_{s_3}}, \frac{\partial f_{\theta}^{-1}}{\partial \theta_{s_4}}\right)
J\left(\frac{\partial^{3}f_{\theta}^{-1}}{\partial  \theta_{s_1}\partial\theta_{s_2}\partial \theta_r} \right).
\end{eqnarray*}
Then, using similar technique from the univariate case, we can show
\begin{itemize}
\item The Gaussian likelihood:
$\Ex\left[\partial \mathcal{L}_{n}(\theta)/\partial \theta_{r}\right] = B_{r,G,n}(\theta) + B_{r,K,n}(\theta)$.

\item The Whittle likelihood:
$\Ex\left[\partial K_{n}(\theta)/\partial \theta_{r}\right] = B_{r,K,n}(\theta)$

\item The boundary corrected Whittle and hybrid Whittle likelihood:
\\$\Ex\left[\partial W_{n}(\theta)/\partial \theta_{r}\right] = \Ex\left[\partial H_{n}(\theta)/\partial \theta_{r}\right] =0$.
\end{itemize}
Substituting the above into (\ref{eq:bias_expansion}) gives the four
difference biases in (\ref{eq:biasMultivar}). Thus we have proved the result.
\hfill $\Box$

%% file: appendix_simulations.tex
\section{Additional Simulations}\label{sec:sim_appendix}

\subsection{Table of results for the AR$(1)$ and MA$(1)$ for a Gaussian time series}\label{sec:AR-MA-Gaussian}

\begin{landscape}
\begin{table}[ht]
    \centering
\scriptsize
  \begin{tabular}{c|rrrrr|rrrrr}
\multirow{2}{*}{\textit{Likelihoods}} & \multicolumn{10}{c}{$\theta$} \\
\cline{2-11}	 
 & 0.1 & 0.3 & 0.5 & 0.7 & 0.9    & 0.1 & 0.3 & 0.5 & 0.7 & 0.9 \\ \hline \hline

& \multicolumn{5}{c}{ \textbf{AR(1)},	 $\{e_{t}\}\sim \mathcal{N}(0,1)$, $n=20$}
& \multicolumn{5}{c}{ \textbf{MA(1)}, $\{e_{t}\}\sim \mathcal{N}(0,1)$, $n=20$} \\ \hline
 Gaussian  &  \color{blue}{-$0.012$}{\scriptsize (0.22)} & \color{blue}{-$0.028$}{\scriptsize (0.21)} & \color{red}{-$0.043$}{\scriptsize (0.19)} & \color{red}{-$0.066$}{\scriptsize (0.18)} & \color{red}{-$0.072$}{\scriptsize (0.14)} 
& \color{red}{$0.010$}{\scriptsize (0.28)} & $0.016${\scriptsize (0.28)} & $0.025${\scriptsize (0.24)} & \color{red}{$0.012$}{\scriptsize (0.21)} & \color{red}{$0.029$}{\scriptsize (0.17)} \\ 

 Whittle  & \color{red}{-$0.015$}{\scriptsize (0.21)} & \color{red}{-$0.041$}{\scriptsize (0.20)} & \color{blue}{-$0.063$}{\scriptsize (0.19)} & -$0.095${\scriptsize (0.18)} & -$0.124${\scriptsize (0.15)} &
 \color{blue}{$0.005$}{\scriptsize (0.29)} & \color{blue}{$0.002$}{\scriptsize (0.28)} & \color{red}{-$0.004$}{\scriptsize (0.24)} & -$0.052${\scriptsize (0.23)} & -$0.152${\scriptsize (0.21)} \\ 

 \color{blue}{Boundary}  & -$0.015${\scriptsize (0.22)} & -$0.037${\scriptsize (0.21)} & -$0.054${\scriptsize (0.19)} & -$0.079${\scriptsize (0.18)} & -$0.103${\scriptsize (0.14)} 
& $0.007${\scriptsize (0.30)} & $0.009${\scriptsize (0.29)} & $0.009${\scriptsize (0.24)} &-$0.022${\scriptsize (0.24)} & -$0.111${\scriptsize (0.20)} \\

 \color{blue}{Hybrid}  & -$0.012${\scriptsize (0.22)} & -$0.030${\scriptsize (0.21)} & -$0.049${\scriptsize (0.19)} & \color{blue}{-$0.072$}{\scriptsize (0.18)} & \color{blue}{-$0.095$}{\scriptsize (0.14)} 
& $0.011${\scriptsize (0.30)} & $0.021${\scriptsize (0.29)} & $0.026${\scriptsize (0.25)} & -$0.007${\scriptsize (0.22)} & \color{blue}{-$0.074$}{\scriptsize (0.17)} \\

 Tapered  & -$0.014${\scriptsize (0.22)} & -$0.036${\scriptsize (0.21)} & -$0.063${\scriptsize (0.19)} & -$0.090${\scriptsize (0.18)} & -$0.117${\scriptsize (0.14)} 
& $0.004${\scriptsize (0.29)} & \color{red}{$0.004$}{\scriptsize (0.28)} & \color{blue}{-$0.006$}{\scriptsize (0.24)} & \color{blue}{-$0.043$}{\scriptsize (0.21)} & -$0.122${\scriptsize (0.18)} \\

 Debiased  & -$0.013${\scriptsize (0.22)} & -$0.033${\scriptsize (0.21)} & -$0.049${\scriptsize (0.19)} & -$0.069${\scriptsize (0.19)} & -$0.085${\scriptsize (0.16)} 
& $0.005${\scriptsize (0.29)} & $0.013${\scriptsize (0.28)} & $0.021${\scriptsize (0.25)} & -$0.005${\scriptsize (0.24)} & -$0.088${\scriptsize (0.21)} \\ \hline \hline

& \multicolumn{5}{c}{ \textbf{AR(1)}, $\{e_{t}\}\sim \mathcal{N}(0,1)$, $n=50$}
& \multicolumn{5}{c}{ \textbf{MA(1)}, $\{e_{t}\}\sim \mathcal{N}(0,1)$, $n=50$} \\ \hline

Gaussian  & \color{blue}{-$0.006$}{\scriptsize (0.14)} &\color{red}{-$0.011$}{\scriptsize (0.14)} & \color{red}{-$0.013$}{\scriptsize (0.12)} & \color{red}{-$0.033$}{\scriptsize (0.11)} & \color{red}{-$0.030	$}{\scriptsize (0.07)} 
& -$0.002${\scriptsize (0.16)} & \color{blue}{$0.008$}{\scriptsize (0.15)} & \color{blue}{$0.017$}{\scriptsize (0.14)} & $0.018${\scriptsize (0.12)} & \color{blue}{$0.014$}{\scriptsize (0.08)} \\ 

 Whittle  & \color{red}{-$0.008$}{\scriptsize (0.14)} & \color{blue}{-$0.016$}{\scriptsize(0.14)} & \color{blue}{-$0.023$}{\scriptsize (0.12)} & -$0.045${\scriptsize (0.11)} & -$0.049${\scriptsize (0.08)} 
& \color{red}{-$0.004$}{\scriptsize (0.15)} & \color{red}{$0.001$}{\scriptsize (0.15)} & \color{red}{$0.001$}{\scriptsize (0.14)} & -$0.020${\scriptsize (0.13)} & -$0.067${\scriptsize (0.11)} \\ 

 \color{blue}{Boundary}  & \color{blue}{-$0.007$}{\scriptsize (0.14)} &-$0.012${\scriptsize (0.14)} & -$0.015${\scriptsize (0.12)} &  \color{blue}{-$0.034$}{\scriptsize (0.11)} & \color{blue}{-$0.036$}{\scriptsize (0.07)} 
& -$0.003${\scriptsize (0.16)} & $0.006${\scriptsize (0.16)} & $0.013${\scriptsize (0.14)} & $0.005${\scriptsize (0.13)} & -$0.026${\scriptsize (0.09)} \\

 \color{blue}{Hybrid}  & -$0.005${\scriptsize (0.14)} & -$0.011${\scriptsize (0.14)} &  -$0.015${\scriptsize (0.13)} & -$0.033${\scriptsize(0.11)} & -$0.035${\scriptsize (0.07)} 
&-$0.001${\scriptsize (0.16)} & $0.010${\scriptsize (0.16)} & $0.015${\scriptsize (0.14)} & \color{blue}{$0.014$}{\scriptsize (0.12)} & \color{red}{-$0.010$}{\scriptsize (0.07)} \\

 Tapered  & -$0.005${\scriptsize (0.14)} & -$0.013${\scriptsize(0.14)} &-$0.018${\scriptsize (0.13)} & -$0.038${\scriptsize (0.11)} & -$0.039${\scriptsize (0.08)} 
& $0${\scriptsize (0.16)} & $0.008${\scriptsize (0.16)} &$0.010${\scriptsize (0.14)} & \color{red}{$0.003$}{\scriptsize (0.12)} & -$0.023${\scriptsize (0.08)} \\
 
Debiased  & -$0.006${\scriptsize (0.14)} &-$0.011${\scriptsize (0.14)} & -$0.015${\scriptsize (0.12)} & -$0.035${\scriptsize (0.11)} & -$0.032${\scriptsize (0.08)} 
& \color{blue}{-$0.002$}{\scriptsize (0.16)} & $0.009${\scriptsize (0.16)} & $0.019${\scriptsize (0.15)} & $0.017${\scriptsize (0.15)} & -$0.011${\scriptsize (0.11)} \\ \hline \hline

& \multicolumn{5}{c}{ \textbf{AR(1)}, $\{e_{t}\}\sim \mathcal{N}(0,1)$, $n=300$}
& \multicolumn{5}{c}{ \textbf{MA(1)}, $\{e_{t}\}\sim \mathcal{N}(0,1)$, $n=300$} \\ \hline

Gaussian  & \color{blue}{$0$}{\scriptsize (0.06)} &\color{red}{-$0.002$}{\scriptsize (0.06)} & \color{blue}{-$0.001$}{\scriptsize (0.05)} & \color{red}{-$0.004$}{\scriptsize (0.04)} & \color{red}{-$0.005$}{\scriptsize (0.03)} 
& $0.002${\scriptsize (0.06)} & \color{red}{$0$}{\scriptsize (0.06)} & \color{blue}{$0.003$}{\scriptsize (0.05)} & \color{red}{$0$}{\scriptsize (0.04)} & \color{blue}{$0.004$}{\scriptsize (0.03)} \\ 

 Whittle  & \color{red}{$0$}{\scriptsize (0.06)} & \color{blue}{-$0.003$}{\scriptsize(0.06)} & \color{red}{-$0.003$}{\scriptsize (0.05)} & -$0.007${\scriptsize (0.04)} & -$0.008${\scriptsize (0.03)} 
& \color{red}{$0.001$}{\scriptsize (0.06)} & \color{blue}{-$0.001$}{\scriptsize (0.06)} & \color{red}{$0$}{\scriptsize (0.05)} & -$0.007${\scriptsize (0.04)} & -$0.020${\scriptsize (0.04)} \\ 

 \color{blue}{Boundary}  & \color{blue}{$0$}{\scriptsize (0.06)} &-$0.002${\scriptsize (0.06)} & \color{blue}{-$0.001$}{\scriptsize (0.05)} &  \color{blue}{-$0.004$}{\scriptsize (0.04)} & \color{blue}{-$0.006$}{\scriptsize (0.03)} 
& \color{blue}{$0.002$}{\scriptsize (0.06)} & $0${\scriptsize (0.06)} & $0.003${\scriptsize (0.05)} & \color{blue}{$0$}{\scriptsize (0.04)} & -$0.002${\scriptsize (0.03)} \\

 \color{blue}{Hybrid}  & $0${\scriptsize (0.06)} & -$0.002${\scriptsize (0.06)} &  -$0.001${\scriptsize (0.05)} & -$0.005${\scriptsize(0.04)} & -$0.006${\scriptsize (0.03)} 
&$0.002${\scriptsize (0.06)} & $0${\scriptsize (0.06)} & $0.004${\scriptsize (0.05)} & $0.001${\scriptsize (0.05)} & \color{red}{$0.003$}{\scriptsize (0.03)} \\

Tapered  & $0${\scriptsize (0.06)} & -$0.002${\scriptsize (0.06)} &  -$0.001${\scriptsize (0.05)} & -$0.005${\scriptsize(0.05)} & -$0.006${\scriptsize (0.03)} 
&$0.002${\scriptsize (0.06)} & $0${\scriptsize (0.06)} & $0.004${\scriptsize (0.05)} & $0.001${\scriptsize (0.05)} & $0.003${\scriptsize (0.03)} \\
 
Debiased  & \color{blue}{$0$}{\scriptsize (0.06)} & -$0.002${\scriptsize (0.06)} &  \color{blue}{-$0.001$}{\scriptsize (0.05)} & -$0.004${\scriptsize(0.04)} & -$0.006${\scriptsize (0.03)} 
&\color{blue}{$0.002$}{\scriptsize (0.06)} & $0${\scriptsize (0.06)} & $0.003${\scriptsize (0.05)} & $0${\scriptsize (0.05)} & $0.009${\scriptsize (0.05)} \\
\end{tabular} 

\caption{\textit{Bias and the standard deviation (in the parentheses) of six different quasi-likelihoods for an AR(1) (left) and MA(1) (right) model for
the standard normal innovations. Length of the time series $n=20, 50$, and $300$. We use {\color{red}red} to denote the smallest RMSE and
{\color{blue}blue} to denote the second smallest RMSE.}}
\label{tab:AR}
\end{table}
\end{landscape}

\subsection{Figures and Table of results for the AR$(1)$ and MA$(1)$  for a non-Gaussian time series}\label{sec:AR-MA-chi}

In this section, we provide figures and table of the results in Section \ref{sec:specified}
when the innovations follow a standardized chi-squared distribution
two degrees of freedom, i.e.  
$\varepsilon_t \sim (\chi^2(2)-2)/2$ 
(this time the asymptotic bias will contain the fourth order cumulant term). 
The results are very similar to the Gaussian innovations.

\begin{figure}[ht]
\begin{center}

\textbf{AR$(1)$ model}

\vspace{1em}

\includegraphics[scale=0.35,page=1]{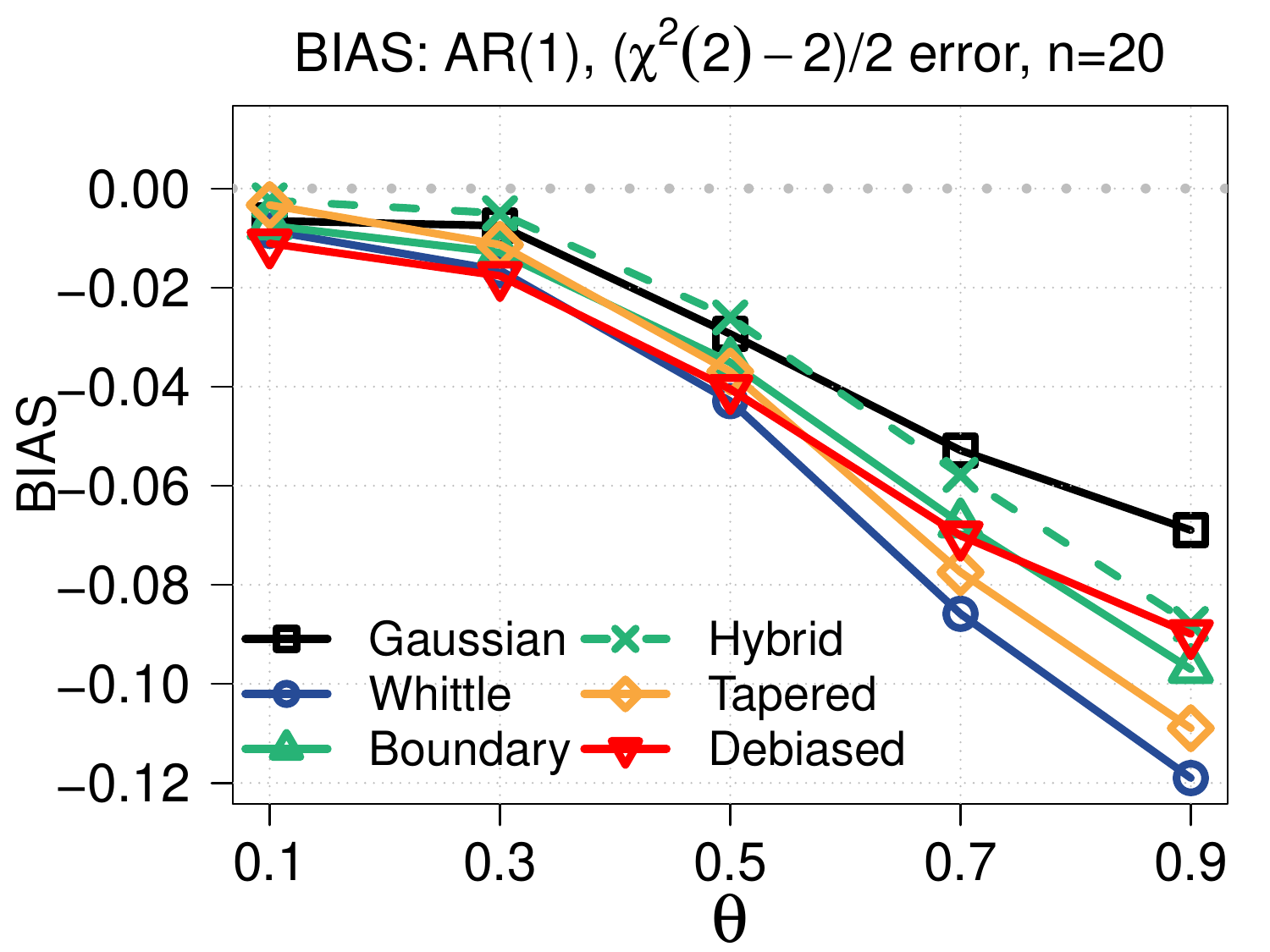}
\includegraphics[scale=0.35,page=3]{plot/AR_chisq2.pdf}
\includegraphics[scale=0.35,page=5]{plot/AR_chisq2.pdf}

\includegraphics[scale=0.35,page=2]{plot/AR_chisq2.pdf}
\includegraphics[scale=0.35,page=4]{plot/AR_chisq2.pdf}
\includegraphics[scale=0.35,page=6]{plot/AR_chisq2.pdf}

\vspace{1em}

\textbf{MA$(1)$ model}

\vspace{1em}

\includegraphics[scale=0.35,page=1]{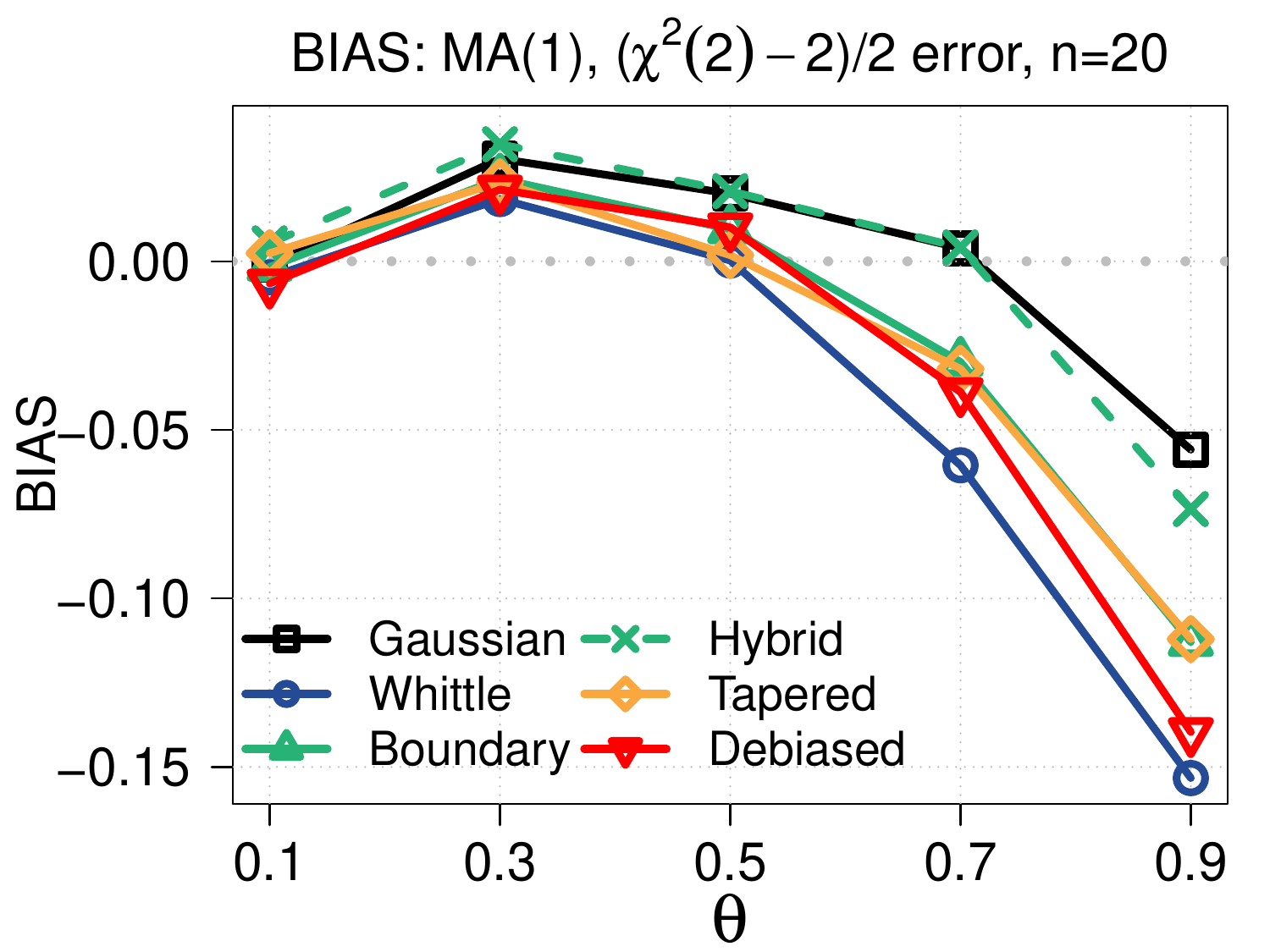}
\includegraphics[scale=0.35,page=3]{plot/MA_chisq2.pdf}
\includegraphics[scale=0.35,page=5]{plot/MA_chisq2.pdf}

\includegraphics[scale=0.35,page=2]{plot/MA_chisq2.pdf}
\includegraphics[scale=0.35,page=4]{plot/MA_chisq2.pdf}
\includegraphics[scale=0.35,page=6]{plot/MA_chisq2.pdf}

\caption{\textit{Bias (first row) and the RMSE (second row) of the parameter estimates for the AR(1) and MA(1) models
where the innovations follow the standardized chi-squared distribution with 2 degrees of freedom.  Length of the time series $n=20$(left), $50$(middle), and $300$(right). } }
\label{fig:ar1ma1.chisq} 
\end{center}
\end{figure}

\begin{landscape}
\begin{table}[ht]
    \centering
\scriptsize
  \begin{tabular}{c|rrrrr|rrrrr}
\multirow{2}{*}{\textit{Likelihoods}} & \multicolumn{10}{c}{$\theta$} \\
\cline{2-11}	 
 & 0.1 & 0.3 & 0.5 & 0.7 & 0.9    & 0.1 & 0.3 & 0.5 & 0.7 & 0.9 \\ \hline \hline

& \multicolumn{5}{c}{ \textbf{AR(1)},	 $\{e_{t}\}\sim (\chi^2(2)-2)/2$, $n=20$}
& \multicolumn{5}{c}{ \textbf{MA(1)}, $\{e_{t}\}\sim (\chi^2(2)-2)/2$, $n=20$} \\ \hline
 Gaussian  &  \color{blue}{-$0.007$}{\scriptsize (0.21)} & \color{blue}{-$0.007$}{\scriptsize (0.20)} & \color{red}{-$0.029$}{\scriptsize (0.19)} & \color{red}{-$0.053$}{\scriptsize (0.17)} & \color{red}{-$0.069$}{\scriptsize (0.13)} 
& -$0.001${\scriptsize (0.28)} & \color{blue}{$0.030$}{\scriptsize (0.25)} & \color{red}{$0.020$}{\scriptsize (0.23)} & \color{red}{$0.004$}{\scriptsize (0.20)} & \color{red}{$0.056$}{\scriptsize (0.17)} \\ 

 Whittle  & -$0.009${\scriptsize (0.21)} & -$0.016${\scriptsize (0.20)} & -$0.043${\scriptsize (0.20)} & -$0.086${\scriptsize (0.18)} & -$0.119${\scriptsize (0.14)} &
 \color{blue}{-$0.005$}{\scriptsize (0.27)} & $0.018${\scriptsize (0.26)} & $0${\scriptsize (0.24)} & -$0.061${\scriptsize (0.22)} & -$0.153${\scriptsize (0.21)} \\ 

 \color{blue}{Boundary}  & -$0.007${\scriptsize (0.22)} & -$0.013${\scriptsize (0.20)} & -$0.035${\scriptsize (0.20)} & -$0.068${\scriptsize (0.18)} & -$0.097${\scriptsize (0.13)} 
& -$0.002${\scriptsize (0.28)} & $0.024${\scriptsize (0.26)} & $0.009${\scriptsize (0.25)} &-$0.030${\scriptsize (0.23)} & -$0.113${\scriptsize (0.20)} \\

 \color{blue}{Hybrid}  & -$0.002${\scriptsize (0.22)} & -$0.005${\scriptsize (0.20)} & -$0.026${\scriptsize (0.20)} & \color{blue}{-$0.058$}{\scriptsize (0.18)} & \color{blue}{-$0.088$}{\scriptsize (0.13)} 
& $0.005${\scriptsize (0.29)} & $0.035${\scriptsize (0.26)} & $0.021${\scriptsize (0.24)} & \color{blue}{$0.004$}{\scriptsize (0.20)} & \color{blue}{-$0.074$}{\scriptsize (0.17)} \\

 Tapered  & -$0.003${\scriptsize (0.21)} & -$0.011${\scriptsize (0.20)} & -$0.037${\scriptsize (0.20)} & -$0.077${\scriptsize (0.18)} & -$0.109${\scriptsize (0.13)} 
& $0.002${\scriptsize (0.28)} & $0.023${\scriptsize (0.25)} & \color{blue}{$0.002$}{\scriptsize (0.23)} & -$0.032${\scriptsize (0.21)} & -$0.112${\scriptsize (0.18)} \\

 Debiased  & \color{red}{-$0.011$}{\scriptsize (0.21)} & \color{red}{-$0.018$}{\scriptsize (0.19)} & \color{blue}{-$0.040$}{\scriptsize (0.20)} & -$0.070${\scriptsize (0.19)} & -$0.090${\scriptsize (0.15)} 
& \color{red}{-$0.007$}{\scriptsize (0.27)} & \color{red}{$0.021$}{\scriptsize (0.25)} & $0.010${\scriptsize (0.24)} & -$0.039${\scriptsize (0.24)} & -$0.140${\scriptsize (0.23)} \\ \hline \hline

& \multicolumn{5}{c}{ \textbf{AR(1)}, $\{e_{t}\}\sim (\chi^2(2)-2)/2$, $n=50$}
& \multicolumn{5}{c}{ \textbf{MA(1)}, $\{e_{t}\}\sim (\chi^2(2)-2)/2$, $n=50$} \\ \hline

Gaussian  & $0.004${\scriptsize (0.13)} &-$0.011${\scriptsize (0.13)} & \color{red}{-$0.012$}{\scriptsize (0.11)} & \color{red}{-$0.031$}{\scriptsize (0.10)} & \color{red}{-$0.029$}{\scriptsize (0.07)} 
& $0.009${\scriptsize (0.15)} & \color{blue}{$0.003$}{\scriptsize (0.15)} & \color{red}{$0.017$}{\scriptsize (0.13)} & $0.014${\scriptsize (0.12)} & \color{red}{$0.010$}{\scriptsize (0.08)} \\ 

 Whittle  & \color{red}{$0.001$}{\scriptsize (0.13)} & \color{blue}{-$0.016$}{\scriptsize(0.13)} & -$0.019${\scriptsize (0.12)} & -$0.044${\scriptsize (0.10)} & -$0.049${\scriptsize (0.07)} 
& \color{red}{$0.005$}{\scriptsize (0.14)} & \color{red}{-$0.004$}{\scriptsize (0.14)} & $0.004${\scriptsize (0.14)} & -$0.020${\scriptsize (0.13)} & -$0.065${\scriptsize (0.12)} \\ 

 \color{blue}{Boundary}  & $0.001${\scriptsize (0.13)} &-$0.013${\scriptsize (0.13)} & -$0.012${\scriptsize (0.12)} &  \color{blue}{-$0.033$}{\scriptsize (0.10)} & -$0.036${\scriptsize (0.07)} 
& $0.006${\scriptsize (0.15)} & $0.001${\scriptsize (0.15)} & $0.015${\scriptsize (0.14)} & $0.001${\scriptsize (0.12)} & -$0.030${\scriptsize (0.10)} \\

 \color{blue}{Hybrid}  & $0.003${\scriptsize (0.13)} & -$0.009${\scriptsize (0.14)} &  -$0.010${\scriptsize (0.12)} & \color{blue}{-$0.032$}{\scriptsize(0.11)} & \color{blue}{-$0.034$}{\scriptsize (0.07)} 
&$0.008${\scriptsize (0.15)} & $0.005${\scriptsize (0.15)} & $0.018${\scriptsize (0.13)} & \color{blue}{$0.010$}{\scriptsize (0.12)} & \color{blue}{-$0.014$}{\scriptsize (0.09)} \\

 Tapered  & $0.003${\scriptsize (0.13)} & -$0.011${\scriptsize(0.14)} &-$0.013${\scriptsize (0.12)} & -$0.036${\scriptsize (0.11)} & -$0.038${\scriptsize (0.07)} 
& $0.007${\scriptsize (0.15)} & $0.004${\scriptsize (0.15)} & \color{blue}{$0.014$}{\scriptsize (0.13)} & \color{red}{$0$}{\scriptsize (0.11)} & -$0.026${\scriptsize (0.08)} \\
 
Debiased  & \color{blue}{$0.002$}{\scriptsize (0.13)} & \color{red}{-$0.013$}{\scriptsize (0.13)} & \color{blue}{-$0.014$}{\scriptsize (0.11)} & -$0.034${\scriptsize (0.11)} & -$0.030${\scriptsize (0.08)} 
& \color{blue}{$0.007$}{\scriptsize (0.15)} & \color{blue}{$0.001$}{\scriptsize (0.15)} & $0.017${\scriptsize (0.14)} & $0.015${\scriptsize (0.14)} & -$0.027${\scriptsize (0.13)} \\ \hline \hline

& \multicolumn{5}{c}{ \textbf{AR(1)}, $\{e_{t}\}\sim (\chi^2(2)-2)/2$, $n=300$}
& \multicolumn{5}{c}{ \textbf{MA(1)}, $\{e_{t}\}\sim (\chi^2(2)-2)/2$, $n=300$} \\ \hline

Gaussian  & $0${\scriptsize (0.06)} &\color{red}{-$0.005$}{\scriptsize (0.05)} & \color{red}{-$0.004$}{\scriptsize (0.05)} & \color{red}{-$0.004$}{\scriptsize (0.04)} & \color{red}{-$0.006$}{\scriptsize (0.03)} 
& $0${\scriptsize (0.06)} & \color{red}{-$0.002$}{\scriptsize (0.05)} & \color{red}{$0$}{\scriptsize (0.05)} & \color{red}{$0.003$}{\scriptsize (0.04)} & \color{blue}{$0.003$}{\scriptsize (0.03)} \\ 

 Whittle  & \color{red}{-$0.001$}{\scriptsize (0.06)} & \color{blue}{-$0.006$}{\scriptsize(0.05)} & -$0.005${\scriptsize (0.05)} & -$0.006${\scriptsize (0.04)} & -$0.009${\scriptsize (0.03)} 
& \color{red}{$0$}{\scriptsize (0.06)} & \color{blue}{-$0.003$}{\scriptsize (0.05)} & -$0.003${\scriptsize (0.05)} & -$0.004${\scriptsize (0.04)} & -$0.018${\scriptsize (0.04)} \\ 

 \color{blue}{Boundary}  & \color{blue}{$0$}{\scriptsize (0.06)} &-$0.005${\scriptsize (0.05)} & \color{blue}{-$0.004$}{\scriptsize (0.05)} &  \color{blue}{-$0.004$}{\scriptsize (0.04)} & \color{blue}{-$0.007$}{\scriptsize (0.03)} 
& \color{blue}{$0$}{\scriptsize (0.06)} & -$0.002${\scriptsize (0.05)} & \color{blue}{$0$}{\scriptsize (0.05)} & \color{blue}{$0.002$}{\scriptsize (0.04)} & -$0.002${\scriptsize (0.03)} \\

 \color{blue}{Hybrid}  & $0${\scriptsize (0.06)} & -$0.006${\scriptsize (0.06)} &  -$0.004${\scriptsize (0.05)} & -$0.004${\scriptsize(0.04)} & -$0.007${\scriptsize (0.03)} 
&$0.001${\scriptsize (0.06)} & -$0.002${\scriptsize (0.06)} & $0${\scriptsize (0.05)} & \color{blue}{$0.003$}{\scriptsize (0.04)} & \color{red}{$0.002$}{\scriptsize (0.03)} \\

Tapered  & $0${\scriptsize (0.06)} & -$0.006${\scriptsize (0.06)} &  -$0.005${\scriptsize (0.05)} & -$0.004${\scriptsize(0.04)} & -$0.007${\scriptsize (0.03)} 
&$0.001${\scriptsize (0.06)} & -$0.002${\scriptsize (0.06)} & $0${\scriptsize (0.05)} & $0.003${\scriptsize (0.04)} & $0.001${\scriptsize (0.03)} \\
 
Debiased  & $0${\scriptsize (0.06)} & -$0.005${\scriptsize (0.05)} &  -$0.004${\scriptsize (0.05)} & -$0.004${\scriptsize(0.04)} & -$0.006${\scriptsize (0.03)} 
&\color{blue}{$0$}{\scriptsize (0.06)} & -$0.002${\scriptsize (0.05)} & $0${\scriptsize (0.05)} & $0.003${\scriptsize (0.05)} & $0.013${\scriptsize (0.05)} \\
\end{tabular} 

\caption{\textit{Bias and the standard deviation (in the parentheses) of six different quasi-likelihoods for an AR(1) (left) and MA(1) (right) model for
the standardized chi-squared innovations. Length of the time series $n=20, 50$, and $300$. We use {\color{red}red} to denote the smallest RMSE and
{\color{blue}blue} to denote the second smallest RMSE.}}
\label{tab:AR.chisq}
\end{table}
\end{landscape}

\subsection{Misspecified model for a non-Gaussian time series}

In this section, we provide figures and table of the results in Section \ref{sec:misspecifiedmodel}
when the innovations follow a standardized chi-squared distribution
two degrees of freedom, i.e.  
$\varepsilon_t \sim (\chi^2(2)-2)/2$. 
The results are given in Tables \ref{tab:arma11.chisq} and \ref{tab:ar2.chisq}.

\begin{table}[ht]
\centering
\small
\begin{tabular}{cc|ccccccc}

$n$ & Parameter & Gaussian & Whittle & {\color{blue}Boundary} & {\color{blue}Hybrid} & Tapered & Debiased \\ \hline \hline
	
\multirow{3}{*}{20} & $\phi$ & {\color{blue} $0.029(0.1)$} & -$0.102(0.16)$ & -$0.032(0.12)$
&{\color{red} -$0.001(0.1)$} & -$0.088(0.13)$ & $0.170(0.12)$ \\ 

& $\psi$ & {\color{blue} $0.066(0.08)$} & -$0.184(0.20)$ & -$0.039(0.15)$ 
& {\color{red}$0.030(0.09)$} & -$0.064(0.12)$ & $0.086(0.09)$ \\ 

& $I_{n}(f;f_\theta)$ & $1.573(0.82)$ & $1.377(3.11)$ & $0.952(0.91)$
& {\color{blue} $1.006(0.84)$} & {\color{red}$0.675(0.63)$} & $2.618(0.84)$\\ \hline

\multirow{3}{*}{50} & $\phi$ & {\color{red}$0.014(0.07)$} & -$0.051(0.10)$ & -$0.004(0.07)$ 
& $0.007(0.07)$ & {\color{blue}-$0.003(0.07)$} & $0.143(0.11)$ \\

& $\psi$ & $0.027(0.06)$ & -$0.118(0.13)$ & -$0.013(0.09)$ 
& {\color{blue} $0.008(0.07)$} & {\color{red} $0.009(0.06)$} & $0.090(0.03)$ \\ 

& $I_{n}(f;f_\theta)$ & $0.342(0.34)$ & $0.478(0.53)$ & $0.298(0.32)$ 
& {\color{blue}$0.230(0.27)$} & {\color{red} $0.222(0.27)$} & $1.158(0.37)$ \\ \hline

\multirow{3}{*}{300} & $\phi$ & {\color{red} $0.001(0.03)$} & -$0.015(0.03)$ & {\color{blue} -$0.002(0.03)$} 
& $0(0.03)$ & -$0.001(0.03)$ & $0.090(0.08)$ \\ 

& $\psi$ & {\color{red}$0.006(0.03)$} & -$0.033(0.05)$ & $0.002(0.03)$ 
& {\color{red}$0.003(0.03)$} & $0.003(0.03)$ & $0.091(0.02)$ \\ 

& $I_{n}(f;f_\theta)$ & $0.029(0.05)$ & $0.067(0.10)$ & $0.034(0.06)$
& {\color{red}$0.027(0.04)$} & {\color{blue}$0.028(0.04)$} & $0.747(0.23)$ \\ \hline

\multicolumn{8}{l}{Best fitting ARMA$(1,1)$ coefficients $\theta = (\phi, \psi)$ and spectral divergence:} \\
\multicolumn{8}{l}{~~$-$ $\theta_{20}=(0.693, 0.845)$, $\theta_{50}=(0.694,0.857)$, $\theta_{300}=(0.696,0.857)$. } \\
\multicolumn{8}{l}{~~$-$ $I_{20}(f; f_{\theta}) = 3.773$, $I_{50}(f; f_{\theta}) = 3.415$, $I_{300}(f; f_{\theta}) = 3.388$.} \\

\end{tabular} 
\caption{\textit{Best fitting (bottom lines) and the bias of estimated
coefficients for six different methods for the ARMA$(3,2)$
misspecified case fitting ARMA$(1,1)$ model for the standardized chi-squared innovations. Standard deviations are in the parentheses. We
use {\color{red}red} to denote the smallest RMSE and 
{\color{blue}blue} to denote the second smallest RMSE.}
}
\label{tab:arma11.chisq}
\end{table}

\begin{table}[ht]
\centering
\footnotesize
\begin{tabular}{cc|ccccccc}

$n$ & Parameter & Gaussian & Whittle & {\color{blue}Boundary} & {\color{blue}Hybrid} & Tapered & Debiased \\ \hline \hline
	
\multirow{3}{*}{20} & $\phi_1$ & {\color{blue}$0.017(0.13)$} & -$0.178(0.23)$ & -$0.047(0.17)$
&{\color{red} -$0.006(0.14)$} & -$0.134(0.15)$ & $0.044(0.14)$ \\ 

& $\phi_2$ & {\color{red} $0.002(0.09)$} & $0.176(0.2)$ & $0.057(0.16)$ 
& $0.023(0.12)$ & $0.135(0.13)$ & {\color{blue}-$0.019(0.13)$} \\ 

& $I_{n}(f;f_\theta)$ & {\color{red}$0.652(0.72)$} & $1.3073(1.46)$ & $0.788(0.85)$ 
& {\color{blue}$0.671(0.8)$} & $0.887(0.97)$ & $0.658(0.81)$\\ \hline

\multirow{3}{*}{50} & $\phi_1$ & {\color{blue}$0.018(0.09)$} & -$0.079(0.12)$ & -$0.010(0.09)$ 
& {\color{red}$0.002(0.09)$} & {\color{blue}-$0.018(0.09)$} & $0.140(0.15)$ \\

& $\phi_2$ & -$0.018(0.06)$ & $0.072(0.11)$ & $0.012(0.07)$ 
& {\color{red} $0.001(0.06)$} & {\color{blue}$0.016(0.06)$} & -$0.1(0.09)$ \\ 

& $I_{n}(f;f_\theta)$ & {\color{red}$0.287(0.36)$} & $0.406(0.52)$ & $0.302(0.39)$ 
& $0.298(0.39)$ & {\color{blue}$0.293(0.38)$} & $0.631(0.7)$ \\ \hline

\multirow{3}{*}{300} & $\phi_1$ & {\color{red} $0.002(0.04)$} & -$0.015(0.04)$ & {\color{blue}-$0.002(0.04)$}
& $0(0.04)$ & -$0.001(0.04)$ & $0.012(0.04)$ \\ 

& $\phi_2$ & -$0.005(0.02)$ & $0.011(0.03)$ & {\color{red}-$0.001(0.02)$} 
& {\color{blue}-$0.001(0.02)$} & -$0.001(0.02)$ & -$0.016(0.04)$ \\ 

& $I_{n}(f;f_\theta)$ & {\color{red}$0.050(0.07)$} & $0.056(0.07)$ & {\color{blue} $0.051(0.07)$}
& $0.052(0.07)$ & $0.054(0.08)$ & $0.061(0.08)$ \\ \hline

\multicolumn{8}{l}{Best fitting AR$(1)$ coefficients $\theta = (\phi_1, \phi_2)$ and spectral divergence:} \\
\multicolumn{8}{l}{~~$-$ $\theta_{20}=(1.367, -0.841)$, $\theta_{50}=(1.364,-0.803)$, $\theta_{300}=(1.365,-0.802)$. } \\
\multicolumn{8}{l}{~~$-$ $I_{20}(f; f_{\theta}) = 2.902$, $I_{50}(f; f_{\theta}) = 2.937$, $I_{300}(f; f_{\theta}) = 2.916$.} \\

\end{tabular} 
\caption{\textit{Same as in Table \ref{tab:ar2.chisq} but fitting an AR(2) model.
}}
\label{tab:ar2.chisq}
\end{table}

\subsection{Comparing the the new likelihoods constructed with the
  predictive DFT with AR$(1)$ coefficients and AIC order selected AR$(p)$ coefficients} \label{sec:fixedP}

In this section we compare the performance of new likelihoods where the order of the AR
model used in the predictive DFT is determined using the AIC with a fixed choice
of order with the AR model (set to $p=1$). We use ARMA$(3,2)$ model considered in
Section \ref{sec:misspecifiedmodel} and fit the the ARMA$(1,1)$ and
AR$(2)$ to the data. We compare the new likelihoods with the Gaussian
likelihood and the Whittle likelihood. The results are given in Tables
 \ref{tab:arma11P} and \ref{tab:arma20P}. 

\begin{table}[ht]
\centering
\small
\begin{tabular}{c|l|ccc}

\multicolumn{2}{c}{} & $\phi$ & $\psi$ & $I_{n}(f;f_\theta)$ \\ \hline \hline

\multicolumn{2}{c}{Best} & $0.694$ & $0.857$ & $3.415$ \\ \hline	
\multirow{6}{*}{Bias}
& Gaussian & $0.012$\scriptsize{(0.07)} & $0.029$\scriptsize{(0.06)} & $0.354$\scriptsize{(0.34)} \\
& Whittle & -$0.054$\scriptsize{(0.09)} & -$0.116$\scriptsize{(0.12)} & $0.457$\scriptsize{(0.46)} \\  \cline{2-5}
& Boundary(AIC) & -$0.006$\scriptsize{(0.07)} & -$0.008$\scriptsize{(0.08)} & $0.292$\scriptsize{(0.3)} \\
& Boundary($p$=1) & -$0.020$\scriptsize{(0.08)} & -$0.045$\scriptsize{(0.09)} & $0.299$\scriptsize{(0.29)} \\  \cline{2-5}
& Hybrid(AIC) & $0.004$\scriptsize{(0.07)} & $0.009$\scriptsize{(0.07)} & $0.235$\scriptsize{(0.28)} \\
& Hybrid($p$=1) & $0.003$\scriptsize{(0.07)} & $0.010$\scriptsize{(0.07)} & $0.261$\scriptsize{(0.3)} \\

\end{tabular} 
\caption{\textit{Best fitting (top row) and the bias of estimated
coefficients for six different methods for the Gaussian ARMA$(3,2)$
misspecified case fitting ARMA$(1,1)$ model. Length of the time series $n$=50. Standard deviations are in the parentheses. 
(AIC): an order $p$ is chosen using AIC; ($p$=1): an order $p$ is set to 1. 
} }
\label{tab:arma11P}
\end{table}

\begin{table}[ht]
\centering
\small
\begin{tabular}{c|l|ccc}

\multicolumn{2}{c}{} & $\phi_1$ & $\phi_2$ & $I_{n}(f;f_\theta)$ \\ \hline \hline

\multicolumn{2}{c}{Best} & $1.364$ & -$0.803$ & $2.937$ \\ \hline	
\multirow{6}{*}{Bias}
& Gaussian & $0.019$\scriptsize{(0.09)} & -$0.024$\scriptsize{(0.06)} & $0.275$\scriptsize{(0.33)} \\
& Whittle & -$0.077$\scriptsize{(0.12)} & $0.066$\scriptsize{(0.1)} & $0.382$\scriptsize{(0.45)} \\  \cline{2-5}
& Boundary(AIC) & -$0.009$\scriptsize{(0.09)} & $0.006$\scriptsize{(0.07)} & $0.283$\scriptsize{(0.37)} \\
& Boundary($p$=1) & -$0.030$\scriptsize{(0.1)} & $0.032$\scriptsize{(0.07)} & $0.295$\scriptsize{(0.35)} \\  \cline{2-5}
& Hybrid(AIC) & $0.003$\scriptsize{(0.09)} & -$0.006$\scriptsize{(0.07)} & $0.283$\scriptsize{(0.37)} \\
& Hybrid($p$=1) & -$0.003$\scriptsize{(0.09)} & $0.003$\scriptsize{(0.06)} & $0.276$\scriptsize{(0.35)} \\

\end{tabular} 
\caption{\textit{Same as in Table \ref{tab:arma11P}, but fitting an AR$(2)$.
} }
\label{tab:arma20P}
\end{table}

\section{Simulations: Estimation for long memory time series}\label{sec:sim-long}

\subsection{Parametric estimation for long memory Gaussian time series}\label{sec:longGaussian}

We conduct some simulations for time series whose spectral
density, $f$, does not satisfies Assumption \ref{assum:A}. We focus on
the ARFIMA$(0,d,0)$ model where
\begin{eqnarray*}
(1-B)^{d}W_t = \varepsilon_t,
\end{eqnarray*} 
$B$ is the backshift operator, $-1/2<d<1/2$ is a fractional
differencing parameter, and $\{\varepsilon_t\}$ is an i.i.d. standard
normal random variable. Let $\Gamma(x)$ denote the gamma function.
The spectral density and autocovariance of the ARFIMA$(0,d,0)$ model (where the variance of
the innovations is set to $\sigma^2=1$) is 
\begin{eqnarray} 
\label{eq:arfima_formula}
f_W(\omega) = (1-e^{-i\omega})^{-2d} = (2\sin(\omega/2))^{-2d}
\quad \text{and} \quad
c_W(k) = \frac{\Gamma(k+d)\Gamma(1-2d)}{\Gamma(k-d+1) \Gamma(1-d) \Gamma(d)}
\end{eqnarray} 
respectively (see \cite{b:gir-12}, Chapter 7.2). Observe that for
$-1/2<d<0$, the $f_{W}(0)=0$, this is called antipersistence. On the other
hand, if $0< d <1/2$, then $f_{W}(0)=\infty$ and $W_{t}$ has long memory.

We generate ARFIMA$(0,d,0)$ models with  $d = -0.4,-0.2,0.2$ and $0.4$ and
Gaussian innovations. We fit both the ARFIMA$(0,d,0)$
model, with $d$ unknown (specified case) and the AR$(2)$ model (with
unknown parameters $\theta = (\phi_1,\phi_2)$)
(misspecified case) to the data. To do so, 
we first demean the time series. We evaluate the (plug-in) 
Gaussian likelihood using the autocovariance function in (\ref{eq:arfima_formula})
and the autocovariance function of AR$(2)$ model.
For the other 5 frequency domain likelihoods, we evaluate the likelihoods
at all the fundamental frequencies with the exception of the zero frequency
$\omega_{n,n}=0$. We fit using the spectral density in (\ref{eq:arfima_formula})
or the spectral density  $f_{\theta}(\omega) = |1-\phi_{1}e^{-i\omega} - \phi_{2}e^{-2i \omega}|^{-2}$ where 
$\theta=(\phi_{1}, \phi_{2})$ (depending on whether the model is
specified or misspecified).
For each simulation, we calculate six different parameter estimators. For the misspecified case,
we also calculate the spectral divergence 
\begin{eqnarray*}
\widetilde{I}_n(f;f_\theta) = \frac{1}{n-1} \sum_{k=1}^{n-1}\left( \frac{f(\omega_{k,n})}{f_\theta(\omega_{k,n})} +
\log f_\theta(\omega_{k,n})  \right),
\end{eqnarray*} 
where we omit the zero frequency.
The best fitting AR(2) model is $\theta^{Best} = \arg \min_{\theta
  \in \Theta} \widetilde{I}_n(f;f_\theta)$. In Tables \ref{tab:FI} and
\ref{tab:misFI} we give a bias and standard deviation
for the parameter estimators for the correctly specified and
misspecified model (in the misspecified case we also give the spectral
divergence). 

\vspace{1em}

\noindent \underline{Correctly specified model} From Table \ref{tab:FI} we observe that the bias of both new
likelihood estimators is consistently the smallest over all sample sizes and all $d$ except for $d=-0.2$. 
The new likelihoods have the
smallest or second smallest RMSE for $n$=300, but not for the small sample sizes (e.g. $n$=20 and 50). This is probably due to increased variation in the new likelihood
estimators caused by the estimation of the AR parameter for the predictive
DFT. Since the time series has a long memory, the AIC is likely to choose a
large order autoregressive order $p$, which will increase the variance
in the estimator (recall that the second order error of the boundary corrected Whittle is $O(p^{3}n^{-3/2})$).  
The (plug-in) Gaussian likelihood has a relatively large bias for all $d$, which 
matches the observations in \cite{p:lib-05}, Table 1. However, it has 
the smallest variance and this results in the smallest RMSE for almost all $n$ when $d$ is negative. 
The debiased Whittle also has a larger bias 
than most of the other estimators. 
However, it has a smaller variance and thus, having the
smallest RMSE for almost all $n$ and positive $d$.

\begin{table}[ht]
    \centering
\small
  \begin{tabular}{c|cccc}
\multirow{2}{*}{\textit{Likelihoods}} & \multicolumn{4}{c}{$d$} \\
\cline{2-5}	 
 & -0.4 & -0.2 & 0.2 & 0.4 \\ \hline \hline

& \multicolumn{4}{c}{$n=20$}\\ \hline
 Gaussian  &  {\color{red}-$0.097$}{\scriptsize (0.23)} &  {\color{red}-$0.148$}{\scriptsize (0.22)} & -$0.240${\scriptsize (0.23)}  & -$0.289${\scriptsize (0.22)}  \\

 Whittle  &  $0.027${\scriptsize (0.3)} & $0.006${\scriptsize (0.28)} & {\color{blue}-$0.008$}{\scriptsize (0.29)}  & {\color{blue}-$0.016$}{\scriptsize (0.29)} \\

 \color{blue}{Boundary} &  $0.014${\scriptsize (0.31)}  & $0${\scriptsize (0.29)} & -$0.005${\scriptsize (0.30)}  & -$0.007${\scriptsize (0.30)} \\

 \color{blue}{Hybrid}  &  $0.009${\scriptsize (0.31)} & -$0.007${\scriptsize (0.3)} & $0.005${\scriptsize (0.30)}  & -$0.001${\scriptsize (0.30)}  \\

 Tapered &  $0.026${\scriptsize (0.3)} & $0${\scriptsize (0.3)} & $0.006${\scriptsize (0.30)}  & $0.003${\scriptsize (0.29)} \\

 Debiased &  {\color{blue}$0.015$}{\scriptsize (0.29)} & {\color{blue}-$0.003$}{\scriptsize (0.27)} & {\color{red}-$0.029$}{\scriptsize (0.26)} & {\color{red}-$0.044$}{\scriptsize (0.27)}  \\ \hline \hline

& \multicolumn{4}{c}{$n=50$}\\ \hline
 Gaussian &  {\color{red}-$0.042$}{\scriptsize (0.13)} &  -$0.073${\scriptsize (0.14)} & -$0.097${\scriptsize (0.14)} & -$0.123${\scriptsize (0.12)}  \\

 Whittle  &  {\color{blue}$0.006$}{\scriptsize (0.15)} & {\color{red}-$0.016$}{\scriptsize (0.15)} & {\color{blue}-$0.013$}{\scriptsize (0.15)} & {\color{blue}-$0.005$}{\scriptsize (0.15)} \\

 \color{blue}{Boundary} &  -$0.005${\scriptsize (0.15)} & -$0.020${\scriptsize (0.16)} & -$0.011${\scriptsize (0.16)} & $0.001${\scriptsize (0.16)} \\

 \color{blue}{Hybrid} &  -$0.011${\scriptsize (0.15)} & -$0.021${\scriptsize (0.15)} & -$0.012${\scriptsize (0.16)}  & $0${\scriptsize (0.16)}  \\

 Tapered &  -$0.007${\scriptsize (0.15)} & -$0.019${\scriptsize (0.16)} & -$0.011${\scriptsize (0.16)} & $0.010${\scriptsize (0.16)} \\

 Debiased&  -$0.008${\scriptsize (0.16)}  & {\color{blue}-$0.020$}{\scriptsize (0.16)} & {\color{red}-$0.019$}{\scriptsize (0.15)} &  {\color{red}-$0.021$}{\scriptsize (0.14)}  \\ \hline \hline

& \multicolumn{4}{c}{$n=300$}\\ \hline
 Gaussian&  {\color{red}-$0.006$}{\scriptsize (0.05)}  &  -$0.013${\scriptsize (0.05)} & -$0.020${\scriptsize (0.05)} &  -$0.083${\scriptsize (0.02)}  \\

 Whittle&  $0.006${\scriptsize (0.05)}  & {\color{red}-$0.001$}{\scriptsize (0.05)} & {\color{blue}-$0.004$}{\scriptsize (0.05)} & {\color{blue}$0.002$}{\scriptsize (0.05)} \\

 \color{blue}{Boundary}&  {\color{blue}$0.002$}{\scriptsize (0.05)}  & {\color{blue}-$0.003$}{\scriptsize (0.05)} & {\color{blue}-$0.003$}{\scriptsize (0.05)}  & {\color{red}$0.002$}{\scriptsize (0.05)} \\

 \color{blue}{Hybrid}&  $0${\scriptsize (0.05)}  & -$0.004${\scriptsize (0.05)} & -$0.004${\scriptsize (0.05)}  & $0.001${\scriptsize (0.05)}  \\

 Tapered&  $0${\scriptsize (0.05)}  & -$0.004${\scriptsize (0.05)} & -$0.004${\scriptsize (0.05)} &  $0.003${\scriptsize (0.05)} \\

 Debiased&  -$0.001${\scriptsize (0.05)}  & -$0.003${\scriptsize (0.05)} & {\color{red}-$0.007$}{\scriptsize (0.05)} &  -$0.060${\scriptsize (0.02)}  \\ \hline \hline

\end{tabular} 

\caption{\textit{Bias and the standard deviation (in the parentheses) of six different quasi-likelihoods for ARFIMA$(0,d,0)$ model for
the standard normal innovations. Length of the time series $n=20, 50$, and $300$. We use {\color{red}red} to denote the smallest RMSE and
{\color{blue}blue} to denote the second smallest RMSE.
}}
\label{tab:FI}
\end{table}


\vspace{1em}

\noindent \underline{Misspecified model} We now compare the estimator when we fit the misspecified AR$(2)$ model to the data.  From Table \ref{tab:misFI}
we observe that the Gaussian likelihood performs uniformly well for all $d$ and $n$, it usually has the smallest bias and RMSE for 
the positive $d$. The tapered Whittle also performs uniformly well for all $d$ and $n$, especially for the negative $d$. In comparison, the new
likelihood estimators do not perform that well as compared with Gaussian and
Whittle likelihood. As mentioned above, this may be due to the
increased variation caused by estimating many AR parameters. However,
it is interesting to note that when $d=0.4$ and $n=300$, the estimated spectral divergence outperforms the Gaussian
likelihood. We leave the theoretical development of the sampling
properties of the new likelihoods and long memory time series for future research.

\begin{table}[ht]
\centering
\scriptsize
\begin{tabular}{ccc|c|cccccc}

\multirow{2}{*}{$d$} & \multirow{2}{*}{$n$} & \multirow{2}{*}{Par.} & \multirow{2}{*}{Best} & 
 \multicolumn{6}{c}{Bias} \\ \cline{5-10}
&&&&
Gaussian & Whittle & {\color{blue}Boundary} & {\color{blue}Hybrid} & Tapered & Debiased \\ \hline \hline

\multirow{9}{*}{-0.4} &
\multirow{3}{*}{20} & $\phi_1$ & -$0.300$ & -$0.028(0.22)$ & {\color{blue}-$0.015(0.22)$} & -$0.022(0.23)$
&-$0.026(0.23)$ & {\color{red}-$0.010(0.21)$} & -$0.026(0.22)$ \\ 

& & $\phi_2$ & -$0.134$ & -$0.067(0.19)$ & {\color{blue}-$0.058(0.19)$} & -$0.064(0.2)$ 
& -$0.068(0.2)$ & {\color{red}-$0.062(0.18)$} & -$0.067(0.19)$ \\ 

& & $I_{n}(f;f_\theta)$ & $1.141$ & $0.103(0.1)$ & $0.099(0.1)$ & $0.110(0.12)$ 
& $0.108(0.11)$ & {\color{red}$0.095(0.1)$} & {\color{blue}$0.106(0.1)$} \\ \cline{2-10}

& \multirow{3}{*}{50} & $\phi_1$ & -$0.319$ & {\color{blue}-$0.006(0.14)$} & {\color{red}$0.02(0.14)$} & -$0.004(0.15)$ 
& -$0.004(0.15)$ & $0.003(0.15)$ & -$0.006(0.15)$ \\

& & $\phi_2$ & -$0.152$ & -$0.028(0.13)$ & {\color{blue}-$0.022(0.13)$} & -$0.027(0.13)$ 
& -$0.029(0.13)$ & {\color{red}-$0.022(0.13)$} & -$0.029(0.13)$ \\ 

& & $I_{n}(f;f_\theta)$ & $1.092$ & {\color{blue}$0.043(0.05)$} & {\color{red}$0.043(0.05)$} & $0.045(0.05)$ 
& $0.046(0.05)$ & $0.043(0.05)$ & $0.045(0.05)$ \\ \cline{2-10}

& \multirow{3}{*}{300} & $\phi_1$ & -$0.331$ & {\color{blue}-$0.003(0.06)$} & {\color{red}-$0.002(0.06)$} & -$0.003(0.06)$
& -$0.002(0.06)$ & -$0.001(0.06)$ & {\color{blue}-$0.003(0.06)$} \\ 

& & $\phi_2$ & -$0.164$ & {\color{blue}-$0.005(0.06)$} & {\color{red}-$0.004(0.05)$} & {\color{blue}-$0.005(0.06)$} 
& -$0.005(0.06)$ & -$0.004(0.06)$ & -$0.005(0.06)$ \\ 

& & $I_{n}(f;f_\theta)$ & $1.062$ & {\color{red}$0.008(0.01)$} & {\color{red}$0.008(0.01)$} & $0.008(0.01)$
& $0.008(0.01)$ & $0.008(0.01)$ & $0.008(0.01)$ \\ \hline \hline

\multirow{9}{*}{-0.2} &
\multirow{3}{*}{20} & $\phi_1$ & -$0.157$ & -$0.027(0.23)$ & {\color{blue}-$0.024(0.23)$} & -$0.028(0.24)$
&-$0.029(0.23)$ & {\color{red}-$0.020(0.22)$} & -$0.028(0.23)$ \\ 

& & $\phi_2$ & -$0.066$ & -$0.084(0.21)$ & {\color{blue}-$0.085(0.2)$} & -$0.088(0.21)$ 
& -$0.089(0.21)$ & {\color{red}-$0.084(0.19)$} & -$0.089(0.21)$ \\ 

& & $I_{n}(f;f_\theta)$ & $1.066$ & $0.107(0.1)$ & $0.104(0.1)$ & $0.111(0.11)$ 
& $0.110(0.11)$ & {\color{red}$0.096(0.1)$} & {\color{blue}$0.108(0.1)$} \\ \cline{2-10}

& \multirow{3}{*}{50} & $\phi_1$ & -$0.170$ & {\color{blue}-$0.020(0.15)$} & {\color{red}-$0.015(0.15)$} & -$0.018(0.15)$ 
& -$0.020(0.15)$ & -$0.016(0.15)$ & -$0.019(0.15)$ \\

& & $\phi_2$ & -$0.079$ & -$0.035(0.14)$ & {\color{blue}-$0.032(0.14)$} & -$0.034(0.14)$ 
& -$0.035(0.14)$ & {\color{red}-$0.031(0.13)$} & -$0.036(0.14)$ \\ 

& & $I_{n}(f;f_\theta)$ & $1.038$ & $0.043(0.04)$ & {\color{blue}$0.043(0.04)$} & $0.045(0.04)$ 
& $0.045(0.04)$ & {\color{red}$0.043(0.04)$} & $0.045(0.04)$ \\ \cline{2-10}

& \multirow{3}{*}{300} & $\phi_1$ & -$0.001$ & {\color{blue} $0(0.06)$} & {\color{red}-$0.001(0.06)$} & $0(0.06)$
& $0.001(0.06)$ & $0.001(0.06)$ & -$0.006(0.06)$ \\ 

& & $\phi_2$ & -$0.088$ & {\color{blue}-$0.007(0.05)$} & {\color{red}-$0.007(0.05)$} & -$0.007(0.05)$ 
& -$0.007(0.06)$ & -$0.007(0.06)$ & -$0.007(0.05)$ \\ 

& & $I_{n}(f;f_\theta)$ & $1.019$ & {\color{red}$0.007(0.01)$} & {\color{red}$0.007(0.01)$} & $0.007(0.01)$
& $0.008(0.01)$ & $0.008(0.01)$ & $0.007(0.01)$ \\ \hline \hline

\multirow{9}{*}{0.2} &
\multirow{3}{*}{20} & $\phi_1$ & $0.167$ & -$0.066(0.24)$ & -$0.072(0.24)$ & -$0.071(0.25)$
& -$0.068(0.25)$ & {\color{red}-$0.080(0.23)$} & {\color{blue}-$0.070(0.24)$} \\ 

& & $\phi_2$ & $0.057$ & -$0.098(0.2)$ & {\color{blue}-$0.106(0.19)$} & -$0.107(0.19)$ 
& -$0.108(0.2)$ & {\color{red}-$0.111(0.18)$} & -$0.105(0.19)$ \\ 

& & $I_{n}(f;f_\theta)$ & $0.938$ & {\color{blue}$0.1(0.11)$} & $0.1(0.11)$ & $0.103(0.12)$ 
& $0.105(0.12)$ & $0.098(0.11)$ & {\color{red}$0.1(0.1)$} \\ \cline{2-10}

& \multirow{3}{*}{50} & $\phi_1$ & $0.186$ & {\color{blue}-$0.025(0.15)$} & {\color{red}-$0.027(0.15)$} & -$0.026(0.15)$ 
& -$0.027(0.15)$ & -$0.034(0.15)$ & -$0.025(0.15)$ \\

& & $\phi_2$ & $0.075$ & -$0.040(0.15)$ & {\color{red}-$0.043(0.15)$} & {\color{blue}-$0.043(0.15)$} 
& -$0.042(0.15)$ & -$0.047(0.15)$ & -$0.042(0.15)$ \\ 

& & $I_{n}(f;f_\theta)$ & $0.971$ & $0.046(0.05)$ & {\color{red}$0.044(0.05)$} & $0.046(0.05)$ 
& $0.048(0.05)$ & $0.047(0.05)$ & {\color{blue}$0.046(0.05)$} \\ \cline{2-10}

& \multirow{3}{*}{300} & $\phi_1$ & $0.208$ & {\color{red}-$0.007(0.06)$} & {\color{blue}-$0.007(0.06)$} & -$0.007(0.06)$ & -$0.007(0.06)$ & -$0.008(0.06)$ & -$0.007(0.06)$ \\ 

& & $\phi_2$ & $0.097$ & -$0.006(0.06)$ & {\color{red}-$0.007(0.06)$} & {\color{blue}-$0.006(0.06)$} & -$0.007(0.07)$ & -$0.008(0.07)$ & -$0.006(0.06)$ \\ 

& & $I_{n}(f;f_\theta)$ & $1.002$ & {\color{red}$0.008(0.01)$} & {\color{blue}$0.008(0.01)$} & {\color{blue}$0.008(0.01)$} & $0.008(0.01)$ & $0.009(0.01)$ & {\color{blue}$0.008(0.01)$ } \\ \hline \hline

\multirow{9}{*}{0.4} &
\multirow{3}{*}{20} & $\phi_1$ & $0.341$ & -$0.072(0.25)$ & {\color{red}-$0.082(0.25)$} & -$0.075(0.26)$ & -$0.074(0.26)$ & -$0.103(0.24)$ & {\color{blue}-$0.077(0.25)$} \\ 

& & $\phi_2$ & $0.094$ & {\color{blue}-$0.116(0.21)$} & -$0.134(0.2)$ & -$0.135(0.2)$ & {\color{blue}-$0.133(0.2)$} & -$0.133(0.19)$ & -$0.130(0.2)$ \\ 

& & $I_{n}(f;f_\theta)$ & $0.877$ & $0.111(0.12)$ & {\color{blue}$0.113(0.12)$} & $0.116(0.12)$ 
& $0.116(0.13)$ & $0.114(0.12)$ & {\color{red}$0.113(0.11)$} \\ \cline{2-10}

& \multirow{3}{*}{50} & $\phi_1$ & $0.378$ & {\color{red}-$0.024(0.15)$} & {\color{blue} -$0.030(0.15)$} & -$0.027(0.15)$ & -$0.029(0.15)$ & -$0.039(0.15)$ & -$0.027(0.15)$ \\

& & $\phi_2$ & $0.129$ & -$0.058(0.15)$ & {\color{red}-$0.069(0.14)$} & {\color{blue}-$0.067(0.15)$} & -$0.064(0.15)$ & -$0.072(0.15)$ & -$0.066(0.15)$ \\ 

& & $I_{n}(f;f_\theta)$ & $0.944$ & {\color{red}$0.051(0.06)$} & $0.053(0.06)$ & $0.053(0.06)$ & {\color{blue}$0.054(0.06)$} & $0.055(0.06)$ & $0.054(0.06)$ \\ \cline{2-10}

& \multirow{3}{*}{300} & $\phi_1$ & $0.428$ & {\color{red}-$0.004(0.06)$} & -$0.004(0.06)$ & {\color{blue}-$0.004(0.06)$} & -$0.007(0.06)$ & -$0.009(0.06)$ & -$0.004(0.06)$ \\ 

& & $\phi_2$ & $0.178$ & -$0.003(0.06)$ & {\color{red}-$0.005(0.06)$} & -$0.004(0.06)$ & -$0.004(0.06)$ & -$0.006(0.06)$ & {\color{blue}-$0.004(0.06)$} \\ 

& & $I_{n}(f;f_\theta)$ & $1.010$ & {\color{red}$0.009(0.01)$} & $0.009(0.01)$ & {\color{red}$0.009(0.01)$} & $0.010(0.01)$ & $0.010(0.01)$ & $0.009(0.01)$ \\

\end{tabular} 
\caption{\textit{ Best fitting and the bias of estimated coefficients using  the six different methods 
for misspecified Gaussian ARFIMA$(0,d,0)$ case fitting AR$(2)$ model. Standard deviations are in the parentheses.
We use {\color{red}red} to denote the smallest RMSE and 
{\color{blue}blue} to denote the second smallest RMSE. 
}}
\label{tab:misFI}
\end{table}

\subsection{Parametric estimation for long memory non-Gaussian time series}

We fit the parametric models described in Appendix \ref{sec:longGaussian}. However, the underlying time series is non-Gaussian and 
generated from the ARFIMA$(0,d,0)$ 
\begin{eqnarray*}
(1-B)^{d}W_t = \varepsilon_t,
\end{eqnarray*} 
where $\{\varepsilon_t\}$ are i.i.d. standardized chi-square random variables with two-degrees of freedom i.e.  
$\varepsilon_t \sim (\chi^2(2)-2)/2$. The results in the specified setting are given in Table 
\ref{tab:FI.chisq} and from the non-specified setting in Table \ref{tab:misFI.chisq}.

\begin{table}[ht]
    \centering
\small
  \begin{tabular}{c|cccc}
\multirow{2}{*}{\textit{Likelihoods}} & \multicolumn{4}{c}{$d$} \\
\cline{2-5}	 
 & -0.4 & -0.2 & 0.2 & 0.4 \\ \hline \hline

& \multicolumn{4}{c}{$n=20$}\\ \hline
 Gaussian& {\color{red}-$0.077$}{\scriptsize (0.22)}  &  {\color{red}-$0.130$}{\scriptsize (0.23)} & {\color{blue}-$0.219$}{\scriptsize (0.21)}  & -$0.277${\scriptsize (0.20)}  \\

 Whittle& $0.089${\scriptsize (0.34)}  & $0.092${\scriptsize (0.36)} & $0.077${\scriptsize (0.32)}  & {\color{blue}$0.058$}{\scriptsize (0.31)} \\

 \color{blue}{Boundary}& $0.077${\scriptsize (0.34)}  & $0.086${\scriptsize (0.37)} & $0.079${\scriptsize (0.33)} & $0.070${\scriptsize (0.32)} \\

 \color{blue}{Hybrid}& $0.077${\scriptsize (0.35)}  & $0.084${\scriptsize (0.37)} & $0.087${\scriptsize (0.33)}  & $0.086${\scriptsize (0.32)}  \\

 Tapered& $0.092${\scriptsize (0.35)}  & $0.096${\scriptsize (0.37)} & $0.089${\scriptsize (0.33)}  & $0.078${\scriptsize (0.31)} \\

 Debiased & {\color{blue}$0.057$}{\scriptsize (0.33)} & {\color{blue}$0.051$}{\scriptsize (0.31)} & {\color{red}$0.009$}{\scriptsize (0.25)}  & {\color{red}-$0.021$}{\scriptsize (0.26)}  \\ \hline \hline

& \multicolumn{4}{c}{$n=50$}\\ \hline
 Gaussian & {\color{red}-$0.047$}{\scriptsize (0.13)} &  {\color{red}-$0.065$}{\scriptsize (0.14)} & -$0.097${\scriptsize (0.13)}  & -$0.130${\scriptsize (0.12)}  \\

 Whittle& {\color{blue}$0.008$}{\scriptsize (0.15)}  & $0.004${\scriptsize (0.16)} & {\color{blue}-$0.001$}{\scriptsize (0.15)}  & {\color{blue}$0.001$}{\scriptsize (0.16)} \\

 \color{blue}{Boundary} & -$0.002${\scriptsize (0.16)}  & $0${\scriptsize (0.16)} & $0.002${\scriptsize (0.15)}  & $0.006${\scriptsize (0.16)} \\

 \color{blue}{Hybrid} & -$0.005${\scriptsize (0.15)} & -$0.001${\scriptsize (0.16)} & $0.004${\scriptsize (0.16)} &  $0.006${\scriptsize (0.16)}  \\

 Tapered & $0.001${\scriptsize (0.15)} & $0.002${\scriptsize (0.16)} & $0.006${\scriptsize (0.16)}  & $0.016${\scriptsize (0.17)} \\

 Debiased & -$0.004${\scriptsize (0.16)} & {\color{blue}-$0.001$}{\scriptsize (0.16)} & {\color{red}-$0.013$}{\scriptsize (0.14)}  & {\color{red}-$0.025$}{\scriptsize (0.14)}  \\ \hline \hline

& \multicolumn{4}{c}{$n=300$}\\ \hline
 Gaussian & {\color{red}-$0.011$}{\scriptsize (0.05)} &  -$0.012${\scriptsize (0.05)} & -$0.017${\scriptsize (0.05)}  & -$0.085${\scriptsize (0.02)}  \\

 Whittle & $0.002${\scriptsize (0.05)} & {\color{red}$0$}{\scriptsize (0.05)} & {\color{blue}-$0.001$}{\scriptsize (0.05)} & {\color{blue}-$0.001$}{\scriptsize (0.05)} \\

 \color{blue}{Boundary} & {\color{blue}-$0.003$}{\scriptsize (0.05)} & {\color{blue}-$0.001$}{\scriptsize (0.05)} & {\color{blue}$0$}{\scriptsize (0.05)} & {\color{red}$0.001$}{\scriptsize (0.05)} \\

 \color{blue}{Hybrid} & -$0.006${\scriptsize (0.05)} & -$0.003${\scriptsize (0.05)} & -$0.001${\scriptsize (0.05)}  & -$0.003${\scriptsize (0.05)}  \\

 Tapered & -$0.006${\scriptsize (0.05)}  & -$0.003${\scriptsize (0.05)} & -$0.001${\scriptsize (0.05)}  & -$0.001${\scriptsize (0.05)} \\

 Debiased & -$0.006${\scriptsize (0.05)} & -$0.002${\scriptsize (0.05)} & {\color{red}-$0.004$}{\scriptsize (0.05)}  & -$0.061${\scriptsize (0.03)}  \\ \hline \hline

\end{tabular} 

\caption{\textit{Same as in Table \ref{tab:FI} but for the chi-squared innovations.
}}
\label{tab:FI.chisq}
\end{table}

\begin{table}[ht]
\centering
\scriptsize
\begin{tabular}{ccc|c|cccccc}

\multirow{2}{*}{$d$} & \multirow{2}{*}{$n$} & \multirow{2}{*}{Par.} & \multirow{2}{*}{Best} & 
 \multicolumn{6}{c}{Bias} \\ \cline{5-10}
&&&&
Gaussian & Whittle & {\color{blue}Boundary} & {\color{blue}Hybrid} & Tapered & Debiased \\ \hline \hline
	
\multirow{9}{*}{-0.4} &
\multirow{3}{*}{20} & $\phi_1$ & -$0.300$ & {\color{blue}-$0.017(0.21)$} & -$0.004(0.21)$ & -$0.012(0.22)$
& -$0.018(0.22)$ & {\color{red}$0.001(0.21)$} & -$0.009(0.22)$ \\ 

& & $\phi_2$ & -$0.134$ & -$0.073(0.18)$ & {\color{blue}-$0.071(0.18)$} & -$0.076(0.18)$ 
& -$0.077(0.18)$ & {\color{red}-$0.064(0.16)$} & -$0.074(0.18)$ \\ 

& & $I_{n}(f;f_\theta)$ & $1.141$ & $0.096(0.11)$ & {\color{red}$0.097(0.11)$} & $0.105(0.11)$ 
& $0.102(0.11)$ & {\color{blue}$0.088(0.11)$} & $0.100(0.11)$ \\ \cline{2-10}

& \multirow{3}{*}{50} & $\phi_1$ & -$0.319$ & {\color{blue}-$0.008(0.14)$} & {\color{red}$0(0.14)$} & -$0.005(0.14)$ 
& -$0.007(0.15)$ &$0.001(0.14)$ & -$0.006(0.14)$ \\

& & $\phi_2$ & -$0.152$ & -$0.035(0.12)$ & {\color{blue}-$0.030(0.12)$} & -$0.035(0.12)$
& -$0.035(0.12)$ & {\color{red}-$0.026(0.12)$} & -$0.035(0.12)$ \\ 

& & $I_{n}(f;f_\theta)$ & $1.092$ & $0.042(0.05)$ & {\color{red}$0.040(0.04)$} & $0.043(0.05)$ 
& $0.043(0.05)$ & {\color{blue}$0.041(0.04)$} & {\color{blue}$0.043(0.05)$} \\ \cline{2-10}

& \multirow{3}{*}{300} & $\phi_1$ & -$0.331$ & {\color{red}-$0.004(0.06)$} & {\color{blue}-$0.002(0.06)$} & -$0.004(0.06)$ & -$0.004(0.06)$ & -$0.003(0.06)$ & -$0.004(0.06)$ \\ 

& & $\phi_2$ & -$0.164$ & {\color{red}-$0.007(0.05)$} & {\color{red}-$0.006(0.05)$} & -$0.007(0.05)$ & -$0.008(0.05)$ & -$0.007(0.05)$ & -$0.007(0.05)$ \\ 

& & $I_{n}(f;f_\theta)$ & $1.062$ & {\color{red}$0.007(0.01)$} & {\color{blue}$0.007(0.01)$} & {\color{blue}$0.007(0.01)$} & {\color{blue}$0.008(0.01)$} & $0.008(0.01)$ & $0.007(0.01)$ \\ \hline \hline

\multirow{9}{*}{-0.2} &
\multirow{3}{*}{20} & $\phi_1$ & -$0.157$ & -$0.036(0.22)$ & -$0.034(0.22)$ & -$0.037(0.23)$
& -$0.039(0.23)$ & {\color{red}-$0.028(0.21)$} & {\color{blue}-$0.034(0.22)$} \\ 

& & $\phi_2$ & -$0.066$ & -$0.083(0.19)$ & {\color{blue}-$0.080(0.18)$} & -$0.082(0.19)$ 
& -$0.083(0.19)$ & {\color{red}-$0.076(0.18)$} & -$0.082(0.19)$ \\ 

& & $I_{n}(f;f_\theta)$ & $1.066$ & $0.099(0.11)$ & {\color{blue}$0.095(0.11)$} & $0.100(0.11)$ 
& $0.101(0.11)$ & {\color{red}$0.086(0.1)$} & $0.097(0.11)$ \\ \cline{2-10}

& \multirow{3}{*}{50} & $\phi_1$ & -$0.170$ & {\color{blue}-$0.018(0.15)$} & {\color{red}-$0.016(0.14)$} & -$0.019(0.15)$ 
& -$0.017(0.15)$ & -$0.012(0.15)$ & -$0.019(0.15)$ \\

& & $\phi_2$ & -$0.079$ & -$0.034(0.13)$ & {\color{blue}-$0.033(0.13)$} & -$0.034(0.13)$
& -$0.035(0.13)$ & {\color{red}-$0.031(0.13)$} & -$0.034(0.13)$ \\ 

& & $I_{n}(f;f_\theta)$ & $1.038$ & $0.042(0.05)$ & {\color{blue}$0.041(0.05)$} & $0.043(0.05)$ 
& $0.043(0.05)$ & {\color{red}$0.041(0.04)$} & $0.043(0.05)$ \\ \cline{2-10}

& \multirow{3}{*}{300} & $\phi_1$ & -$0.179$ & {\color{red}-$0.001(0.06)$} & {\color{red}$0(0.06)$} & -$0.001(0.06)$ & $0(0.06)$ & $0(0.06)$ & -$0.001(0.06)$ \\ 

& & $\phi_2$ & -$0.088$ & {\color{red}-$0.008(0.05)$} & {\color{red}-$0.007(0.05)$} & -$0.008(0.06)$ & -$0.008(0.06)$ & -$0.007(0.06)$ & -$0.008(0.06)$ \\ 

& & $I_{n}(f;f_\theta)$ & $1.019$ & {\color{red}$0.007(0.01)$} & {\color{red}$0.007(0.01)$} & {\color{red}$0.007(0.01)$} & $0.007(0.01)$ & $0.007(0.01)$ & {\color{red}$0.007(0.01)$} \\ \hline \hline

\multirow{9}{*}{0.2} &
\multirow{3}{*}{20} & $\phi_1$ & $0.167$ & -$0.053(0.23)$ & -$0.064(0.22)$ & -$0.062(0.23)$
& -$0.055(0.23)$ & {\color{blue}-$0.071(0.21)$} & {\color{red}-$0.066(0.21)$} \\ 

& & $\phi_2$ & $0.057$ & -$0.091(0.2)$ & -$0.100(0.2)$ & -$0.100(0.2)$ 
& -$0.101(0.2)$ & {\color{red}-$0.099(0.19)$} & {\color{blue}-$0.095(0.19)$} \\ 

& & $I_{n}(f;f_\theta)$ & $0.938$ & $0.094(0.1)$ & $0.091(0.09)$ & $0.094(0.1)$ 
& $0.097(0.1)$ & {\color{blue}$0.086(0.09)$} & {\color{red}$0.086(0.08)$} \\ \cline{2-10}

& \multirow{3}{*}{50} & $\phi_1$ & $0.186$ & {\color{red}-$0.023(0.15)$} & {\color{blue}-$0.026(0.15)$} & -$0.025(0.15)$ 
& -$0.024(0.15)$ & -$0.030(0.15)$ & -$0.024(0.15)$ \\

& & $\phi_2$ & $0.075$ & -$0.044(0.14)$ & {\color{blue}-$0.051(0.14)$} & -$0.051(0.14)$
& -$0.047(0.14)$ & {\color{red}-$0.051(0.13)$} & -$0.050(0.14)$ \\ 

& & $I_{n}(f;f_\theta)$ & $0.971$ & {\color{red}$0.043(0.04)$} & {\color{red}$0.042(0.04)$} & $0.043(0.04)$ 
& $0.044(0.04)$ & {\color{blue}$0.042(0.04)$} & $0.043(0.04)$ \\ \cline{2-10}

& \multirow{3}{*}{300} & $\phi_1$ & $0.208$ & {\color{red}-$0.003(0.06)$} & {\color{blue}-$0.004(0.06)$} & -$0.003(0.06)$ & -$0.003(0.06)$ & -$0.005(0.06)$ & -$0.003(0.06)$ \\ 

& & $\phi_2$ & $0.097$ & {\color{blue}-$0.008(0.06)$} & {\color{red}-$0.009(0.06)$} & -$0.008(0.06)$ & -$0.010(0.07)$ & -$0.011(0.06)$ & -$0.008(0.06)$ \\ 

& & $I_{n}(f;f_\theta)$ & $1.002$ & {\color{red}$0.007(0.01)$} & {\color{red}$0.007(0.01)$} & {\color{red}$0.007(0.01)$} & $0.007(0.01)$ & $0.008(0.01)$ & {\color{red}$0.007(0.01)$ } \\ \hline \hline

\multirow{9}{*}{0.4} &
\multirow{3}{*}{20} & $\phi_1$ & $0.341$ & -$0.073(0.24)$ & {\color{blue}-$0.080(0.23)$} & -$0.074(0.24)$ & -$0.070(0.25)$ & -$0.103(0.23)$ & {\color{red}-$0.085(0.23)$} \\ 

& & $\phi_2$ & $0.094$ & -$0.106(0.21)$ & -$0.129(0.19)$ & -$0.128(0.2)$ & -$0.126(0.2)$ & {\color{red}-$0.122(0.19)$} & {\color{blue}-$0.123(0.19)$} \\ 

& & $I_{n}(f;f_\theta)$ & $0.877$ & {\color{red}$0.103(0.11)$} & {\color{blue}$0.103(0.11)$} & $0.107(0.12)$ 
& $0.108(0.12)$ & $0.103(0.11)$ & $0.101(0.11)$ \\ \cline{2-10}

& \multirow{3}{*}{50} & $\phi_1$ & $0.378$ & {\color{red}-$0.021(0.14)$} & {\color{red} -$0.027(0.14)$} & -$0.023(0.15)$ & -$0.026(0.15)$ & -$0.037(0.15)$ & -$0.026(0.14)$ \\

& & $\phi_2$ & $0.129$ & {\color{red}-$0.044(0.14)$} & {\color{blue}-$0.054(0.14)$} & -$0.051(0.14)$ & -$0.049(0.14)$ & -$0.059(0.14)$ & -$0.051(0.14)$ \\ 

& & $I_{n}(f;f_\theta)$ & $0.944$ & $0.044(0.05)$ & {\color{blue}$0.045(0.05)$} & $0.046(0.05)$ & $0.046(0.05)$ & $0.047(0.05)$ & {\color{red}$0.045(0.05)$} \\ \cline{2-10}

& \multirow{3}{*}{300} & $\phi_1$ & $0.428$ & {\color{blue}-$0.003(0.06)$} & {\color{blue}-$0.003(0.06)$} & {\color{red}-$0.003(0.06)$} & -$0.004(0.06)$ & -$0.006(0.06)$ & -$0.002(0.06)$ \\ 

& & $\phi_2$ & $0.178$ & {\color{blue}-$0.004(0.06)$} & {\color{red}-$0.006(0.06)$} & {\color{blue}-$0.004(0.06)$} & -$0.006(0.06)$ & -$0.008(0.06)$ & {\color{blue}-$0.005(0.06)$} \\ 

& & $I_{n}(f;f_\theta)$ & $1.010$ & $0.010(0.01)$ & {\color{red}$0.009(0.01)$} & {\color{red}$0.009(0.01)$} & $0.010(0.01)$ & $0.010(0.01)$ & {\color{red}$0.009(0.01)$} \\

\end{tabular} 
\caption{\textit{Best fitting and the bias of estimated coefficients using the six different methods 
for misspecified ARFIMA$(0,d,0)$ case fitting AR$(2)$ model for the chi-squared innovations. 
Standard deviations are in the parentheses.
We use {\color{red}red} to denote the smallest RMSE and 
{\color{blue}blue} to denote the second smallest RMSE.
}}
\label{tab:misFI.chisq}
\end{table}

\subsection{Semi-parametric estimation for Gaussian time series}\label{sec:LWGaussian}

Suppose the time series $\{X_t\}$ has a spectral density $f(\cdot)$ with $\lim_{\omega \rightarrow 0+} f(\omega) \sim C\omega^{-2d}$ for some
$d \in (-1/2, 1/2)$. The local Whittle (LW)  estimator is an estimation method for estimating $d$ without using assuming any parametric structure on $d$. It was first proposed in 
\cite{p:kun-87,p:rob-95,p:che-03}, see also \cite{b:gir-12}, Chapter 8). The LW estimator is defined as
$\widehat{d} = \arg\min R(d)$ where
\begin{eqnarray} \label{eq:LW}
R(d) = \log \left( M^{-1} \sum_{k=1}^{M} \frac{|J_n(\omega_{k,n})|^2}{\omega_{k,n}^{-2d}} \right)
- \frac{2d}{M} \sum_{k=1}^{M} \log \omega_{k,n},
\end{eqnarray} 
and $M=M(n)$ is an integer such that $M^{-1} + M/n \rightarrow 0$ as $n \rightarrow \infty$.
The objective function can be viewed as ``locally'' fitting a spectral density of form  $f_{\theta}(\omega) = C\omega^{-2d}$ where $\theta = (C,d)$ using
the Whittle likelihood. 

Since $\widetilde{J}_n (\omega_{k,n} ; f) \overline{J_n(\omega_{k,n})}$ is an unbiased estimator of 
true spectral density $f(\omega_{k,n})$, it is possible that replacing the periodogram with 
the (feasible) complete periodogram my lead to a better estimator of $d$. Based on this 
we define the (feasible) hybrid LW criterion,
\begin{eqnarray*}
Q(d) = \log \left( M^{-1} \sum_{k=1}^{M} \frac{\widetilde{J}_n (\omega_{k,n} ; \widehat{f}_p) \overline{J_{n,\underline{h}_n}(\omega_{k,n})}}{\omega_{k,n}^{-2d}} \right) - \frac{2d}{M} \sum_{k=1}^{M} \log \omega_{k,n}.
\end{eqnarray*} 
In a special case that the data taper $h_{t,n} \equiv 1$, we call it the boundary corrected LW criterion.

To empirically assess the validity of the above estimation scheme, 
we generate a Gaussian ARFIMA$(0,d,0)$ model from Section \ref{sec:longGaussian} for $d=-0.4, -0.2, 0.2$ and $0.4$ and evaulate the LW, tapered LW (using tapered DFT in (\ref{eq:LW})), boundary corrected LW, and hybrid LW. We set $M \approx n^{0.65}$ where $n$ is a length of the time series and we use Tukey taper with 10\% of the taper on each end of the time series. For each simulation, we obtain four different LW estimators.

Table \ref{tab:LW} summarizes the bias and standard deviation (in the parentheses) of LW estimators. We observe that the bondary corrected LW has a smaller bias than the regular Local Whittle likelihood except when $d=-0.2$ and $n=50, 300$. However, the standard error tends to be larger (this is probably because of the additional error caused by estimating the AR$(p)$ parameters in the new likelihoods). Despite the larger standard error, in terms of RMSE, the 
 boundary corrected LW (or hybrid) tends to have overall at least the second smallest RSME for most $d$ and $n$.

\begin{table}[h]
    \centering
\small
  \begin{tabular}{c|cccc}
\multirow{2}{*}{\textit{Local likelihoods}} & \multicolumn{4}{c}{$d$} \\
\cline{2-5}	 
 & -0.4 & -0.2 & 0.2 & 0.4 \\ \hline \hline

& \multicolumn{4}{c}{$n=20$}\\ \hline

 Whittle& {\color{red}$0.283$}{\scriptsize (0.46)}  & {\color{red}$0.111$}{\scriptsize (0.5)} & {\color{red}-$0.075$}{\scriptsize (0.48)}  & -$0.226${\scriptsize (0.43)} \\

 \color{blue}{Boundary}& $0.282${\scriptsize (0.46)}  & {\color{blue}$0.109$}{\scriptsize (0.51)} & {\color{blue}-$0.075$}{\scriptsize (0.48)} & -$0.220${\scriptsize (0.44)} \\

 \color{blue}{Hybrid}& $0.275${\scriptsize (0.46)}  & $0.115${\scriptsize (0.52)} & -$0.067${\scriptsize (0.49)}  & {\color{blue}-$0.209$}{\scriptsize (0.44)}  \\

 Tapered& {\color{blue}$0.279$}{\scriptsize (0.46)}  & $0.121${\scriptsize (0.52)} & -$0.068${\scriptsize (0.49)}  & {\color{red}-$0.204$}{\scriptsize (0.43)} \\ \hline \hline

& \multicolumn{4}{c}{$n=50$}\\ \hline
 Whittle& $0.060${\scriptsize (0.26)}  & -$0.056${\scriptsize (0.33)} & {\color{red}-$0.089$}{\scriptsize (0.37)}  & -$0.109${\scriptsize (0.32)} \\

 \color{blue}{Boundary} & $0.045${\scriptsize (0.26)}  & {\color{blue}-$0.063$}{\scriptsize (0.34)} & {\color{blue}-$0.088$}{\scriptsize (0.38)}  & {\color{blue}-$0.106$}{\scriptsize (0.32)} \\

 \color{blue}{Hybrid} & {\color{blue}$0.033$}{\scriptsize (0.25)} & -$0.069${\scriptsize (0.34)} & -$0.090${\scriptsize (0.38)} &  -$0.110${\scriptsize (0.32)}  \\

 Tapered & {\color{red}$0.035$}{\scriptsize (0.25)} & {\color{red}-$0.068$}{\scriptsize (0.34)} & -$0.085${\scriptsize (0.38)}  & {\color{red}-$0.085$}{\scriptsize (0.31)} \\ \hline \hline

& \multicolumn{4}{c}{$n=300$}\\ \hline
 Whittle & {\color{red}$0.056$}{\scriptsize (0.12)} & {\color{red}-$0.014$}{\scriptsize (0.11)} & {\color{red}-$0.010$}{\scriptsize (0.12)} & {\color{blue}$0.004$}{\scriptsize (0.11)} \\

 \color{blue}{Boundary} & {\color{blue}$0.052$}{\scriptsize (0.12)} & {\color{blue}-$0.017$}{\scriptsize (0.11)} & {\color{blue}-$0.009$}{\scriptsize (0.12)} & {\color{red}$0.003$}{\scriptsize (0.11)} \\

 \color{blue}{Hybrid} & $0.054${\scriptsize (0.12)} & -$0.018${\scriptsize (0.12)} & -$0.011${\scriptsize (0.12)}  & -$0.002${\scriptsize (0.11)}  \\

 Tapered & $0.057${\scriptsize (0.12)}  & -$0.018${\scriptsize (0.12)} & -$0.011${\scriptsize (0.12)}  & $0.003${\scriptsize (0.12)} \\ \hline \hline

\end{tabular} 

\caption{\textit{Bias and the standard deviation (in the parentheses) of four different Local Whittle estimators for ARFIMA$(0,d,0)$ model for
the standard normal innovations. Length of the time series $n=20, 50$, and $300$. We use {\color{red}red} to denote the smallest RMSE and
{\color{blue}blue} to denote the second smallest RMSE.
}}
\label{tab:LW}
\end{table}

\subsection{Semi-parametric estimation for long memory non-Gaussian time series}

Once again we consider the semi-parametric Local Whittle estimator described in Appendix \ref{sec:LWGaussian}. However, this time we assess the estimation scheme for non-Gaussian time series. We generate the ARFIMA$(0,d,0)$ 
\begin{eqnarray*}
(1-B)^{d}W_t = \varepsilon_t,
\end{eqnarray*} 
where $\{\varepsilon_t\}$ are i.i.d. standardarized chi-square random variables with 
two-degrees of freedom i.e.  $\varepsilon_t \sim (\chi^2(2)-2)/2$. 
The results are summarized in Table \ref{tab:LW.chisq}. 

\begin{table}[ht]
    \centering
\small
  \begin{tabular}{c|cccc}
\multirow{2}{*}{\textit{Local likelihoods}} & \multicolumn{4}{c}{$d$} \\
\cline{2-5}	 
 & -0.4 & -0.2 & 0.2 & 0.4 \\ \hline \hline

& \multicolumn{4}{c}{$n=20$}\\ \hline

 Whittle& $0.354${\scriptsize (0.52)}  & {\color{blue}$0.235$}{\scriptsize (0.54)} & {\color{red}$0.002$}{\scriptsize (0.51)}  & -$0.154${\scriptsize (0.41)} \\

 \color{blue}{Boundary}& {\color{blue}$0.349$}{\scriptsize (0.52)}  & {\color{red}$0.231$}{\scriptsize (0.54)} & $0.003${\scriptsize (0.51)} & -$0.146${\scriptsize (0.41)} \\

 \color{blue}{Hybrid}& $0.347${\scriptsize (0.53)}  & $0.229${\scriptsize (0.55)} & $0.016${\scriptsize (0.52)}  & {\color{red}-$0.133$}{\scriptsize (0.4)}  \\

 Tapered& {\color{red}$0.351$}{\scriptsize (0.52)}  & $0.229${\scriptsize (0.54)} & {\color{blue}$0.023$}{\scriptsize (0.51)}  & {\color{blue}-$0.144$}{\scriptsize (0.41)} \\ \hline \hline

& \multicolumn{4}{c}{$n=50$}\\ \hline
 Whittle& $0.019${\scriptsize (0.25)}  & {\color{red}-$0.080$}{\scriptsize (0.34)} & {\color{red}-$0.125$}{\scriptsize (0.38)}  & -$0.149${\scriptsize (0.33)} \\

 \color{blue}{Boundary} & $0.007${\scriptsize (0.25)}  & {\color{blue}-$0.086$}{\scriptsize (0.35)} & {\color{blue}-$0.123$}{\scriptsize (0.39)}  & {\color{blue}-$0.146$}{\scriptsize (0.33)} \\

 \color{blue}{Hybrid} & {\color{red}-$0.003$}{\scriptsize (0.24)} & -$0.098${\scriptsize (0.35)} & -$0.128${\scriptsize (0.39)} &  -$0.153${\scriptsize (0.34)}  \\

 Tapered & {\color{blue}$0.001$}{\scriptsize (0.24)} & -$0.100${\scriptsize (0.35)} & -$0.125${\scriptsize (0.39)}  & {\color{red}-$0.133$}{\scriptsize (0.33)} \\ \hline \hline

& \multicolumn{4}{c}{$n=300$}\\ \hline
 Whittle & {\color{red}$0.104$}{\scriptsize (0.14)} & {\color{red}-$0.006$}{\scriptsize (0.14)} & {\color{red}-$0.015$}{\scriptsize (0.14)} & {\color{blue}-$0.011$}{\scriptsize (0.13)} \\

 \color{blue}{Boundary} & $0.101${\scriptsize (0.15)} & {\color{blue}-$0.008$}{\scriptsize (0.14)} & {\color{blue}-$0.014$}{\scriptsize (0.14)} & {\color{red}-$0.012$}{\scriptsize (0.13)} \\

 \color{blue}{Hybrid} & {\color{blue}$0.100$}{\scriptsize (0.15)} & -$0.011${\scriptsize (0.15)} & -$0.017${\scriptsize (0.14)}  & -$0.020${\scriptsize (0.13)}  \\

 Tapered & $0.104${\scriptsize (0.15)}  & -$0.012${\scriptsize (0.15)} & -$0.017${\scriptsize (0.15)}  & -$0.016${\scriptsize (0.13)} \\ \hline \hline

\end{tabular} 

\caption{\textit{Same as in Table \ref{tab:LW} but for the chi-square innovations.
}}
\label{tab:LW.chisq}
\end{table}

\section{Simulations: Alternative methods for estimating the
  predictive DFT} \label{sec:alternative}


As pointed out by the referees, using the Yule-Walker estimator to
estimate the prediction coefficients in the predictive DFT may in
certain situations be problematic. We discuss the issues and potential
solutions below. 

The first issue is that Yule-Walker
estimator suffers a finite sample bias, especially when the spectral density has a root close to the 
unit circle (see, e.g., \cite{p:tjo-83}). One remedy to reduce the
bias is via data tapering (\cite{p:dah-88} and \cite{p:zha-92}).
Therefore, we define the \textbf{b}oundary \textbf{c}orrected Whittle likelihood using \textbf{t}apered \textbf{Y}ule-\textbf{W}alker (BC-tYW) replace
$\widehat{f}_p$ with $\widetilde{f}_p$ in (\ref{eq:Winf})
where
$\widetilde{f}_p$ is a spectral density of AR$(p)$ process where the AR coefficients are estimated using 
Yule-Walker with tapered time series. In the simulations we use the Tukey taper with 
$d=n/10$ and select the order $p$ using the AIC.

The second issue is if the underlying time series is complicated in
the sense that the underlying AR representation has multiple roots. Then fitting a large order
AR$(p)$ model may result in a loss of efficiency. As an alternative,
we consider a fully nonparametric estimator of $\widehat{J}_n(\omega;f)$ based on the estimated spectral density function. 
To do so, 
we recall from Section \ref{sec:firstorder} the first order approximation of $\widehat{J}_n(\omega;f)$ is $\widehat{J}_{\infty,n}(\omega;f)$ where
\begin{eqnarray*}
\widehat{J}_{\infty,n}(\omega_{};f) 
&=& 
\frac{n^{-1/2}}{\phi(\omega;f)} \sum_{t=1}^{n}X_{t}\phi_{t}^{\infty}(\omega;f)
+e^{i(n+1)\omega}
\frac{n^{-1/2}}{ \overline{\phi(\omega;f)}} \sum_{t=1}^{n}X_{n+1-t}\overline{\phi_{t}^{\infty}(\omega;f)} \\
&=&
\frac{\psi(\omega;f)}{\sqrt{n}} \sum_{t=1}^{n}X_{t} \sum_{s=0}^{\infty} \phi_{s+t}(f)e^{-is\omega}
+e^{i(n+1)\omega}
\frac{\overline{\psi(\omega;f)}}{\sqrt{n}} \sum_{t=1}^{n}X_{n+1-t} \sum_{s=0}^{\infty} \phi_{s+t}(f)e^{is\omega}
,
\end{eqnarray*} 
where 
$\psi(\omega;f) = \sum_{j=0}^{\infty} \psi_j(f) e^{-ij \omega}$ be an MA transfer function.
Our goal is to estimate $\psi(\omega;f)$ and $\{\phi_j(f)\}$ based on the observed time series.
We use the method proposed in Section 2.2. of \cite{p:kra-18}. 
We first start from the well known  Szeg{\"o}'s  identity 
\begin{eqnarray*}
\log f(\cdot) = \log \sigma^2
|\psi(\cdot;f)|^2 = \log \sigma^2 + \log \psi(\cdot;f) + 
\log \overline{\psi(\cdot;f)}.
\end{eqnarray*}
Next, let $\alpha_{k}(f)$ be the $k$-th Fourier coefficient of $\log f$,
i.e., $\alpha_{k}(f) = (2\pi)^{-1} \int_{-\pi}^{\pi} \log f(\lambda) e^{-ik \lambda} d\lambda$. Then,
since $\log f$ is real, $\alpha_{-k}(f) = \overline{\alpha_k(f)}$.
Plug in the expansion of $\log f$ to the above identity gives 
\begin{eqnarray*}
\log \psi(\omega;f) = 
 \sum_{j=1}^{\infty} \alpha_j(f) e^{-ij\omega}.
\end{eqnarray*} 
Using above identity, we estimator $\psi(\cdot;f)$. let $\widehat{f}$
be a spectral density estimator and let $\widehat{\alpha}_k$ be the estimated
$k$-th Fourier coefficient of $\log \widehat{f}$. Then define
\begin{eqnarray*}
\widehat{\psi}(\omega;\widehat{f}) = \exp \left( \sum_{j=1}^{M} \widehat{\alpha}_j e^{-ij\omega} \right)
\end{eqnarray*} 
for some large enough $M$.
To estimate the AR$(\infty)$ coefficients
we use the recursive formula in equation (2.7) in  \cite{p:kra-18},
\begin{eqnarray*}
\widehat{\phi}_{k+1} = -\sum_{j=0}^{k} \left( 1-\frac{j}{k+1} \right) \widehat{\alpha}_{k+1-j} \widehat{\phi}_{j}
\qquad k=0,1,...,M-1
\end{eqnarray*} 
where $\widehat{\phi}_0=-1$.
Based on this a nonparametric estimator of $\widehat{J}_{n}(\omega_{};f)$ is
\begin{eqnarray*}
\widehat{J}_n(\omega;\widehat{f})
=
\frac{\widehat{\psi}(\omega;\widehat{f})}{\sqrt{n}} \sum_{t=1}^{n \wedge M}X_{t} \sum_{s=0}^{M-t} \widehat{\phi}_{s+t}e^{-is\omega}
+e^{i(n+1)\omega}
\frac{\overline{\widehat{\psi}(\omega;\widehat{f})}}{\sqrt{n}} \sum_{t=1}^{n \wedge M}X_{n+1-t} \sum_{s=0}^{M -t} \widehat{\phi}_{s+t}e^{is\omega}
\end{eqnarray*} 
where $n \wedge M = \min(n,M)$.
In the simulations we estimate $\widehat{f}$ using \texttt{iospecden} function in R (smoothing with infinite order Flat-top kernel)
 and set $M$=30. 

By replacing $\widehat{J}_n(\omega;f)$ 
with its nonparametric estimator $\widehat{J}_n(\omega;\widehat{f})$
in (\ref{eq:Winf}) leads us to define a new feasible criterion which we call the  
\textbf{b}oundary \textbf{c}orrected Whittle likelihood using
\textbf{N}on\textbf{p}arametric  estimation (BC-NP). 

\subsection{Alternative methods for estimating the predictive DFT results for a Gaussian time series}

To access the performance of all the different likelihoods (with different 
estimates of the predictive DFT), we generate the AR$(8)$ model 
\begin{eqnarray*}
 U_{t} = \phi_{U}(B)\varepsilon_{t}
\end{eqnarray*}
where $\{\varepsilon_{t}\}$ are i.i.d. normal random variables, 
\begin{eqnarray}
\label{eq:phiU}
\phi_U(z) = \prod_{j=1}^{4} (1-r_j e^{i \lambda_j}z) (1-r_j e^{-i \lambda_j}z) = 1-\sum_{j=1}^{8}\phi_j z^{j}
\end{eqnarray} 
$\underline{r} = (r_1,r_2,r_3,r_4) = (0.95,0.95,0.95,0.95)$ and
$\underline{\lambda} = (\lambda_1, \lambda_2, \lambda_3,\lambda_4) = (0.5,1,2,2.5)$. 
We observe that corresponding spectral density $f_U(\omega) = |\phi_U(e^{-i\omega})|^{-2}$
has pronounced peaks at $\omega=0.5,1,1.5$ and $2$. For all the simulations below we use $n=100$.

For each simulation, we fit AR$(8)$ model, evaluate six likelihoods from the previous sections plus two likelihoods (BC-tYW and BC-NP),
and calculate the parameter estimators. Table \ref{tab:AR8} summarizes the bias and standard derivation of the estimators and the last row is 
an average $\ell_2$-distance between the true and estimator scaled with $n$. The Gaussian likelihood has the smallest bias and the smallest RMSE. 
As mentioned in Section \ref{sec:specified}, our methods still need to estimate AR coefficients which has an additional error of order $O(p^3n^{-3/2})$
and it could potentially increase the bias compared to the Gaussian likelihood.
The boundary corrected Whittle and hybrid Whittle have smaller bias than the Whittle, tapered, and debiased Whittle. Especially, the hybrid Whittle usually has the second smallest RMSE.

\begin{table}[ht]
\centering
\scriptsize
\begin{tabular}{c|cccccccc}

 \multirow{2}{*}{Par.} & 
 \multicolumn{8}{c}{Bias} \\ \cline{2-9}
&
Gaussian & Whittle & {\color{blue}Boundary} & {\color{blue}Hybrid} & Tapered & Debiased
& {\color{blue} BC-tYW} & {\color{blue} BC-NP} \\ \hline \hline
	
$\phi_1 (0.381)$ & {\color{red}-$0.008(0.08)$} & -$0.025(0.09)$ & -$0.009(0.08)$ & -$0.006(0.09)$ & -$0.012(0.09)$ & -$0.008(0.09)$ & {\color{blue}-$0.008(0.08)$} & -$0.005(0.12)$ \\ 

$\phi_2 (\text{-}0.294)$ & {\color{red}$0.002(0.09)$} & $0.024(0.1)$ & $0.005(0.09)$ & {\color{blue}$0.002(0.09)$} & $0.010(0.09)$ & $0.003(0.1)$ & $0.003(0.09)$ & $0.002(0.13)$ \\ 

$\phi_3 (0.315)$ & {\color{red}-$0.009(0.08)$} & -$0.038(0.09)$ & -$0.011(0.09)$ & {\color{blue}-$0.009(0.09)$} & -$0.023(0.09)$ & -$0.010(0.09)$ & -$0.009(0.09)$ & -$0.010(0.12)$ \\ 

$\phi_4 (\text{-}0.963)$ & {\color{red}$0.031(0.09)$} & $0.108(0.1)$ & $0.042(0.09)$ & {\color{blue}$0.034(0.09)$} & $0.075(0.09)$ & $0.043(0.1)$ & $0.037(0.09)$ & $0.076(0.12)$ \\ 

$\phi_5 (0.285)$ & {\color{red}-$0.015(0.08)$} & -$0.049(0.09)$ & -$0.020(0.09)$ & -$0.016(0.08)$ & {\color{blue}-$0.029(0.08)$} & -$0.017(0.1)$ & -$0.018(0.09)$ & -$0.022(0.12)$ \\ 

$\phi_6 (\text{-}0.240)$ & {\color{red}$0.010(0.08)$} & $0.040(0.09)$ & $0.014(0.09)$ & $0.010(0.09)$ & {\color{blue}$0.024(0.08)$} & $0.012(0.1)$ & $0.011(0.09)$ & $0.022(0.11)$ \\ 

$\phi_7 (0.280)$ & {\color{red}-$0.017(0.08)$} & -$0.053(0.09)$ & -$0.021(0.09)$ & -{\color{blue}$0.020(0.09)$} & -$0.039(0.08)$ & -$0.022(0.09)$ & -$0.020(0.09)$ & -$0.027(0.1)$ \\ 

$\phi_8 (\text{-}0.663)$ & {\color{red}$0.049(0.08)$} & $0.116(0.08)$ & $0.059(0.08)$ & {\color{blue}$0.055(0.08)$} & $0.096(0.08)$ & $0.061(0.09)$ & $0.056(0.08)$ &$0.101(0.1)$ \\

$n\|\underline{\phi}-\widehat{\underline{\phi}}\|_2$  & {\color{red}$6.466$} & $18.607$ & $8.029$ & {\color{blue}$7.085$} & $13.611$ & $8.164$ & $7.470$ & $13.280$ \\ 

\end{tabular} 
\caption{\textit{ Bias and the standard deviation (in the parenthesis) of eight different quasi-likelihoods for the Gaussian AR$(8)$ model. Length of time series $n$=100.
True AR coefficients are in the parenthesis of the first column. 
}}
\label{tab:AR8}
\end{table}

Bear in mind that neither of the two new criteria uses a hybrid method (tapering on the actual DFT), the
BC-tYW significantly reduces the bias than the boundary corrected Whittle and it is comparable with the hybrid Whittle. 
This gives some credence to the referee's claim that the bias due to the Yule-Walker estimation can be alleviated using tapered Yule-Walker estimation.
Whereas, BC-NP reduces the bias for the first few coefficients but overall, has a larger bias than the boundary corrected Whittle. Also, the standard deviation of BC-NP is quite large than other methods. We suspect that the nonparametric estimator $\widehat{J}(\omega;\widehat{f})$ is
sensitive to the choice of the tuning parameters (e.g. bandwidth, kernel function, etc). Moreover, since the true model follows a finite autoregressive process, other methods (boundary corrected Whittle, BC-tYW, and hybrid Whittle) have an advantage over the nonparametric method. Therefore, by choosing appropriate tuning parameters under certain underlying process (e.g., seasonal ARMA model) can improve the estimators, and this will be 
investigated in future research.

\subsection{Alternative methods for estimating the predictive DFT results for a non-Gaussian time series}

This time we assess the different estimation schemes for non-Gaussian time series. We generate the same AR$(8)$ model as above with
\begin{eqnarray*}
 V_{t} = \phi_{U}(B)\varepsilon_{t}
\end{eqnarray*}
where $\{\varepsilon_t\}$ are i.i.d. standardarized chi-square random variables with 
two-degrees of freedom i.e.  $\varepsilon_t \sim (\chi^2(2)-2)/2$ and 
$\phi_{U}(z)$ is defined as in (\ref{eq:phiU}). 
For each simulation, we fit AR$(8)$ model, evaluate six likelihoods from the previous sections plus two likelihoods (BC-tYW and BC-NP), and calculate the parameter estimators.
The results are summarized in  Table \ref{tab:AR8.chisq}.

\begin{table}[ht]
\centering
\scriptsize 
\begin{tabular}{c|cccccccc}

 \multirow{2}{*}{Par.} & 
 \multicolumn{8}{c}{Bias} \\ \cline{2-9}
&
Gaussian & Whittle & {\color{blue}Boundary} & {\color{blue}Hybrid} & Tapered & Debiased
& {\color{blue} BC-tYW} & {\color{blue} BC-NP} \\ \hline \hline
	
$\phi_1 (0.381)$ & {\color{red}$0.001(0.08)$}& -$0.013(0.09)$ & -$0.002(0.09)$ & $0.001(0.09)$ & -$0.003(0.09)$ & $0.004(0.09)$ & {\color{blue}$0(0.09)$} & $0.001(0.12)$ \\ 

$\phi_2 (\text{-}0.294)$ & {\color{red}-$0.001(0.09)$} & $0.014(0.1)$ & -$0.001(0.09)$ & {\color{blue} -$0.002(0.09)$} & $0.006(0.09)$ & -$0.008(0.11)$ & -$0.002(0.09)$ & -$0.010(0.13)$ \\ 

$\phi_3 (0.315)$ & {\color{red}-$0.004(0.09)$} & -$0.027(0.1)$ & -$0.005(0.09)$ & {\color{blue} -$0.003(0.09)$} & -$0.015(0.09)$ & $0(0.1)$ & -$0.003(0.09)$ & -$0.005(0.12)$ \\ 

$\phi_4 (\text{-}0.963)$ & {\color{red}$0.034(0.09)$} & $0.097(0.09)$ & $0.040(0.09)$ & {\color{blue}$0.034(0.09)$} & $0.073(0.09)$ & $0.038(0.11)$ & $0.036(0.09)$ & $0.068(0.12)$ \\ 

$\phi_5 (0.285)$ & {\color{red}-$0.007(0.09)$} & -$0.032(0.09)$ & -$0.009(0.09)$ & -$0.005(0.09)$ & {\color{blue}-$0.018(0.09)$} & -$0.004(0.1)$ & -$0.007(0.09)$ & -$0.005(0.12)$ \\ 

$\phi_6 (\text{-}0.240)$ & {\color{red}$0.007(0.09)$} & $0.029(0.09)$ & $0.009(0.09)$ & $0.006(0.09)$ & {\color{blue}$0.018(0.09)$} & $0.003(0.1)$ & $0.007(0.09)$ & $0.006(0.12)$ \\ 

$\phi_7 (0.280)$ & {\color{red}-$0.019(0.08)$} & -$0.047(0.09)$ & -$0.021(0.09)$ & {\color{blue}-$0.018(0.09)$} & -$0.034(0.09)$ & -$0.020(0.1)$ & -$0.019(0.09)$ & -$0.026(0.11)$ \\ 

$\phi_8 (\text{-}0.663)$ & {\color{red}$0.058(0.08)$} & $0.114(0.08)$ & $0.062(0.09)$ & {\color{blue}$0.059(0.09)$} & $0.098(0.08)$ & $0.065(0.1)$ & $0.060(0.08)$ &$0.107(0.1)$ \\

$n\|\underline{\phi}-\widehat{\underline{\phi}}\|_2$  & {\color{red}$7.006$} & $16.607$ & $7.728$ & {\color{blue}$7.107$} & $13.054$ & $7.889$ & $7.319$ & $13.001$ \\ 

\end{tabular} 
\caption{\textit{Bias and the standard deviation (in the parenthesis) of eight different quasi-likelihoods for the AR$(8)$ model for the standardized chi-squared innovations. Length of time series $n$=100.
True AR coefficients are in the parenthesis of the first column. We use {\color{red}red} to denote the smallest RMSE and 
{\color{blue}blue} to denote the second smallest RMSE.
}}
\label{tab:AR8.chisq}
\end{table}

%% file: whittle_likelihood_revision.bbl
\begin{thebibliography}{57}
\providecommand{\natexlab}[1]{#1}
\providecommand{\url}[1]{\texttt{#1}}
\expandafter\ifx\csname urlstyle\endcsname\relax
  \providecommand{\doi}[1]{doi: #1}\else
  \providecommand{\doi}{doi: \begingroup \urlstyle{rm}\Url}\fi

\bibitem[Abadir et~al.(2007)Abadir, Distaso, and Giraitis]{p:abd-07}
K.~M. Abadir, W.~Distaso, and L.~Giraitis.
\newblock Nonstationarity-extended local {W}hittle estimation.
\newblock \emph{J. Econometrics}, 141\penalty0 (2):\penalty0 1353--1384, 2007.

\bibitem[Bartlett(1953)]{p:bar-53}
M.~S. Bartlett.
\newblock Approximate confidence intervals. {II}. {M}ore than one unknown
  parameter.
\newblock \emph{Biometrika}, 40\penalty0 (3/4):\penalty0 306--317, 1953.

\bibitem[Baxter(1962)]{p:bax-62}
G.~Baxter.
\newblock An asymptotic result for the finite predictor.
\newblock \emph{Math. Scand.}, 10:\penalty0 137--144, 1962.

\bibitem[Baxter(1963)]{p:bax-63}
G.~Baxter.
\newblock A norm inequality for a ``finite-section'' {W}iener-{H}opf equation.
\newblock \emph{Illinois J. Math.}, 7\penalty0 (1):\penalty0 97--103, 1963.

\bibitem[Bhansali(1982)]{p:bha-82}
R.~J. Bhansali.
\newblock The evaluation of certain quadratic forms occurring in autoregressive
  model fitting.
\newblock \emph{Ann. Statist.}, 10\penalty0 (1):\penalty0 121--131, 1982.

\bibitem[Bhansali(1996)]{p:bha-96}
R.~J. Bhansali.
\newblock Asymptotically efficient autoregressive model selection for multistep
  prediction.
\newblock \emph{Ann. Inst. Statist. Math.}, 48\penalty0 (3):\penalty0 577--602,
  1996.

\bibitem[B{\"o}ttcher and Silbermann(2013)]{b:bot-13}
Albrecht B{\"o}ttcher and Bernd Silbermann.
\newblock \emph{Analysis of Toeplitz operators}.
\newblock Springer Science \& Business Media, 2013.

\bibitem[Brillinger(2001)]{b:bri-01}
David~R. Brillinger.
\newblock \emph{Time series: Data Analysis and theory}, volume~36 of
  \emph{Classics Appl. Math.}
\newblock SIAM, Philadelphia, PA, 2001.

\bibitem[Chen and Hurvich(2003)]{p:che-03}
W.~W. Chen and C.~M. Hurvich.
\newblock Semiparametric estimation of multivariate fractional cointegration.
\newblock \emph{J. Amer. Statist. Assoc.}, 98\penalty0 (463):\penalty0
  629--642, 2003.

\bibitem[Cheng and Pourahmadi(1993)]{p:che-pou-93}
R.~Cheng and M.~Pourahmadi.
\newblock Baxter's inequality and convergence of finite predictors of
  multivariate stochastic processes.
\newblock \emph{Probab. Theory Related Fields}, 95:\penalty0 115--124, 1993.

\bibitem[Choudhuri et~al.(2004)Choudhuri, Ghosal, and Roy]{p:cho-04}
N.~Choudhuri, S.~Ghosal, and A.~Roy.
\newblock Bayesian estimation of the spectral density of a time series.
\newblock \emph{J. Amer. Statist. Assoc.}, 99\penalty0 (468):\penalty0
  1050--1059, 2004.

\bibitem[Coursol and Dacunha-Castelle(1982)]{p:cou-82}
J.~Coursol and D.~Dacunha-Castelle.
\newblock Remarques sur l'approximation de la vraisemblance d'un processus
  {G}aussien stationnaire.
\newblock \emph{Teor. Veroyatnost. i Primenen. (Theory of Probability and its
  Applications)}, 27\penalty0 (1):\penalty0 155--160, 1982.

\bibitem[Cox and Snell(1968)]{p:cox-68}
D.~R. Cox and E.~J. Snell.
\newblock A general definition of residuals.
\newblock \emph{J. R. Stat. Soc. Ser. B. Stat. Methodol.}, 30\penalty0
  (2):\penalty0 248--265, 1968.

\bibitem[Dahlhaus(1983)]{p:dah-83}
R.~Dahlhaus.
\newblock Spectral analysis with tapered data.
\newblock \emph{J. Time Series Anal.}, 4\penalty0 (3):\penalty0 163--175, 1983.

\bibitem[Dahlhaus(1988)]{p:dah-88}
R.~Dahlhaus.
\newblock Small sample effects in time series analysis: a new asymptotic theory
  and a new estimate.
\newblock \emph{Ann. Statist.}, 16\penalty0 (2):\penalty0 808--841, 1988.

\bibitem[Dahlhaus(1990)]{p:dah-90}
R.~Dahlhaus.
\newblock Nonparametric high resolution spectral estimation.
\newblock \emph{Probability Theory and Related Fields}, 85\penalty0
  (2):\penalty0 147--180, 1990.

\bibitem[Dahlhaus(2000)]{p:dah-00}
R.~Dahlhaus.
\newblock A likelihood approximation for locally stationary processes.
\newblock \emph{Ann. Statist.}, 28\penalty0 (6):\penalty0 1762--1794, 2000.

\bibitem[Dahlhaus and K\"{u}nsch(1987)]{p:dah-87}
R.~Dahlhaus and H.~K\"{u}nsch.
\newblock Edge effects and efficient parameter estimation for stationary random
  fields.
\newblock \emph{Biometrika}, 74\penalty0 (4):\penalty0 877--882, 1987.

\bibitem[Das et~al.(2020)Das, Subba~Rao, and Yang]{p:dsy-20}
S.~Das, S.~Subba~Rao, and J.~Yang.
\newblock Spectral methods for small sample time series: A complete periodogram
  approach.
\newblock \emph{arXiv preprint arXiv:2007.00363}, 2020.

\bibitem[Fox and Taqqu(1986)]{p:fox-taq-86}
R.~Fox and M.~S. Taqqu.
\newblock Large-sample properties of parameter estimates for strongly dependent
  stationary {G}aussian time series.
\newblock \emph{Ann. Statist.}, 14\penalty0 (2):\penalty0 517--532, 1986.

\bibitem[Galbraith and Galbraith(1974)]{p:gal-74}
R.~F. Galbraith and J.~I. Galbraith.
\newblock On the inverses of some patterned matrices arising in the theory of
  stationary time series.
\newblock \emph{J. Appl. Probab.}, 11\penalty0 (1):\penalty0 63--71, 1974.

\bibitem[Giraitis and Robinson(2001)]{p:gir-01}
L.~Giraitis and P.~M. Robinson.
\newblock Whittle estimation of {ARCH} models.
\newblock \emph{Econometric Theory}, 17\penalty0 (3):\penalty0 608--631, 2001.

\bibitem[Giraitis et~al.(2012)Giraitis, Koul, and Surgailis]{b:gir-12}
Liudas Giraitis, Hira~L. Koul, and Donatas Surgailis.
\newblock \emph{Large sample inference for long memory processes}.
\newblock Imperial College Press, London, 2012.

\bibitem[Hurvich and Chen(2000)]{p:che-00}
C.~M. Hurvich and W.~W. Chen.
\newblock An efficient taper for potentially overdifferenced long-memory time
  series.
\newblock \emph{J. Time Series Anal.}, 21\penalty0 (2):\penalty0 155--180,
  2000.

\bibitem[Ing and Wei(2005)]{p:ing-wei-05}
C.-K. Ing and C.-Z. Wei.
\newblock Order selection for same-realization predictions in autoregressive
  processes.
\newblock \emph{Ann. Statist.}, 33\penalty0 (5):\penalty0 2423--2474, 2005.

\bibitem[Inoue and Kasahara(2006)]{p:ino-06}
A.~Inoue and Y.~Kasahara.
\newblock Explicit representation of finite predictor coefficients and its
  applications.
\newblock \emph{Ann. Statist.}, 34\penalty0 (2):\penalty0 973--993, 2006.

\bibitem[Inoue et~al.(2018)Inoue, Kasahara, and Pourahmadi]{p:ino-18}
A.~Inoue, Y.~Kasahara, and M.~Pourahmadi.
\newblock Baxter's inequality for finite predictor coefficients of multivariate
  long-memory stationary processes.
\newblock \emph{Bernoulli}, 24\penalty0 (2):\penalty0 1202--1232, 2018.

\bibitem[Kasahara et~al.(2009)Kasahara, Pourahmadi, and Inoue]{p:kas-09}
Y.~Kasahara, M.~Pourahmadi, and A.~Inoue.
\newblock Duals of random vectors and processes with applications to prediction
  problems with missing values.
\newblock \emph{Statist. Probab. Lett.}, 79\penalty0 (14):\penalty0 1637--1646,
  2009.

\bibitem[Kirch et~al.(2019)Kirch, Edwards, Meier, and Meyer]{p:kir-19}
C.~Kirch, M.~C. Edwards, A.~Meier, and R.~Meyer.
\newblock Beyond {W}hittle: Nonparametric correction of a parametric likelihood
  with a focus on bayesian time series analysis.
\newblock \emph{Bayesian Anal.}, 14\penalty0 (4):\penalty0 1037--1073, 2019.

\bibitem[Kley et~al.(2019)Kley, Preu{\ss}, and Fryzlewicz]{p:kle-19}
T.~Kley, P.~Preu{\ss}, and P.~Fryzlewicz.
\newblock Predictive, finite-sample model choice for time series under
  stationarity and non-stationarity.
\newblock \emph{Electron. J. Stat.}, 13\penalty0 (2):\penalty0 3710--3774,
  2019.

\bibitem[Krampe et~al.(2018)Krampe, Kreiss, and Paparoditis]{p:kra-18}
J.~Krampe, J.-P. Kreiss, and E.~Paparoditis.
\newblock Estimated {W}old representation and spectral-density-driven bootstrap
  for time series.
\newblock \emph{J. R. Stat. Soc. Ser. B. Stat. Methodol.}, 80:\penalty0
  703--726, 2018.

\bibitem[Kreiss et~al.(2011)Kreiss, Paparoditis, and Politis]{p:kre-11}
J.-P. Kreiss, E.~Paparoditis, and D.~N. Politis.
\newblock On the range of validity of the autoregressive sieve bootstrap.
\newblock \emph{Ann. Statist.}, 39\penalty0 (4):\penalty0 2103--2130, 2011.

\bibitem[K{\"u}nsch(1987)]{p:kun-87}
H.~R. K{\"u}nsch.
\newblock Statistical aspects of self-similar processes.
\newblock In \emph{Proceedings of the 1st World Congress of the Bernoulli
  Society}, volume~1, pages 67--74. VNU Science Press, 1987.

\bibitem[Lahiri(2003)]{p:lah-03}
S.~N. Lahiri.
\newblock A necessary and sufficient condition for asymptotic independence of
  discrete {F}ourier transforms under short- and long-range dependence.
\newblock \emph{Ann. Statist.}, 31\penalty0 (2):\penalty0 613--641, 2003.

\bibitem[Lieberman(2005)]{p:lib-05}
O.~Lieberman.
\newblock On plug-in estimation of long memory models.
\newblock \emph{Econometric Theory}, 21\penalty0 (2):\penalty0 431--454, 2005.

\bibitem[Meyer et~al.(2017)Meyer, Jentsch, and Kreiss]{p:kre-17}
M.~Meyer, C.~Jentsch, and J.-P. Kreiss.
\newblock Baxter's inequality and sieve bootstrap for random fields.
\newblock \emph{Bernoulli}, 23\penalty0 (4B):\penalty0 2988--3020, 2017.

\bibitem[Panaretos and Tavakoli(2013)]{p:pan-13}
V.~M. Panaretos and S.~Tavakoli.
\newblock Fourier analysis of stationary time series in function space.
\newblock \emph{Ann. Statist.}, 41\penalty0 (2):\penalty0 568--603, 2013.

\bibitem[Parzen(1983)]{p:par-83}
Emanuel Parzen.
\newblock Autoregressive spectral estimation.
\newblock In \emph{Time series in the frequency domain}, volume~3 of
  \emph{Handbook of Statist.}, pages 221--247. North-Holland, Amsterdam, 1983.

\bibitem[Pourahmadi(2001)]{b:pou-01}
Mohsen Pourahmadi.
\newblock \emph{Foundations of time series analysis and prediction theory}.
\newblock Wiley Series in Probability and Statistics: Applied Probability and
  Statistics. Wiley-Interscience, New York, 2001.

\bibitem[Priestley(1981)]{b:pri-81}
Maurice~B. Priestley.
\newblock \emph{Spectral analysis and time series. {V}ol. 2}.
\newblock Academic Press, Inc. [Harcourt Brace Jovanovich, Publishers],
  London-New York, 1981.
\newblock Multivariate series, prediction and control, Probability and
  Mathematical Statistics.

\bibitem[Robinson(1995)]{p:rob-95}
P.~M. Robinson.
\newblock Gaussian semiparametric estimation of long range dependence.
\newblock \emph{Ann. Statist.}, 23\penalty0 (5):\penalty0 1630--1661, 1995.

\bibitem[Shaman(1975)]{p:sha-75}
P.~Shaman.
\newblock An approximate inverse for the covariance matrix of moving average
  and autoregressive processes.
\newblock \emph{Ann. Statist.}, 3\penalty0 (2):\penalty0 532--538, 1975.

\bibitem[Shaman(1976)]{p:sha-76}
P.~Shaman.
\newblock Approximations for stationary covariance matrices and their inverses
  with application to {ARMA} models.
\newblock \emph{Ann. Statist.}, 4\penalty0 (2):\penalty0 292--301, 1976.

\bibitem[Shaman and Stine(1988)]{p:sha-sti-88}
P.~Shaman and R.~A. Stine.
\newblock The bias of autoregressive coefficient estimators.
\newblock \emph{J. Amer. Statist. Assoc.}, 83\penalty0 (403):\penalty0
  842--848, 1988.

\bibitem[Shao and Wu(2007)]{p:sha-wu-07}
X.~Shao and W.~B. Wu.
\newblock Local whittle estimation of fractional integration for nonlinear
  processes.
\newblock \emph{Econometric Theory}, 23\penalty0 (5):\penalty0 899--929, 2007.

\bibitem[Siddiqui(1958)]{p:sid-58}
M.~M. Siddiqui.
\newblock On the inversion of the sample covariance matrix in a stationary
  autoregressive process.
\newblock \emph{Ann. Math. Statist.}, 29\penalty0 (2):\penalty0 585--588, 1958.

\bibitem[Subba~Rao(2018)]{p:sub-18}
S.~Subba~Rao.
\newblock Orthogonal samples for estimators in time series.
\newblock \emph{J. Time Series Anal.}, 39:\penalty0 313--337, 2018.

\bibitem[Sykulski et~al.(2019)Sykulski, Olhede, Guillaumin, Lilly, and
  Early]{p:olh-19}
A.~M. Sykulski, S.~C. Olhede, A.~P. Guillaumin, J.~M. Lilly, and J.~J. Early.
\newblock The debiased {W}hittle likelihood.
\newblock \emph{Biometrika}, 106\penalty0 (2):\penalty0 251--266, 2019.

\bibitem[Szeg\"o(1921)]{p:sze-21}
G.~Szeg\"o.
\newblock {\"U}ber die randwerte einer analytischen funktion.
\newblock \emph{Math. Ann.}, 84:\penalty0 232--244, 1921.

\bibitem[Tanaka(1984)]{p:tan-84}
K.~Tanaka.
\newblock An asymptotic expansion associated with the maximum likelihood
  estimators in {ARMA} models.
\newblock \emph{J. R. Stat. Soc. Ser. B. Stat. Methodol.}, 46\penalty0
  (1):\penalty0 58--67, 1984.

\bibitem[Taniguchi(1983)]{p:tan-83}
M.~Taniguchi.
\newblock On the second order asymptotic efficiency of estimators of gaussian
  {ARMA} processes.
\newblock \emph{Ann. Statist.}, 11:\penalty0 157--169, 1983.

\bibitem[Tj{\o}stheim and Paulsen(1983)]{p:tjo-83}
D.~Tj{\o}stheim and J.~Paulsen.
\newblock Bias of some commonly-used time series estimates.
\newblock \emph{Biometrika}, 70\penalty0 (2):\penalty0 389--399, 1983.

\bibitem[van Delft and Eichler(2020)]{p:del-eic-19}
A.~van Delft and M.~Eichler.
\newblock A note on {H}erglotz's theorem for time series on functional spaces.
\newblock \emph{Stochastic Process. Appl.}, 130\penalty0 (6):\penalty0
  3687--3710, 2020.

\bibitem[Walker(1964)]{p:wal-64}
A.~M. Walker.
\newblock Asymptotic properties of least-squares estimates of parameters of the
  spectrum of a stationary non-deterministic time-series.
\newblock \emph{J. Aust. Math. Soc.}, 4:\penalty0 363--384, 1964.

\bibitem[Whittle(1953)]{p:whi-53}
P.~Whittle.
\newblock The analysis of multiple stationary time series.
\newblock \emph{J. R. Stat. Soc. Ser. B. Stat. Methodol.}, 15:\penalty0
  125--139, 1953.

\bibitem[Whittle(1951)]{p:whi-51}
Peter Whittle.
\newblock \emph{Hypothesis {T}esting in {T}ime {S}eries {A}nalysis}.
\newblock Thesis, Uppsala University, 1951.

\bibitem[Zhang(1992)]{p:zha-92}
H.-C. Zhang.
\newblock Reduction of the asymptotic bias of autoregressive and spectral
  estimators by tapering.
\newblock \emph{J. Time Series Anal.}, 13\penalty0 (5):\penalty0 451--469,
  1992.

\end{thebibliography}
